%% file: m2an170020.tex
\documentclass{m2an}
\usepackage{ifpdf}
\usepackage{ifthen}

\usepackage[latin1]{inputenc}
\usepackage[T1]{fontenc}

\usepackage{array}
\usepackage{tabularx}
\usepackage{multirow}
\usepackage{booktabs}

\usepackage{graphicx}
\graphicspath{{figures/}}
\usepackage{color}
\usepackage[dvipsnames]{xcolor}
\usepackage{float}
\usepackage{caption}
\usepackage{subcaption}
\usepackage{rotating}

\usepackage{tablefootnote}

\usepackage{pgf,pgfplots,pgfplotstable}
\pgfplotsset{compat=newest}
\pgfplotsset{
xlabel near ticks,
ylabel near ticks,
label style={font=\small},
tick label style={font=\footnotesize},
legend style={font=\small},
xticklabel style={/pgf/number format/set thousands separator={\,}},
yticklabel style={/pgf/number format/set thousands separator={\,}},
tick scale binop=\times,
try min ticks=6,
legend pos=outer north east
}
\usepgfplotslibrary{external} 
\usepgfplotslibrary{units}

\usepackage{tikz,tkz-base,xstring}
\usetikzlibrary{shapes,decorations,shadows}
\usetikzlibrary{decorations.pathreplacing}
\usetikzlibrary{decorations.pathmorphing}
\usetikzlibrary{decorations.shapes}
\usetikzlibrary{decorations.text}
\usetikzlibrary{decorations.footprints}
\usetikzlibrary{decorations.fractals}
\usetikzlibrary{fadings}
\usetikzlibrary{patterns}
\usetikzlibrary{calc}
\usetikzlibrary{shapes.geometric}
\usetikzlibrary{shapes.gates.logic.IEC}
\usetikzlibrary{shapes.gates.logic.US}
\usetikzlibrary{fit,chains}
\usetikzlibrary{positioning}
\usepgflibrary{shapes}
\usetikzlibrary{scopes}
\usetikzlibrary{arrows}
\usetikzlibrary{arrows.meta}
\usetikzlibrary{backgrounds}
\usetikzlibrary{intersections}
\usetikzlibrary{matrix}
\usetikzlibrary{pgfplots.units}
\usetikzlibrary{trees}
\usetikzlibrary{hobby}

\newlength\figureheight
\newlength\figurewidth
\usepgfplotslibrary{units}

\pgfdeclaredecoration{penciline}{initial}{
\state{initial}[width=+\pgfdecoratedinputsegmentremainingdistance,auto corner on length=1mm]{
\pgfpathcurveto
{% from
\pgfqpoint{\pgfdecoratedinputsegmentremainingdistance}{\pgfdecorationsegmentamplitude}}{% control
\pgfpointadd{\pgfqpoint{\pgfdecoratedinputsegmentremainingdistance}{0pt}}{
\pgfqpoint{-\pgfdecorationsegmentaspect\pgfdecoratedinputsegmentremainingdistance}{\pgfmathresult\pgfdecorationsegmentamplitude}}}
{% to
\pgfpointadd{\pgfpointdecoratedinputsegmentlast}{\pgfpoint{1pt}{1pt}}}
}
\state{final}{}
}

\usepackage[all]{xy}

\usepackage{amsmath,amssymb,amsfonts,amsthm,bm}
\usepackage{textcomp,wasysym,eurosym,stmaryrd}
\usepackage{mathrsfs}
\usepackage{mathtools}
\usepackage{pifont}
\usepackage{upgreek}
\usepackage{centernot}
\usepackage{algorithm,algorithmic}

\usepackage[pdftex,colorlinks,citecolor=blue,filecolor=blue,linkcolor=blue,urlcolor=blue]{hyperref}

\usepackage{cite}
\usepackage{breakcites}

\usepackage{xspace}
\usepackage{enumitem}
\setlist[itemize]{label={--}}
\newtheorem{assumption}[thrm]{Assumption}

\def\Ac{\mathcal{A}}

\def\Dc{\mathcal{D}}
\def\Ec{\mathcal{E}}
\def\Fc{\mathcal{F}}

\def\Ic{\mathcal{I}}

\def\Mc{\mathcal{M}}
\def\Nc{\mathcal{N}}
\def\Oc{\mathcal{O}}

\def\Sc{\mathcal{S}}
\def\Tc{\mathcal{T}}
\def\Uc{\mathcal{U}}
\def\Vc{\mathcal{V}}
\def\Wc{\mathcal{W}}

\def\xE{\mathbb{E}}
\def\xV{\mathbb{V}}

\newcommand{\forallin}[2]{\forall #1 \in #2}
\newcommand{\existsin}[2]{\exists #1 \in #2}
\newcommand{\abs}[1]{\lvert#1\rvert}

\newcommand{\norm}[1]{\lVert#1\rVert}

\newcommand{\set}[1]{\{#1\}}

\newcommand{\setst}[2]{\{#1\mathrel{:}#2\}}

\newcommand{\scalprod}[2]{#1\cdot#2}
\newcommand{\scalproda}[2]{\langle#1,#2\rangle}

\newcommand{\restrictto}{\mathclose{}|\mathopen{}}

\def\xd{\rm d}

\DeclareMathOperator*{\argmax}{arg\,max}

\DeclareMathOperator{\spann}{span}
\newcommand{\spanset}[1]{\ensuremath\spann\set{#1}}
\newcommand{\spansetst}[2]{\ensuremath\spann\setst{#1}{#2}}
\DeclareMathOperator{\meas}{meas}
\DeclareMathOperator{\dist}{dist}
\DeclareMathOperator{\ceil}{ceil}

\newcommand{\interval}[4]{\mathopen{#1}#2 \mathclose{}\mathpunct{},#3 \mathclose{#4}}
\newcommand{\intervalcc}[2]{\interval{[}{#1}{#2}{]}}

\renewcommand{\(}{\left(}
\renewcommand{\)}{\right)}

\renewcommand{\geq}{\geqslant}
\renewcommand{\leq}{\leqslant}
\let\oldtimes\times
\renewcommand{\times}{\!\oldtimes\!}

\def\eg{\emph{e.g.\/}}
\def\aposteriori{\emph{a posteriori\/}}

\newcommand{\algorithmicbreak}{\textbf{break}}
\newcommand{\BREAK}{\STATE \algorithmicbreak}

\begin{document}

\title{A multiscale method for semi-linear elliptic equations with localized uncertainties and non-linearities}\thanks{The financial support of the French National Research Agency - Agence Nationale de la Recherche (ANR) - under Grant ICARE ANR-12-MONU-0002 is acknowledged by the authors.}
\runningtitle{Multiscale method for semi-linear elliptic stochastic equations}

\author{Anthony Nouy}\address{Centrale Nantes, LMJL UMR CNRS 6629, 1 rue de la No\"e, BP 92101, 44321 Nantes Cedex 3, France ;
\email{anthony.nouy@ec-nantes.fr}}
\author{Florent Pled}\address{Universit\'e Paris-Est, Laboratoire Mod\'elisation et Simulation Multi Echelle, MSME UMR 8208 CNRS, 5 bd Descartes, 77454 Marne-la-Vall\'ee, France ;
\email{florent.pled@univ-paris-est.fr}}
\date{\today}

\begin{abstract}
A multiscale numerical method is proposed for the solution of semi-linear elliptic stochastic partial differential equations with localized uncertainties and non-linearities, the uncertainties being modeled by a set of random parameters. It relies on a domain decomposition method which introduces several subdomains of interest (called patches) containing the different sources of uncertainties and non-linearities. 
An iterative algorithm is then introduced, which requires the solution of a sequence of linear global problems (with deterministic operators and uncertain right-hand sides), and non-linear local problems (with uncertain operators and/or right-hand sides) over the patches. Non-linear local problems are solved using an adaptive sampling-based least-squares method for the construction of sparse polynomial approximations of local solutions as functions of the random parameters. Consistency, convergence and robustness of the algorithm are proved under general assumptions on the semi-linear elliptic operator. A convergence acceleration technique (Aitken's dynamic relaxation) is also introduced to speed up the convergence of the algorithm. The performances of the proposed method are illustrated through numerical experiments carried out on a stationary non-linear diffusion-reaction problem.
\end{abstract}
\maketitle

Uncertainty quantification has become a topical issue in computational sciences and engineering. Numerous methods have been proposed to propagate uncertainties through models governed by partial differential equations (see \eg{} \cite{Nou09,Xiu09,LeMai10}). While these methods have reached a certain degree of maturity and become nowadays widespread, a major concern has emerged for multiscale models where uncertainties occur at various scales.

Several numerical methods dedicated to deterministic multiscale models have been extended to the stochastic framework. For multiscale problems with global sources of uncertainties, spectral stochastic methods have been combined with deterministic multiscale methods, \eg{} the multiscale finite element method (FEM) 
\cite{Hou97%,Hou99,Efe09
}, the variational multiscale method 
\cite{%Hug95,
Hug98} or the heterogeneous multiscale method 
\cite{E03%,E07
}, leading to the so-called multiscale stochastic FEM 
\cite{Xu07} and its variants \cite{Nar05,Aso06,Gan07bis,Dos08,Xu09,Gan09,Gin10,%Ma11,Jia12,
LeBri14}. These methods are well adapted to global uncertainties and are proved to be efficient when assuming small random fluctuations and scale separation. Meanwhile, traditional substructuring and domain decomposition methods have been introduced for stochastic monoscale models \cite{Jin07,Zha08,Sar09%,Gho09
} and recently extended to multiscale models \cite{Gan11,Whe11} in order to benefit from scalable parallel algorithms available in the deterministic framework. These methods are also well adapted to problems where the uncertainties are scattered in the whole domain.

The present work focuses on non-linear stochastic multiscale models where localized sources of uncertainties and non-linearities may occur in some regions of interest. Concurrent approaches, initially developed in the deterministic framework, have been proposed to couple numerical models. %exhibiting specific local features at various scales. 
First of all, mono-model approaches currently rely on adaptive remeshing 
strategies \cite{Ver96,Ste97%,Dus07
} or enrichment techniques (\eg{} the eXtended FEM 
\cite{Bel99bis,Moe99} or the Generalized FEM \cite{Str00bis}) and generally require high computational resources or specific (intrusive) implementations. Conversely, multi-model approaches based on patches have a high potential to manage complex multiscale problems by operating a separation of scales. The separation of scales allows to capture the local features of multiscale solutions at a micro scale (local level) while keeping a simplified global description at a macro scale (global level). Several multiscale coupling methods have been developed within the deterministic framework and some have been extended to the stochastic framework. They distinguish themselves by the way of defining and constructing the coupling operator between global and local models. First, superposition methods, such as the method of finite element patches \cite{Lio99,%Bre01,Glo03,
Glo05%,Loz07,Kam07
} and the method of harmonic patches \cite{He07}, consist in adding a fine local correction to a coarse global solution. Second, surface coupling methods include the Chimera-Schwarz method \cite{Ste83,%Ben83,Lio88,
Bre01%,Hec09,Loz11
}, the Semi-Schwarz method \cite{Pir09}, the Semi-Schwarz-Lagrange method \cite{%Whi91,
Gen09,Loz10,Gen11,Hag12} and the local multigrid method \cite{%Par90,
Gra08,Ran09,Pas11,Pas13}. Both multiscale superposition and surface coupling methods are based on global-local iterative algorithms originally developed for domain decomposition methods or multigrid methods. Nevertheless, the former can be interpreted as a local model refinement technique, while the latter can be seen as a local model substitution technique. Third, volume coupling methods, such as the Arlequin method \cite{Dhia98%,Dhia01,Dhia05
}, belong to the class of overlapping domain decomposition methods and require the definition of a coupling zone between the different models. Among all these multi-model approaches, few have been explored in the stochastic framework. The Arlequin (volume coupling) method has been applied to deterministic-stochastic coupling in \cite{Cha08b,Cot11} for homogenization purposes. Besides, the Semi-Schwarz-Lagrange (surface coupling) method has been recently extended to linear stochastic multiscale models with localized sources of uncertainties in \cite{Che13a}.

This work extends \cite{Che13a} to a class of non-linear stochastic multiscale models. Alternative multiscale approaches have been recently proposed to handle non-linear elliptic problems. For semi-linear elliptic equations, we refer to the variational multiscale method proposed in \cite{Hen14}, while for other non-linear elliptic equations, we refer to the multiscale FEM presented in \cite{Efe04,Efe09}, the variational multiscale method developed in \cite{Nor10} or the heterogeneous multiscale method proposed in \cite{Hen15}. In the present work, a dedicated multiscale method based on a domain decomposition is proposed to exploit the localized side of uncertainties and non-linearities. It relies on a global-local iterative algorithm which requires the solution of a sequence of linear global problems (with deterministic operators and uncertain right-hand sides) at a macro scale and non-linear local problems (with uncertain operators and right-hand sides) at a micro scale (over patches). Appropriate approximation spaces and solvers can be considered to solve both types of problems efficiently. This multiscale approach then appears to be flexible and non-intrusive in the sense that it requires no modification of both global and local models and solvers, which makes possible the use of stand-alone codes. The main motivation is the deployment and transfer of methods towards complex large-scale industrial applications \cite{All11}. Besides, different types of uncertainties can be considered in the non-linear local models. They may be associated with some variabilities of the operator but also of the geometry, the source terms or the boundary conditions. 

The remainder of the paper is structured as follows. Section~\ref{sec:formulation} introduces the initial formulation of the semi-linear elliptic stochastic partial differential equation with localized uncertainties and non-linearities and states suitable assumptions. Section~\ref{sec:formulationpatch} presents the global-local (two-scale) formulation with patches containing localized variabilities and non-linearities. A global-local iterative algorithm is then introduced and analyzed in Section~\ref{sec:algo}. Consistency, convergence and robustness properties are deduced from the assumptions introduced in Section~\ref{sec:formulation}. Subsequently, the computational aspects associated with the solution of both global and local problems are detailed in Section~\ref{sec:computation}. In particular, the stochastic local problems are solved using sampling-based (non-intrusive) approaches and working sets algorithms proposed in \cite{Chk13} for the adaptive construction of sparse polynomial approximations of local solutions. Finally, the proposed method is illustrated through numerical examples in Section~\ref{sec:results}. 
\section{Problem statement}\label{sec:formulation}

Let $\xi$ denote a set of real-valued random variables modeling the different sources of uncertainties (on the operator, geometry, source terms and boundary conditions). We assume that $\xi$ takes values in a set $\Xi$ 
 and we let $\mu$ be the probability law of $\xi$.
We consider the following semi-linear second-order stochastic partial differential equation
\begin{subequations}\label{strongformulation}
\begin{equation}\label{initialpb}
-\nabla \cdot B(u,\nabla u;x,\xi) + C(u,\nabla u;x,\xi) = f(x,\xi) \quad \text{for } x \in \Omega(\xi),
\end{equation}
where $\Omega(\xi)$ is an uncertain domain of $\xR^d$ 
with sufficiently smooth (\eg{} Lipschitz) boundary $\partial \Omega(\xi)$.
Here $B(\cdot,\cdot;x,\xi) \colon \xR \times \xR^d \to \xR^d$ and $C(\cdot,\cdot;x,\xi) \colon \xR \times \xR^d \to \xR$, and $f(\cdot,\xi) \colon \Omega(\xi) \to \xR$ is a given source term. For a given value of $\xi$, the solution $u(\cdot,\xi)$ is a function from $\Omega(\xi)$ to $\xR$.
We supply equation \eqref{initialpb} with the following Dirichlet and Neumann boundary conditions 
\begin{alignat}{2}
&u = 0 \quad &&\text{on } \Gamma_D(\xi),\\
&B(u,\nabla u;x,\xi) \cdot n = g(x,\xi) \quad &&\text{on } \Gamma_N(\xi),
\end{alignat}
\end{subequations}
where $\Gamma_D(\xi)$ and $\Gamma_N(\xi)$ are disjoint and complementary parts of $\partial \Omega(\xi)$ such that $\overline{\Gamma_D(\xi) \cup \Gamma_N(\xi)} = \partial \Omega(\xi)$ and $\Gamma_D(\xi) \cap \Gamma_N(\xi) = \emptyset$, and $\meas(\Gamma_D(\xi)) \neq 0$. $g(\cdot,\xi) \colon \Gamma_N(\xi) \to \xR$ is a prescribed normal flux on $\Gamma_N(\xi)$, and $n$ is the unit outward normal to $\Gamma_N(\xi)$.

\begin{xmpl}[Non-linear diffusion-reaction equation]\label{ex:diffusionreactionpb}
As a model example, we consider a non-linear diffusion-reaction equation \eqref{initialpb} in dimension $d \leq 3$, with 
\begin{equation*}
B(u,\nabla u;x,\xi) = K(x,\xi) \nabla u \quad \text{and} \quad C(u,\nabla u;x,\xi) = R(x,\xi) u^3,
\end{equation*}
where $K$ and $R$ are respectively the diffusion and reaction coefficients. This example will serve as a guideline. Such a semi-linear second-order stochastic partial differential equation describes transport phenomena at equilibrium such as steady-state diffusion-reaction processes arising from mathematical models in population dynamics \cite{Baz69,Aro75,%Aro78,
Bra04} as well as in chemical kinetics  (kinetics/dynamics of autocatalytic chemical reactions% such as a third-order autocatalytic reaction
) \cite{Han82,Hua96,Edw02,Kae02,Spa03,Lec03,Kop08,Sah13}.
\end{xmpl}

\subsection{Localized uncertainties and non-linearities}

We consider that non-linearities and uncertainties on operator and geometry only affect a given subdomain of interest $\Lambda_{\star} \subset \Omega$. For the sake of simplicity, we also consider that uncertainties on the right-hand side are localized in $\Lambda_{\star}$.

First, the subdomain $\Lambda_{\star}$ may depend on $\xi$ while the complementary subdomain $\Omega \setminus \Lambda_{\star}$ is supposed independent of $\xi$, which means that geometrical uncertainties are contained in $\Lambda_{\star}$.
The boundary $\partial \Lambda_{\star}$ of $\Lambda_{\star}$ contains the possible uncertainties of the boundary $\partial \Omega$ of domain $\Omega$.

Also, $B$ and $C$ are supposed linear and independent of $\xi$ outside $\Lambda_{\star}$.
More precisely, we suppose that $B$ can be split into a linear part $B_L$ (such that $u \mapsto B_L(u,\nabla u;x,\xi)$ is linear) and a non-linear part $B_N$, such that $B=B_L+B_N$. The same decomposition is introduced for $C = C_L+C_N$. Then, we consider
that $B$ and $C$ are such that
\begin{equation}\label{localizationuncertainties}
B_L(\cdot,\cdot;x,\xi) = B_L(\cdot,\cdot;x) \quad \text{and} \quad C_L(\cdot,\cdot;x,\xi) = C_L(\cdot,\cdot;x) \quad \text{for } x \in \Omega \setminus \Lambda_{\star},
\end{equation}
and 
\begin{equation}\label{localizationnonlinearities}
B_N(\cdot,\cdot;x,\xi) = C_N(\cdot,\cdot;x,\xi) = 0 \quad \text{for } x \in \Omega \setminus \Lambda_{\star}.
\end{equation}
Also, the prescribed source term $f$ and normal flux $g$ are such that
\begin{equation}\label{localizationuncertaintiesrhs}
f(x,\xi) = f(x) \quad \text{and} \quad g(x,\xi) = g(x) \quad \text{for } x \in \Omega \setminus \Lambda_{\star}.
\end{equation}

\begin{xmpl}
In Example~\ref{ex:diffusionreactionpb}, we consider that diffusion coefficient $K$ and reaction parameter $R$ are such that $K(x,\xi)=K(x)$ and $R(x,\xi)=0$ for $x \in \Omega \setminus \Lambda_{\star}$. A practical application is the chlorite-thiosulfate autocatalytic reaction \cite{Kop08} or the iodate-arsenous acid autocatalytic reaction %(evolution of a propagating front of iodide in an iodate-arsenous acid solution) 
\cite{Han82,Hua96,Edw02,Spa03,Lec03,Sah13} described by a third-order reaction-diffusion equation for the iodide concentration $u$ in the absence of convection, with localized and random molecular diffusion coefficient $K$ and reaction rate kinetic coefficient $R$.
\end{xmpl}

In the following, a function $v(x,\xi)$ of two variables defined for $\xi \in \Xi$ and $x \in D(\xi)$, with $D(\xi)$ a parametrized domain of $\xR^d$, will be equivalently considered as a function $v(\xi)$ defined on $D(\xi)$. For the sake of readability, we will often omit the dependence on $\xi$ for geometrical domains and for function spaces defined on these domains.

\subsection{Assumptions}

Here $\Oc$ denotes a subset of $\Omega$. For a function $v \in \xHone(\Oc)$, we denote
\begin{equation*}
\abs{v}_{\xHone(\Oc)} = \norm{\nabla v}_{\xLtwo(\Oc)} \quad \text{ and } \quad \norm{v}_{\xHone(\Oc)}^2 = \abs{v}_{\xHone(\Oc)}^2 + \norm{v}_{\xLtwo(\Oc)}^2.
\end{equation*}

\subsubsection{Assumptions on the source terms}

For a function $v$ defined on $\Oc$, we introduce the linear form
\begin{equation}
\ell_{\Oc}(v;\xi) = \int_{\Oc} f(\cdot,\xi) v + \int_{\Gamma_N \cap \partial \Oc} g(\cdot,\xi) v,
\end{equation}
and we assume that $f$ and $g$ are such that the following assumption holds.
\begin{assumption}[Properties of linear form $\ell_{\Oc}$]\label{assdata}
We assume that the linear form $\ell_{\Oc}(\cdot;\xi) \colon \xHone(\Oc) \to \xR$ is almost surely continuous, that means there exists a random variable $\kappa(\xi)>0$ such that it holds 
\begin{equation}\label{continuity_l}
\abs{\ell_{\Oc}(v;\xi)} \leq \kappa(\xi) \norm{v}_{\xHone(\Oc)} \quad \forallin{v}{\xHone(\Oc)},
\end{equation}
and we further assume that $\kappa \in \xLn{p}_{\mu}(\Xi)$ for some $2 \leq p \leq +\infty$.
\end{assumption}

\subsubsection{Assumptions on the differential operator}

For functions $u,v$ defined on $\Oc$, we introduce the semi-linear form
\begin{equation}
d_{\Oc}(u,v;\xi) = \int_{\Oc} B(u,\nabla u;\cdot,\xi) \cdot \nabla v + \int_{\Oc} C(u,\nabla u;\cdot,\xi) v,
\end{equation}
which can be written as
\begin{equation*}
d_{\Oc}(u,v;\xi) = a_{\Oc}(u,v;\xi) + n_{\Oc}(u,v;\xi),
\end{equation*}
where $a_{\Oc}(\cdot,\cdot;\xi)$ is a bilinear form and $n_{\Oc}(\cdot,\cdot;\xi)$ is a semi-linear form, respectively defined by
\begin{align*}
a_{\Oc}(u,v;\xi) &= \int_{\Oc} B_L(u,\nabla u;\cdot,\xi) \cdot \nabla v + \int_{\Oc} C_L(u,\nabla u;\cdot,\xi) v,\\
n_{\Oc}(u,v;\xi) &= \int_{\Oc} B_N(u,\nabla u;\cdot,\xi) \cdot \nabla v + \int_{\Oc} C_N(u,\nabla u;\cdot,\xi) v.
\end{align*}
We make the following assumptions.
\begin{assumption}[Properties of bilinear form $a_{\Oc}$]\label{assa}
We assume that the bilinear form $a_{\Oc}(\cdot,\cdot;\xi) \colon \xHone(\Oc) \times \xHone(\Oc) \to \xR$ is 
such that there exist constants $0 < \alpha_a \leq \beta_a < +\infty$ such that it holds almost surely
\begin{alignat}{2}
a_{\Oc}(v,v;\xi) &\geq \alpha_a \abs{v}^2_{\xHone(\Oc)} \quad &&\forallin{v}{\xHone(\Oc)},\label{coercivity_a}\\
\abs{a_{\Oc}(u,v;\xi)} &\leq \beta_a \norm{u}_{\xHone(\Oc)} \norm{v}_{\xHone(\Oc)} \quad &&\forallin{u,v}{\xHone(\Oc)},\label{continuity_a}
\end{alignat}
and we further assume that $\alpha_a$ and $\beta_a$ are independent of $\xi$ and $\Oc$.
\end{assumption}

\begin{assumption}[Properties of semi-linear form $n_{\Oc}$]\label{assn}
We assume that the semi-linear form $n_{\Oc}(\cdot,\cdot;\xi) \colon \xHone(\Oc) \times \xHone(\Oc) \to \xR$ is almost surely continuous with respect to the second variable and radially continuous with respect to the first variable, that means for all $u,v \in \xHone(\Oc)$, the map $t \mapsto n_{\Oc}(u+tv,v;\xi)$ is almost surely continuous. We also assume that $n_{\Oc}(\cdot,\cdot;\xi)$ is almost surely monotone in the first variable, that means
\begin{equation}\label{monotonicity_n}
n_{\Oc}(u,u-v;\xi) - n_{\Oc}(v,u-v;\xi) \geq 0 \quad \forallin{u,v}{\xHone(\Oc)}
\end{equation}
holds almost surely.
Finally, we assume that $n_{\Oc}(\cdot,\cdot;\xi)$ satisfies almost surely
\begin{equation}\label{zerocondition_n}
n_{\Oc}(0,v;\xi) = 0 \quad \forallin{v}{\xHone(\Oc)}.
\end{equation}
\end{assumption}

\begin{xmpl}
Concerning Example~\ref{ex:diffusionreactionpb}, assumption~\ref{assa} on $a_{\Oc}$ is satisfied if the diffusion coefficient $K$ is 
such that $0 < K_{\mathrm{inf}} \abs{\zeta}^2 \leq \scalprod{K(x,\xi) \zeta}{\zeta} \leq K_{\mathrm{sup}} \abs{\zeta}^2 < +\infty$ for all $\zeta \in \xR^d$ holds almost surely and almost everywhere on $\Oc$, where $K_{\mathrm{inf}}$ and $K_{\mathrm{sup}}$ are some strictly positive constants independent of $\xi$ and independent of the considered subdomain $\Oc \subset \Omega$. Also, assumption~\ref{assn} on $n_{\Oc}$ is satisfied if the reaction coefficient $R$ is such that $0 \leq R(x,\xi) \leq R_{\mathrm{sup}} < +\infty$ holds almost surely and almost everywhere on $\Oc$, where $R_{\mathrm{sup}}$ is a strictly positive constant independent of $\xi$ and independent of the considered subdomain $\Oc \subset \Omega$. That means assumptions~\ref{assa} and \ref{assn} on $a_{\Oc}$ and $n_{\Oc}$ are satisfied if the diffusion coefficient $K$ is almost surely uniformly bounded and elliptic and the reaction parameter $R$ is almost surely non negative and uniformly bounded.
\end{xmpl}

\subsubsection{Assumption on the geometry}

We suppose that the considered domains have sufficiently smooth boundary (\eg{} Lipschitz).
For a subset $\Ec \subset \partial \Oc$ with non zero measure, we denote by $\xHn{1/2}(\Ec)$ the space of traces on $\Ec$ of functions in $\xHone(\Oc)$. We recall that we have
\begin{equation}\label{generalizedpoincare}
\norm{v}_{\xHone(\Oc)} \leq C_{\Oc,\Ec} \( \abs{v}_{\xHone(\Oc)} + \norm{v}_{\xHn{1/2}(\Ec)} \),
\end{equation}
for all $v \in \xHone(\Oc)$, with a constant $C_{\Oc,\Ec}$ depending only on $\Oc$ and $\Ec$ (see \cite[Theorem 7.3.13]{Atk09}).
\begin{assumption}\label{assgeometry}
For any considered domains $\Oc$ and $\Ec \subset \partial \Oc$, we assume that the constant $C_{\Oc,\Ec}$ involved in \eqref{generalizedpoincare} is independent of $\xi$.
\end{assumption}
Assumption~\ref{assgeometry} is obviously satisfied if the domains $\Oc$ and $\Ec$ are independent of $\xi$. In the case of uncertain domains $\Oc(\xi)$ and $\Ec(\xi)$, assumption~\ref{assgeometry} implies some restrictions on the dependence of the geometry on the parameters $\xi$. Let us describe a typical situation where the uncertain domain is described through a parametrized mapping defined on a fixed domain.
Assume that there exist domains $\Oc_0$ and $\Ec_0 \subset \partial \Oc_0$ independent of $\xi$, and 
a parametrized diffeomorphism $\phi(\cdot;\xi) \colon \overline{\Oc_0} \to \overline{\Oc(\xi)}$ such that $\phi(\Oc_0;\xi) = \Oc(\xi)$ and $\phi(\Ec_0;\xi) = \Ec(\xi)$. Then it can be proved that assumption~\ref{assgeometry} is satisfied if the mapping $\phi(\cdot;\xi)$ satisfies almost surely
\begin{equation*}
\alpha_{\phi} \abs{\zeta} \leq \abs{\nabla \phi(x_0;\xi) \zeta} \leq \beta_{\phi} \abs{\zeta} \quad \forallin{\zeta}{\xR^d}, \quad \forallin{x_0}{\Oc_0},
\end{equation*}
with constants $\alpha_{\phi}$ and $\beta_{\phi}$ independent of $\xi$. That means assumption~\ref{assgeometry} on the geometry is satisfied if we consider a fixed (deterministic) domain $\Oc_0$ and a random geometrical transformation $\phi(\cdot;\xi)$ mapping the fixed domain $\overline{\Oc_0}$ to the parameter-dependent domain $\overline{\Oc(\xi)}$ such that the singular values of the jacobian matrix $\nabla \phi(\cdot;\xi)$ are almost surely uniformly bounded (from above and from below).

\section{Global-local formulation with patch}\label{sec:formulationpatch}

\subsection{Domain decomposition: introduction of a patch}

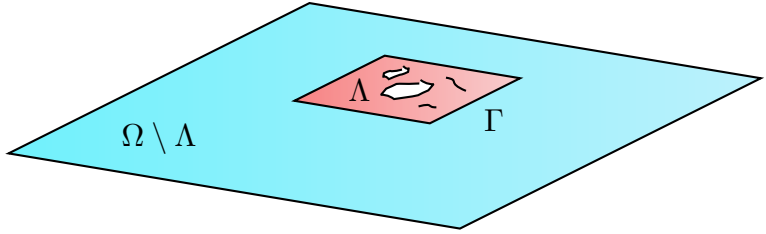
\begin{figure}[h!]
\centering
\tikzsetnextfilename{partition}
\input{partition}
\caption{Representation of interface $\Gamma$ between patch $\Lambda$ and complementary subdomain $\Omega \setminus \Lambda$}\label{partition}
\end{figure}

We introduce a subdomain $\Lambda \subset \Omega$, hereafter called a patch, such that $\Lambda_{\star} \subset \Lambda$, and such that $\Omega \setminus \Lambda$ is independent of $\xi$.
This yields the following partition of domain $\Omega(\xi)$:
\begin{equation*}
\Omega(\xi) = (\Omega \setminus \Lambda) \cup \Lambda(\xi).
\end{equation*}
The patch $\Lambda$ is chosen such that
\begin{equation}\label{localization_delta}
\dist(\Lambda_{\star},\Omega \setminus \Lambda) > \delta,
\end{equation}
that means uncertainties on operator, geometry and right-hand side affect a region in $\Lambda$ whose distance to $\Omega \setminus \Lambda$ is greater than $\delta$.
We assume that the patch $\Lambda$ has a sufficiently smooth boundary (\eg{} Lipschitz). We denote by
\begin{equation*}
\Gamma = \partial \Lambda \cap \partial(\Omega \setminus \Lambda)
\end{equation*}
the deterministic interface between the patch $\Lambda$ and the exterior subdomain $\Omega \setminus \Lambda$ (see Figure~\ref{partition}).

We denote by $U(\xi)$ and $w(\xi)$ the restrictions of $u(\xi)$ to subdomains $\Omega \setminus \Lambda$ and $\Lambda$, respectively. For $U(\xi) \in \xHone(\Omega \setminus \Lambda)$ and $w(\xi) \in \xHone(\Lambda)$, we denote by $U(\xi)_{\restrictto{\Gamma}}$ and $w(\xi)_{\restrictto{\Gamma}}$ in $\xHn{1/2}(\Gamma)$ the traces on $\Gamma$ of $U(\xi)$ and $w(\xi)$, respectively.
A weak continuity condition is enforced on interface $\Gamma$ by imposing 
\begin{equation}\label{weakcontinuitycondition}
b_{\Gamma}(\delta \lambda,U(\xi)_{\restrictto{\Gamma}}) = b_{\Gamma}(\delta \lambda,w(\xi)_{\restrictto{\Gamma}}) \quad \forallin{\delta \lambda}{\xHn{1/2}(\Gamma)^{\ast}},
\end{equation}
where $b_{\Gamma}$ denotes the duality pairing between $\xHn{1/2}(\Gamma)$ and its topological dual $\xHn{1/2}(\Gamma)^{\ast}$.
In the following, we denote by $\Mc = \xHn{1/2}(\Gamma)^{\ast}$ and $\Vert \cdot \Vert_\Mc = \Vert \cdot \Vert_{\xHn{1/2}(\Gamma)^{\ast}}$. We use the same notation $b_{\Gamma}$ for the bilinear form $b_{\Gamma} \colon \Mc \times \xHn{1/2}(\Gamma) \to \xR$ and its extension to $\Mc \times \xHone(\Omega \setminus \Lambda)$ (resp. $\Mc \times \xHone(\Lambda)$) defined using the trace operator from $\xHone(\Omega \setminus \Lambda)$ (resp. $\xHone(\Lambda)$) to $\xHn{1/2}(\Gamma)$.

\begin{rmrk}\label{multiplepatches}
The proposed multiscale approach can be naturally extended to the case where the sources of uncertainties and possible non-linearities are localized in several non-overlapping local subdomains of interest (or patches). The patch $\Lambda$ and the interface $\Gamma$ can then be respectively interpreted as the disjoint union of $Q$ patches $\set{\Lambda_q}_{q=1}^{Q}$ and $Q$ interfaces $\set{\Gamma_q}_{q=1}^{Q}$, 
where $\Gamma_q = \partial \Lambda_q \cap \partial(\Omega \setminus \Lambda)$ is the deterministic interface between the patch $\Lambda_q$ and the exterior subdomain $\Omega \setminus \Lambda$.
\end{rmrk}

\subsubsection{Weak formulation}

We introduce the Hilbert spaces
\begin{align*}
\Uc &= \setst{U \in \xHone(\Omega \setminus \Lambda)}{U=0 \; \text{on} \; \Gamma_D \cap \partial (\Omega \setminus \Lambda)}, \\
\Wc &= \setst{w \in \xHone(\Lambda)}{w=0 \; \text{on} \; \Gamma_D \cap \partial \Lambda},
\end{align*}
equipped with norms $\Vert \cdot \Vert_{\Uc} = \Vert \cdot \Vert_{\xHone(\Omega \setminus \Lambda)}$ and $\Vert \cdot \Vert_{\Wc} = \Vert \cdot \Vert_{\xHone(\Lambda)}$, respectively.
Also, we introduce the Hilbert space 
\begin{equation*}
\widehat{\Vc} = \setst{u \colon \Omega \to \xR}{u_{\restrictto{\Omega \setminus \Lambda}} \in \Uc \text{ and } u_{\restrictto{\Lambda}} \in \Wc}
\end{equation*}
equipped with the norm $\Vert \cdot \Vert_{\Vc}$ defined by 
\begin{equation*}
\norm{u}_{\Vc}^2 = \norm{u_{\restrictto{\Omega \setminus \Lambda}}}_{\xHone(\Omega \setminus \Lambda)}^2 + \norm{u_{\restrictto{\Lambda}}}_{\xHone(\Lambda)}^2,
\end{equation*}
and the closed linear subspace 
\begin{equation*}
\Vc = \setst{u \in \widehat{\Vc}}{b_{\Gamma}(\delta \lambda, u_{\restrictto{\Omega \setminus \Lambda}}) = b_{\Gamma}(\delta \lambda, u_{\restrictto{\Lambda}}) \; \text{for all } \delta \lambda \in \Mc},
\end{equation*}
which is a Hilbert space when equipped with norm $\Vert \cdot \Vert_{\Vc}$.

\begin{lmm}\label{lmm:Cv}
There exists a constant $C_{\Vc}$ such that
$\abs{v}_{\Vc} \leq \norm{v}_{\Vc} \leq C_{\Vc} \abs{v}_{\Vc} \quad \forallin{v}{\Vc},
$ with
\begin{equation*}
\abs{v}_{\Vc}^2 = \abs{v_{\restrictto{\Omega \setminus \Lambda}}}_{\xHone(\Omega \setminus \Lambda)}^2 + \abs{v_{\restrictto{\Lambda}}}_{\xHone(\Lambda)}^2.
\end{equation*}
Under assumption~\ref{assgeometry}, $C_{\Vc}$ is independent of $\xi$.
\end{lmm}
\begin{proof}
See Section~\ref{sec:lmm_Cv} in Appendix~\ref{sec:appendix}.
\end{proof}

In the following, for a given Hilbert space $H$ (possibly dependent on $\xi$) equipped with a norm $\Vert \cdot \Vert_{H}$, we denote by $H^{\Xi}$ the space 
$H^{\Xi} \coloneqq \set{v \colon \xi \in \Xi \mapsto v(\xi) \in H(\xi)}$, 
and we identify functions in $H^{\Xi}$ that are equal almost surely. 
We denote by 
$
\xLn{p}_{\mu}(\Xi;H) = \setst{v \in H^{\Xi}}{\xE(\norm{v(\xi)}_{H(\xi)}^p) < +\infty},
$
where $\xE$ denotes the mathematical expectation defined by $\xE(a(\xi)) = \int_{\Xi} a(\xi) \mu(\xd\xi)$.

We consider the following weak formulation of the problem: find $u \in \Vc^{\Xi}$ such that it holds almost surely
\begin{equation}\label{weakformulation}
d_{\Omega}(u(\xi),\delta u ; \xi) = \ell_{\Omega}(\delta u ; \xi)\quad \forallin{\delta u}{\Vc}.
\end{equation}

\begin{thrm}\label{thrm:wellposedness_weakformulation}
Under assumptions~\ref{assdata}, \ref{assa} and \ref{assn}, problem \eqref{weakformulation} is well-posed, that means for almost all $\xi \in \Xi$, it admits a unique solution $u(\xi) \in \Vc$ and the application that maps $\ell_{\Omega}(\cdot;\xi)$ to this solution $u(\xi)$ is Lipschitz continuous with Lipschitz constant ${C_{\Vc}^2}/{\alpha_a}$.
Moreover, under assumption~\ref{assgeometry}, the solution
$
u \in \xLn{p}_{\mu}(\Xi;\Vc),
$
with exponent $p$ defined in assumption~\ref{assdata}.
\end{thrm}
\begin{proof}
See Section~\ref{sec:thrm_wellposedness_weakformulation} in Appendix~\ref{sec:appendix}.
\end{proof}

\subsubsection{Reformulation using a Lagrange multiplier}

From \eqref{localizationuncertainties} and \eqref{localizationnonlinearities}, we have that
\begin{align*}
d_{\Omega}(u(\xi),\delta u;\xi) &= a_{\Omega \setminus \Lambda}(U(\xi),\delta U) + a_{\Lambda}(w(\xi),\delta w;\xi) + n_{\Lambda}(w(\xi),\delta w;\xi),\\
\ell_{\Omega}(\delta u;\xi) &= \ell_{\Omega \setminus \Lambda}(\delta U) + \ell_{\Lambda}(\delta w;\xi),
\end{align*}
for all $\delta u \colon \Omega \to \xR$ such that $\delta u_{\restrictto{\Omega \setminus \Lambda}} = \delta U$ and $\delta u_{\restrictto{\Lambda}} = \delta w$. 
A formulation equivalent to \eqref{weakformulation} can be written as follows:
find $(U,w,\lambda) \in \Uc^{\Xi} \times \Wc^{\Xi} \times \Mc^{\Xi}$ such that it satisfies almost surely
\begin{subequations}\label{globallocalformulation}
\begin{align}
&a_{\Omega \setminus \Lambda}(U(\xi),\delta U) + b_{\Gamma}(\lambda(\xi),\delta U) = \ell_{\Omega \setminus \Lambda}(\delta U),\label{complementarypart}\\
&a_{\Lambda}(w(\xi),\delta w;\xi) + n_{\Lambda}(w(\xi),\delta w;\xi) - b_{\Gamma}(\lambda(\xi),\delta w) = \ell_{\Lambda}(\delta w;\xi),\label{patchpart}\\
&b_{\Gamma}(\delta \lambda,U(\xi)) - b_{\Gamma}(\delta \lambda,w(\xi)) = 0,\label{interfacepart}
\end{align}
\end{subequations}
for all $(\delta U,\delta w,\delta \lambda) \in \Uc \times \Wc \times \Mc$, where $\lambda$ represents the Lagrange multiplier allowing to ensure the weak continuity condition \eqref{weakcontinuitycondition} of solution $u$ across interface $\Gamma$.

\begin{thrm}\label{thrm:existence_globallocalformulation}
Under assumptions~\ref{assdata}, \ref{assa} and \ref{assn}, problem \eqref{globallocalformulation} admits a unique solution $(U(\xi),w(\xi),\lambda(\xi)) \in \Uc \times \Wc \times \Mc$ for almost all $\xi \in \Xi$. Moreover, under assumption~\ref{assgeometry},
 $U \in \xLn{p}_{\mu}(\Xi;\Uc)$, $w \in \xLn{p}_{\mu}(\Xi;\Wc)$ and $\lambda \in \xLn{p}_{\mu}(\Xi;\Mc)$, with exponent $p$ defined in assumption~\ref{assdata}.
\end{thrm}
\begin{proof}
See Section~\ref{sec:thrm_existence_globallocalformulation} in Appendix~\ref{sec:appendix}.
\end{proof}

\subsection{Reformulation with extended domain: introduction of a fictitious patch}

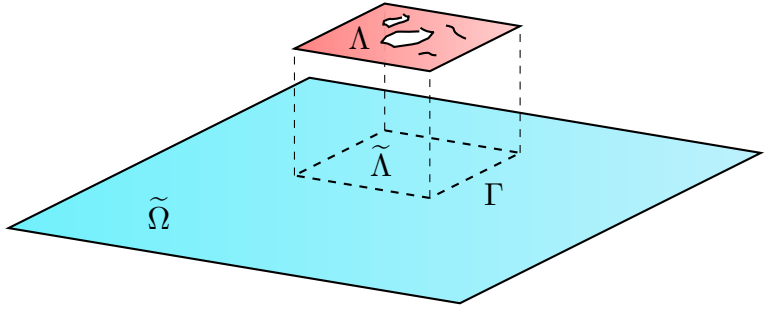
\begin{figure}[h!]
\centering
\tikzsetnextfilename{fictitious_patch}
\input{fictitious_patch}
\caption{Representation of fictitious domain $\widetilde{\Omega}$, fictitious patch $\widetilde{\Lambda}$, real patch $\Lambda$ and interface $\Gamma$}\label{fictitious_patch}
\end{figure}

Let us now introduce a deterministic fictitious patch $\widetilde{\Lambda} \supset \Lambda$ such that $\Gamma \subset \partial \widetilde{\Lambda}$ and define the corresponding deterministic fictitious domain $\widetilde{\Omega} \supset \Omega$ such that $\widetilde{\Omega} = (\Omega \setminus \Lambda) \cup \widetilde{\Lambda}$ and $\widetilde{\Omega} \setminus \widetilde{\Lambda} = \Omega \setminus \Lambda$ (see Figure~\ref{fictitious_patch}). Note that in the case where the patch $\Lambda$ does not contain any geometrical details (\ie{} no internal boundary such as holes, cracks, \etc), we simply have $\widetilde{\Lambda} = \Lambda$ and $\widetilde{\Omega} = \Omega$.

We now consider an extension of global solution $U$ from subdomain $\Omega \setminus \Lambda$ to fictitious domain $\widetilde{\Omega}$. We introduce the new Hilbert space $\widetilde{\Uc} = \setst{U \in \xHone(\widetilde{\Omega})}{U=0 \; \text{on} \; \Gamma_D \cap \partial \widetilde{\Omega}}$ equipped with the norm 
$\Vert \cdot \Vert_{\widetilde{\Uc}}= \Vert \cdot \Vert_{\xHone(\widetilde{\Omega})}$. We then define a new bilinear form $c_{\widetilde{\Omega}} \colon \widetilde{\Uc} \times \widetilde{\Uc} \to \xR$ as the following extension of $a_{\Omega \setminus \Lambda} \colon \Uc \times \Uc \to \xR$ to $\widetilde{\Uc} \times \widetilde{\Uc}$: for all $U,V \in \widetilde{\Uc}$,
\begin{equation}\label{extension}
c_{\widetilde{\Omega}}(U,V) = a_{\Omega \setminus \Lambda}(U,V) + c_{\widetilde{\Lambda}}(U,V),
\end{equation}
where, for a subdomain $\Oc \subset \widetilde{\Omega}$, $c_{\Oc}$ is a bilinear form defined by
\begin{equation*}
c_{\Oc}(U,V) = \int_{\Oc} \widetilde{B}_L(U,\nabla U;\cdot) \cdot \nabla V + \int_{\Oc} \widetilde{C}_L(U,\nabla U;\cdot) V,
\end{equation*}
where $\widetilde{B}_L(\cdot,\cdot;x,\xi) \colon \xR \times \xR^d \to \xR^d$ and $\widetilde{C}_L(\cdot,\cdot;x,\xi) \colon \xR \times \xR^d \to \xR$ are such that
\begin{equation*}
\widetilde{B}_L(\cdot,\cdot;x) = B_L(\cdot,\cdot;x) \quad \text{and} \quad \widetilde{C}_L(\cdot,\cdot;x) = C_L(\cdot,\cdot;x) \quad \text{for } x \in \Omega \setminus \Lambda.
\end{equation*}
We make the following assumption. Here $\Oc$ denotes a subset of $\widetilde{\Omega}$.
\begin{assumption}[Properties of bilinear form $c_{\Oc}$]\label{assc}
We assume that {the} bilinear form $c_{\Oc} \colon {\xHone(\Oc)} \times {\xHone(\Oc)} \to \xR$ is symmetric and 
such that there exist constants $0<\alpha_c \leq \beta_c < +\infty$ such that
\begin{alignat}{2}
c_{\Oc}(V,V) &\geq \alpha_c \abs{V}^2_{\xHone(\Oc)} \quad &&\forallin{V}{\xHone(\Oc)},\label{coercivity_c}\\
\abs{c_{\Oc}(U,V)} &\leq \beta_c \norm{U}_{\xHone(\Oc)} \norm{V}_{\xHone(\Oc)} \quad &&\forallin{U,V}{\xHone(\Oc)},\label{continuity_c}
\end{alignat}
and we further assume that $\alpha_c$ and $\beta_c$ are independent of $\xi$ and {$\Oc$}.
\end{assumption}

\begin{xmpl}
In Example~\ref{ex:diffusionreactionpb}, $\widetilde{B}_L$ and $\widetilde{C}_L$ can be respectively defined by
\begin{equation*}
\widetilde{B}_L(U,\nabla U;x) = \widetilde{K}(x) \nabla U \quad \text{and} \quad \widetilde{C}_L(U,\nabla U;x) = 0,
\end{equation*}
where $\widetilde{K}$ is a fictitious diffusion coefficient such that $\widetilde{K}(x) = K(x)$ for $x \in \Omega \setminus \Lambda$. Assumption~\ref{assc} on $c_{\widetilde{\Omega}}$ is satisfied if the fictitious diffusion coefficient $\widetilde{K}$ is uniformly bounded and elliptic on $\widetilde{\Omega}$, that means condition $0 < \widetilde{K}_{\mathrm{inf}} \abs{\zeta}^2 \leq \scalprod{\widetilde{K}(x) \zeta}{\zeta} \leq \widetilde{K}_{\mathrm{sup}} \abs{\zeta}^2 < +\infty$ for all $\zeta \in \xR^d$ holds almost everywhere on ${\widetilde{\Omega}}$, where $\widetilde{K}_{\mathrm{inf}}$ and $\widetilde{K}_{\mathrm{sup}}$ are some strictly positive constants.
\end{xmpl}

Afterwards, a reformulation of the global-local problem \eqref{globallocalformulation} reads: find $(U,w,\lambda) \in \widetilde{\Uc}^{\Xi} \times \Wc^{\Xi} \times \Mc^{\Xi}$ such that it satisfies almost surely
\begin{subequations}\label{globallocalreformulation}
\begin{align}
&c_{\widetilde{\Omega}}(U(\xi),\delta U) - c_{\widetilde{\Lambda}}(U(\xi),\delta U) + b_{\Gamma}(\lambda(\xi),\delta U) = \ell_{\Omega \setminus \Lambda}(\delta U),\label{globalpart}\\
&a_{\Lambda}(w(\xi),\delta w;\xi) + n_{\Lambda}(w(\xi),\delta w;\xi) - b_{\Gamma}(\lambda(\xi),\delta w) = \ell_{\Lambda}(\delta w;\xi),\label{localpart}\\
&b_{\Gamma}(\delta \lambda,U(\xi)) - b_{\Gamma}(\delta \lambda,w(\xi)) = 0,\label{couplingpart}
\end{align}
\end{subequations}
for all $(\delta U,\delta w,\delta \lambda) \in \widetilde{\Uc} \times \Wc \times \Mc$. Let us here mention that problem \eqref{globallocalreformulation} admits infinitely many solutions $(U,w,\lambda)$ that only differ by the value of global solution $U$ in fictitious patch $\widetilde{\Lambda}$. A particular solution can be uniquely defined by defining the value of $U$ in $\widetilde{\Lambda}$ as a particular extension of the value of $U$ on interface $\Gamma$. The global-local iterative algorithm presented in the next section will be proven to converge to a solution $(U,w,\lambda)$ in a subspace of $\widetilde{\Uc}^{\Xi} \times \Wc^{\Xi} \times \Mc^{\Xi}$ corresponding to a particular definition of the extension.

\section{Global-local iterative algorithm}\label{sec:algo}

We now introduce and analyze an iterative algorithm to solve problem \eqref{globallocalreformulation}.

\subsection{Description of the algorithm}\label{sec:descriptionalgo}
We initialize the algorithm with $U^0=w^0=\lambda^0=0$. Then, at iteration $k \geq 1$, $(U^k,w^k,\lambda^k) \in \widetilde{\Uc}^{\Xi} \times \Wc^{\Xi} \times \Mc^{\Xi}$ is defined by three steps (global step, relaxation step and local step), described below.
 
\subsubsection{Global step}\label{sec:globalstep}

We first define $\widehat{U}^{k} \in \widetilde{\Uc}^{\Xi}$ such that it satisfies almost surely
\begin{equation}\label{globalpb}
c_{\widetilde{\Omega}}(\widehat{U}^k(\xi),\delta U) = c_{\widetilde{\Lambda}}(U^{k-1}(\xi),\delta U) - b_{\Gamma}(\lambda^{k-1}(\xi),\delta U) + \ell_{\Omega \setminus \Lambda}(\delta U)
\end{equation}
for all $\delta U \in \widetilde{\Uc}$.
The computation of $\widehat{U}^{k} \in \widetilde{\Uc}^{\Xi}$ thus requires the solution of a linear problem defined on fictitious domain $\widetilde{\Omega}$ with a deterministic operator and an uncertain right-hand side (involving Lagrange multiplier $\lambda^{k-1}$ on interface $\Gamma$ and global iterate $U^{k-1}$ in fictitious patch $\widetilde{\Lambda}$ at previous iteration $k-1$).

\begin{rmrk}\label{determinsiticfictitiousoperators}
Although $\widetilde{B}_L$ and $\widetilde{C}_L$ could \apriori{} be chosen arbitrarily (uncertain or deterministic) on $\widetilde{\Lambda}$, a convenient choice consists in taking for $\widetilde{B}_L$ and $\widetilde{C}_L$ parameter-independent functions, \ie{} $\widetilde{B}_L(\cdot,\cdot;x,\xi) = \widetilde{B}_L(\cdot,\cdot;x)$ and $\widetilde{C}_L(\cdot,\cdot;x,\xi) = \widetilde{C}_L(\cdot,\cdot;x)$ for $x \in \widetilde{\Lambda}$, which allows to preserve a linear global problem with deterministic linear operator throughout iterations.
Also, a natural choice consists in taking for $\widetilde{B}_L$ and $\widetilde{C}_L$ over $\widetilde{\Lambda}$ the mean value of the corresponding linear functions $B_L$ and $C_L$ over $\Lambda$, \ie{} $\widetilde{B}_L(\cdot,\cdot;x) = \xE \(B_L(\cdot,\cdot;x,\xi)\)$ and $\widetilde{C}_L(\cdot,\cdot;x) = \xE \(C_L(\cdot,\cdot;x,\xi)\)$ for $x \in \widetilde{\Lambda}$.
Besides, choosing parameter-dependent functions $\widetilde{B}_L$ and $\widetilde{C}_L$ on $\widetilde{\Lambda}$ could allow to accelerate the convergence of the algorithm (see Remark~\ref{choiceforconvergence}). Another possible choice would consist in taking for $\widetilde{B}_L$ and $\widetilde{C}_L$ over $\widetilde{\Lambda}$ the tangent linear functions to the corresponding semi-linear functions $B=B_L+B_N$ and $C=C_L+C_N$ over $\Lambda$.
\end{rmrk}

\begin{rmrk}
Assume that $\partial \widetilde{\Lambda} = \Gamma \cup (\partial \widetilde{\Lambda} \cap \Gamma_D)$. By using Green's formula in the definition of $c_{\widetilde{\Lambda}}$, the global problem \eqref{globalpb} can be reformulated as
\begin{equation}\label{reformulatedglobalpb}
c_{\widetilde{\Omega}}(\widehat{U}^k(\xi),\delta U) = - b_{\Gamma}(\mu^{k-1}(\xi) + \lambda^{k-1}(\xi),\delta U) + \ell_{\Omega \setminus \Lambda}(\delta U) + \ell_{\widetilde{\Lambda}}(\delta U;\xi)
\end{equation}
for all $\delta U \in \widetilde{\Uc}$, where $\mu^{k-1}(\xi) \in \Mc$ is defined by the following expression (interpreted in a weak sense):
$\mu^{k-1}(\xi) = - \widetilde{B}_L(U^{k-1}{(\xi)},\nabla U^{k-1}{(\xi)};\cdot) \cdot n \quad \text{on } \Gamma,
$ with $n$ the unit normal to $\Gamma$ pointing outward $\widetilde{\Lambda}$, and where $\ell_{\widetilde{\Lambda}}(\cdot;\xi)$ is a linear form defined by
$
\ell_{\widetilde{\Lambda}}(V;\xi) = - \int_{\widetilde{\Lambda}} \nabla \cdot \widetilde{B}_L(U^{k-1}(\xi),\nabla U^{k-1}(\xi);\cdot) V + \int_{\widetilde{\Lambda}} \widetilde{C}_L(U^{k-1}(\xi),\nabla U^{k-1}(\xi);\cdot) V.
$
The quantity $\mu^{k-1} + \lambda^{k-1}$ is seen as a flux discontinuity on the interface $\Gamma$ between global and local models. The iterative algorithm can then be interpreted as a modified Newton method (with constant linear global operator) formulated on the flux equilibrium over interface $\Gamma$ (interpreted in a weak sense) \cite{Gen09,Gen11}.
\end{rmrk}

\subsubsection{Relaxation step}\label{sec:relaxationstep}

We then define $U^{k} \in \widetilde{\Uc}^{\Xi}$ by
\begin{equation}\label{relaxation}
U^{k}(\xi) 
= \rho_k \widehat{U}^k(\xi) + (1-\rho_k) U^{k-1}(\xi),
\end{equation}
where $\rho_k>0$ is a relaxation parameter (possibly depending on $\xi$) chosen sufficiently small to ensure convergence (see convergence analysis in Section~\ref{sec:convergence}).
Relaxation parameter $\rho_k$ may have a significant impact on the convergence and stability properties of the algorithm. Practical choices for $\rho_k$ will be discussed in Section~\ref{sec:relaxation}.

\subsubsection{Local step}\label{sec:localstep}
We finally define $(w^k,\lambda^k) \in \Wc^{\Xi} \times \Mc^{\Xi}$ such that it satisfies almost surely
\begin{subequations}\label{localpb}
\begin{align}
&a_{\Lambda}(w^k(\xi),\delta w;\xi) + n_{\Lambda}(w^k(\xi),\delta w;\xi) - b_{\Gamma}(\lambda^k(\xi),\delta w) = \ell_{\Lambda}(\delta w;\xi),\label{localpb1}\\
&b_{\Gamma}(\delta \lambda,w^k(\xi)) = b_{\Gamma}(\delta \lambda,U^k(\xi)),\label{localpb2}
\end{align}
\end{subequations}
for all $(\delta w,\delta \lambda) \in \Wc \times \Mc$.
The computation of $(w^k,\lambda^k) \in \Wc^{\Xi} \times \Mc^{\Xi}$ thus requires the solution of a non-linear problem defined on patch $\Lambda$ with uncertain operator and right-hand side (involving global iterate $U^k$ at current iteration $k$ as a boundary data). The Lagrange multiplier $\lambda^k$ allows to enforce the weak continuity conditions on interface $\Gamma$ between local iterate $w^k$ and global iterate $U^k$, which corresponds to non-homogeneous Dirichlet boundary conditions imposed on an external boundary $\Gamma$ of patch $\Lambda$ in the local computation. Recall that, contrary to the global step, the local step takes into account possible non-linearities and uncertainties in the operator, as well as possible uncertainties in the geometry of the domain.

\begin{rmrk}\label{singlefieldlocalpb}
The local problem \eqref{localpb} can be reformulated as a single-field problem by noting $w^k(\xi) = \widetilde{w}^k(\xi) + z^k(\xi)$, where $\widetilde{w}^k(\xi) \in \Wc$ is an extension of global iterate $U^k(\xi)$ from interface $\Gamma$ to patch $\Lambda$ such that 
$b_{\Gamma}(\delta \lambda,\widetilde{w}^k(\xi)) = b_{\Gamma}(\delta \lambda,U^k(\xi))
$
 for all $\delta \lambda \in \Mc$, and $z^k(\xi) \in \Wc_0 = \setst{z \in \Wc}{z=0 \; \text{on} \; \Gamma}$.
Local problem then consists in computing $z^k \in \Wc^{\Xi}_0$ such that it satisfies almost surely
$
a_{\Lambda}(\widetilde{w}^k(\xi)+z^k(\xi),\delta z;\xi) + n_{\Lambda}(\widetilde{w}^k(\xi)+z^k(\xi),\delta z;\xi) = \ell_{\Lambda}(\delta z;\xi) 
$ for all $\delta z \in \Wc_0$.
The Lagrange multiplier $\lambda^k \in \Mc^{\Xi}$ is then determined \aposteriori{} from \eqref{localpb1}.
\end{rmrk}

\begin{rmrk}\label{localpbgeomvariability}
When the patch $\Lambda$ contains geometrical variabilities, the local stochastic problem~\eqref{localpb} can be reformulated on a fixed (deterministic) domain by using either a random mapping technique \cite{Xiu06,Tar06%,Moh11
} or a fictitious domain method \cite{Canu07,Nou08a,Nou11b,Che13a}. Both techniques do not require any remeshing procedure and allow to handle complex geometries. The former approach consists in introducing a suitable random mapping between a random domain and a reference deterministic domain, thus transforming a (deterministic or stochastic) partial differential equation defined on a random domain into a stochastic partial differential equation (with uncertain operator and right-hand side depending on the mapping and its derivatives) defined on a deterministic domain. The latter approach consists in embedding a random domain into a deterministic fictitious domain and considering a prolongation of the solution on the fictitious domain. In the present context of multiscale problems with localized geometrical variabilities, different fictitious domain formulations (depending on the type of boundary conditions) have been proposed in \cite{Che13a}.
\end{rmrk}

\begin{rmrk}\label{multiplepatcheslocalpb}
Following Remark~\ref{multiplepatches}, in the case of $Q$ non-overlapping patches $\set{\Lambda_q}_{q=1}^{Q}$, the local step consists in solving $Q$ independent non-linear local problems defined on each of the patches $\Lambda_q$. The solution of such uncoupled problems can be performed independently on each patch $\Lambda_q$ in a fully parallel way.
\end{rmrk}

\subsection{Analysis of the algorithm}\label{sec:analysisalgo}
\subsubsection{Consistency}\label{sec:consistency}

Let $(U,w,\lambda) \in \Uc^{\Xi} \times \Wc^{\Xi} \times \Mc^{\Xi}$ denote the solution of the initial problem \eqref{globallocalformulation}. We now introduce the closed linear subspace $\widetilde{\Uc}_{\star}$ of $\widetilde{\Uc}$ defined by
\begin{equation*}
\widetilde{\Uc}_{\star} = \setst{V \in \widetilde{\Uc}}{c_{\widetilde{\Lambda}}(V,\delta U) = 0 \text{ for all } \delta U \in \xHone_0(\widetilde{\Lambda})},
\end{equation*}
where $\xHone_0(\widetilde{\Lambda})$ is considered as the subset of functions of $\widetilde{\Uc}$ which are zero on $\widetilde{\Omega} \setminus \widetilde{\Lambda}$. 
For any function $V \in \Uc$, there exists a unique extension $\widetilde{V} \in \widetilde{\Uc}_{\star}$ such that $\widetilde{V}=V$ on $\Omega \setminus \Lambda$ and the restriction of $\widetilde{V}$ to $\widetilde{\Lambda}$ is uniquely defined by the trace of $V$ on the interface $\Gamma$. 
Then, we also denote by $U(\xi) \in \widetilde{\Uc}_{\star}$ the unique extension to $\widetilde{\Omega}$ of the global solution $U(\xi) \in \Uc$.

\begin{lmm}\label{lmm:Uk_Ustar}
All global iterates $U^k(\xi)$ belong to the subspace $\widetilde{\Uc}_{\star}$.
\end{lmm}

\begin{proof}
Considering test functions $\delta U \in \xHone_0(\widetilde{\Lambda})$ in global problem~\eqref{globalpb}, we obtain that for all $k \geq 1$,
$c_{\widetilde{\Lambda}}(\widehat{U}^k(\xi),\delta U) = c_{\widetilde{\Lambda}}(U^{k-1}(\xi),\delta U)$ for all $\delta U \in \xHone_0(\widetilde{\Lambda}) $. Then, using \eqref{relaxation}, we have that $c_{\widetilde{\Lambda}}({U}^k(\xi),\delta U) = c_{\widetilde{\Lambda}}(U^{k-1}(\xi),\delta U)$ for all $\delta U \in \xHone_0(\widetilde{\Lambda})$. Since $U^0=0$, we obtain by induction that all global iterates $U^k(\xi)$ belong to $\widetilde{\Uc}_{\star}$.
\end{proof}

We then derive the following consistency result.

\begin{thrm}[Consistency]\label{thrm:consistency}
If the sequence $\set{(U^k(\xi),w^k(\xi),\lambda^k(\xi))}_{k \in \xN}$ strongly converges to an element $(\widetilde{U}(\xi),w(\xi),\lambda(\xi))$ in $\widetilde{\Uc} \times \Wc \times \Mc$, then $(\widetilde{U}(\xi)_{\restrictto{\Omega \setminus \Lambda}},w(\xi),\lambda(\xi)) \in \Uc \times \Wc \times \Mc$ is the unique solution $(U(\xi),w(\xi),\lambda(\xi))$ of problem \eqref{globallocalformulation}. Also, the limit $\widetilde{U}(\xi)$ is the unique extension of $U(\xi)$ to $\widetilde{\Uc}_{\star}$.
\end{thrm}

\begin{proof}
Taking the limit with $k$ in \eqref{globalpb}, \eqref{relaxation} and \eqref{localpb}, we obtain that $(\widetilde{U}(\xi),w(\xi),\lambda(\xi))$ satisfies problem \eqref{globallocalreformulation}, and therefore $(\widetilde{U}(\xi)_{\restrictto{\Omega \setminus \Lambda}},w(\xi),\lambda(\xi)) \in \Uc \times \Wc \times \Mc$ is the unique solution of problem \eqref{globallocalformulation}. Then, as all global iterates $U^k(\xi)$ belong to the closed linear subspace $\widetilde{\Uc}_{\star}$ of $\widetilde{\Uc}$ (see Lemma~\ref{lmm:Uk_Ustar}), the limit $\widetilde{U}(\xi)$ also belongs to $\widetilde{\Uc}_{\star}$.
\end{proof}

Note that problem \eqref{globallocalreformulation} is well-posed in $\widetilde{\Uc}_{\star} \times \Wc \times \Mc$ and admits $(U(\xi),w(\xi),\lambda(\xi)) \in \widetilde{\Uc}_{\star} \times \Wc \times \Mc$ as its unique solution. The algorithm can then be analyzed in the subspace $\widetilde{\Uc}_{\star}$ of $\widetilde{\Uc}$ and we have the following useful result which proves that $\Vert \cdot \Vert_{\widetilde{\Uc}}$ defines a norm equivalent to $\Vert \cdot \Vert_{\Uc}$ on $\widetilde{\Uc}_{\star}$.

\begin{lmm}\label{lmm:equiv_norm_Ustar}
The norms $\Vert \cdot \Vert_{\Uc}$ and $\Vert \cdot \Vert_{\widetilde{\Uc}}$ are equivalent on $\widetilde{\Uc}_{\star}$, with 
$\norm{V}_{\Uc} \leq \norm{V}_{\widetilde{\Uc}} \leq C_{\widetilde{\Uc}} \norm{V}_{\Uc}$ for all $V\in {\widetilde{\Uc}_{\star}}$, with a constant $C_{\widetilde{\Uc}}$ independent of $\xi$.
\end{lmm}

\begin{proof}
See Section~\ref{sec:lmm_equiv_norm_Ustar} in Appendix~\ref{sec:appendix}.
\end{proof}

From Theorem~\ref{thrm:existence_globallocalformulation} and Lemma~\ref{lmm:equiv_norm_Ustar}, we directly deduce the following property.
\begin{crllr}\label{U_Lp}
The extended global solution $U$ is in $\xLn{p}_{\mu}(\Xi;\widetilde{\Uc})$, with exponent $p$ defined in assumption~\ref{assdata}.
\end{crllr}

\subsubsection{Convergence}\label{sec:convergence}

We now prove the convergence of the sequence $\set{(U^k(\xi),w^k(\xi),\lambda^k(\xi))}_{k \in \xN}$ to the exact solution $(U(\xi),w(\xi),\lambda(\xi))$ in $\widetilde{\Uc}_{\star} \times \Wc \times \Mc$. The global problem~\eqref{globalpb} being linear, the solution $\widehat{U}^k \in \widetilde{\Uc}^{\Xi}$ can be written as
\begin{equation*}
\widehat{U}^k(\xi) = \overline{U} + \Upsilon(U^{k-1}(\xi)) + \Phi(\lambda^{k-1}(\xi)),
\end{equation*}
where $\Upsilon \colon \widetilde{\Uc} \to \widetilde{\Uc}$ and $\Phi \colon \Mc \to \widetilde{\Uc}$ are linear mappings. Mapping $\Upsilon$ is such that for $V \in \widetilde{\Uc}$, $\Upsilon(V) \in \widetilde{\Uc}$ is the unique solution of 
\begin{subequations}\label{globalpb_linearmap}
\begin{equation}\label{globalpbUpsilon}
c_{\widetilde{\Omega}}(\Upsilon(V),\delta U) = c_{\widetilde{\Lambda}}(V,\delta U) \quad \forallin{\delta U}{\widetilde{\Uc}}.
\end{equation}
Similarly, mapping $\Phi$ is such that for $\beta \in \Mc$, $\Phi(\beta) \in \widetilde{\Uc}$ is the unique solution of 
\begin{equation}\label{globalpbPhi}
c_{\widetilde{\Omega}}(\Phi(\beta),\delta U) = - b_{\Gamma}(\beta,\delta U) \quad \forallin{\delta U}{\widetilde{\Uc}}.
\end{equation}
\end{subequations}
Lastly, $\overline{U} \in \widetilde{\Uc}$ is the unique solution of 
\begin{equation*}
c_{\widetilde{\Omega}}(\overline{U},\delta U) = \ell_{\Omega \setminus \Lambda}(\delta U) \quad \forallin{\delta U}{\widetilde{\Uc}}.
\end{equation*}
The solution $(w^k,\lambda^k) \in \Wc^{\Xi} \times \Mc^{\Xi}$ of the local problem~\eqref{localpb} can be written as
\begin{equation*}
w^k(\xi) = \Theta(U^k(\xi);\xi) \quad \text{and} \quad \lambda^k(\xi) = \Psi(U^k(\xi);\xi),
\end{equation*}
where $\Theta(\cdot;\xi) \colon \widetilde{\Uc} \to \Wc$ and $\Psi(\cdot;\xi) \colon \widetilde{\Uc} \to \Mc$ are non-linear mappings. Mappings $\Theta$ and $\Psi$ are such that for $V \in \widetilde{\Uc}$, $(\Theta(V;\xi),\Psi(V;\xi)) \in \Wc \times \Mc$ is the solution of
\begin{subequations}\label{localpb_nonlinearmap}
\begin{align}
&a_{\Lambda}(\Theta(V;\xi),\delta w;\xi) + n_{\Lambda}(\Theta(V;\xi),\delta w;\xi) - b_{\Gamma}(\Psi(V;\xi),\delta w) = \ell_{\Lambda}(\delta w;\xi) & \forallin{\delta w}{\Wc},\label{localpbTheta}\\
&b_{\Gamma}(\delta \lambda,\Theta(V;\xi)) = b_{\Gamma}(\delta \lambda,V) & \forallin{\delta \lambda}{\Mc}.\label{localpbPsi}
\end{align}
\end{subequations}
Consequently, the algorithm generates a sequence $\set{(U^k,w^k,\lambda^k)}_{k \in \xN}$ by applying the following iterative scheme:
\begin{subequations}\label{iterativescheme}
\begin{align}
U^{k}(\xi) &= \rho_k \(\overline{U} + \Upsilon(U^{k-1}(\xi)) + \Phi(\lambda^{k-1}(\xi))\) + (1-\rho_k) U^{k-1}(\xi),\\
w^k(\xi) &= \Theta(U^k(\xi);\xi),\\
\lambda^k(\xi) & = \Psi(U^k(\xi);\xi).
\end{align}
\end{subequations}

\begin{lmm}\label{lmm:continuity_linearmap}
The linear mappings $\Upsilon \colon \widetilde{\Uc} \to \widetilde{\Uc}$ and $\Phi \colon \Mc \to \widetilde{\Uc}$ defined in \eqref{globalpb_linearmap} are continuous, with respective continuity constants $\beta_{\Upsilon}$ and $\beta_{\Phi}$ independent of $\xi$.
\end{lmm}

\begin{proof}
See Section~\ref{sec:lmm_continuity_linearmap} in Appendix~\ref{sec:appendix}.
\end{proof}

\begin{lmm}\label{lmm:continuity_nonlinearmap}
The non-linear mappings $\Theta(\cdot;\xi) \colon \widetilde{\Uc} \to \Wc$ and $\Psi(\cdot;\xi) \colon \widetilde{\Uc} \to \Mc$ defined in \eqref{localpb_nonlinearmap} are Lipschitz continuous, with respective Lipschitz constants $\beta_{\Theta}$ and $\beta_{\Psi}$ independent of $\xi$.
\end{lmm}

\begin{proof}
See Section~\ref{sec:lmm_continuity_nonlinearmap} in Appendix~\ref{sec:appendix}.
\end{proof}

Let us now define the errors at a given iteration of the algorithm. At the global level, the error at iteration $k$ is
\begin{align*}
\widehat{U}^{k}(\xi)-U(\xi) &= \Upsilon(U^{k-1}(\xi)) - \Upsilon(U(\xi)) + \Phi(\lambda^{k-1}(\xi)) - \Phi(\lambda(\xi)) \\
&= \Upsilon(U^{k-1}(\xi) - U(\xi)) + \Phi(\Psi(U^{k-1}(\xi);\xi) - \Psi(U(\xi);\xi)), \\
&= U^{k-1}(\xi) - U(\xi) - (A(U^{k-1}(\xi);\xi) -A(U(\xi);\xi)),
\end{align*}
and 
\begin{align*}
U^{k}(\xi)-U(\xi) &= \rho_k (\widehat{U}^{k}(\xi)-U(\xi)) + (1-\rho_k) (U^{k-1}(\xi)-U(\xi)) \\
&= U^{k-1}(\xi)-U(\xi) - \rho_k (A(U^{k-1}(\xi);\xi) -A(U(\xi);\xi)),
\end{align*}
where $A(\cdot;\xi) \colon \widetilde{\Uc} \to \widetilde{\Uc}$ is the non-linear mapping defined by
\begin{equation}\label{mapA}
A(V;\xi) = V - \Upsilon(V) - \Phi(\Psi(V;\xi)).
\end{equation}
From the definitions of $\Upsilon(V) \in \widetilde{\Uc}$ and $\Phi(\Psi(V;\xi)) \in \widetilde{\Uc}$ and from \eqref{extension}, we deduce that mapping $A$ is such that for $V \in \widetilde{\Uc}$, $A(V;\xi) \in \widetilde{\Uc}$ is the solution of
\begin{equation}\label{globalpbA}
c_{\widetilde{\Omega}}(A(V;\xi),\delta U) = a_{\Omega \setminus \Lambda}(V,\delta U) + b_{\Gamma}(\Psi(V;\xi),\delta U) \quad \forallin{\delta U}{\widetilde{\Uc}}.
\end{equation}
Note that the non-linear nature of map $A$ is inherited from that of map $\Psi$. At the local level, the error at iteration $k$ writes 
\begin{subequations}\label{local_error}
\begin{align}
w^{k}(\xi)-w(\xi) &= \Theta(U^{k}(\xi);\xi) - \Theta(U(\xi);\xi), \\
\lambda^{k}(\xi)-\lambda(\xi) &= \Psi(U^{k}(\xi);\xi) - \Psi(U(\xi);\xi).
\end{align}
\end{subequations}
Given that non-linear map $\Theta(\cdot;\xi)$ (resp. $\Psi(\cdot;\xi)$) is Lipschitz continuous, the almost sure convergence of the local sequence $\set{w^k(\xi)}_{k \in \xN}$ (resp. $\set{\lambda^k(\xi)}_{k \in \xN}$) to $w(\xi)$ (resp. $\lambda(\xi)$) in $\Wc$ (resp. $\Mc$) can be directly obtained from that of the global sequence $\set{U^k(\xi)}_{k \in \xN}$ to $U(\xi)$ in $\widetilde{\Uc}$. Recalling that the exact global solution $U(\xi)$ as well as all global iterates $U^k(\xi)$ belong to subspace $\widetilde{\Uc}_{\star} \subset \widetilde{\Uc}$, one can restrict the convergence analysis to that of the sequence $\set{U^k(\xi)}_{k \in \xN}$ to $U(\xi)$ in the subspace $\widetilde{\Uc}_{\star}$. Let $C_{\widetilde{\Omega}} \colon \widetilde{\Uc} \to \widetilde{\Uc}$ be the linear map defined for all $U,V \in \widetilde{\Uc}$ by $\scalproda{C_{\widetilde{\Omega}}(U)}{V}_{\widetilde{\Uc}} = c_{\widetilde{\Omega}}(U,V)$.

\begin{lmm}\label{lmm:propD}
The non-linear mapping $D(\cdot;\xi) \colon \widetilde{\Uc}_{\star} \to \widetilde{\Uc}$ defined by $D(\cdot;\xi) = C_{\widetilde{\Omega}}(A(\cdot;\xi))$ with $A(\cdot;\xi)$ defined in \eqref{mapA}, is Lipschitz continuous and strongly monotone, with Lipschitz constant $\beta_D$ and strong monotonicity constant $\alpha_D$ both independent of $\xi$.
\end{lmm}
\begin{proof}
See Section~\ref{sec:lmm_propD} in Appendix~\ref{sec:appendix}.
\end{proof}
One iteration of the algorithm can be written as
$$U^{k}(\xi) = B_{\rho_k}(U^{k-1}(\xi);\xi), $$ 
where $B_{\rho_k}(\cdot;\xi) \colon \widetilde{\Uc} \to \widetilde{\Uc}$ is the non-linear iteration map defined by $B_{\rho_k}(V;\xi) = \rho_k \overline{U} + V - \rho_k A(V;\xi)$.
We finally derive the following convergence result.
\begin{thrm}[Convergence]\label{thrm:convergence}
Assume that the sequence of relaxation parameters $\set{\rho_k}_{k \in \xN}$ is such that 
\begin{equation}\label{convergencecondition}
0<\rho_{\mathrm{inf}} \leq \rho_k \leq \rho_{\mathrm{sup}}<+\infty,
\end{equation}
for some strictly positive constants $\rho_{\mathrm{inf}}$ and $\rho_{\mathrm{sup}}$ independent of $\xi$ and $k$. Then, for $\rho_{\mathrm{sup}}$ sufficiently small, the sequence $\set{(U^k(\xi),w^k(\xi),\lambda^k(\xi))}_{k \in \xN}$ converges almost surely to the unique solution $(U(\xi),w(\xi),\lambda(\xi))$ of problem \eqref{globallocalreformulation} in $\widetilde{\Uc}_{\star} \times \Wc \times \Mc$. Also, the sequence $\set{U^k}_{k \in \xN}$ (resp. $\set{w^k}_{k \in \xN}$ and $\set{\lambda^k}_{k \in \xN}$) converges to $U$ (resp. $w$ and $\lambda$) in $\xLn{p}_{\mu}(\Xi;\widetilde{\Uc})$ (resp. $\xLn{p}_{\mu}(\Xi;\Wc)$ and $\xLn{p}_{\mu}(\Xi;\Mc)$).
\end{thrm}

\begin{proof}
See Section~\ref{sec:thrm_convergence} in Appendix~\ref{sec:appendix}.
\end{proof}

\begin{rmrk}\label{choiceforconvergence}
As long as \eqref{convergencecondition} is satisfied, the algorithm converges to the exact solution whatever the choice of relaxation parameter $\rho_k$ (bounded from above by a sufficiently small $\rho_{\mathrm{sup}}$ and from below by any value $\rho_{\mathrm{inf}}>0$ to ensure the convergence) and fictitious operators $\widetilde{B}_L$ and $\widetilde{C}_L$ (satisfying assumption~\ref{assc} on $c_{\widetilde{\Omega}}$). Nevertheless, these choices may have a significant influence on the convergence properties of the algorithm. Note that $\widetilde{B}_L$ and $\widetilde{C}_L$ play the role of preconditioners for the iterative algorithm.
\end{rmrk}

\begin{rmrk}
If $\Lambda = \widetilde{\Lambda}$ and if $B=B_L=\widetilde{B}_L$ and $C=C_L=\widetilde{C}_L$, that means if $n_{\Lambda} = 0$ and $c_{\widetilde \Lambda} = a_\Lambda$ (\ie{} in the linear case with the same geometry and bilinear form over $\Lambda$ and $\widetilde{\Lambda}$), then $A(\cdot;\xi)$ is such that $A(U;\xi)-A(V;\xi) = U-V$ for all $U,V\in \widetilde{\Uc}_{\star}$ and $B_{\rho_k}(\cdot;\xi)$ is such that $B_{\rho_k}(U;\xi)-B_{\rho_k}(V;\xi) = (1-\rho_k) (U-V)$ for all $U,V\in \widetilde{\Uc}_{\star}$, so that the convergence of the algorithm is achieved if the condition $0 < \rho_{\mathrm{inf}} \leq \rho_k \leq \rho_{\mathrm{sup}} < 2$ is fulfilled, with a convergence in two iterations for a fixed relaxation parameter $\rho_k =1$.
\end{rmrk}

\subsubsection{Robustness with respect to approximations}\label{sec:robustness}

Let us now consider that some approximations are introduced in the different steps of the global-local iterative algorithm.

Such approximations arise when using finite element approximations of spatial functions (see Section~\ref{sec:FEapprox}), approximations of parameter-dependent functions (see Sections~\ref{sec:approx-stochastic} and \ref{sec:leastsquares}), or when resorting to the use of unconverged iterative solvers for the approximate solution of either global or local problems with a certain prescribed accuracy. For example, the solution of non-linear local problems may be performed by means of classical non-linear solvers, such as Newton-type iterative solvers, leading to approximate local solutions. 

We now analyze the sensitivity of the global-local iterative algorithm with respect to these approximations at the different steps of the algorithm. Due to these approximations, the algorithm, initially defined by the unperturbed iterative scheme \eqref{iterativescheme}, generates a sequence $\set{(U_{\varepsilon}^k,w_{\varepsilon}^k,\lambda_{\varepsilon}^k)}_{k \in \xN}$ defined by the following perturbed iterative scheme:
\begin{subequations}\label{perturbediterativescheme}
\begin{align}
U_{\varepsilon}^{k}(\xi) &= \rho_k \(\overline{U}_{\varepsilon} + \Upsilon_{\varepsilon}(U_{\varepsilon}^{k-1}(\xi)) + \Phi_{\varepsilon}(\lambda_{\varepsilon}^{k-1}(\xi))\) + (1-\rho_k) U_{\varepsilon}^{k-1}(\xi),\\
w_{\varepsilon}^k(\xi) & = \Theta_{\varepsilon}(U_{\varepsilon}^k(\xi);\xi), \quad \lambda_{\varepsilon}^k(\xi) = \Psi_{\varepsilon}(U_{\varepsilon}^k(\xi);\xi),
\end{align}
\end{subequations}
where $\Upsilon_{\varepsilon}$ and $\Phi_{\varepsilon}$ (resp. $\Theta_{\varepsilon}$ and $\Psi_{\varepsilon}$) are approximations of linear maps $\Upsilon$ and $\Phi$ (resp. non-linear maps $\Theta$ and $\Psi$). Similarly, $\overline{U}_{\varepsilon}$ represents an approximation of $\overline{U}$.
At the global level, the unperturbed global iterate $U^{k} \in \widetilde{\Uc}_{\star}^{\Xi}$ at iteration $k$ satisfies
$
U^k(\xi) = B_{\rho_k}(U_{\varepsilon}^{k-1}(\xi);\xi).
$ 
The approximate (or perturbed) global iterate $U_{\varepsilon}^k$ at iteration $k$ is assumed to belong to $\widetilde{\Uc}_{\star}^{\Xi}$ and satisfies
$ 
U_{\varepsilon}^k(\xi) = B_{\rho_k}^{\varepsilon}(U_{\varepsilon}^{k-1}(\xi);\xi),
$ 
where $B_{\rho_k}^{\varepsilon}(\cdot;\xi)$ denotes an approximation of iteration map $B_{\rho_k}(\cdot;\xi)$ defined by $B_{\rho_k}^{\varepsilon}(V;\xi)= \rho_k \overline{U}_{\varepsilon} + V - \rho_k A_{\varepsilon}(V;\xi)$, with $A_{\varepsilon}(V;\xi) = V - \Upsilon_{\varepsilon}(V;\xi) - \Phi_{\varepsilon}(\Psi_{\varepsilon}(V;\xi);\xi)$.

We assume that approximations are controlled in a $\xLn{p}_{\mu}$-norm (typically $p=2$ or $p=\infty$), that means the approximation error at iteration $k$, 
$ 
U_{\varepsilon}^k(\xi) - U^k(\xi) = B_{\rho_k}^{\varepsilon}(U_{\varepsilon}^{k-1}(\xi);\xi) - B_{\rho_k}(U_{\varepsilon}^{k-1}(\xi);\xi)
$, should satisfy 
\begin{equation*}
\norm{U_{\varepsilon}^k - U^k}_{\xLn{p}_{\mu}(\Xi;\widetilde{\Uc})} \leq \varepsilon \norm{U}_{\xLn{p}_{\mu}(\Xi;\widetilde{\Uc})} + \varepsilon^{\ast} \norm{U_{\varepsilon}^{k-1} - U}_{\xLn{p}_{\mu}(\Xi;\widetilde{\Uc})},
\end{equation*}
where $\varepsilon$ conveys an absolute error with respect to the solution norm $\norm{U}_{\xLn{p}_{\mu}(\Xi;\widetilde{\Uc})}$, while $\varepsilon^{\ast}$ conveys an approximation error controlled relatively to the solution error in $\xLn{p}_{\mu}$-norm $\norm{U_{\varepsilon}^{k-1} - U}_{\xLn{p}_{\mu}(\Xi;\widetilde{\Uc})}$ at previous iteration $k-1$. In practice, $\varepsilon$ (resp. $\varepsilon^{\ast}$) is related to the user-specified tolerance (prescribed to the iterative solver) for the precision of the residual norms associated with global problem \eqref{globalpb} and local problem \eqref{localpb} formulated on the current iterates $U^k$ and $(w^k,\lambda^k)$ (resp. on the current increments $\delta U^k = U^k - U_{\varepsilon}^{k-1}$ and $(\delta w^k, \delta \lambda^k) = (w^k - w_{\varepsilon}^{k-1},\lambda^k - \lambda_{\varepsilon}^{k-1}$)) (see \cite[Section~3.5]{Che13a} for further details). We then provide for a robustness result relative to both types of errors.

\begin{thrm}[Robustness]\label{thrm:robustness}
Suppose that the set of iteration maps $\set{B_{\rho_k}(\cdot;\xi)}_{k \geq 1}$ is uniformly contractive on $\widetilde{\Uc}_{\star}$, that means 
\begin{equation*}
\norm{B_{\rho_k}(V;\xi) - B_{\rho_k}(W;\xi)}_{\widetilde{\Uc}} \leq \rho_B \norm{V-W}_{\widetilde{\Uc}},
\end{equation*}
for all $V, W \in \widetilde{\Uc}_{\star}$, with a contractivity constant $\rho_B < 1$ independent of $\xi$.
Further assume that the set of perturbed iteration maps $\set{B_{\rho_k}^{\varepsilon}(\cdot;\xi)}_{k \geq 1}$ is such that for all $V$ in a $\delta$-neighborhood $\Vc_{\delta}$ of the exact global solution $U$, defined by
$\Vc_{\delta} = \setst{V \in \xLn{p}_{\mu}(\Xi;\widetilde{\Uc}_{\star})}{\norm{V-U}_{\xLn{p}_{\mu}(\Xi;\widetilde{\Uc})} < \delta \norm{U}_{\xLn{p}_{\mu}(\Xi;\widetilde{\Uc})}}$, we have 
\begin{equation*}
\norm{B_{\rho_k}^{\varepsilon}(V) - B_{\rho_k}(V)}_{\xLn{p}_{\mu}(\Xi;\widetilde{\Uc})} \leq \varepsilon \norm{U}_{\xLn{p}_{\mu}(\Xi;\widetilde{\Uc})} + \varepsilon^{\ast} \norm{V-U}_{\xLn{p}_{\mu}(\Xi;\widetilde{\Uc})},
\end{equation*}
for some given tolerances $0 \leq \varepsilon^{\ast} < 1 - \rho_B$ and $0 \leq \varepsilon \leq \delta (1 - \rho_B - \varepsilon^{\ast})$.
Then, if the initial iterate $U_{\varepsilon}^0=0 \in \Vc_{\delta}$, the approximate sequence $\set{U_{\varepsilon}^k}_{k \in \xN}$ is such that
\begin{equation}\label{robustness_Lp}
\limsup_{k \to +\infty} \norm{U_{\varepsilon}^k - U}_{\xLn{p}_{\mu}(\Xi;\widetilde{\Uc})} \leq \gamma(\varepsilon,\varepsilon^{\ast}) \norm{U}_{\xLn{p}_{\mu}(\Xi;\widetilde{\Uc})},
\end{equation}
with $\gamma(\varepsilon,\varepsilon^{\ast}) = \frac{\varepsilon}{1 - (\rho_B + \varepsilon^{\ast})} \to 0$ as $\varepsilon \to 0$, and tends to a neighborhood of $U$ in $\xLn{p}_{\mu}(\Xi;\widetilde{\Uc})$ whose size is proportional to $\gamma(\varepsilon,\varepsilon^{\ast})$.
\end{thrm}
\begin{proof}
See Section~\ref{sec:thrm_robustness} in Appendix~\ref{sec:appendix}.
\end{proof}

Finally, the approximate sequence $\set{U_{\varepsilon}^k}_{k \in \xN}$ (resp. $\set{w_{\varepsilon}^k}_{k \in \xN}$ and $\set{\lambda_{\varepsilon}^k}_{k \in \xN}$) generated by the perturbed iterative scheme \eqref{perturbediterativescheme} converges in $\xLn{p}_{\mu}(\Xi;\widetilde{\Uc})$ (resp. $\xLn{p}_{\mu}(\Xi;\Wc)$ and $\xLn{p}_{\mu}(\Xi;\Mc)$) to a neighborhood of the exact solution $U$ (resp. $w$ and $\lambda$). Therefore, the proposed global-local iterative algorithm exhibits robustness properties with respect to possible perturbations such as numerical approximations, either controlled with relative precision $\varepsilon^*$ or absolute precision $\varepsilon$, which is an essential feature from a numerical point of view.

\begin{rmrk}
The convergence rate of the perturbed algorithm \eqref{perturbediterativescheme}, which is $\rho_B + \varepsilon^*$ (see the proof of 
Theorem~\ref{thrm:robustness} in Section~\ref{sec:thrm_robustness} in Appendix~\ref{sec:appendix}), depends on the perturbations with controlled relative precision $\varepsilon^*$ but not on perturbations with controlled absolute precision $\varepsilon$.  
\end{rmrk}

\begin{rmrk}
Under the more restrictive assumption that the set of perturbed iteration maps $\set{B_{\rho_k}^{\varepsilon}(\cdot;\xi)}_{k \geq 1}$ is such that for all $V$ in a $\delta$-neighborhood $\Vc_{\delta}(\xi)$ of the exact global solution $U(\xi)$, defined by
$\Vc_{\delta}(\xi) = \setst{V \in \widetilde{\Uc}_{\star}}{\norm{V-U(\xi)}_{\widetilde{\Uc}} < \delta \norm{U(\xi)}_{\widetilde{\Uc}}}$, we have almost surely
\begin{equation*}
\norm{B_{\rho_k}^{\varepsilon}(V;\xi) - B_{\rho_k}(V;\xi)}_{\widetilde{\Uc}} \leq \varepsilon \norm{U(\xi)}_{\widetilde{\Uc}} + \varepsilon^{\ast} \norm{V-U(\xi)}_{\widetilde{\Uc}},
\end{equation*}
with $0 \leq \varepsilon^{\ast} < 1 - \rho_B$ and $0 \leq \varepsilon \leq \delta (1 - \rho_B - \varepsilon^{\ast})$, then we can prove that the approximate sequence $\set{(U_{\varepsilon}^k(\xi),w_{\varepsilon}^k(\xi),\lambda_{\varepsilon}^k(\xi)}_{k \in \xN}$ generated by the perturbed iterative scheme \eqref{perturbediterativescheme} converges almost surely to a neighborhood of the exact solution $(U(\xi),w(\xi),\lambda(\xi))$. 
\end{rmrk}

\section{Computational aspects}\label{sec:computation}

In this section, we address computational aspects related to the proposed global-local iterative algorithm.

\subsection{Finite element approximations at spatial level}\label{sec:FEapprox}

At the spatial level, we employ a standard Galerkin finite element method by introducing finite-dimensional approximation spaces $\widetilde{\Uc}_H \subset \widetilde{\Uc}$ and $\Wc_h \subset \Wc$ with dimensions $n_U$ and $n_w$, respectively. We denote by $\Tc_H(\widetilde{\Omega})$ (resp. $\Tc_h(\Lambda)$) the finite element mesh of fictitious domain $\widetilde{\Omega}$ (resp. patch $\Lambda$) composed of elements of maximum size $H$ (resp. $h$). For the sake of simplicity, we further make the following assumptions:
\begin{itemize}
	\item domain $\widetilde{\Omega}$ and patch $\Lambda$ are exactly covered by global mesh $\Tc_H(\widetilde{\Omega})$ and local mesh $\Tc_h(\Lambda)$, respectively;
	\item global mesh $\Tc_H(\widetilde{\Omega})$ is partitioned into two submeshes $\Tc_H(\Omega \setminus \Lambda)$ and $\Tc_H(\widetilde{\Lambda})$ (associated with subdomain $\Omega \setminus \Lambda$ and fictitious patch $\widetilde{\Lambda}$, respectively) such that $\Tc_H(\widetilde{\Omega}) = \Tc_H(\Omega \setminus \Lambda) \cup \Tc_H(\widetilde{\Lambda})$, that means interface $\Gamma$ coincides with the intersection of boundaries of both submeshes $\Tc_H(\Omega \setminus \Lambda)$ and $\Tc_H(\widetilde{\Lambda})$ and therefore does not cut any element of global mesh $\Tc_H(\widetilde{\Omega})$.
\end{itemize}
Both meshes $\Tc_H(\widetilde{\Omega})$ and $\Tc_h(\Lambda)$ are \apriori{} not compatible at interface $\Gamma$, that means they may not match on interface $\Gamma$. Note that interface $\Gamma$ is part of the boundary of meshes $\Tc_h(\Lambda)$, $\Tc_H(\Omega \setminus \Lambda)$ and $\Tc_H(\widetilde{\Lambda})$. We now introduce an approximation space $\Mc_h \subset \Mc$ with dimension $n_{\lambda}$. In the general case of non-matching meshes and geometries, where both meshes $\Tc_H(\widetilde{\Omega})$ and $\Tc_h(\Lambda)$ do not match and even do not share a common interface $\Gamma$, one should pay attention to the construction of a suitable approximation space $\Mc_h$ of Lagrange multipliers% satisfying discrete inf-sup conditions for bilinear form $b_{\Gamma}$ \cite{Bab73,Bre74,Bre91} in approximation space $\Wc_h \times \Mc_h$
. The interested reader can refer to \cite{Bel99,Woh01,Kim01} for further information on the construction of appropriate Lagrange multiplier spaces using mortar (non-conforming) finite elements. In our particular case where both meshes $\Tc_H(\Omega \setminus \Lambda)$ and $\Tc_h(\Lambda)$ share a common interface $\Gamma$, a natural choice consists in taking for $\Mc_h$ a finite-dimensional subspace of trace space $\xHn{1/2}(\Gamma)$, so that $\Mc_h \subset \xHn{1/2}(\Gamma) \subset \xLtwo(\Gamma) \subset \xHn{1/2}(\Gamma)^{\ast} = \Mc$. If interface $\Gamma$ does not present any boundary, a convenient choice consists in taking the trace of $\Wc_h$ on interface $\Gamma$ for the practical construction of $\Mc_h$%, which leads to a continuous mortar approximation space
. Otherwise, if interface $\Gamma$ has a boundary, an alternative choice consists in taking a subspace of the trace of $\Wc_h$ on interface $\Gamma$ 
(see \cite{Woh01,Kim01}). The interested reader can refer to \cite{Bel99,Woh01,Hag12} for details about the properties of trace spaces and mortar projection operators.

For a given Hilbert space $H$ (possibly dependent on $\xi$), we denote by $H_n$ a finite element approximation subspace of $H$ spanned by basis functions $\set{\varphi_i}_{i \in \Ic}$ and with dimension $n = \#\Ic$. A function $v \in H_n$ can then be identified with a vector $\mathbf{v} = (v_i)_{i \in \Ic} \in \xR^n$ such that $v = \sum_{i \in \Ic} v_i \varphi_i$. Similarly, an element $v \in (H_n)^{\Xi}$ can be identified with a random vector $\mathbf{v} = (v_i)_{i \in \Ic} \in (\xR^n)^{\Xi}$ such that $v(\xi) = \sum_{i \in \Ic} v_i(\xi) \varphi_i$. For the bilinear forms $c_{\Oc}$ and $a_{\Oc}$, semi-linear form $n_{\Oc}$ and linear form $\ell_{\Oc}$, we introduce the finite element matrices $\mathbf{C}_{\Oc}$ and $\mathbf{A}_{\Oc}(\xi)$, discretized random non-linear map $\mathbf{N}_{\Oc}(\cdot;\xi)$ and finite element random vector $\mathbf{l}_{\Oc}(\xi)$, respectively defined for a subdomain $\Oc \subset \widetilde{\Omega}$ by
\begin{alignat*}{2}
c_{\Oc}(u,v) &= \mathbf{v}^T \mathbf{C}_{\Oc} \mathbf{u}, \quad a_{\Oc}(u,v;\xi) &&= \mathbf{v}^T \mathbf{A}_{\Oc}(\xi) \mathbf{u},\\
n_{\Oc}(u,v;\xi) &= \mathbf{v}^T \mathbf{N}_{\Oc}(\mathbf{u};\xi), \quad \ell_{\Oc}(v;\xi) &&= \mathbf{v}^T \mathbf{l}_{\Oc}(\xi).
\end{alignat*}
As the coupling bilinear form $b_{\Gamma}$ is defined on the two distinct subspaces $\Mc_h \times \widetilde{\Uc}_H$ and $\Mc_h \times \Wc_h$, we introduce two finite element matrices $\widetilde{\mathbf{B}}_{\Gamma}$ and $\mathbf{B}_{\Gamma}$ defined by
\begin{equation*}
b_{\Gamma}(\lambda,v) = \mathbf{v}^T \widetilde{\mathbf{B}}_{\Gamma} \boldsymbol{\lambda} \ \text{for } v \in \widetilde{\Uc}_H, \quad \text{and} \quad 
b_{\Gamma}(\lambda,v) = 
\mathbf{v}^T \mathbf{B}_{\Gamma} \boldsymbol{\lambda} \ \text{for } v \in \Wc_h.
\end{equation*}

In the discrete setting, an approximation of global problem \eqref{globalpb} reads: find $\widehat{U}^{k} \in \widetilde{\Uc}_H^{\Xi}$ satisfying \eqref{globalpb} for all $\delta U \in \widetilde{\Uc}_H$. In an algebraic setting, it boils down to solving the following system of linear algebraic equations:
\begin{equation}\label{algebraicglobalpb}
\mathbf{C}_{\widetilde{\Omega}} \mathbf{\widehat{U}}^k(\xi) = \mathbf{C}_{\widetilde{\Lambda}} \mathbf{U}^{k-1}(\xi) - \widetilde{\mathbf{B}}_{\Gamma} \boldsymbol{\lambda}^{k-1}(\xi) + \mathbf{l}_{\Omega \setminus \Lambda}.
\end{equation}
An approximation of local problem \eqref{localpb} reads: find $(w^k,\lambda^k) \in \Wc_h^{\Xi} \times \Mc_h^{\Xi}$ such that it satisfies almost surely \eqref{localpb} for all $(\delta w,\delta \lambda) \in \Wc_h \times \Mc_h$. In an algebraic setting, it comes down to solving the following system of non-linear algebraic equations:
\begin{subequations}\label{algebraiclocalpb}
\begin{align}
&\mathbf{A}_{\Lambda}(\xi) \mathbf{w}^k(\xi) + \mathbf{N}_{\Lambda}(\mathbf{w}^k(\xi);\xi) - \mathbf{B}_{\Gamma} \boldsymbol{\lambda}^k(\xi) = \mathbf{l}_{\Lambda}(\xi),\\
&\mathbf{B}^T_{\Gamma} \mathbf{w}^k(\xi) = \mathbf{\widetilde{B}}^T_{\Gamma} \mathbf{U}^k(\xi).
\end{align}
\end{subequations}

\begin{rmrk}
The convergence properties of the algorithm may be affected by the choice of spatial approximation spaces. The discretization errors can be viewed as additional perturbations occurring at both global and local steps of the algorithm. The impact of these perturbations on the behavior of the algorithm is addressed through Theorem~\ref{thrm:robustness}.
\end{rmrk}

\subsection{Approximations at stochastic level}\label{sec:approx-stochastic}

At the stochastic level, we introduce a basis $\set{\psi_{\alpha}}_{\alpha\in \Fc}$ of $\xLtwo_{\mu}(\Xi)$ (typically a polynomial basis) and we consider approximation spaces $\Sc_\Ac = \spanset{\psi_{\alpha}}_{\alpha\in \Ac}$, where 
$\Ac$ is a finite subset of $\Fc = \xN^m$ which is a partially ordered set such that for $\alpha, \beta \in \Fc$, $\alpha\leq \beta \iff \alpha_i\leq \beta_i$ for all $i\in \set{1,\dots,m}$. 
Then a function $v = \sum_{\alpha\in \Ac} v_\alpha \psi_{\alpha} \in \Sc_\Ac$ is identified with the vector of its coefficients $(v_{\alpha})_{\alpha \in \Ac} \in \xR^{\#\Ac}$ on the basis of $\Sc_\Ac$.

At the global level, suppose that finite element random vectors associated with $U^{k-1}$ and $\lambda^{k-1}$ are respectively given by $\mathbf{U}^{k-1}(\xi) = \sum_{\alpha \in \Ac} \mathbf{U}_{\alpha}^{k-1} \psi_{\alpha}(\xi)$ and $\boldsymbol{\lambda}^{k-1}(\xi) = \sum_{\alpha \in \Ac} \boldsymbol{\lambda}_{\alpha}^{k-1} \psi_{\alpha}(\xi)$.
Then, finite element random vectors associated with $\widehat{U}^k$ and $U^k$ admit the expansions 
$\mathbf{\widehat{U}}^{k}(\xi) = \sum_{\alpha \in \Ac} \mathbf{\widehat{U}}_{\alpha}^{k} \psi_{\alpha}(\xi)$ 
and $\mathbf{U}^{k}(\xi) = \sum_{\alpha \in \Ac} \mathbf{U}_{\alpha}^{k} \psi_{\alpha}(\xi)$, respectively, 
 with 
\begin{equation*}
\mathbf{\widehat{U}}^k_{\alpha} = \mathbf{C}_{\widetilde{\Omega}}^{-1} \( \mathbf{C}_{\widetilde{\Lambda}} \mathbf{U}_{\alpha}^{k-1} - \widetilde{\mathbf{B}}_{\Gamma} \boldsymbol{\lambda}_{\alpha}^{k-1} + \mathbf{l}_{\Omega \setminus \Lambda} \xE(\psi_{\alpha}(\xi)) \) \quad \text{and} \quad \mathbf{U}^k_{\alpha} = \rho_k \mathbf{\widehat{U}}^k_{\alpha} + (1-\rho_k) \mathbf{U}^{k-1}_{\alpha}.
\end{equation*}
It is worthy noticing that since $\widetilde{B}_L$ and $\widetilde{C}_L$ are chosen deterministic on $\widetilde{\Lambda}$ (see Remark~\ref{determinsiticfictitiousoperators}), then $\mathbf{C}_{\widetilde{\Omega}}$ and $\mathbf{C}_{\widetilde{\Lambda}}$ are independent of $\xi$, and the set of expansion coefficients $\set{\mathbf{\widehat{U}}^k_{\alpha}}_{\alpha \in \Ac}$ are simply obtained by solving a system of only $\#\Ac$ uncoupled linear algebraic equations with the same deterministic global finite element matrix $\mathbf{C}_{\widetilde{\Omega}}$, which can be factorized only once for all at initialization of the iterative procedure. The solution of such uncoupled global problems can be performed in parallel using traditional solvers available in standard deterministic finite element codes.

\begin{rmrk}
Convergence acceleration techniques based on Quasi-Newton or Newton update formulas have been proposed in \cite{Gen09,Duv16} within the deterministic framework and rely on successive corrections of the global finite element matrix $\mathbf{C}_{\widetilde{\Omega}}$ at each global step of the iterative procedure in order to improve the convergence rate of the algorithm. In the present stochastic framework, this would yield to a parameter-dependent matrix $\mathbf{C}_{\widetilde{\Omega}}$, unless using deterministic approximations of the successive corrections of $\mathbf{C}_{\widetilde{\Omega}}$.
\end{rmrk}

At the local level, approximations of finite element random vectors $\mathbf{w}^k(\xi)$ and $\boldsymbol{\lambda}^k(\xi)$ associated with local iterates $w^k$ and $\lambda^k$, respectively, are searched under the form
$
\mathbf{w}^k(\xi) \approx \sum_{\alpha \in \Ac} \mathbf{w}_{\alpha}^k \psi_{\alpha}(\xi) $ {and} $ \boldsymbol{\lambda}^k(\xi) \approx \sum_{\alpha \in \Ac} \boldsymbol{\lambda}_{\alpha}^k \psi_{\alpha}(\xi).
$ 
The determination of these approximations through Galerkin projection methods \cite{Mat05,Nou09} requires the solution of a large system of $\#\Ac$ coupled non-linear algebraic equations whose computational cost and memory storage requirements may be prohibitive and whose implementation may be cumbersome as it usually requires a modification (or at least an adaptation) of existing deterministic codes. Note however that non intrusive (or weakly intrusive) implementations of Galerkin methods can be introduced \cite{Gir14,Gir15}.

Here, for computing approximations of local iterates, we rather rely on an adaptive least-squares method which uses evaluations of the solution of \eqref{localpb} at some samples $\set{\xi^l}_{l=1}^N$ of random variables $\xi$. These evaluations are obtained by $N$ calls to an existing non-linear deterministic solver, \ie{} without requiring any modification of the underlying deterministic computer code. For computing the solution $(\mathbf{w}^k(\xi^l),\boldsymbol{\lambda}^k(\xi^l))$ of \eqref{algebraiclocalpb} for $\xi = \xi^l$, 
 we employ a Newton-type iterative algorithm with some prescribed tolerance. Note that the resulting numerical error can be viewed as an additional perturbation occurring at each local step of the iterative algorithm, whose impact on the behavior of the algorithm is analyzed in Section~\ref{sec:robustness}. 
 The adaptive least-squares method is described in Section~\ref{sec:leastsquares}.

\subsection{Adaptive least-squares method for sparse polynomial approximation}\label{sec:leastsquares}

Here we describe an adaptive least-squares method for sparse approximation of a random vector $\mathbf{u} \in \xR^n \otimes \xLtwo_{\mu}(\Xi)$. We assume $\Xi \subset \xR^m$ ($m<\infty$) and we consider an orthonormal tensor product basis 
$\set{\psi_{\alpha}(\xi) = \prod_{i= 1}^m \psi^{(i)}_{\alpha_i}(\xi_i)}_{\alpha \in \Fc}$ of $\xLtwo_{\mu}(\Xi)$, where $\psi^{(i)}_{k}$ is a univariate polynomial of degree $k$. For a given subset $\Ac \subset \Fc$, we define the corresponding polynomial space $\Sc_\Ac = \spanset{\psi_{\alpha}}_{\alpha\in \Ac}$. A subset $\Ac$ is called monotone (or lower or downward closed) if it is such that $(\beta \in \Ac \quad \text{and} \quad \alpha \leq \beta) \implies \alpha \in \Ac$. If $\Ac$ is monotone, the subspace $\Sc_\Ac$ coincides with the polynomial space $\mathbb{P}_{\Ac} = \spansetst{\xi_1^{\alpha_1}\dots \xi_m^{\alpha_m}}{\alpha \in \Ac}$ whatever the choice of univariate polynomial bases. 
Note that univariate polynomial bases could be replaced by other hierarchical bases (such as wavelet bases) with which we can expect accurate sparse approximations of the random vector. 

 \subsubsection{Approximation in a given subspace}
 
 For a given subset $\Ac$, a least-squares approximation $\mathbf{v}$ of $\mathbf{u}$ in $\xR^n\otimes \Sc_{\Ac}$ can be written as $\mathbf{v}(\xi) = \sum_{\alpha \in \Ac} \mathbf{v}_\alpha \psi_{\alpha}(\xi)$, where the set of coefficients $\mathbf{V} = (\mathbf{v}_\alpha)_{\alpha\in \Ac} \in \xR^{n\times \#\Ac}$ is solution of 
\begin{align}\label{minpb}
\min_{(\mathbf{v}_\alpha)_{\alpha\in \Ac}} \sum_{l=1}^N \norm{\mathbf{u}(\xi^l) - \sum_{\alpha \in \Ac} \mathbf{v}_\alpha \psi_{\alpha}(\xi^l)}_2^2.
%= \min_{\mathbf{V}} \norm{\mathbf{Y} - \boldsymbol{\Psi}\mathbf{V}^T}_F^2.
\end{align}
Assuming $N\geq \#\Ac$ and $ \boldsymbol{\Psi}^T\boldsymbol{\Psi}$ invertible, this yields
$
\mathbf{V}^T = (\boldsymbol{\Psi}^T\boldsymbol{\Psi})^{-1}\boldsymbol{\Psi}^T\mathbf{Y},
$ 
where $\boldsymbol{\Psi} = (\psi_{\alpha}(\xi^l))_{1\leq l\leq N,\alpha\in\Ac} \in \xR^{N \times {\#\Ac}}$ and $\mathbf{Y} = (\mathbf{u}(\xi^l))_{1\leq l\leq N} \in \xR^{N \times n}$. 
The stability of the least-squares approximation is related to the properties of the random matrix $\boldsymbol{\Psi}^T \boldsymbol{\Psi}$. Some theoretical results can be found in \cite{%Coh13b,Mig14,
Chk15b} and the references therein. In practice, for a given set $\Ac$, the stability of the least-squares approximation can be improved by increasing the number of samples. 

The approximation error can be estimated \aposteriori{} using cross-validation techniques which are classical statistical methods for computing error estimates based on a random partitioning of the available sample set into two subsets, the training set (or learning set) and the test set (or validation set). In the $k$-fold cross-validation procedure, the sample set $\upxi = \set{\xi^l}_{l=1}^{N}$ is randomly partitioned into $k$ disjoint and complementary sample subsets $\set{\upxi_s}_{s=1}^k$ of nearly equal size. Each subset $\upxi_s$ is in turn retained as the test set, while the remaining $k-1$ subsets gathered in $\upchi_s = \upxi \setminus \upxi_s$ are used as the training set. An approximation $\mathbf{v}$ is computed independently for each training set $\upchi_s$ and tested against the corresponding validation set $\upxi_s$ in order to assess its accuracy. The cross-validation error is estimated for each of the $k$ training sets $\upchi_s$ and then averaged over the $k$ sets. Such a cross-validation technique requires $k$ additional calls to the least-squares solver and thus may be computationally demanding. In practice, the vector of cross-validation error estimates $\boldsymbol{\varepsilon} = (\varepsilon_i)_{i\in \Ic}$ can be directly obtained from the approximate random vector $\mathbf{v} = (v_i)_{i\in \Ic}$ (computed using the available sample set $\upxi = \set{\xi^l}_{l=1}^{N}$) using the Bartlett matrix inversion formula \cite{Bar51} (a special case of well-known Sherman-Morrison-Woodbury formula) without any additional call to the least-squares solver. % \cite{She49,She50,Hag89}. 
The leave-one-out cross-validation procedure is a special case of $k$-fold cross-validation procedure where the number of folds $k$ is equal to the number of samples $N$. Note that the $k$-fold cross-validation technique depends on the chosen partition, contrary to the leave-one-out cross-validation technique. In the present work, we use the fast leave-one-out cross-validation procedure presented in \cite{Caw04} and summarized in Algorithm~\ref{LOOCValgo} {(see Section~\ref{sec:LOOCValgo} in Appendix~\ref{sec:appendix})} to assess the accuracy of $\mathbf{v}$.

\subsubsection{Working set strategy for adaptive approximation}

Now, we introduce a working set strategy for the 
construction of a sequence of approximation spaces 
$(\Sc_{\Ac_j})_{j\geq 1}$, where $(\Ac_j)_{j\geq 1}$ is an increasing sequence of monotone sets. 
Given $\Ac_j$, we define $\Ac_{j+1} = \Ac_j \cup \Nc_j$, where $\Nc_j$ is selected 
in a set of candidate multi-indices in $\Fc\setminus \Ac_j$. A natural approach consists in choosing for 
$\Nc_j$ a subset of the margin 
$\Mc_j = \Mc(\Ac_j)$ of $\Ac_j$, where the margin of a monotone set $\Ac$ is defined by 
\begin{equation*}
\Mc(\Ac) = \setst{\alpha \not \in \Ac}{\existsin{i}{\set{1,\dots,m}} \text{ such that } \alpha_i \neq 0 \implies \alpha - e_i \in \Ac},
\end{equation*}
where $(e_i)_j =\delta_{ij}$ is the Kronecker delta for $i,j\in \xN^{\ast}$.
A strategy for the selection of $\Nc_j$ (referred to as bulk search strategy in \cite{Chk13}) consists in computing 
a least-squares approximation $\mathbf{v}(\xi) = \sum_{\alpha \in \Ac_j\cup \Mc_j} \mathbf{v}_{\alpha} \psi_{\alpha}(\xi)$ associated with the augmented approximation space $\Sc_{\Ac_j\cup \Mc_j}$, and then in selecting 
a subset $\Nc_j$ such that $e(\Nc_j)\geq \theta e(\Mc_j)$, where $\theta \in \intervalcc{0}{1}$ is a parameter and where for a given set $\Nc$, $e(\Nc) = \sum_{\alpha \in \Nc} \norm{\mathbf{v}_\alpha}^2_2$ corresponds to the contribution of 
 coefficients $(\mathbf{v}_\alpha)_{\alpha \in \Nc}$ to the $\xLtwo_{\mu}$-norm of $\mathbf{v}$. 
Note that the construction of an optimal (smallest) monotone subset $\Nc_j$ in the margin $\Mc_j$ of $\Ac_j$ by a fast algorithm is still an open question. 
\begin{rmrk}
A practical choice for constructing a monotone set $\Nc_j$ is to consider the smallest subset in the margin $\Mc_j$ such that $e(\Nc_j)\geq \theta e(\Mc_j)$ and which contains the multi-indices $\alpha$ corresponding to the largest elements in the monotone envelope\footnote{For a given set $\Ac$, the monotone envelope (also called monotone majorant) $(\mathfrak{v}_{\alpha})_{\alpha \in \Ac}$ of a bounded sequence $(\norm{\mathbf{v}_\alpha}_2)_{\alpha \in \Ac}$ is defined by $\mathfrak{v}_{\alpha} = \max_{\beta \in \Ac, \beta \geq \alpha} \norm{\mathbf{v}_{\beta}}_2$ for $\alpha \in \Ac$.} $(\mathfrak{v}_{\alpha})_{\alpha \in \Mc_j}$ of the bounded sequence $(\norm{\mathbf{v}_\alpha}_2)_{\alpha \in \Mc_j}$. 
\end{rmrk}
Also, as the cardinality of the margin $\Mc(\Ac_j)$ may become prohibitively large in high dimension $m$, an alternative strategy consists in considering for $\Mc_j$ the reduced margin $\Mc_{\text{red}}(\Ac_j)$ of $\Ac_j$, where 
\begin{equation*}
\Mc_\text{red}(\Ac) = \setst{\alpha \not \in \Ac}{\forallin{i}{\set{1,\dots,m}} \text{ such that } \alpha_i \neq 0 \implies \alpha - e_i \in \Ac}.
\end{equation*}
The additional set $\Nc_j$ is then defined as the smallest non-empty subset of the reduced margin $\Mc_j$
of $\Ac_j$ such that $e(\Nc_j)\geq \theta e(\Mc_j)$, which is a monotone set by construction. Therefore, $\Ac_{j+1} = \Ac_j \cup \Nc_j$, as a union of monotone sets, is a monotone set. For $\theta = 1$, the selected subset $\Nc_j = \setst{\alpha\in \Mc_j}{\norm{\mathbf{v}_\alpha}_2 \neq 0}$. 
For $\theta=0$, $\Nc_j$ is one arbitrary element of $\setst{\alpha\in \Mc_j}{\alpha=\argmax_{\alpha\in \Mc_j}\norm{\mathbf{v}_\alpha}_2}$, and the strategy corresponds 
to the largest neighbor strategy proposed in \cite{Chk13}. In the numerical experiments, we will consider the strategy with a parameter $\theta =0.5$.

\subsubsection{Adaptive strategy}

In order to reach a desired accuracy, we finally propose an algorithm with adaptive random sampling and an adaptive selection of the approximation space $\Sc_\Ac$ with the working set strategy presented above. The algorithm is summarized in Algorithm~\ref{LSalgo} {(see Section~\ref{sec:LSalgo} in Appendix~\ref{sec:appendix})}. The convergence, stagnation and overfitting criteria in Algorithm~\ref{LSalgo} are respectively defined by
\begin{equation*}
\norm{\boldsymbol{\varepsilon}}_2 \leq \varepsilon_{\text{cv}}, \quad \frac{\norm{\boldsymbol{\varepsilon} - \boldsymbol{\varepsilon}^{\text{prev}}}_2}{\norm{\boldsymbol{\varepsilon}}_2} \leq \varepsilon_{\text{stagn}} \quad \text{and} \quad \frac{\norm{\boldsymbol{\varepsilon}}_2}{\norm{\boldsymbol{\varepsilon}^{\text{prev}}}_2} > 1+\varepsilon_{\text{overfit}},
\end{equation*}
where $\boldsymbol{\varepsilon}$ (resp. $\boldsymbol{\varepsilon}^{\text{prev}}$) is the vector of cross-validation error estimates computed at current (resp. previous) iteration, and $\varepsilon_{\text{cv}}$ (resp. $\varepsilon_{\text{stagn}}$ and $\varepsilon_{\text{overfit}}$) is the convergence (resp. stagnation and overfitting) threshold.

For computing sparse approximations of the solution $(\mathbf{w}^k(\xi),\boldsymbol{\lambda}^k(\xi))$ of non-linear local problem \eqref{algebraiclocalpb} at iteration $k$ of the algorithm, we apply the adaptive least-squares method to $\mathbf{u}(\xi) = (\mathbf{w}^k(\xi),\boldsymbol{\lambda}^k(\xi))$. In the adaptive sparse least-squares solver presented in Algorithm~\ref{LSalgo}, the adaptive random sampling step is performed by computing $N_{\text{add}}$ additional samples $(\mathbf{w}^k(\xi^{N+l}),\boldsymbol{\lambda}^k(\xi^{N+l}))_{1\leq l\leq N_{\text{add}}}$ of random vectors $(\mathbf{w}^k(\xi),\boldsymbol{\lambda}^k(\xi))$ (thus solving $N_{\text{add}}$ deterministic non-linear local problems) until the vectors of cross-validation error estimates $\boldsymbol{\varepsilon}^{\mathbf{w}}$ and $\boldsymbol{\varepsilon}^{\boldsymbol{\lambda}}$ for both random vectors $\mathbf{w}^k(\xi)$ and $\boldsymbol{\lambda}^k(\xi)$ satisfy a stagnation criterium. Then, the working set strategy is applied to $\mathbf{w}^k(\xi)$ and $\boldsymbol{\lambda}^k(\xi)$ separately to construct a sparse polynomial approximation space $\Sc_{\Ac}$ dedicated to each random vector.

\begin{rmrk}
Following Remark~\ref{localpbgeomvariability}, in the case of a patch $\Lambda$ containing geometrical variabilities, remeshing techniques could also be used in conjunction with sampling-based approaches. Indeed, the algorithm requires the computation of local iterates $\lambda^k$ (defined on a deterministic interface $\Gamma$) but not of local iterates $w^k$ (defined on the uncertain domain $\Lambda(\xi)$). An approximation of $\lambda^k$ can therefore be computed from samples of the local iterate $\lambda^k(\xi)$, where each new sample %$\lambda^k(\xi^l)$ 
involves a specific mesh %$\Tc_h(\Lambda(\xi^l))$ 
of $\Lambda(\xi)$. After convergence of the algorithm (at final iteration $k$), samples of the local solution $w^k(\xi)$ can be obtained by using remeshing techniques for each sample, from which we can deduce samples of quantities of interest expressed as functionals of $w^k(\xi)$. 
An explicit approximation of $w^k$ in terms of $\xi$ can also be obtained by using approaches based on reformulations of the local problem on a fixed deterministic domain (such as the random mapping techniques or the fictitious domain approaches briefly described in Remark~\ref{localpbgeomvariability}).
\end{rmrk}

\subsection{Relaxation step}\label{sec:relaxation}

The relaxation step may affect the convergence rate of the iterative algorithm as it can be interpreted as a line-search step of a non-linear solver.

%\subsubsection{Fixed relaxation parameter}\label{sec:fixedrelaxation}

The simplest method is to choose a fixed relaxation parameter $\rho$ throughout iterations. A large relaxation parameter may speed up the convergence but can lead to a divergence of the algorithm, while a small relaxation parameter ensures the convergence but triggers more iterations in return. The computation of an optimal relaxation parameter leading to an optimal convergence rate of the algorithm is not obvious in the non-linear framework. Even in the linear case, the optimal fixed value of relaxation parameter $\rho$ is problem-dependent and not known \apriori{}.
%
%\begin{rmrk}
%If $B=B_L$ and $C=C_L$ and if the mapping $A(\cdot;\xi)$ given in \eqref{mapA} is linear, symmetric and positive definite on $\widetilde{\Uc}_{\star}$, then an optimal convergent rate of the global-local iterative algorithm can be achieved by computing the optimal relaxation parameter $\rho_{\mathrm{opt}}$ defined by $\rho_{\mathrm{opt}} = \frac{2}{\lambda_{\mathrm{min}} + \lambda_{\mathrm{max}}}$ \cite{Xu92},
%where $\lambda_{\mathrm{min}}$ and $\lambda_{\mathrm{max}}$ are the lowest and largest eigenvalues of $A$ (see \cite{Che13a}).
%\end{rmrk}

%\subsubsection{Dynamic relaxation: Aitken's acceleration}\label{sec:dynamicrelaxation}

The Aitken's Delta-Squared method \cite{Iro69,Mac86} is a convergence acceleration technique which allows improving the current solution by using information gained at two previous iterations. The current global iterate $U^{k}$ is then obtained from the two pairs $(\widehat{U}^{k},U^{k-1})$ and $(\widehat{U}^{k-1},U^{k-2})$ and defined as
\begin{equation*}
U^{k}(\xi) = \widehat{U}^{k}(\xi) - \frac{\scalproda{\delta^{k}(\xi) - \delta^{k-1}(\xi)}{\widehat{U}^{k}(\xi) - \widehat{U}^{k-1}(\xi)}_{\widetilde{\Uc}}}{\norm{\delta^{k}(\xi) - \delta^{k-1}(\xi)}^2_{\widetilde{\Uc}}} \delta^{k}(\xi),
\end{equation*}
where $\delta^{k}(\xi) = \widehat{U}^{k}(\xi) - U^{k-1}(\xi)$ is the difference between the current global solution $\widehat{U}^{k}(\xi)$ and the previous global iterate $U^{k-1}(\xi)$. Then, using \eqref{relaxation}, the relaxation parameter $\rho_k$ is dynamically updated and defined by
\begin{equation}\label{Aitkenrelaxation}
\rho_k = -\rho_{k-1} \frac{\scalproda{\delta^{k}(\xi) - \delta^{k-1}(\xi)}{\delta^{k-1}(\xi)}_{\widetilde{\Uc}}}{\norm{\delta^{k}(\xi) - \delta^{k-1}(\xi)}^2_{\widetilde{\Uc}}}.
\end{equation}
As the Aitken's recursive formula \eqref{Aitkenrelaxation} requires two iterations of the algorithm, the first two values $\rho_1$ and $\rho_2$ are commonly set to $1$. For the subsequent iterations, the relaxation parameter $\rho_k$ can de defined as
\begin{equation*}
\rho_k = T_{\intervalcc{\rho_{\mathrm{inf}}}{\rho_{\mathrm{sup}}}}\(-\rho_{k-1} \frac{\scalproda{\delta^{k}(\xi) - \delta^{k-1}(\xi)}{\delta^{k-1}(\xi)}_{\widetilde{\Uc}}}{\norm{\delta^{k}(\xi) - \delta^{k-1}(\xi)}^2_{\widetilde{\Uc}}}\)
\end{equation*}
where $T_{\intervalcc{\rho_{\mathrm{inf}}}{\rho_{\mathrm{sup}}}}(\rho)$ is the projection of $\rho$ on the interval $\intervalcc{\rho_{\mathrm{inf}}}{\rho_{\mathrm{sup}}}$, 
which allows to ensure the convergence of the algorithm (see convergence condition \eqref{convergencecondition}).

Such a convergence acceleration technique is very simple to implement and computationally cheap. Also, the Aitken's acceleration method has been successfully applied to relaxation-based fixed-point algorithms in the context of fluid-structure interaction \cite{Kut08} and multiscale coupling \cite{Liu14,Duv16} problems. It has been proved to be particularly efficient with good convergence properties at low cost, compared with other relaxation methods such as the steepest descent method.

\section{Numerical results}\label{sec:results}

In order to demonstrate the efficiency and the robustness of the proposed method, we present different numerical experiments for a %different two-dimensional numerical examples are considered. 
%All the computations presented in this section were performed using the Parallel Computing toolbox of Matlab software on $1$ compute node containing $4$ $8$-core Intel processors ($12$ GB of memory per core) for a total of $32$ cores and $384$ GB of memory. Parallel computations have been performed on $32$ cores of a single compute node. 
stationary non-linear diffusion-reaction equation defined on a deterministic rectangular (two-dimensional) domain $\Omega = (0,2) \times (0,16) \subset \xR^2$. This equation is complemented with deterministic homogeneous Dirichlet boundary conditions $u=0$ applied on the entire boundary $\Gamma_D = \partial \Omega$. A deterministic volumetric source term $f = 1$ is imposed on the whole domain $\Omega$. The only sources of uncertainties come from diffusion coefficient $K(x,\xi)$ and reaction parameter $R(x,\xi)$ which are input random fields depending on a set of random variables $\xi \in \Xi$. The variabilities are assumed to be confined in $Q = 8$ patches $\set{\Lambda_q}_{q=1}^{Q}$ distributed along the $y$ (vertical) axis. Each patch $\Lambda_q$ is a square subdomain 
%$\Lambda_q = ({c_q}_x-\frac{L_q}{2},{c_q}_x+\frac{L_q}{2}) \times ({c_q}_y-\frac{L_q}{2},{c_q}_y+\frac{L_q}{2})$ with side length $L_q=1$ and center $c_q=(1,2q-1)$, so that 
$\Lambda_q=(0.5,1.5)\times(2q-1.5,2q-0.5)$. None of the patches present geometrical details. Domain $\Omega$ and the set of $Q$ patches $\set{\Lambda_q}_{q=1}^Q$ are illustrated on Figure~\ref{plotdomainglobalpatches_iso}.

\setlength\figureheight{0.25\textheight}
\setlength\figurewidth{0.25\textwidth}
\begin{figure}[h!]
\centering
\definecolor{mycolor1}{rgb}{1.00000,1.00000,0.00000}
\definecolor{mycolor2}{rgb}{1.00000,0.00000,1.00000}
\definecolor{mycolor3}{rgb}{0.00000,1.00000,1.00000}
\begin{subfigure}[t]{\figurewidth}
	\centering
	\tikzsetnextfilename{domain_global_patches_iso}
	\input{domain_global_patches}
	\caption{{Domain and patches}}\label{plotdomainglobalpatches_iso}
\end{subfigure}%\hfill
%\begin{subfigure}[t]{\figurewidth}
%	\centering
%	\tikzsetnextfilename{domain_global_patches_iso}
%	\input{domain_global_patches}
%	\caption{Isotropic case: $\gamma_q = 1$}\label{plotdomainglobalpatches_iso}
%\end{subfigure}%\hfill
%\begin{subfigure}[t]{1.5\figurewidth}
%	\centering
%	\tikzsetnextfilename{domain_global_patches_aniso}
%	\input{domain_global_patches}
%	\caption{Anisotropic case: $\gamma_q = 1-0.1(q+1)$}\label{plotdomainglobalpatches_aniso}
%\end{subfigure}%\hfill
\begin{subfigure}[t]{\figurewidth}
	\centering
	\includegraphics[height=\figureheight]{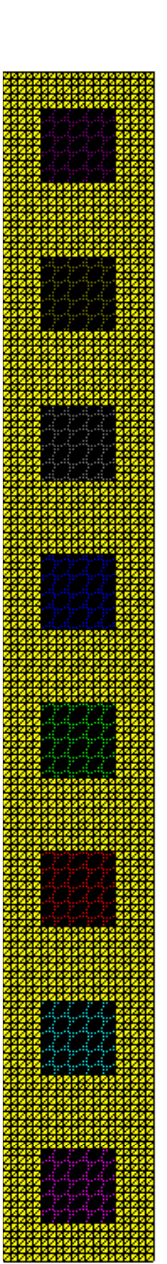}
	\caption{Nested triangulations}\label{plotmeshglobalpatches_iso}
\end{subfigure}%\hfill
\begin{subfigure}[t]{\figurewidth}
	\centering
	\includegraphics[height=0.5\figureheight]{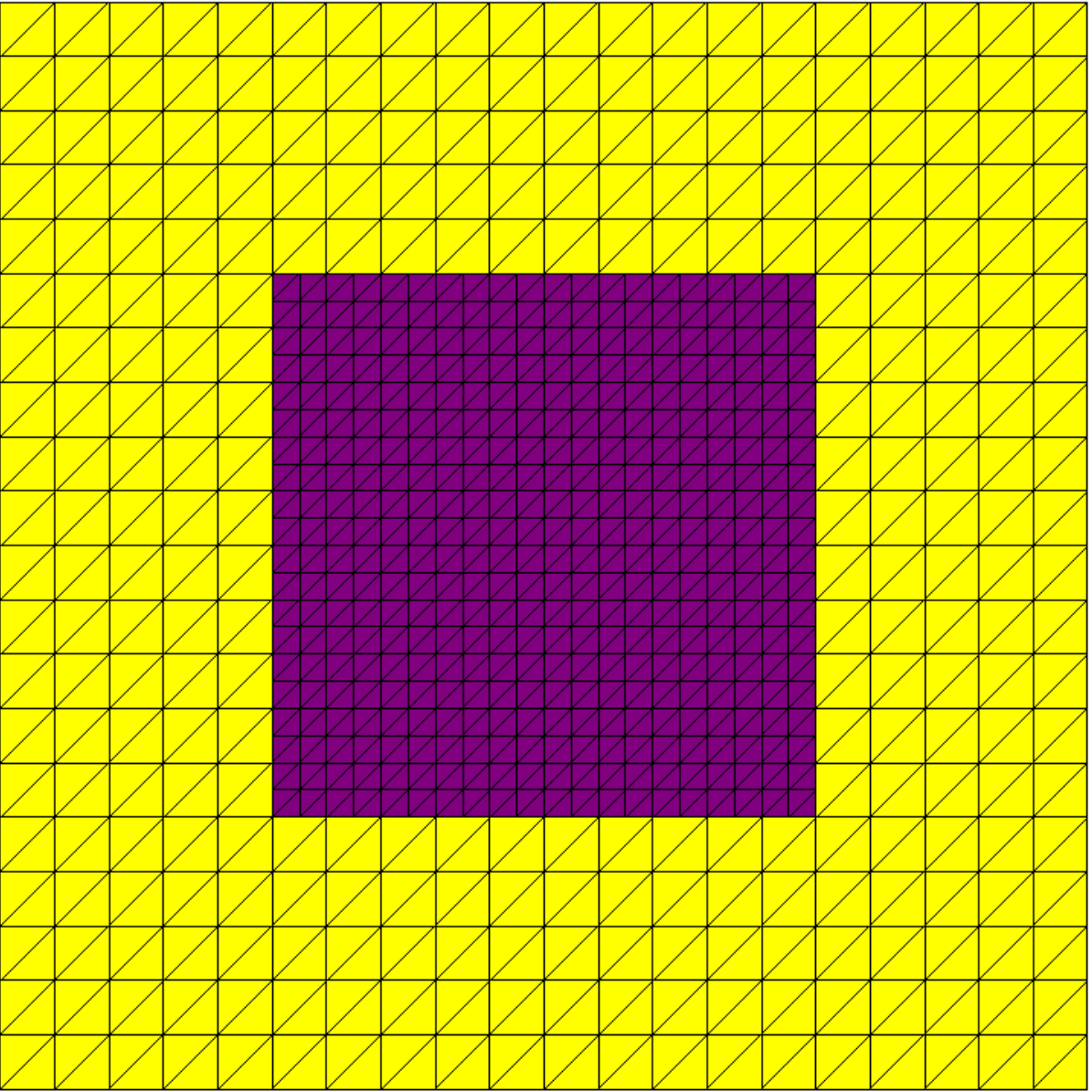}
	\caption{Zoom around patch $\Lambda_1$}\label{plotmeshglobalpatcheszoom_iso}
\end{subfigure}%\hfill
\ref{domainpatches_iso}
\caption{(\subref{plotdomainglobalpatches_iso}) Domain $\Omega$ partitioned into $Q = 8$ patches $\set{\Lambda_q}_{q=1}^{Q}$ and the complementary subdomain $\Omega \setminus \Lambda$ with $\Lambda = \bigcup_{q=1}^{Q} \Lambda_q$, 
%(\subref{plotdomainglobalpatches_iso}) isotropic case with the same weights $\gamma_q$ in every patch $\Lambda_q$ and (\subref{plotdomainglobalpatches_aniso}) anisotropic case with different weights $\gamma_q$ in every patch $\Lambda_q$, 
and (\subref{plotmeshglobalpatches_iso}) associated nested global and local finite element meshes $\Tc_h(\Lambda)$ and $\Tc_H(\Omega\setminus \Lambda)$ with (\subref{plotmeshglobalpatcheszoom_iso}) zoom around patch $\Lambda_1$}
\end{figure}

The solution $u$ satisfies almost surely
\begin{alignat*}{2}
-\nabla \cdot (K(x,\xi) \nabla u) + R(x,\xi) u^3 = f \quad \text{on } \Omega,
\quad u = 0 \quad \text{on } \Gamma_D = \partial \Omega,
\end{alignat*}
where random diffusion coefficient $K$ and random reaction parameter $R$ are such that
\begin{align*}
K(x,\xi) &= \begin{cases}
K_0 = 1 & \text{for } x \in \Omega \setminus \Lambda \\
K_q(x,\xi_{2q-1}) = 1 + \gamma_q\xi_{2q-1} \chi_q(x) & \text{for } x \in \Lambda_q, \ \text{for all } q \in \set{1,\dots,Q}
\end{cases},\\
R(x,\xi) &= \begin{cases}
0 & \text{for } x \in \Omega \setminus \Lambda \\
R_q(x,\xi_{2q}) = \gamma_q\xi_{2q} \chi_q(x) & \text{for } x \in \Lambda_q, \ \text{for all } q \in \set{1,\dots,Q}
\end{cases},
\end{align*}
with $\chi_q(x)$ the indicator function of subdomain $\Lambda^{\star}_q = (0.75,1.25) 
\times(2q-1.25,2q-0.75) \subset \Lambda_q$ for all $q \in \set{1,\dots,Q}$, and where the weights $\gamma_q$ are real coefficients in $(0,1)$ whose values define a level of uncertainty in the patch $\Lambda_q$. We consider two different situations: (i) an isotropic case, for which all weights $\gamma_q = 1$; (ii) an anisotropic case, for which the weights $\gamma_q = 1-0.1(q+1)$. Random diffusion coefficient $K$ and reaction parameter $R$ are parametrized by a set $\xi = (\xi_i)_{i=1}^m$ of $m = 2Q = 16$ real-valued random variables $\xi_i$ assumed to be mutually independent and uniformly distributed on $ (0,1)$. The parameter space is then the hypercube $\Xi = (0,1)^m \subset \xR^m$ endowed with the uniform probability measure. Each patch $\Lambda_q$ is characterized by a random diffusion coefficient $K_q$ (parametrized by $\xi_{2q-1}$) and a random reaction parameter $R_q$ (parametrized by $\xi_{2q}$). The {material} properties of patch $\Lambda_q$ therefore depends on $2$ real-valued random variables $\xi_{2q-1}$ and $\xi_{2q}$. The ranges of variations of $K_q$ and $R_q$ in each patch $\Lambda_q$ are respectively $(1,2)$ and $(0,1)$ for the isotropic case, and $(1,1+\gamma_q)$ and $(0,\gamma_q)$ for the anisotropic case. In this example, all uncertain material parameters $K_q$ and $R_q$ depend on the random variables $\xi$ in an affine manner, so do the bilinear forms $a_{\Lambda_q}$ and semi-linear forms $n_{\Lambda_q}$ for all $q \in \set{1,\dots,Q}$. As the Hilbert space $\Vc$ is assumed to be deterministic, the solution $u \in \xLn{p}_{\mu}(\Xi;\Vc)$ belongs to the tensor product vector space $\Vc \otimes 
\xLn{p}_{\mu}(\Xi)$. Given that Dirichlet boundary conditions are applied on the whole boundary $\partial \Omega$, and the source term $f$ is deterministic on domain $\Omega$, the linear form $\ell_{\Omega}$ is independent of $\xi$. Lastly, as domain $\Omega$ does not contain any geometrical defects, $\widetilde{\Omega} = \Omega$ in this example.

\subsection{Approximation spaces}

At the spatial level, we introduce nested finite element approximation spaces $\widetilde{\Uc}_H \subset \widetilde{\Uc}$ and $\Wc^q_h \subset \Wc^q$ for all $q \in \set{1,\dots,Q}$ (see Figures~\ref{plotmeshglobalpatches_iso} and \ref{plotmeshglobalpatcheszoom_iso}). The coarse global mesh $\Tc_H(\widetilde{\Omega})$ is a regular triangulation of fictitious domain $\widetilde{\Omega}$ which is composed of $3$-nodes linear triangular elements with uniform element size $H = 0.1$. It thus comprises $3\,381$ nodes and $6\,400$ elements. For every $q \in \set{1,\dots,Q}$, the fine local mesh $\Tc_h(\Lambda_q)$ is a regular triangulation of patch $\Lambda_q$ which is composed of $3$-nodes linear triangular elements with uniform element size $h_q = 0.05$. It is thus made of $441$ nodes and $800$ elements. Every fine local mesh $\Tc_h(\Lambda_q)$ corresponds to a uniform refinement of the corresponding coarse local mesh $\Tc_H(\widetilde{\Lambda}_q)$. 
The resulting spatial approximation spaces $\widetilde{\Uc}_H$ and $\Wc^q_h$ for all $q \in \set{1,\dots,Q}$, have dimensions $n_U = \xdim(\widetilde{\Uc}_H) = 3\,381$ and $n_{w_q} = \xdim(\Wc^q_h) = 441$, respectively.

At the stochastic level, we adaptively build a multidimensional polynomial approximation space $\Sc_\Ac $ spanned by generalized polynomial chaos basis $\set{\psi_{\alpha}}_{\alpha \in \Ac}$ (multidimensional Legendre polynomials) by using the adaptive sparse least-squares solver described in Algorithm~\ref{LSalgo}. 
At each iteration of the iterative algorithm, the linear global problem \eqref{globalpb} defined on fictitious domain $\widetilde{\Omega}$ is solved exactly (at the machine precision using a direct solver) as it involves a deterministic operator, while the $Q$ non-linear local problems \eqref{localpb} defined on the $Q$ patches $\set{\Lambda_q}_{q=1}^Q$ are solved using the adaptive sampling-based least-squares method described in Section~\ref{sec:leastsquares}. For each patch $\Lambda_q$, the $N_q$ deterministic non-linear local problems \eqref{algebraiclocalpb} associated with the $N_q$ samples $\set{\xi^l}_{l=1}^{N_q}$ are partially solved using a tangent-Newton iterative algorithm with a prescribed tolerance set to $\epsilon = 10^{-12}$. In our application case, one deterministic non-linear local problem typically requires only few iterations (less than $5$) to reach this stopping criterion. All the local computations have been performed in parallel on $32$ cores of a single computer node. 
The sample set $\set{\xi^l}_{l=1}^{N_q}$ and the approximation spaces $\Sc_{\Ac_q}$ are sequentially enriched (independently for each patch)
in order to control the accuracy of the local approximations $(w^k_q,\lambda^k_q)$. The stagnation and overfitting thresholds in Algorithm~\ref{LSalgo} are both set to $\varepsilon_{\text{stagn}} = \varepsilon_{\text{overfit}} = 10^{-1}$. An initial sample set of size $N=1$ is used with a sampling factor $p_{\text{add}}=0.1$ (percentage of additional samples) and a parameter $\theta=0.5$ for a good trade-off between computational efficiency and stability of local solutions $(w^k_q,\lambda^k_q)$.

\subsection{Convergence analysis}

The accuracy of global approximations $U^k$ is measured in $\xLtwo_{\mu}$-norm with respect to a global reference solution $U^{\text{ref}}$ using the relative error indicator $\varepsilon_{\Xi,\Omega\setminus\Lambda}$ defined as
\begin{equation*}%\label{errorindicator}
\varepsilon_{\Xi,\Omega\setminus\Lambda}(U^k;U^{\text{ref}}) = \frac{\norm{U^k-U^{\text{ref}}}_{\xLtwo_{\mu}(\Xi;\xLtwo(\Omega\setminus\Lambda))}}{\norm{U^{\text{ref}}}_{\xLtwo_{\mu}(\Xi;\xLtwo(\Omega\setminus\Lambda))}}, \quad 
\text{with} \quad \norm{U}^2_{\xLtwo_{\mu}(\Xi;\xLtwo(\Omega\setminus\Lambda))} = \xE(\norm{U(\xi)}^2_{\xLtwo(\Omega\setminus\Lambda)}).
\end{equation*}
The reference solution $(U^{\text{ref}},w^{\text{ref}},\lambda^{\text{ref}})$ is obtained by directly solving the full-scale coupled problem \eqref{globallocalformulation} using the adaptive sparse least-squares method described in Section~\ref{sec:leastsquares}. Following Theorem~\ref{thrm:consistency}, the global reference solution $U^{\text{ref}}$ is the restriction to subdomain $\Omega\setminus \Lambda$ of the limit $U$ of the sequence of global iterates $U^k$. At the spatial level, the global (resp. local) reference solution $U^{\text{ref}}$ (resp. $(w_q^{\text{ref}},\lambda_q^{\text{ref}})$) is discretized using the same finite element mesh as the global (resp. local) approximations $U^k$ (resp. $(w^k_q,\lambda^k_q)$). At the stochastic level, the number of samples $N^{\text{ref}}$ and the approximation spaces $\Sc_{\Ac^{\text{ref}}}$ are controlled by using the leave-one-out cross-validation procedure presented in Algorithm~\ref{LOOCValgo}. The prescribed tolerance for the convergence of Algorithm~\ref{LSalgo} is set to $\varepsilon^{\text{ref}}_{\text{cv}} = 10^{-6}$. The resulting sample size is $N^{\text{ref}} = 795$ for the isotropic case and $N^{\text{ref}} = 491$ for the anisotropic case. The partial polynomial degrees $p_i^{\text{ref}}$ (in each random variable $\xi_i$) and the dimension $\#\Ac^{\text{ref}}$ of approximation spaces $\Sc_{\Ac^{\text{ref}}}$ for global and local reference solutions $U^{\text{ref}}$ and $(w_q^{\text{ref}},\lambda_q^{\text{ref}})$ are reported in Table~\ref{degree_dim_ref} for both isotropic and anisotropic cases. The local reference solution $(w_q^{\text{ref}},\lambda_q^{\text{ref}})$ mainly depends on the random variables $\xi_{2q-1}$ and $\xi_{2q}$ associated with the corresponding patch $\Lambda_q$, as well as on the random variables confined in the surrounding patches. Note that, in the anisotropic setting, the stronger the variabilities in the material properties within a patch $\Lambda_q$ are, the more the reference local solution $(w_q^{\text{ref}},\lambda_q^{\text{ref}})$ is sensitive to the random variables $\xi_{2q-1}$ and $\xi_{2q}$ associated with $\Lambda_q$.

\begin{table}[h!]
\footnotesize{
\begin{subtable}[t]{\textwidth}
\centering
\begin{tabular}{cccccccccccccccccc}
\toprule
& $p^{\text{ref}}_1$ & $p^{\text{ref}}_2$ & $p^{\text{ref}}_3$ & $p^{\text{ref}}_4$ & $p^{\text{ref}}_5$ & $p^{\text{ref}}_6$ & $p^{\text{ref}}_7$ & $p^{\text{ref}}_8$ & $p^{\text{ref}}_9$ & $p^{\text{ref}}_{10}$ & $p^{\text{ref}}_{11}$ & $p^{\text{ref}}_{12}$ & $p^{\text{ref}}_{13}$ & $p^{\text{ref}}_{14}$ & $p^{\text{ref}}_{15}$ & $p^{\text{ref}}_{16}$ & $\#\Ac^{\text{ref}}$\\
\midrule
\multicolumn{1}{c}{$U^{\text{ref}}$} & $7$ & $3$ & $6$ & $4$ & $6$ & $4$ & $6$ & $4$ & $6$ & $4$ & $6$ & $4$ & $6$ & $4$ & $7$ & $3$ & $384$\\
\midrule
\multicolumn{1}{c}{$w^{\text{ref}}_1$} & \textcolor{red}{$7$} & \textcolor{red}{$3$} & $5$ & $3$ & $3$ & $2$ & $2$ & $2$ & $1$ & $1$ & $0$ & $0$ & $0$ & $0$ & $0$ & $0$ & $85$\\
\multicolumn{1}{c}{$\lambda^{\text{ref}}_1$} & \textcolor{red}{$8$} & \textcolor{red}{$4$} & $6$ & $4$ & $5$ & $3$ & $3$ & $2$ & $2$ & $1$ & $0$ & $1$ & $0$ & $0$ & $0$ & $0$ & $186$\\
\midrule
\multicolumn{1}{c}{$w^{\text{ref}}_2$} & $6$ & $3$ & \textcolor{red}{$7$} & \textcolor{red}{$4$} & $5$ & $3$ & $3$ & $2$ & $2$ & $2$ & $1$ & $1$ & $0$ & $0$ & $0$ & $0$ & $142$\\
\multicolumn{1}{c}{$\lambda^{\text{ref}}_2$} & $6$ & $3$ & \textcolor{red}{$7$} & \textcolor{red}{$4$} & $5$ & $3$ & $4$ & $2$ & $2$ & $2$ & $1$ & $1$ & $0$ & $0$ & $0$ & $0$ & $170$\\
\midrule
\multicolumn{1}{c}{$w^{\text{ref}}_3$} & $4$ & $2$ & $5$ & $3$ & \textcolor{red}{$7$} & \textcolor{red}{$4$} & $5$ & $3$ & $3$ & $2$ & $2$ & $1$ & $1$ & $1$ & $0$ & $0$ & $150$\\
\multicolumn{1}{c}{$\lambda^{\text{ref}}_3$} & $5$ & $2$ & $5$ & $3$ & \textcolor{red}{$7$} & \textcolor{red}{$4$} & $5$ & $3$ & $4$ & $3$ & $2$ & $2$ & $1$ & $1$ & $1$ & $0$ & $215$\\
\midrule
\multicolumn{1}{c}{$w^{\text{ref}}_4$} & $3$ & $1$ & $3$ & $2$ & $5$ & $3$ & \textcolor{red}{$7$} & \textcolor{red}{$4$} & $5$ & $3$ & $3$ & $2$ & $2$ & $2$ & $1$ & $1$ & $168$\\
\multicolumn{1}{c}{$\lambda^{\text{ref}}_4$} & $3$ & $2$ & $4$ & $3$ & $5$ & $3$ & \textcolor{red}{$7$} & \textcolor{red}{$4$} & $5$ & $3$ & $4$ & $3$ & $3$ & $2$ & $2$ & $1$ & $235$\\
\midrule
\multicolumn{1}{c}{$w^{\text{ref}}_5$} & $1$ & $1$ & $2$ & $1$ & $3$ & $2$ & $5$ & $3$ & \textcolor{red}{$7$} & \textcolor{red}{$4$} & $5$ & $3$ & $3$ & $2$ & $3$ & $1$ & $167$\\
\multicolumn{1}{c}{$\lambda^{\text{ref}}_5$} & $2$ & $1$ & $3$ & $2$ & $4$ & $3$ & $5$ & $3$ & \textcolor{red}{$7$} & \textcolor{red}{$4$} & $5$ & $3$ & $4$ & $2$ & $3$ & $2$ & $234$\\
\midrule
\multicolumn{1}{c}{$w^{\text{ref}}_6$} & $1$ & $0$ & $1$ & $1$ & $2$ & $1$ & $3$ & $2$ & $5$ & $3$ & \textcolor{red}{$7$} & \textcolor{red}{$3$} & $5$ & $3$ & $4$ & $2$ & $149$\\
\multicolumn{1}{c}{$\lambda^{\text{ref}}_6$} & $1$ & $0$ & $1$ & $1$ & $3$ & $2$ & $4$ & $3$ & $6$ & $3$ & \textcolor{red}{$7$} & \textcolor{red}{$4$} & $6$ & $3$ & $5$ & $2$ & $239$\\
\midrule
\multicolumn{1}{c}{$w^{\text{ref}}_7$} & $0$ & $0$ & $0$ & $0$ & $0$ & $1$ & $2$ & $2$ & $3$ & $2$ & $5$ & $3$ & \textcolor{red}{$8$} & \textcolor{red}{$4$} & $5$ & $3$ & $164$\\
\multicolumn{1}{c}{$\lambda^{\text{ref}}_7$} & $0$ & $0$ & $0$ & $0$ & $1$ & $1$ & $2$ & $2$ & $4$ & $2$ & $5$ & $3$ & \textcolor{red}{$7$} & \textcolor{red}{$4$} & $7$ & $3$ & $174$\\
\midrule
\multicolumn{1}{c}{$w^{\text{ref}}_8$} & $0$ & $0$ & $0$ & $0$ & $0$ & $1$ & $1$ & $1$ & $2$ & $1$ & $3$ & $2$ & $5$ & $3$ & \textcolor{red}{$8$} & \textcolor{red}{$3$} & $91$\\
\multicolumn{1}{c}{$\lambda^{\text{ref}}_8$} & $0$ & $1$ & $0$ & $0$ & $1$ & $1$ & $2$ & $2$ & $3$ & $2$ & $5$ & $3$ & $6$ & $4$ & \textcolor{red}{$8$} & \textcolor{red}{$4$} & $202$\\
\bottomrule
\end{tabular}
\caption{Isotropic case}\label{degree_dim_ref_iso}
\end{subtable}
\\
\begin{subtable}[t]{\textwidth}
\centering
\begin{tabular}{cccccccccccccccccc}
\toprule
& $p^{\text{ref}}_1$ & $p^{\text{ref}}_2$ & $p^{\text{ref}}_3$ & $p^{\text{ref}}_4$ & $p^{\text{ref}}_5$ & $p^{\text{ref}}_6$ & $p^{\text{ref}}_7$ & $p^{\text{ref}}_8$ & $p^{\text{ref}}_9$ & $p^{\text{ref}}_{10}$ & $p^{\text{ref}}_{11}$ & $p^{\text{ref}}_{12}$ & $p^{\text{ref}}_{13}$ & $p^{\text{ref}}_{14}$ & $p^{\text{ref}}_{15}$ & $p^{\text{ref}}_{16}$ & $\#\Ac^{\text{ref}}$\\
\midrule
\multicolumn{1}{c}{$U^{\text{ref}}$} & $6$ & $3$ & $5$ & $3$ & $5$ & $3$ & $4$ & $3$ & $4$ & $3$ & $4$ & $3$ & $3$ & $2$ & $3$ & $2$ & $173$\\
\midrule
\multicolumn{1}{c}{$w^{\text{ref}}_1$} & \textcolor{red}{$7$} & \textcolor{red}{$3$} & $4$ & $3$ & $2$ & $2$ & $1$ & $1$ & $0$ & $1$ & $0$ & $0$ & $0$ & $0$ & $0$ & $0$ & $57$\\
\multicolumn{1}{c}{$\lambda^{\text{ref}}_1$} & \textcolor{red}{$7$} & \textcolor{red}{$3$} & $5$ & $3$ & $4$ & $2$ & $2$ & $2$ & $1$ & $1$ & $0$ & $0$ & $0$ & $0$ & $0$ & $0$ & $129$\\
\midrule
\multicolumn{1}{c}{$w^{\text{ref}}_2$} & $5$ & $2$ & \textcolor{red}{$6$} & \textcolor{red}{$3$} & $4$ & $3$ & $3$ & $2$ & $1$ & $1$ & $0$ & $1$ & $0$ & $0$ & $0$ & $0$ & $95$\\
\multicolumn{1}{c}{$\lambda^{\text{ref}}_2$} & $5$ & $3$ & \textcolor{red}{$6$} & \textcolor{red}{$3$} & $4$ & $3$ & $3$ & $2$ & $2$ & $1$ & $1$ & $1$ & $0$ & $0$ & $0$ & $0$ & $115$\\
\midrule
\multicolumn{1}{c}{$w^{\text{ref}}_3$} & $4$ & $2$ & $4$ & $3$ & \textcolor{red}{$6$} & \textcolor{red}{$3$} & $4$ & $2$ & $2$ & $2$ & $1$ & $1$ & $0$ & $0$ & $0$ & $0$ & $106$\\
\multicolumn{1}{c}{$\lambda^{\text{ref}}_3$} & $4$ & $2$ & $5$ & $3$ & \textcolor{red}{$6$} & \textcolor{red}{$3$} & $4$ & $3$ & $3$ & $2$ & $1$ & $1$ & $0$ & $1$ & $0$ & $0$ & $127$\\
\midrule
\multicolumn{1}{c}{$w^{\text{ref}}_4$} & $2$ & $1$ & $3$ & $2$ & $4$ & $2$ & \textcolor{red}{$5$} & \textcolor{red}{$3$} & $3$ & $2$ & $2$ & $2$ & $1$ & $1$ & $0$ & $0$ & $84$\\
\multicolumn{1}{c}{$\lambda^{\text{ref}}_4$} & $3$ & $1$ & $3$ & $2$ & $4$ & $3$ & \textcolor{red}{$5$} & \textcolor{red}{$3$} & $4$ & $3$ & $2$ & $2$ & $1$ & $1$ & $1$ & $0$ & $112$\\
\midrule
\multicolumn{1}{c}{$w^{\text{ref}}_5$} & $1$ & $1$ & $1$ & $1$ & $2$ & $2$ & $3$ & $2$ & \textcolor{red}{$5$} & \textcolor{red}{$3$} & $3$ & $2$ & $2$ & $1$ & $1$ & $1$ & $68$\\
\multicolumn{1}{c}{$\lambda^{\text{ref}}_5$} & $2$ & $1$ & $2$ & $1$ & $3$ & $2$ & $4$ & $3$ & \textcolor{red}{$5$} & \textcolor{red}{$3$} & $3$ & $2$ & $2$ & $2$ & $1$ & $1$ & $94$\\
\midrule
\multicolumn{1}{c}{$w^{\text{ref}}_6$} & $0$ & $0$ & $0$ & $1$ & $1$ & $1$ & $2$ & $2$ & $3$ & $2$ & \textcolor{red}{$4$} & \textcolor{red}{$3$} & $3$ & $2$ & $2$ & $1$ & $54$\\
\multicolumn{1}{c}{$\lambda^{\text{ref}}_6$} & $0$ & $0$ & $1$ & $1$ & $2$ & $2$ & $3$ & $2$ & $4$ & $3$ & \textcolor{red}{$5$} & \textcolor{red}{$3$} & $3$ & $2$ & $2$ & $1$ & $84$\\
\midrule
\multicolumn{1}{c}{$w^{\text{ref}}_7$} & $0$ & $0$ & $0$ & $0$ & $0$ & $0$ & $1$ & $1$ & $2$ & $2$ & $3$ & $2$ & \textcolor{red}{$4$} & \textcolor{red}{$2$} & $2$ & $1$ & $40$\\
\multicolumn{1}{c}{$\lambda^{\text{ref}}_7$} & $0$ & $0$ & $0$ & $0$ & $1$ & $1$ & $2$ & $1$ & $3$ & $2$ & $3$ & $2$ & \textcolor{red}{$4$} & \textcolor{red}{$3$} & $3$ & $2$ & $55$\\
\midrule
\multicolumn{1}{c}{$w^{\text{ref}}_8$} & $0$ & $0$ & $0$ & $0$ & $0$ & $0$ & $0$ & $0$ & $1$ & $1$ & $2$ & $1$ & $2$ & $2$ & \textcolor{red}{$3$} & \textcolor{red}{$2$} & $23$\\
\multicolumn{1}{c}{$\lambda^{\text{ref}}_8$} & $0$ & $0$ & $0$ & $0$ & $0$ & $0$ & $1$ & $1$ & $2$ & $1$ & $3$ & $2$ & $3$ & $2$ & \textcolor{red}{$3$} & \textcolor{red}{$2$} & $41$\\
\bottomrule
\end{tabular}
\caption{Anisotropic case}\label{degree_dim_ref_aniso}
\end{subtable}
\caption{Partial polynomial degrees $p^{\text{ref}}_i$ (in each random variable $\xi_i$), $i \in \set{1,\dots,16}$, and dimension $\#\Ac^{\text{ref}}$ of approximation spaces $\Sc_{\Ac^{\text{ref}}}$ for global and local reference solutions $U^{\text{ref}}$ and $(w^{\text{ref}}_q,\lambda^{\text{ref}}_q)$, $q \in \set{1,\dots,8}$, where values displayed in red correspond to random variables $\xi_{2q-1}$ and $\xi_{2q}$ associated with patch $\Lambda_q$}\label{degree_dim_ref}
}
\end{table}

\begin{figure}[h!]
\centering
\tikzsetnextfilename{error_rho_iso}
\input{error_rho}
\caption{Isotropic problem: evolution of error indicator $\varepsilon_{\Xi,\Omega\setminus\Lambda}$ with respect to iteration number $k$ for different fixed relaxation parameters $\rho$ and for Aitken's dynamic relaxation parameter $\rho_k$}\label{ploterrorrho_iso}
\end{figure}
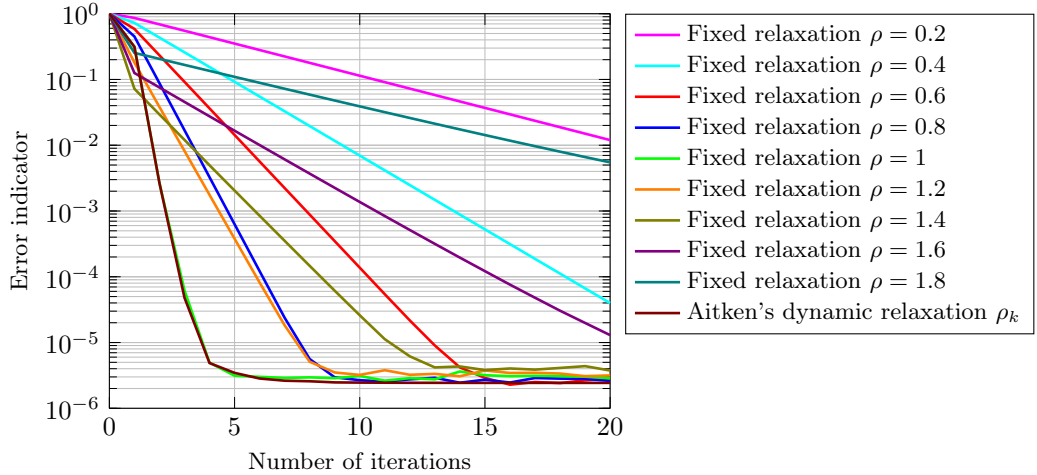

We first consider a fixed cross-validation tolerance $\varepsilon_{\text{cv}} = 10^{-3}$ in Algorithm~\ref{LSalgo} for the accuracy of local solutions $(w^k_q,\lambda^k_q)$ and we study the influence of relaxation parameter $\rho_k$ on the convergence of the global-local iterative algorithm. Figure~\ref{ploterrorrho_iso} represents the evolution of relative error indicator $\varepsilon_{\Xi,\Omega\setminus\Lambda}$ with respect to the number of iterations $k$ for different fixed values of relaxation parameter $\rho \in \set{0.2,0.4,0.6,0.8,1,1.2,1.4,1.6,1.8}$ and for a relaxation parameter $\rho_k$ dynamically updated through the Aitken's Delta-Squared method presented in Section~\ref{sec:relaxation} in the isotropic case. As expected, the relaxation parameter has a strong influence on the convergence rate of the iterative algorithm. The Aitken's Delta-Squared acceleration technique provides similar results as those obtained with an optimal fixed relaxation parameter without any additional computational cost. The relative error indicator decreases sharply ($\varepsilon_{\Xi,\Omega\setminus\Lambda} = 5.10^{-5}$ after only $k=3$ iterations) and reaches a plateau $\varepsilon_{\Xi,\Omega\setminus\Lambda} = 3.10^{-6}$ for $k\geq 5$. Similar results can be obtained for the anisotropic case.

We now use the Aitken's dynamic relaxation and we investigate the influence of the prescribed tolerance $\varepsilon_{\text{cv}}$ for cross-validation in Algorithm~\ref{LSalgo} on the performances of the global-local iterative algorithm in terms of accuracy and computational efficiency. 
Figure~\ref{ploterrorcputimetol} shows the evolution of relative error indicator $\varepsilon_{\Xi,\Omega\setminus\Lambda}$ and computational cost per iteration as functions of the number of iterations $k$ of the global-local iterative algorithm for different cross-validation tolerances $\varepsilon_{\text{cv}} \in \set{10^{-2}, 10^{-3}, 10^{-4}, 10^{-5}}$ for both problems. The iterative algorithm converges quite fast until the relative error indicator $\varepsilon_{\Xi,\Omega\setminus\Lambda}$ stabilizes around a value smaller than the tolerance $\varepsilon_{\text{cv}}$ imposed to Algorithm~\ref{LSalgo} for cross-validation. The cross-validation threshold $\varepsilon_{\text{cv}}$ can then be seen as the level of a perturbation occurring at each local step of the iterative algorithm and having an impact on its convergence properties. The accuracy of multiscale solution $u=(U,w_1,\dots,w_Q)$ obtained at convergence of the global-local iterative algorithm is then controlled by the cross-validation tolerance $\varepsilon_{\text{cv}}$ prescribed to Algorithm~\ref{LSalgo} at each local stage. Note that using a relatively high cross-validation tolerance $\varepsilon_{\text{cv}}=10^{-2}$ allows to reach a rather small precision $\varepsilon_{\Xi,\Omega\setminus\Lambda} = 6.10^{-5}$ (resp. $4.10^{-5}$) after only $k=3$ iterations for the isotropic (resp. anisotropic) case. 
Figure~\ref{plotrhotol} displays the evolution of Aitken's dynamic relaxation parameter $\rho_k$ with respect to iteration number $k$ for the aforementioned cross-validation tolerances $\varepsilon_{\text{cv}} = 10^{-2}, \dots, 10^{-5}$ for both isotropic and anisotropic cases. As mentioned in \cite{Kut08}, relaxation parameter $\rho_k$ varies but does not seem to follow a definite pattern during the iterations of the algorithm. In our numerical experiments, $\rho_k$ varies slightly around the value $1$ during the first iterations, then decreases toward zero as soon as the relative error indicator $\varepsilon_{\Xi,\Omega\setminus\Lambda}$ stagnates around a certain value depending on the prescribed tolerance $\varepsilon_{\text{cv}}$ for cross-validation. 
Figure~\ref{plotnbsamplesdimstochasticbasis4tol} shows the evolutions of the sample size $N_q$ and the dimension $\#\Ac_q$ of the approximation space for local solution $w^k_q$ and Lagrange multiplier $\lambda^k_q$ within patch $\Lambda_q$, for $q=4$, as functions of the number of iterations $k$ for cross-validation threshold values $\varepsilon_{\text{cv}}$ varying from $10^{-2}$ to $10^{-5}$ for both isotropic and anisotropic cases. The number of samples $N_q$ and the dimension $\#\Ac_q$ of the polynomial spaces increase during the first iterations and then stagnate around a certain value which is higher as the prescribed tolerance $\varepsilon_{\text{cv}}$ for cross-validation is lower. The sample sizes and the dimensions of approximation spaces are higher for the isotropic case than for the anisotropic case. Note that the dimension of the approximation space obtained for Lagrange multiplier $\lambda_q$ is higher than the one for local solution $w_q$.

\setlength\figureheight{0.2\textheight}
\setlength\figurewidth{0.4\textwidth}
\begin{figure}[h!]
\centering
\definecolor{mycolor1}{rgb}{1.00000,0.00000,1.00000}%
\definecolor{mycolor2}{rgb}{0.00000,1.00000,1.00000}%
\begin{subfigure}[t]{\figurewidth}
	\centering
	\tikzsetnextfilename{error_tol_iso}
	\input{error_tol_no_legend}
	\caption{Isotropic case}\label{ploterrortol_iso}
\end{subfigure}\hfill
\begin{subfigure}[t]{1.5\figurewidth}
	\centering
	\tikzsetnextfilename{error_tol_aniso}
	\input{error_tol}
	\caption{Anisotropic case}\label{ploterrortol_aniso}
\end{subfigure}
\begin{subfigure}[t]{\figurewidth}
	\centering
	\tikzsetnextfilename{cpu_time_iso}
	\input{cpu_time_no_legend}
	\caption{Isotropic case}\label{plotcputime_iso}
\end{subfigure}\hfill
\begin{subfigure}[t]{1.5\figurewidth}
	\centering
	\tikzsetnextfilename{cpu_time_aniso}
	\input{cpu_time}
	\caption{Anisotropic case}\label{plotcputime_aniso}
\end{subfigure}
\caption{Evolutions of (\subref{ploterrortol_iso})-(\subref{ploterrortol_aniso}) error indicator $\varepsilon_{\Xi,\Omega\setminus\Lambda}$ and (\subref{plotcputime_iso})-(\subref{plotcputime_aniso}) CPU time per iteration with respect to iteration number $k$ for different cross-validation tolerances $\varepsilon_{\text{cv}}$}\label{ploterrorcputimetol}
\end{figure}
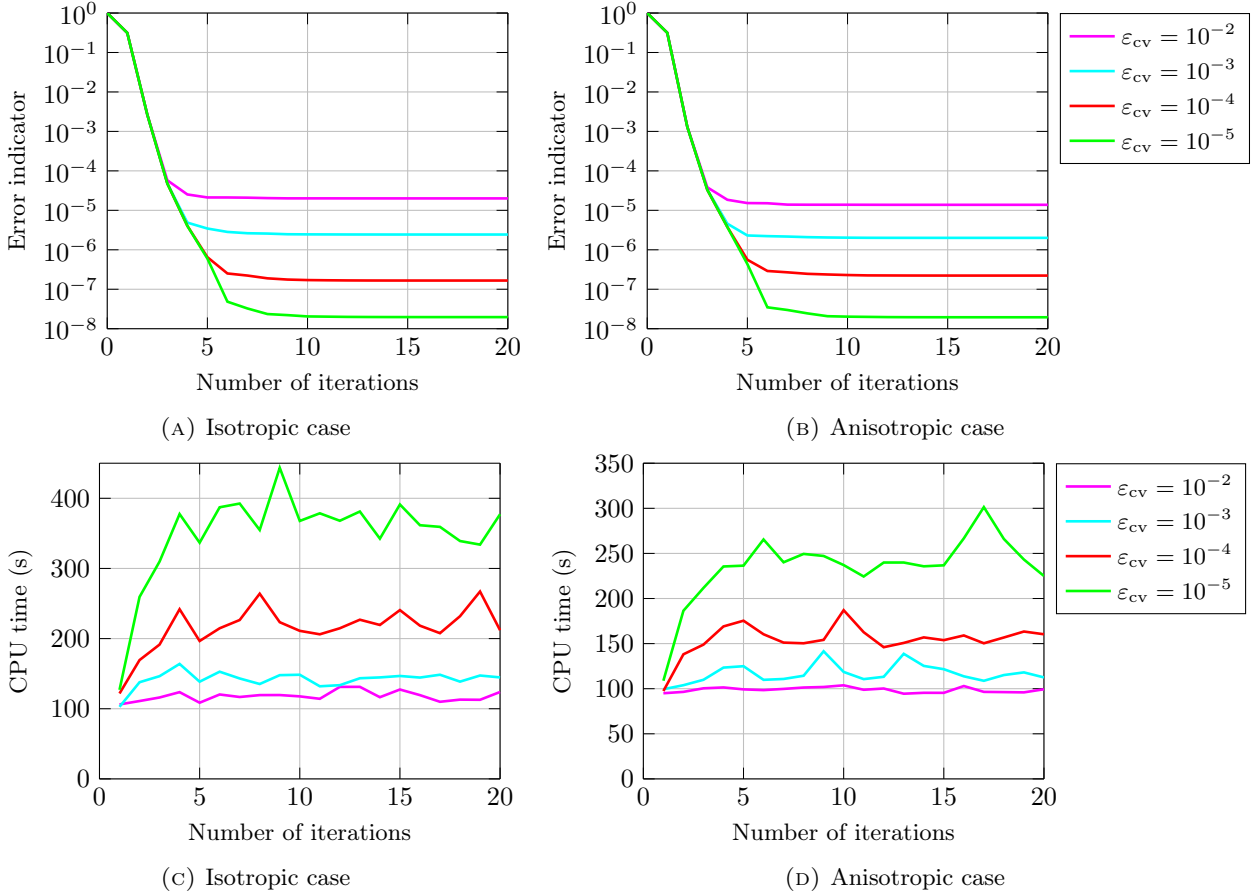

\setlength\figureheight{0.2\textheight}
\setlength\figurewidth{0.4\textwidth}
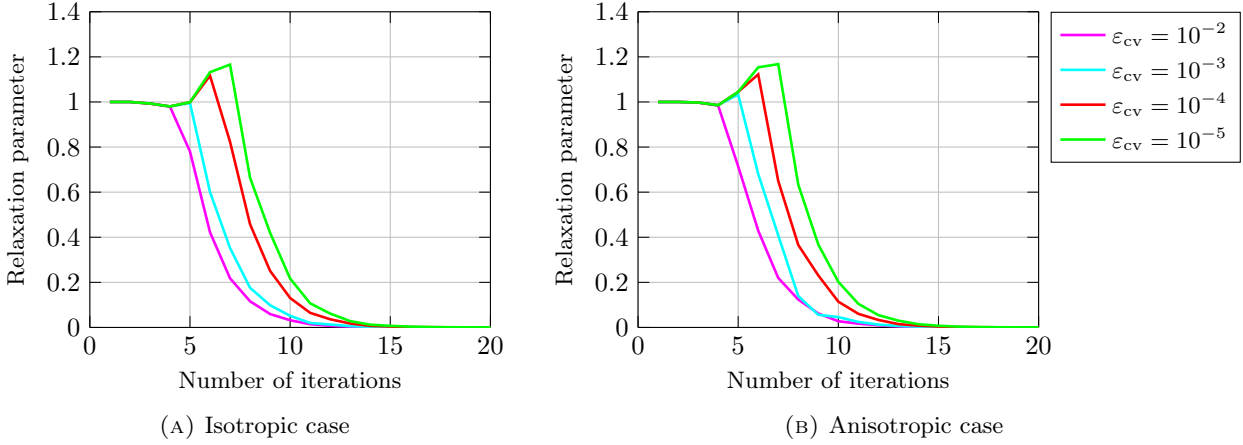
\begin{figure}[h!]
\centering
\definecolor{mycolor1}{rgb}{1.00000,0.00000,1.00000}%
\definecolor{mycolor2}{rgb}{0.00000,1.00000,1.00000}%
\begin{subfigure}[t]{\figurewidth}
	\centering
	\tikzsetnextfilename{relaxation_parameter_iso}
	\input{relaxation_parameter_no_legend}
	\caption{Isotropic case}\label{plotrhotol_iso}
\end{subfigure}\hfill
\begin{subfigure}[t]{1.5\figurewidth}
	\centering
	\tikzsetnextfilename{relaxation_parameter_aniso}
	\input{relaxation_parameter}
	\caption{Anisotropic case}\label{plotrhotol_aniso}
\end{subfigure}
\caption{Evolutions of Aitken's dynamic relaxation parameter $\rho_k$ with respect to iteration number $k$ for different cross-validation tolerances $\varepsilon_{\text{cv}}$}\label{plotrhotol}
\end{figure}

\setlength\figureheight{0.2\textheight}
\setlength\figurewidth{0.4\textwidth}
\begin{figure}[h!]
\centering
\definecolor{mycolor1}{rgb}{1.00000,0.00000,1.00000}%
\definecolor{mycolor2}{rgb}{0.00000,1.00000,1.00000}%
\begin{subfigure}[t]{\figurewidth}
	\centering
	\tikzsetnextfilename{nb_samples_4_tol_iso}
	\input{nb_samples_4_tol_no_legend}
	\caption{Isotropic case}\label{plotnbsamples4tol_iso}
\end{subfigure}\hfill
\begin{subfigure}[t]{1.5\figurewidth}
	\centering
	\tikzsetnextfilename{nb_samples_4_tol_aniso}
	\input{nb_samples_4_tol}
	\caption{Anisotropic case}\label{plotnbsamples4tol_aniso}
\end{subfigure}
\begin{subfigure}[t]{\figurewidth}
	\centering
	\tikzsetnextfilename{dim_stochastic_basis_4_tol_iso}
	\input{dim_stochastic_basis_4_tol_no_legend}
	\caption{Isotropic case}\label{plotdimstochasticbasis4tol_iso}
\end{subfigure}\hfill
\begin{subfigure}[t]{1.5\figurewidth}
	\centering
	\tikzsetnextfilename{dim_stochastic_basis_4_tol_aniso}
	\input{dim_stochastic_basis_4_tol}
	\caption{Anisotropic case}\label{plotdimstochasticbasis4tol_aniso}
\end{subfigure}
\caption{Evolutions of (\subref{plotnbsamples4tol_iso})-(\subref{plotnbsamples4tol_aniso}) the number of samples $N_4$ and (\subref{plotdimstochasticbasis4tol_iso})-(\subref{plotdimstochasticbasis4tol_aniso}) the dimension $\#\Ac_4$ of the approximation basis for local solution $w^k_4$ (solid lines) and Lagrange multiplier $\lambda^k_4$ (dashed lines) within patch $\Lambda_4$ with respect to iteration number $k$ for different cross-validation tolerances $\varepsilon_{\text{cv}}$}\label{plotnbsamplesdimstochasticbasis4tol}
\end{figure}
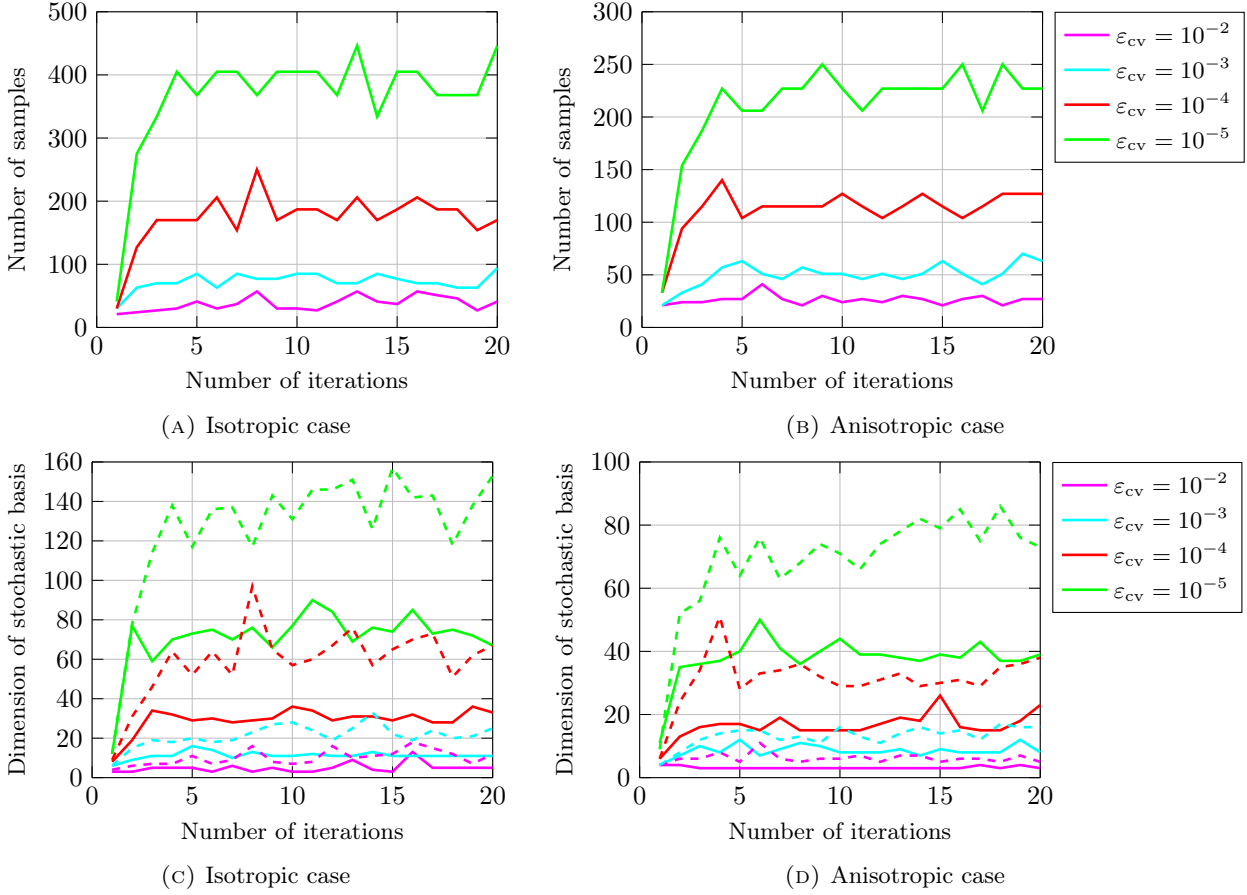

In order to illustrate the capability of the adaptive least-squares solver given in Algorithm~\ref{LSalgo} to capture sparse high-dimensional polynomial approximations of local solutions, Table~\ref{degree_dim} shows the partial polynomial degrees $p_i$ with respect to each random variable $\xi_i$, $i\in\set{1,\dots,m}$, and the dimension $\#\Ac$ of approximation spaces $\Sc_\Ac$ for global and local solutions $U$ and $(w^q,\lambda^q)$, $q\in\set{1,\dots,Q}$, obtained at convergence of the global-local iterative algorithm and using a fixed cross-validation tolerance $\varepsilon_{\text{cv}} = 10^{-3}$ for the convergence of Algorithm~\ref{LSalgo}. We observe that the use of the adaptive sparse least-squares solver allows to detect sparsity in local solutions $(w^q,\lambda^q)$. Indeed, Algorithm~\ref{LSalgo} gives local solutions $(w^q,\lambda^q)$ with a very low effective dimensionality in so far as they are mainly dependent on only few random variables, especially the random variables $\xi_{2q-1}$ and $\xi_{2q}$ associated with patch $\Lambda_q$.

\begin{table}[h!]
\footnotesize{
\begin{subtable}[t]{\textwidth}
\centering
\begin{tabular}{cccccccccccccccccc}
\toprule
& $p_1$ & $p_2$ & $p_3$ & $p_4$ & $p_5$ & $p_6$ & $p_7$ & $p_8$ & $p_9$ & $p_{10}$ & $p_{11}$ & $p_{12}$ & $p_{13}$ & $p_{14}$ & $p_{15}$ & $p_{16}$ & $\#\Ac$\\
\midrule
\multicolumn{1}{c}{$U$} & $5$ & $2$ & $4$ & $3$ & $4$ & $2$ & $4$ & $3$ & $4$ & $3$ & $4$ & $3$ & $4$ & $2$ & $5$ & $2$ & $192$\\
\midrule
\multicolumn{1}{c}{$w_1$} & \textcolor{red}{$3$} & \textcolor{red}{$1$} & $1$ & $1$ & $1$ & $0$ & $0$ & $0$ & $0$ & $0$ & $0$ & $0$ & $0$ & $0$ & $0$ & $0$ & $9$\\
\multicolumn{1}{c}{$\lambda_1$} & \textcolor{red}{$4$} & \textcolor{red}{$2$} & $2$ & $2$ & $1$ & $1$ & $0$ & $0$ & $0$ & $0$ & $1$ & $0$ & $0$ & $0$ & $0$ & $0$ & $21$\\
\midrule
\multicolumn{1}{c}{$w_2$} & $2$ & $1$ & \textcolor{red}{$3$} & \textcolor{red}{$2$} & $1$ & $1$ & $0$ & $0$ & $0$ & $0$ & $0$ & $0$ & $0$ & $0$ & $0$ & $0$ & $13$\\
\multicolumn{1}{c}{$\lambda_2$} & $3$ & $1$ & \textcolor{red}{$3$} & \textcolor{red}{$2$} & $2$ & $2$ & $0$ & $1$ & $0$ & $0$ & $0$ & $0$ & $0$ & $0$ & $0$ & $0$ & $19$\\
\midrule
\multicolumn{1}{c}{$w_3$} & $0$ & $0$ & $1$ & $1$ & \textcolor{red}{$3$} & \textcolor{red}{$2$} & $1$ & $1$ & $0$ & $0$ & $0$ & $0$ & $0$ & $0$ & $0$ & $0$ & $11$\\
\multicolumn{1}{c}{$\lambda_3$} & $1$ & $1$ & $2$ & $1$ & \textcolor{red}{$3$} & \textcolor{red}{$2$} & $2$ & $1$ & $0$ & $1$ & $0$ & $0$ & $0$ & $0$ & $1$ & $1$ & $21$\\
\midrule
\multicolumn{1}{c}{$w_4$} & $0$ & $0$ & $0$ & $0$ & $1$ & $1$ & \textcolor{red}{$3$} & \textcolor{red}{$2$} & $1$ & $1$ & $0$ & $0$ & $0$ & $0$ & $0$ & $0$ & $11$\\
\multicolumn{1}{c}{$\lambda_4$} & $0$ & $0$ & $0$ & $1$ & $2$ & $2$ & \textcolor{red}{$4$} & \textcolor{red}{$2$} & $2$ & $1$ & $1$ & $1$ & $0$ & $1$ & $1$ & $0$ & $25$\\
\midrule
\multicolumn{1}{c}{$w_5$} & $0$ & $0$ & $0$ & $0$ & $0$ & $0$ & $1$ & $1$ & \textcolor{red}{$3$} & \textcolor{red}{$2$} & $1$ & $1$ & $0$ & $0$ & $0$ & $0$ & $11$\\
\multicolumn{1}{c}{$\lambda_5$} & $0$ & $0$ & $0$ & $0$ & $0$ & $1$ & $1$ & $2$ & \textcolor{red}{$3$} & \textcolor{red}{$2$} & $2$ & $1$ & $1$ & $1$ & $0$ & $0$ & $20$\\
\midrule
\multicolumn{1}{c}{$w_6$} & $0$ & $0$ & $0$ & $0$ & $0$ & $0$ & $0$ & $0$ & $1$ & $1$ & \textcolor{red}{$3$} & \textcolor{red}{$2$} & $1$ & $1$ & $0$ & $0$ & $11$\\
\multicolumn{1}{c}{$\lambda_6$} &$0$ & $0$ & $1$ & $0$ & $0$ & $0$ & $1$ & $1$ & $2$ & $2$ & \textcolor{red}{$4$} & \textcolor{red}{$2$} & $2$ & $2$ & $1$ & $1$ & $28$\\
\midrule
\multicolumn{1}{c}{$w_7$} & $0$ & $0$ & $0$ & $0$ & $0$ & $0$ & $0$ & $0$ & $0$ & $0$ & $1$ & $1$ & \textcolor{red}{$3$} & \textcolor{red}{$2$} & $2$ & $1$ & $13$\\
\multicolumn{1}{c}{$\lambda_7$} & $0$ & $0$ & $0$ & $0$ & $0$ & $0$ & $0$ & $0$ & $1$ & $1$ & $2$ & $2$ & \textcolor{red}{$3$} & \textcolor{red}{$2$} & $3$ & $1$ & $20$\\
\midrule
\multicolumn{1}{c}{$w_8$} & $0$ & $1$ & $1$ & $0$ & $0$ & $0$ & $0$ & $0$ & $1$ & $0$ & $1$ & $0$ & $1$ & $1$ & \textcolor{red}{$3$} & \textcolor{red}{$1$} & $12$\\
\multicolumn{1}{c}{$\lambda_8$} & $0$ & $1$ & $1$ & $0$ & $0$ & $0$ & $0$ & $0$ & $1$ & $1$ & $1$ & $1$ & $3$ & $2$ & \textcolor{red}{$5$} & \textcolor{red}{$2$} & $26$\\
\bottomrule
\end{tabular}
\caption{Isotropic case}\label{degree_dim_iso}
\end{subtable}
\\
\begin{subtable}[t]{\textwidth}
\centering
\begin{tabular}{cccccccccccccccccc}
\toprule
& $p_1$ & $p_2$ & $p_3$ & $p_4$ & $p_5$ & $p_6$ & $p_7$ & $p_8$ & $p_9$ & $p_{10}$ & $p_{11}$ & $p_{12}$ & $p_{13}$ & $p_{14}$ & $p_{15}$ & $p_{16}$ & $\#\Ac$\\
\midrule
\multicolumn{1}{c}{$U$} & $5$ & $3$ & $3$ & $2$ & $3$ & $2$ & $3$ & $2$ & $3$ & $2$ & $3$ & $1$ & $2$ & $3$ & $2$ & $1$ & $112$\\
\midrule
\multicolumn{1}{c}{$w_1$} & \textcolor{red}{$3$} & \textcolor{red}{$1$} & $1$ & $1$ & $0$ & $0$ & $0$ & $0$ & $0$ & $0$ & $0$ & $0$ & $0$ & $0$ & $0$ & $0$ & $8$\\
\multicolumn{1}{c}{$\lambda_1$} & \textcolor{red}{$4$} & \textcolor{red}{$2$} & $2$ & $1$ & $1$ & $1$ & $1$ & $0$ & $0$ & $0$ & $0$ & $0$ & $0$ & $0$ & $0$ & $0$ & $18$\\
\midrule
\multicolumn{1}{c}{$w_2$} & $1$ & $1$ & \textcolor{red}{$2$} & \textcolor{red}{$1$} & $1$ & $1$ & $0$ & $0$ & $0$ & $0$ & $0$ & $0$ & $0$ & $0$ & $0$ & $0$ & $9$\\
\multicolumn{1}{c}{$\lambda_2$} & $2$ & $1$ & \textcolor{red}{$3$} & \textcolor{red}{$2$} & $1$ & $1$ & $0$ & $1$ & $0$ & $0$ & $0$ & $0$ & $0$ & $0$ & $0$ & $0$ & $15$\\
\midrule
\multicolumn{1}{c}{$w_3$} & $1$ & $0$ & $1$ & $1$ & \textcolor{red}{$2$} & \textcolor{red}{$1$} & $1$ & $2$ & $0$ & $0$ & $0$ & $0$ & $0$ & $0$ & $0$ & $0$ & $13$\\
\multicolumn{1}{c}{$\lambda_3$} & $1$ & $1$ & $2$ & $1$ & \textcolor{red}{$3$} & \textcolor{red}{$2$} & $1$ & $2$ & $0$ & $1$ & $0$ & $0$ & $0$ & $0$ & $0$ & $0$ & $18$\\
\midrule
\multicolumn{1}{c}{$w_4$} & $0$ & $0$ & $0$ & $0$ & $1$ & $1$ & \textcolor{red}{$2$} & \textcolor{red}{$1$} & $1$ & $1$ & $0$ & $0$ & $0$ & $0$ & $0$ & $0$ & $8$\\
\multicolumn{1}{c}{$\lambda_4$} & $0$ & $0$ & $0$ & $1$ & $1$ & $1$ & \textcolor{red}{$2$} & \textcolor{red}{$2$} & $1$ & $2$ & $0$ & $1$ & $0$ & $0$ & $1$ & $0$ & $16$\\
\midrule
\multicolumn{1}{c}{$w_5$} & $0$ & $0$ & $0$ & $0$ & $0$ & $0$ & $1$ & $1$ & \textcolor{red}{$2$} & \textcolor{red}{$1$} & $0$ & $1$ & $0$ & $0$ & $0$ & $0$ & $7$\\
\multicolumn{1}{c}{$\lambda_5$} & $0$ & $0$ & $0$ & $0$ & $0$ & $1$ & $1$ & $1$ & \textcolor{red}{$2$} & \textcolor{red}{$2$} & $1$ & $1$ & $0$ & $0$ & $0$ & $0$ & $11$\\
\midrule
\multicolumn{1}{c}{$w_6$} & $0$ & $0$ & $0$ & $0$ & $0$ & $0$ & $0$ & $0$ & $1$ & $1$ & \textcolor{red}{$2$} & \textcolor{red}{$1$} & $0$ & $1$ & $0$ & $0$ & $7$\\
\multicolumn{1}{c}{$\lambda_6$} & $0$ & $0$ & $0$ & $0$ & $0$ & $0$ & $0$ & $0$ & $1$ & $1$ & \textcolor{red}{$2$} & \textcolor{red}{$1$} & $1$ & $1$ & $0$ & $0$ & $8$\\
\midrule
\multicolumn{1}{c}{$w_7$} & $0$ & $0$ & $0$ & $0$ & $0$ & $0$ & $0$ & $0$ & $0$ & $0$ & $0$ & $1$ & \textcolor{red}{$1$} & \textcolor{red}{$1$} & $1$ & $0$ & $5$\\
\multicolumn{1}{c}{$\lambda_7$} & $0$ & $0$ & $0$ & $0$ & $0$ & $0$ & $0$ & $0$ & $0$ & $0$ & $1$ & $1$ & \textcolor{red}{$2$} & \textcolor{red}{$1$} & $1$ & $1$ & $8$\\
\midrule
\multicolumn{1}{c}{$w_8$} & $0$ & $0$ & $0$ & $0$ & $0$ & $0$ & $0$ & $0$ & $0$ & $0$ & $0$ & $0$ & $0$ & $1$ & \textcolor{red}{$1$} & \textcolor{red}{$1$} & $4$\\
\multicolumn{1}{c}{$\lambda_8$} & $0$ & $0$ & $0$ & $0$ & $0$ & $0$ & $0$ & $0$ & $0$ & $0$ & $0$ & $1$ & $1$ & $1$ & \textcolor{red}{$2$} & \textcolor{red}{$1$} & $7$\\
\bottomrule
\end{tabular}
\caption{Anisotropic case}\label{degree_dim_aniso}
\end{subtable}
\caption{Partial polynomial degrees $p_i$ (in each random variable $\xi_i$), $i \in \set{1,\dots,16}$, and dimension $\#\Ac$ of approximation spaces $\Sc_{\Ac}$ for global and local solutions $U$ and $(w_q,\lambda_q)$, $q \in \set{1,\dots,8}$, where values displayed in red correspond to random variables $\xi_{2q-1}$ and $\xi_{2q}$ associated with patch $\Lambda_q$}\label{degree_dim}
}
\end{table}

\subsection{Illustration of quantities of interest}

We now look at the effect of input uncertainties in material properties, namely diffusion coefficient $K$ and reaction parameter $R$, on the variability of the solution. We apply the Aitken's acceleration technique to the relaxation step of the global-local iterative algorithm and we set the cross-validation tolerance to $\varepsilon_{\text{cv}} = 10^{-3}$ in Algorithm~\ref{LSalgo}. We then consider the multiscale solution $u=(U,w_1,\dots,w_Q)$ obtained at final iteration of the algorithm. Figures~\ref{fig:mean_var_global_local_multiscale_solution_iso} and~\ref{fig:mean_var_global_local_multiscale_solution_aniso} show the mean and variance of global solution $U$ and local solutions $w_q$ as well as that of multiscale solution $u$, for isotropic and anisotropic cases, respectively. In the anisotropic case, the variability in the material properties of patch $\Lambda_q$ decreases with $q$ due to the anisotropy introduced in the weights $\gamma_q$. The highest spatial contributions to the variance $\xV(u)$ of multiscale solution $u$ are fully captured by the local solution $w_q$ within every patch $\Lambda_q$ and localized in the first patches in the anisotropic case. 

\setlength\figureheight{0.2\textheight}
\setlength\figurewidth{0.2\textwidth}
\begin{figure}[h!]
\centering
\begin{subfigure}[t]{\textwidth}
\centering
\begin{subfigure}[t]{1.5\figurewidth}
	\centering
	\includegraphics[height=\figureheight]{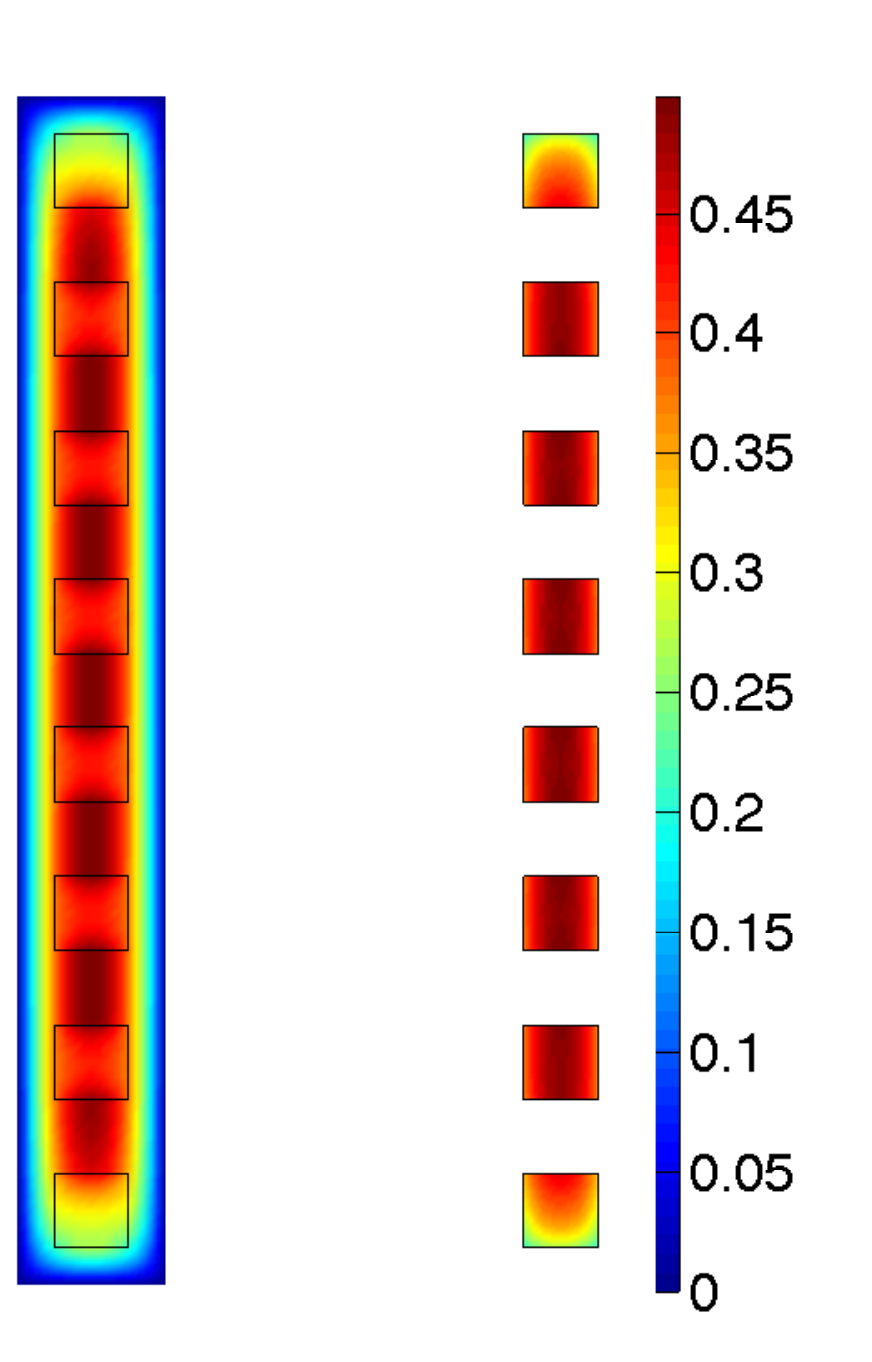}
	\caption*{$\xE(U)$ (left) and $\xE(w_q)$ (right)}%\label{subfig:mean_global_local_solution_iso}
\end{subfigure}\hfill
\begin{subfigure}[t]{\figurewidth}
	\centering
	\includegraphics[height=\figureheight]{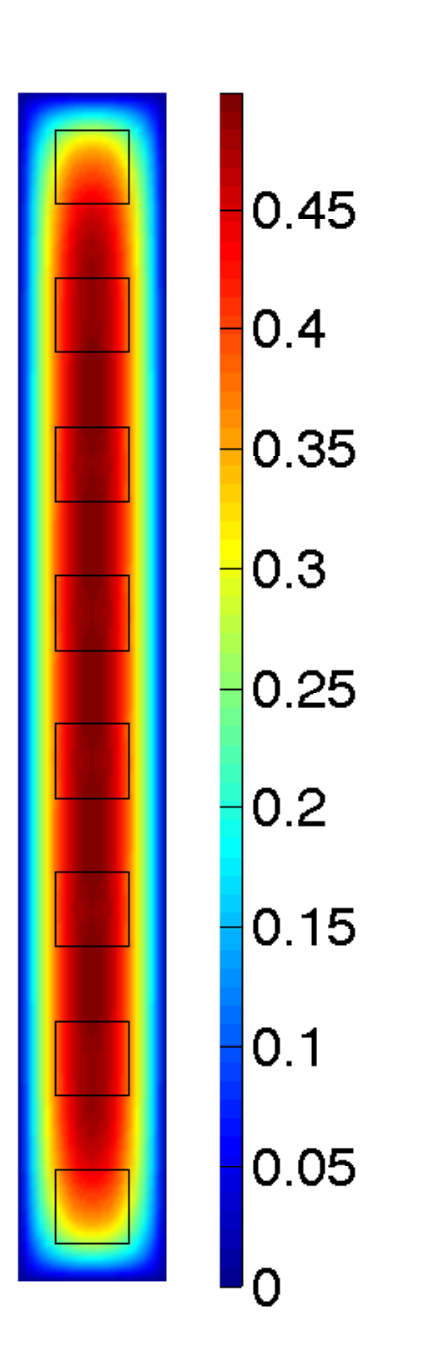}
	\caption*{$\xE(u)$}%\label{subfig:mean_multiscale_solution_iso}
\end{subfigure}\hfill
\begin{subfigure}[t]{1.5\figurewidth}
	\centering
	\includegraphics[height=\figureheight]{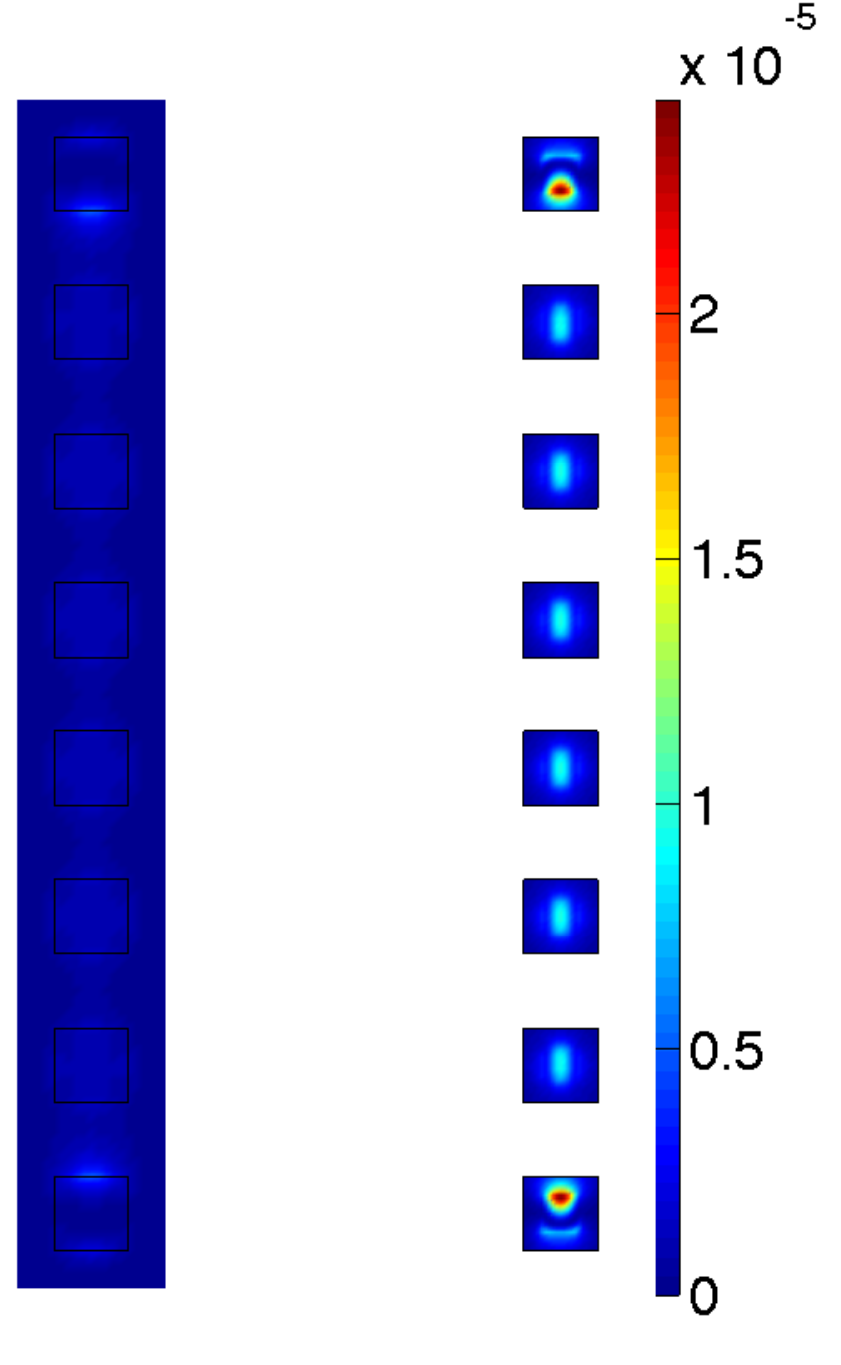}
	\caption*{$\xV(U)$ (left) and $\xV(w_q)$ (right)}%\label{subfig:var_global_local_solution_iso}
\end{subfigure}\hfill
\begin{subfigure}[t]{\figurewidth}
	\centering
	\includegraphics[height=\figureheight]{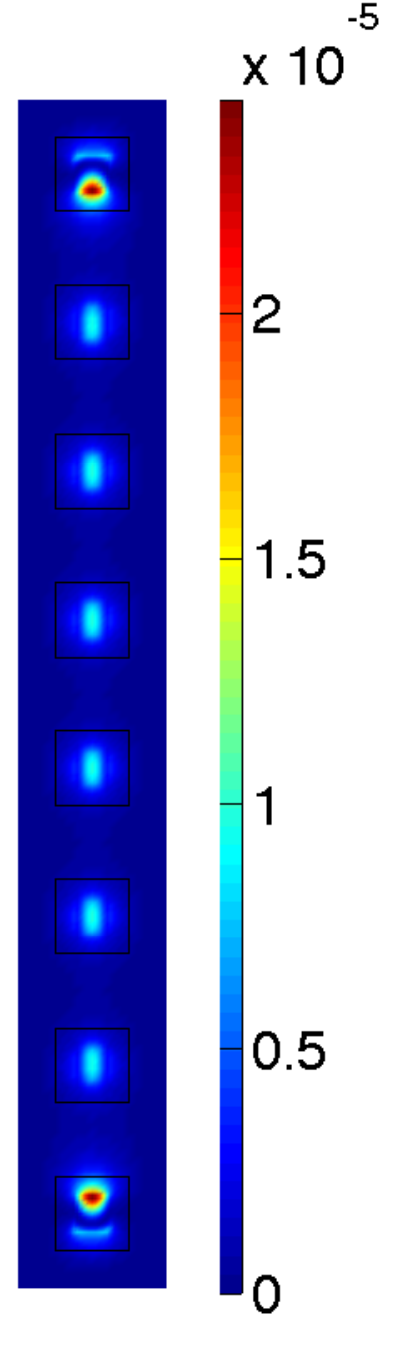}
	\caption*{$\xV(u)$}%\label{subfig:var_multiscale_solution_iso}
\end{subfigure}
\caption{Isotropic case}\label{fig:mean_var_global_local_multiscale_solution_iso}
\end{subfigure}

\begin{subfigure}[t]{\textwidth}
\centering
\begin{subfigure}[t]{1.5\figurewidth}
	\centering
	\includegraphics[height=\figureheight]{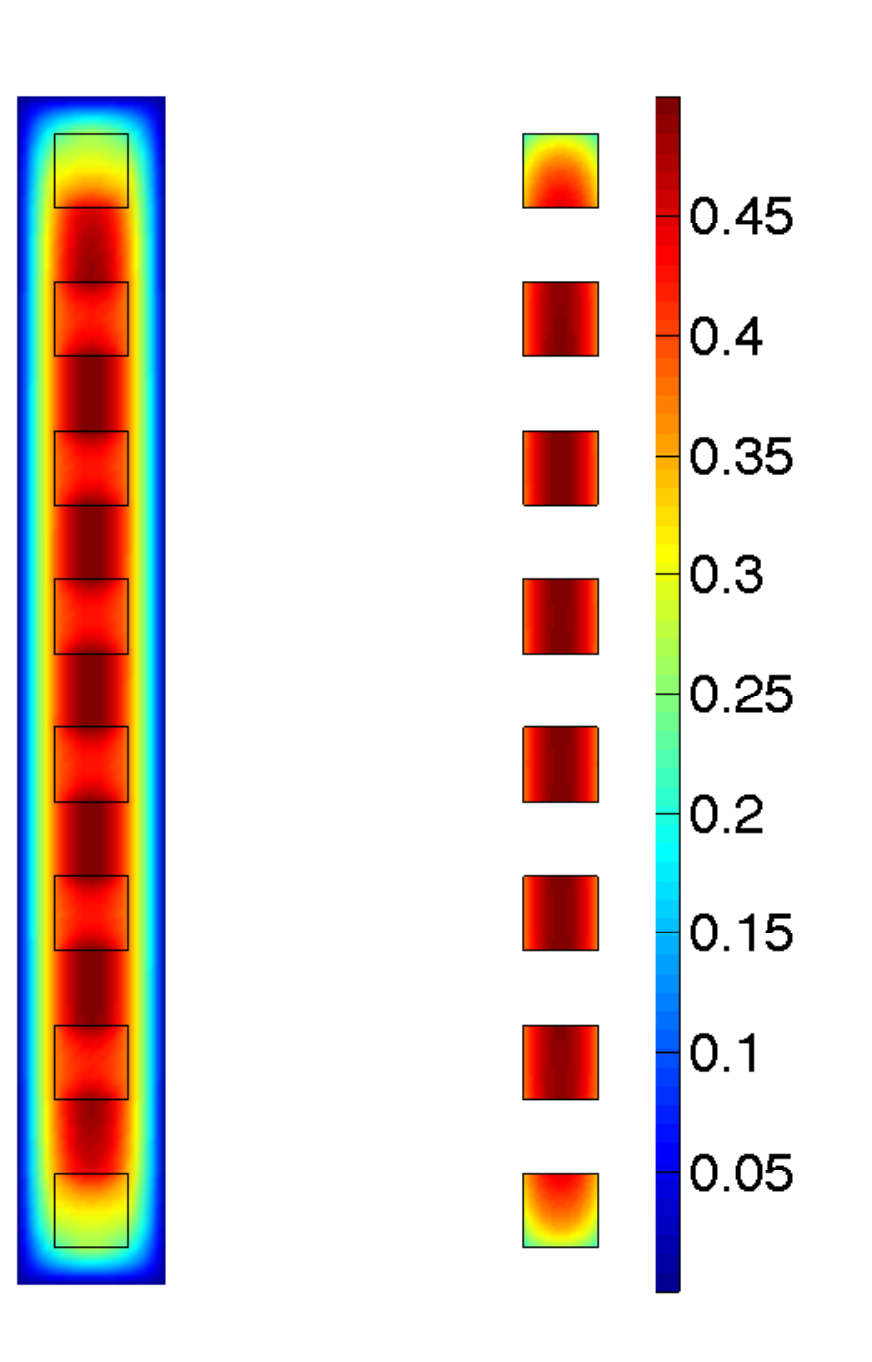}
	\caption*{$\xE(U)$ (left) and $\xE(w_q)$ (right)}%\label{subfig:mean_global_local_solution_aniso}
\end{subfigure}\hfill
\begin{subfigure}[t]{\figurewidth}
	\centering
	\includegraphics[height=\figureheight]{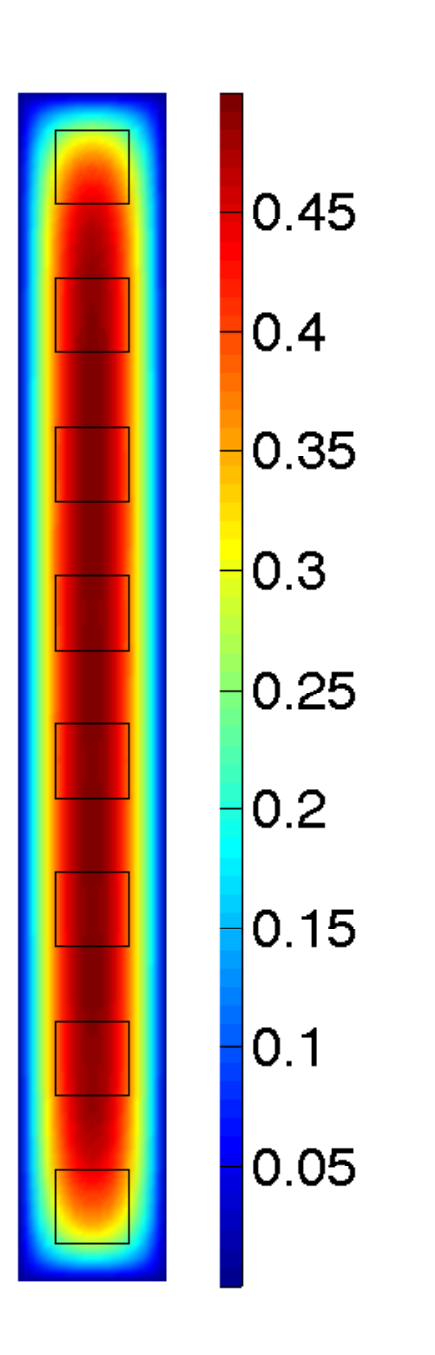}
	\caption*{$\xE(u)$}%\label{subfig:mean_multiscale_solution_aniso}
\end{subfigure}\hfill
\begin{subfigure}[t]{1.5\figurewidth}
	\centering
	\includegraphics[height=\figureheight]{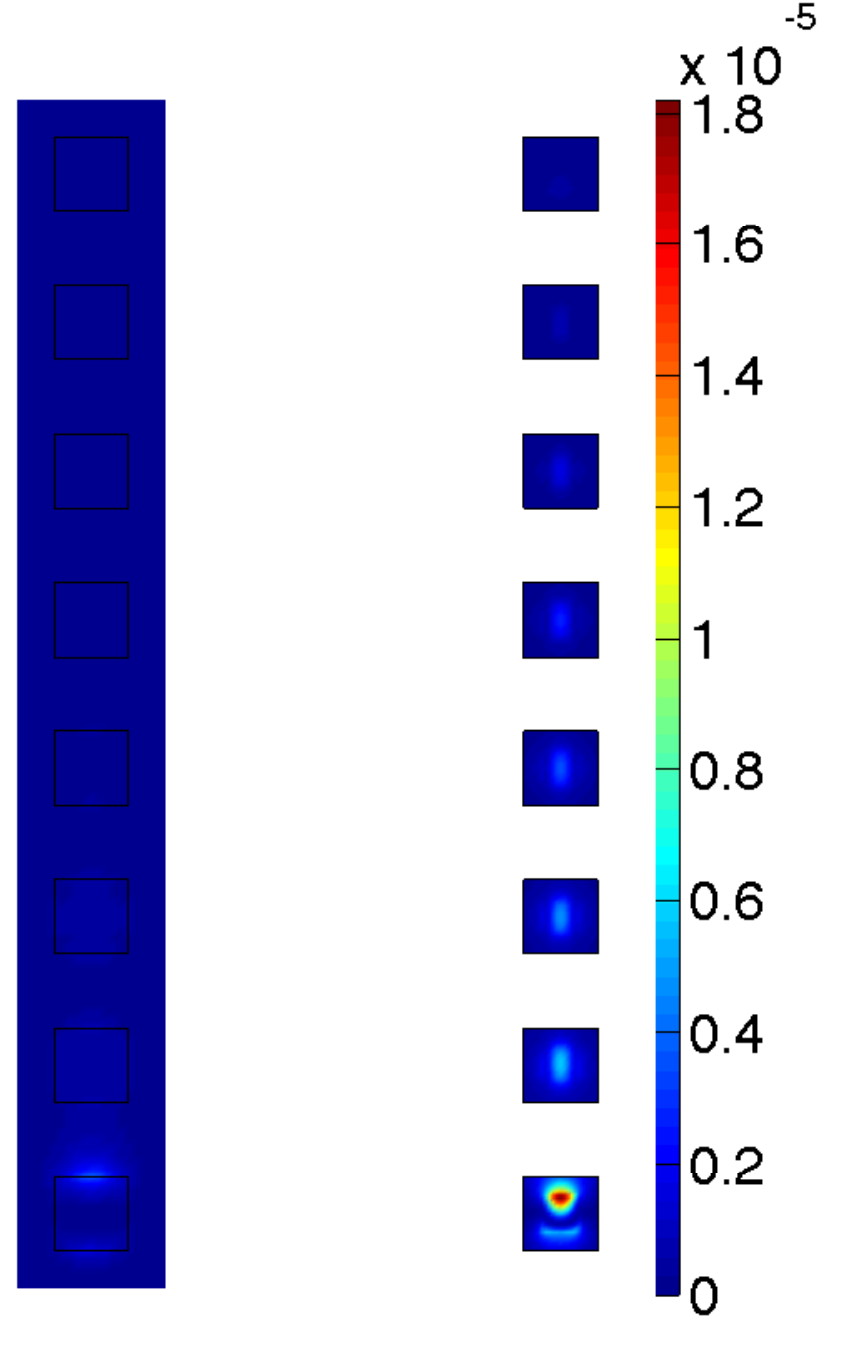}
	\caption*{$\xV(U)$ (left) and $\xV(w_q)$ (right)}%\label{subfig:var_global_local_solution_aniso}
\end{subfigure}\hfill
\begin{subfigure}[t]{\figurewidth}
	\centering
	\includegraphics[height=\figureheight]{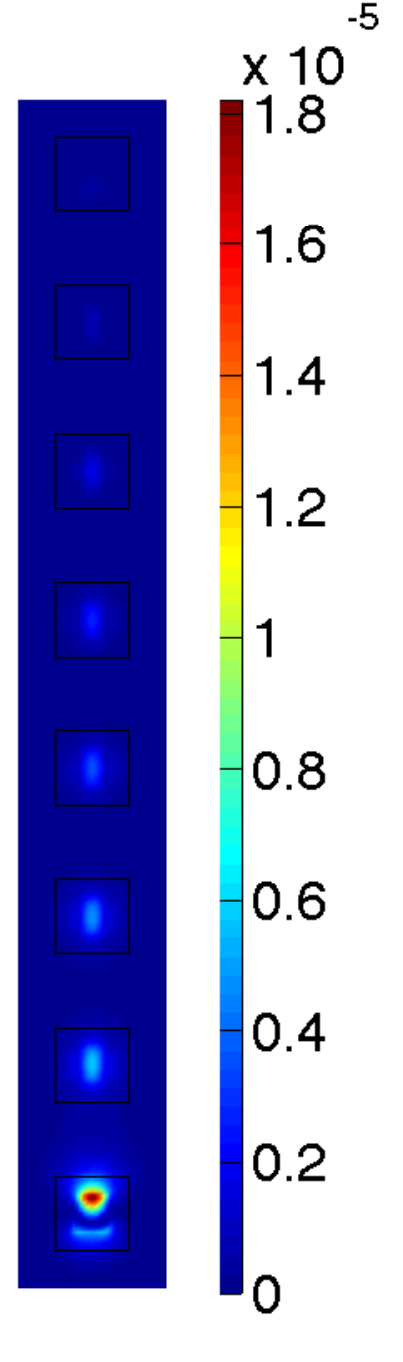}
	\caption*{$\xV(u)$}%\label{subfig:var_multiscale_solution_aniso}
\end{subfigure}
\caption{Anisotropic case}\label{fig:mean_var_global_local_multiscale_solution_aniso}
\end{subfigure}
\caption{Mean $\xE$ and variance $\xV$ of global solution $U$, local solutions $w_q$ and multiscale solution $u$}\label{fig:mean_var_global_local_multiscale_solution}
\end{figure}

In order to quantify the relative impact of each input random variable $\xi_i$ on the variability of solution $u$, we introduce the following global sensitivity indices:
\begin{equation*} \widetilde{S}_i(u) = \frac{\xV(\xE(u(x,\xi)\vert\xi_i))}{\max_{\substack{x \in \Omega}}(\xV(u(x,\xi)))},
\end{equation*}
where $\xE(u(x,\xi)\vert\xi_i)$ is the conditional expectation of solution $u$ with respect to random variable $\xi_i$. 
$\widetilde{S}_i(u)$ is a sensitivity index which reflects the zone of influence of a random variable $\xi_i$ (associated with patch $\Lambda_q$ for $i\in\set{2q-1,2q}$) on the variability of solution $u$. Note that global sensitivity indices $\widetilde{S}_i(u)$ can be straightforwardly computed from the expansion of $u(x,\xi)$ on an orthonormal polynomial basis (see \cite{Sud08}). Figure \ref{fig:sensitivity_indices_multiscale_solution} shows the spatial distributions of sensitivity indices $\widetilde{S}_i(u)$, computed at final iteration of the global-local iterative algorithm for all $i\in\set{1,\dots,m}$. We observe that random variables $\xi_{2q-1}$ and $\xi_{2q}$ have only a local influence within the corresponding patch $\Lambda_q$ on the variance of solution $u$. In the anisotropic case, the magnitude of sensitivity indices $\widetilde{S}_i(u)$ again reflects the highest input variabilities in the material properties within the first patches. Note that the patches $\Lambda_q$ are sufficiently large to capture the main effects of the input uncertainties on local solutions $w_q$, $q\in\set{1,\dots,Q}$.

\setlength\figureheight{0.16\textheight}
\setlength\figurewidth{0.12\textwidth}

\begin{figure}[h!]
\centering
\begin{subfigure}[t]{\textwidth}
\centering
\begin{subfigure}[t]{\figurewidth}
	\centering
	\includegraphics[height=\figureheight]{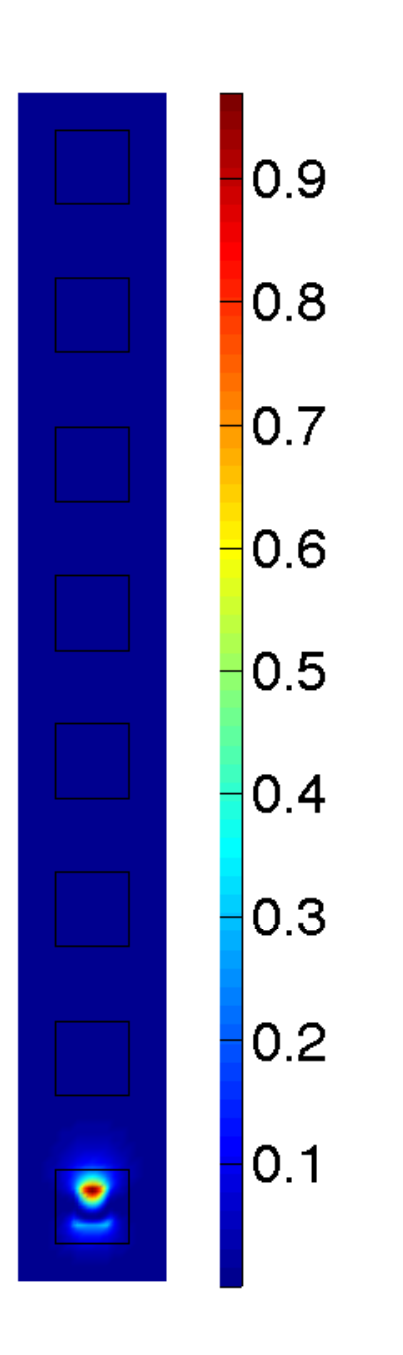}
	\caption*{$\widetilde{S}_1(u)$}%\label{subfig1:sensitivity_indices_multiscale_solution_var_iso}
\end{subfigure}\hfill
\begin{subfigure}[t]{\figurewidth}
	\centering
	\includegraphics[height=\figureheight]{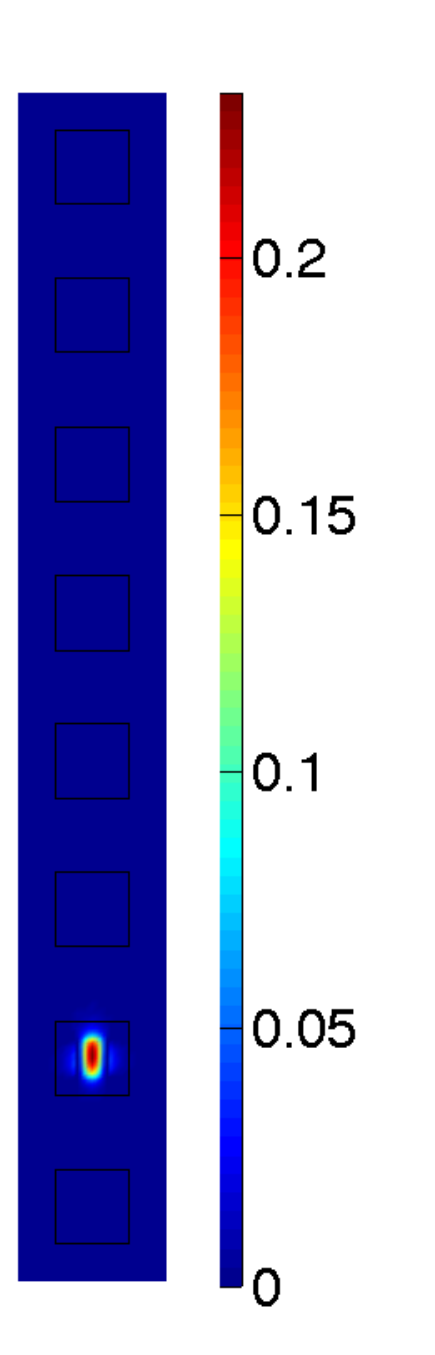}
	\caption*{$\widetilde{S}_3(u)$}%\label{subfig3:sensitivity_indices_multiscale_solution_var_iso}
\end{subfigure}\hfill
\begin{subfigure}[t]{\figurewidth}
	\centering
	\includegraphics[height=\figureheight]{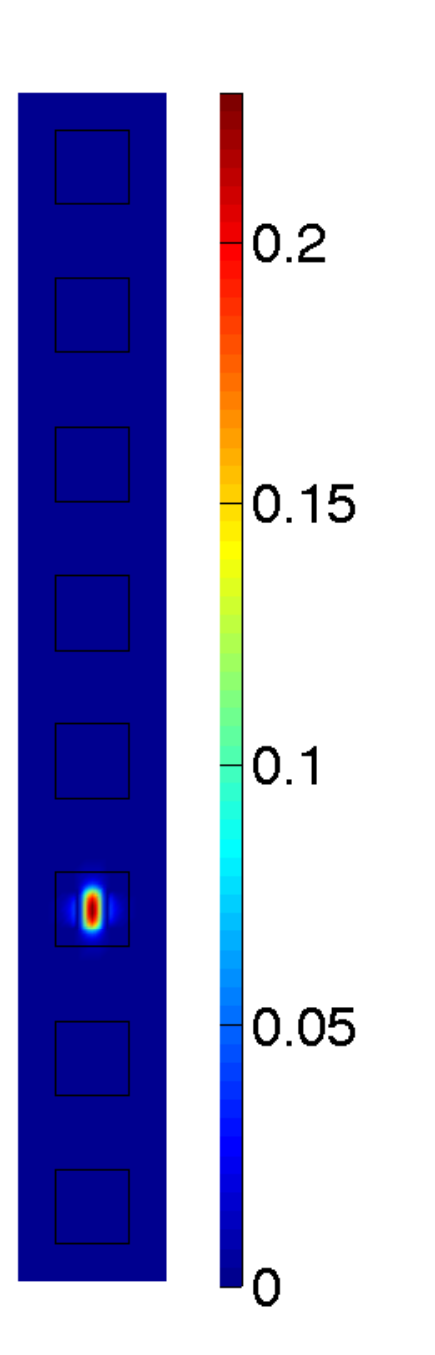}
	\caption*{$\widetilde{S}_5(u)$}%\label{subfig5:sensitivity_indices_multiscale_solution_var_iso}
\end{subfigure}\hfill
\begin{subfigure}[t]{\figurewidth}
	\centering
	\includegraphics[height=\figureheight]{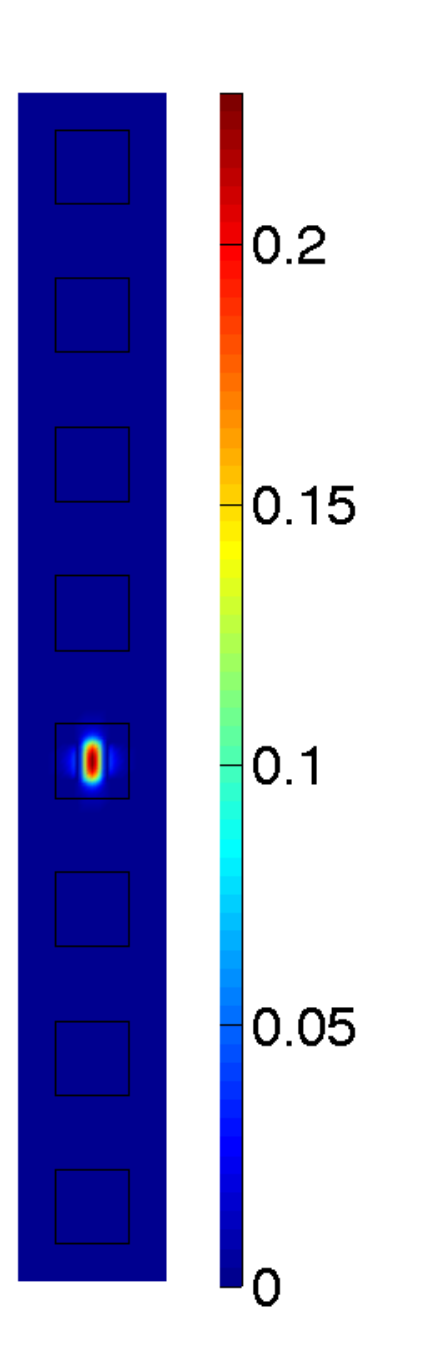}
	\caption*{$\widetilde{S}_7(u)$}%\label{subfig7:sensitivity_indices_multiscale_solution_var_iso}
\end{subfigure}\hfill
\begin{subfigure}[t]{\figurewidth}
	\centering
	\includegraphics[height=\figureheight]{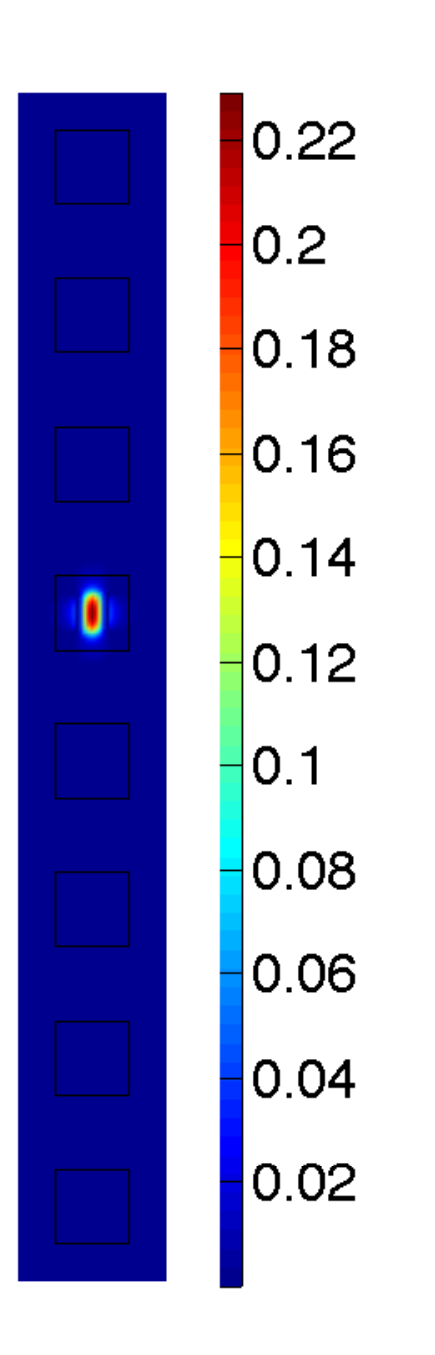}
	\caption*{$\widetilde{S}_9(u)$}%\label{subfig9:sensitivity_indices_multiscale_solution_var_iso}
\end{subfigure}\hfill
\begin{subfigure}[t]{\figurewidth}
	\centering
	\includegraphics[height=\figureheight]{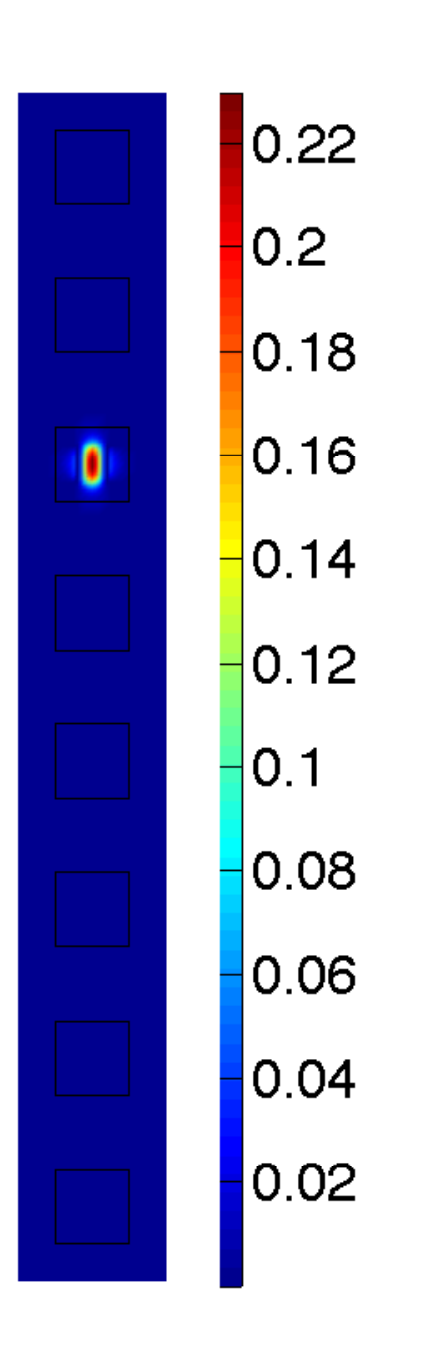}
	\caption*{$\widetilde{S}_{11}(u)$}%\label{subfig11:sensitivity_indices_multiscale_solution_var_iso}
\end{subfigure}\hfill
\begin{subfigure}[t]{\figurewidth}
	\centering
	\includegraphics[height=\figureheight]{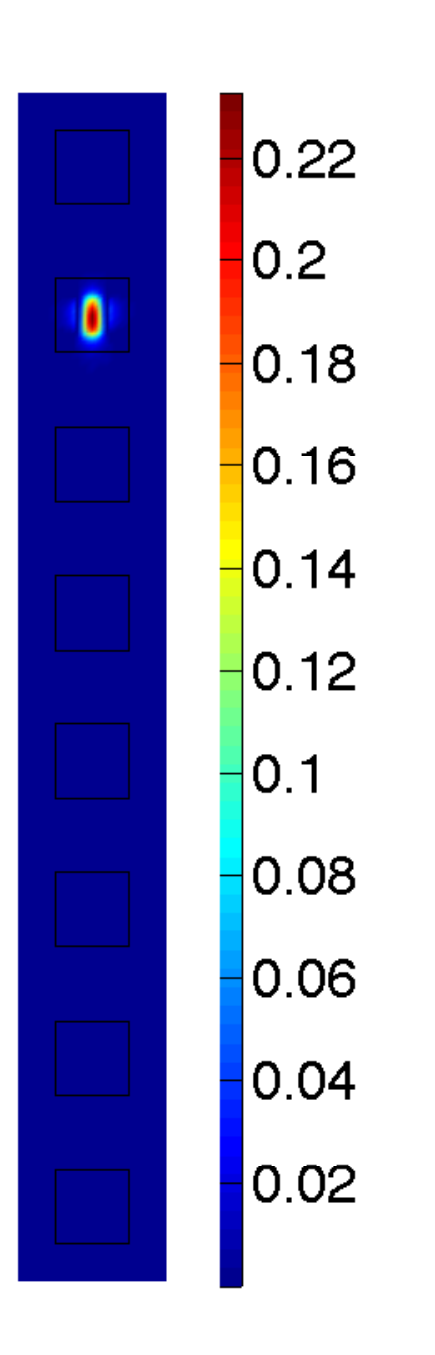}
	\caption*{$\widetilde{S}_{13}(u)$}%\label{subfig13:sensitivity_indices_multiscale_solution_var_iso}
\end{subfigure}\hfill
\begin{subfigure}[t]{\figurewidth}
	\centering
	\includegraphics[height=\figureheight]{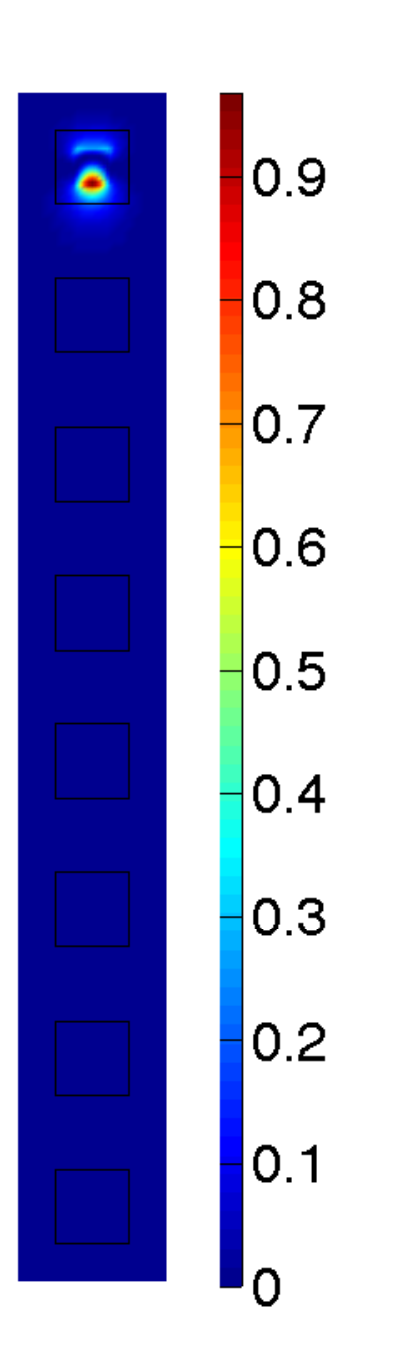}
	\caption*{$\widetilde{S}_{15}(u)$}%\label{subfig15:sensitivity_indices_multiscale_solution_var_iso}
\end{subfigure}

\begin{subfigure}[t]{\figurewidth}
	\centering
	\includegraphics[height=\figureheight]{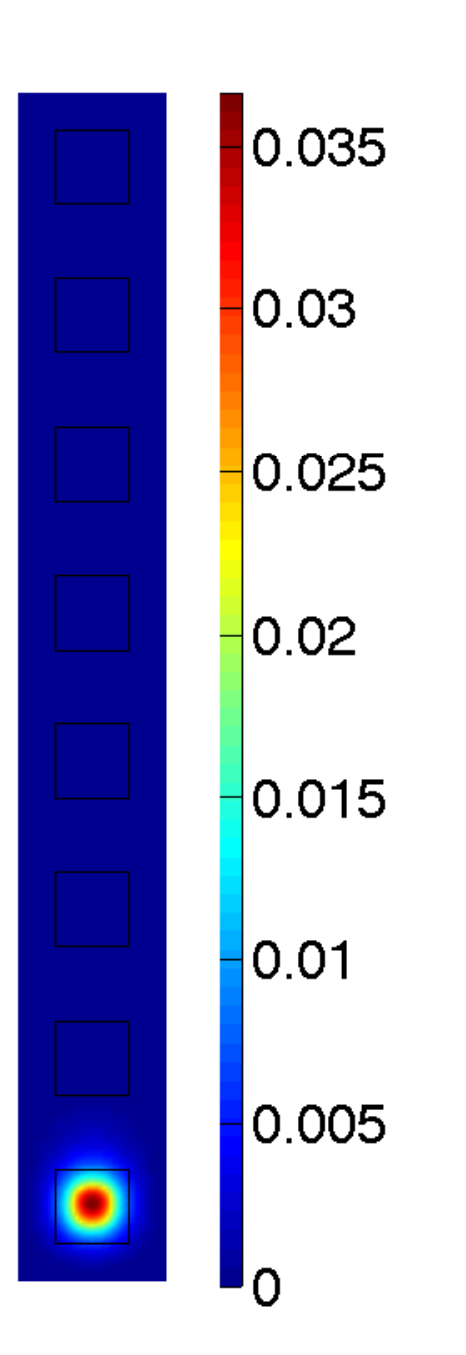}
	\caption*{$\widetilde{S}_2(u)$}%\label{subfig2:sensitivity_indices_multiscale_solution_var_iso}
\end{subfigure}\hfill
\begin{subfigure}[t]{\figurewidth}
	\centering
	\includegraphics[height=\figureheight]{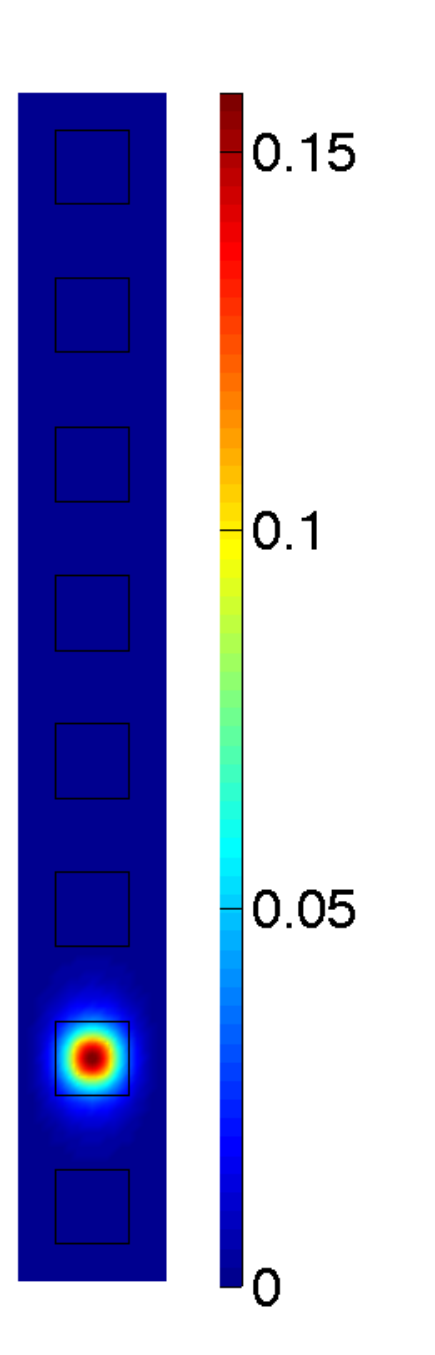}
	\caption*{$\widetilde{S}_4(u)$}%\label{subfig4:sensitivity_indices_multiscale_solution_var_iso}
\end{subfigure}\hfill
\begin{subfigure}[t]{\figurewidth}
	\centering
	\includegraphics[height=\figureheight]{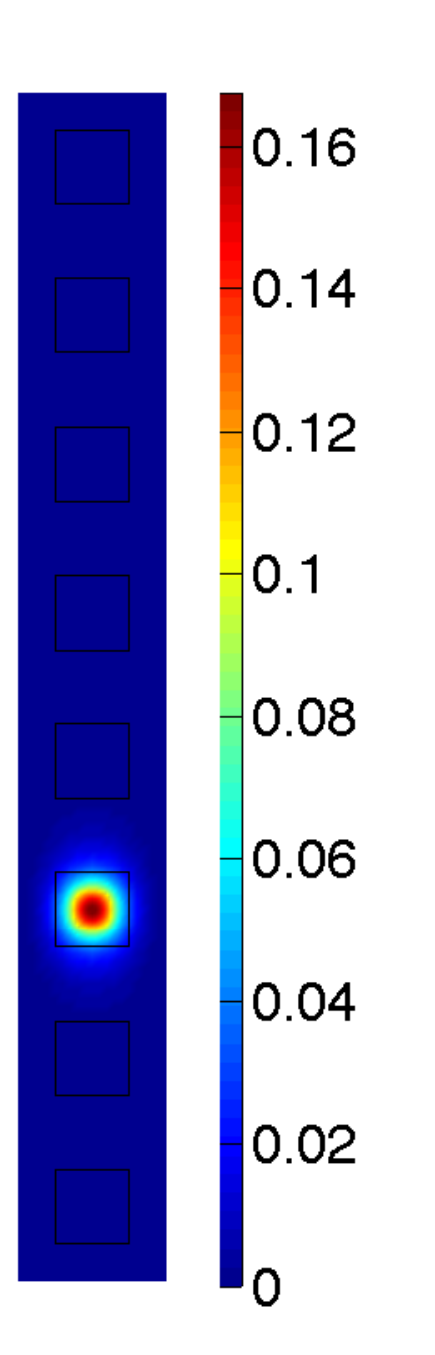}
	\caption*{$\widetilde{S}_6(u)$}%\label{subfig6:sensitivity_indices_multiscale_solution_var_iso}
\end{subfigure}\hfill
\begin{subfigure}[t]{\figurewidth}
	\centering
	\includegraphics[height=\figureheight]{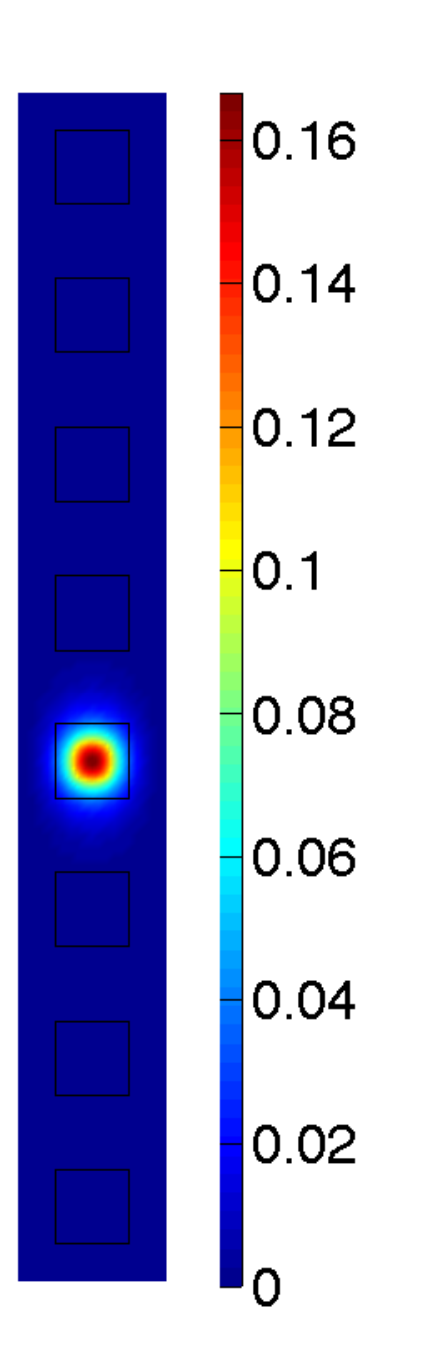}
	\caption*{$\widetilde{S}_8(u)$}%\label{subfig8:sensitivity_indices_multiscale_solution_var_iso}
\end{subfigure}\hfill
\begin{subfigure}[t]{\figurewidth}
	\centering
	\includegraphics[height=\figureheight]{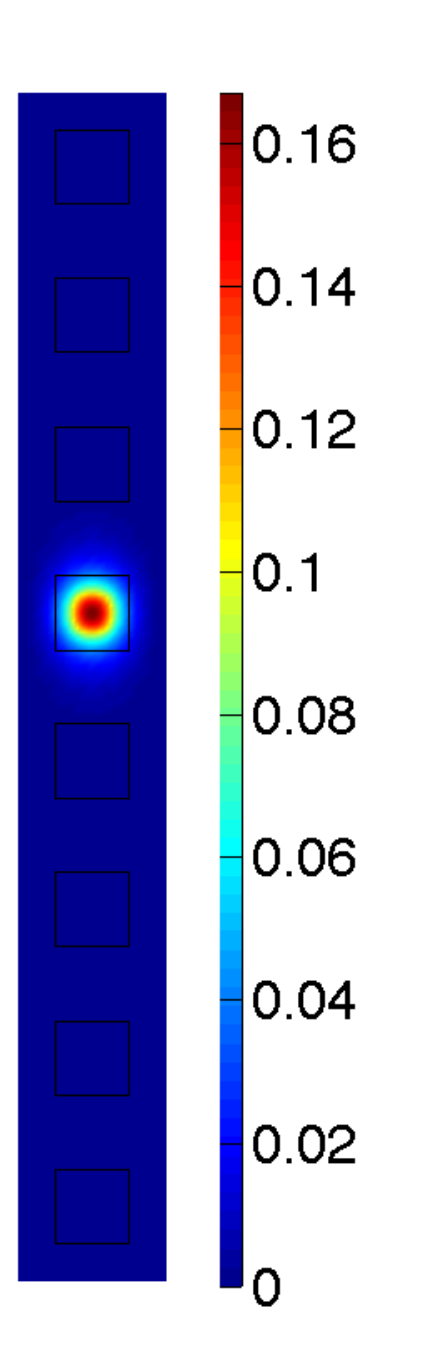}
	\caption*{$\widetilde{S}_{10}(u)$}%\label{subfig10:sensitivity_indices_multiscale_solution_var_iso}
\end{subfigure}\hfill
\begin{subfigure}[t]{\figurewidth}
	\centering
	\includegraphics[height=\figureheight]{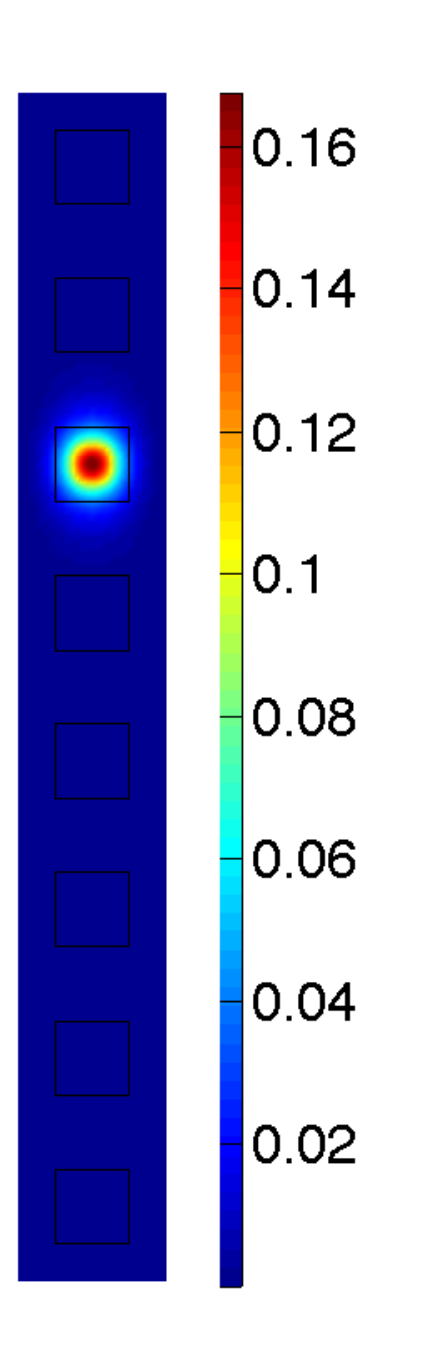}
	\caption*{$\widetilde{S}_{12}(u)$}%\label{subfig12:sensitivity_indices_multiscale_solution_var_iso}
\end{subfigure}\hfill
\begin{subfigure}[t]{\figurewidth}
	\centering
	\includegraphics[height=\figureheight]{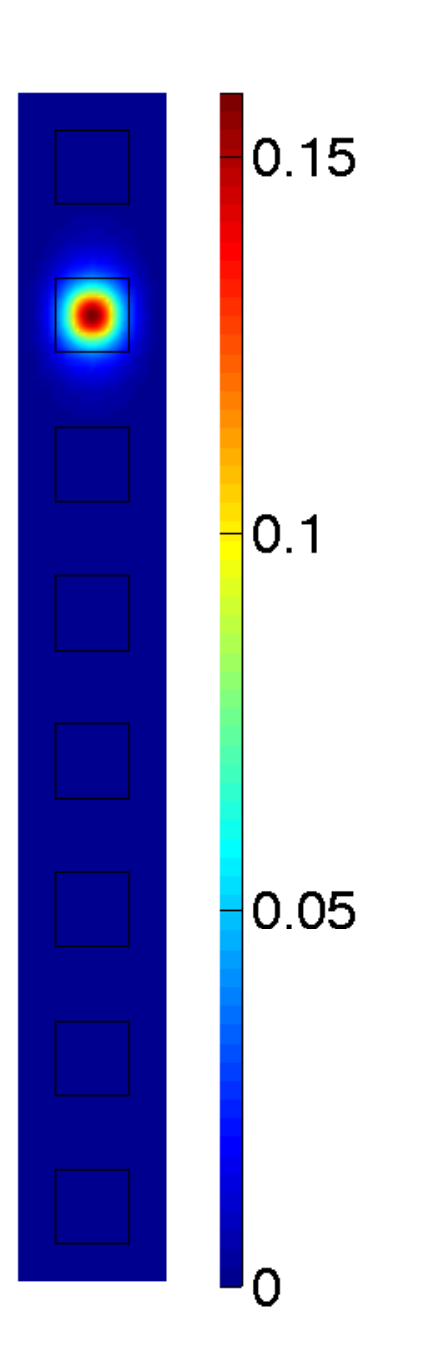}
	\caption*{$\widetilde{S}_{14}(u)$}%\label{subfig14:sensitivity_indices_multiscale_solution_var_iso}
\end{subfigure}\hfill
\begin{subfigure}[t]{\figurewidth}
	\centering
	\includegraphics[height=\figureheight]{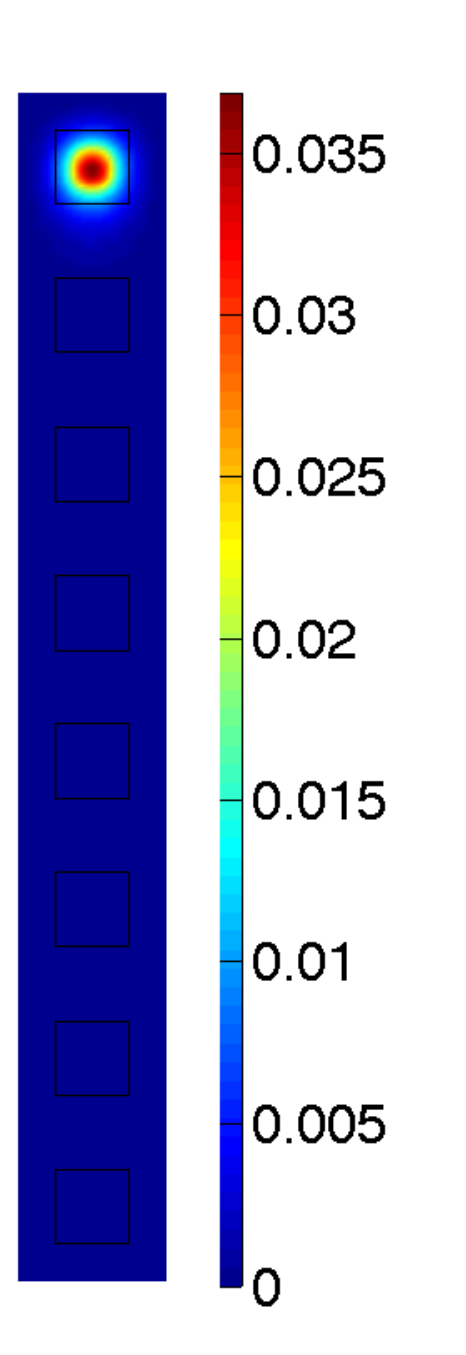}
	\caption*{$\widetilde{S}_{16}(u)$}%\label{subfig16:sensitivity_indices_multiscale_solution_var_iso}
\end{subfigure}
\caption{Isotropic case}\label{fig:sensitivity_indices_multiscale_solution_iso}
\end{subfigure}

\begin{subfigure}[t]{\textwidth}
\centering
\begin{subfigure}[t]{\figurewidth}
	\centering
	\includegraphics[height=\figureheight]{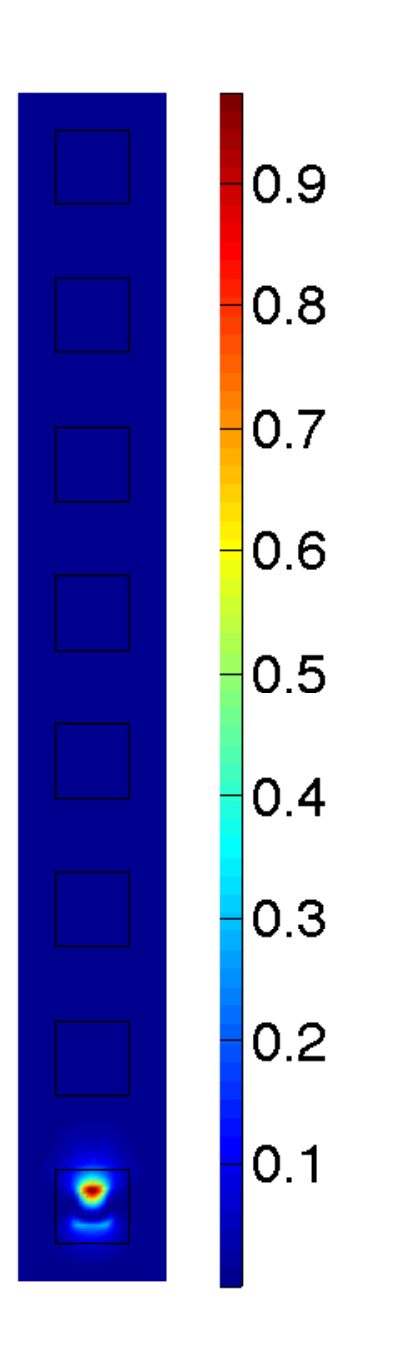}
	\caption*{$\widetilde{S}_1(u)$}%\label{subfig1:sensitivity_indices_multiscale_solution_var_aniso}
\end{subfigure}\hfill
\begin{subfigure}[t]{\figurewidth}
	\centering
	\includegraphics[height=\figureheight]{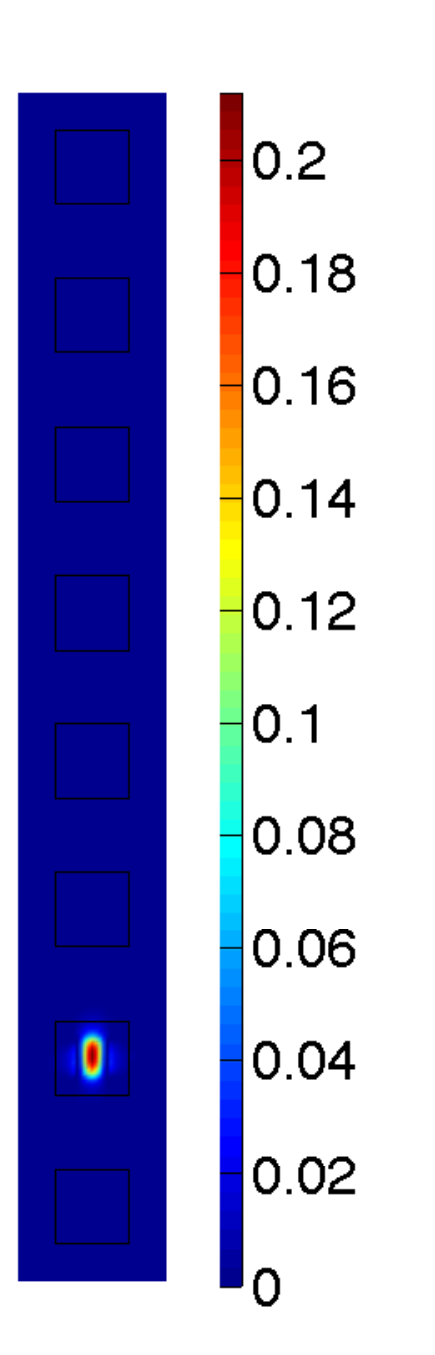}
	\caption*{$\widetilde{S}_3(u)$}%\label{subfig3:sensitivity_indices_multiscale_solution_var_aniso}
\end{subfigure}\hfill
\begin{subfigure}[t]{\figurewidth}
	\centering
	\includegraphics[height=\figureheight]{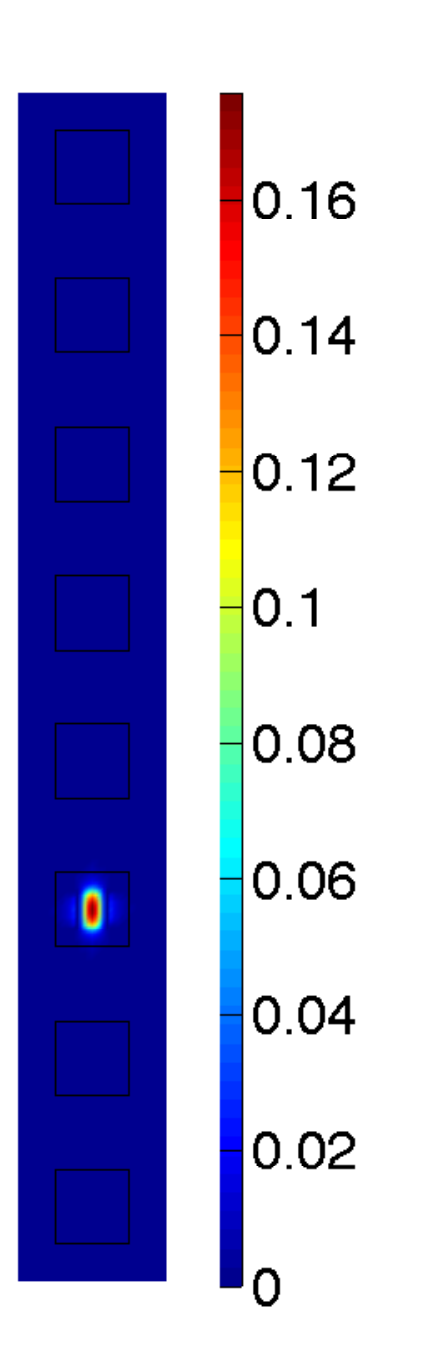}
	\caption*{$\widetilde{S}_5(u)$}%\label{subfig5:sensitivity_indices_multiscale_solution_var_aniso}
\end{subfigure}\hfill
\begin{subfigure}[t]{\figurewidth}
	\centering
	\includegraphics[height=\figureheight]{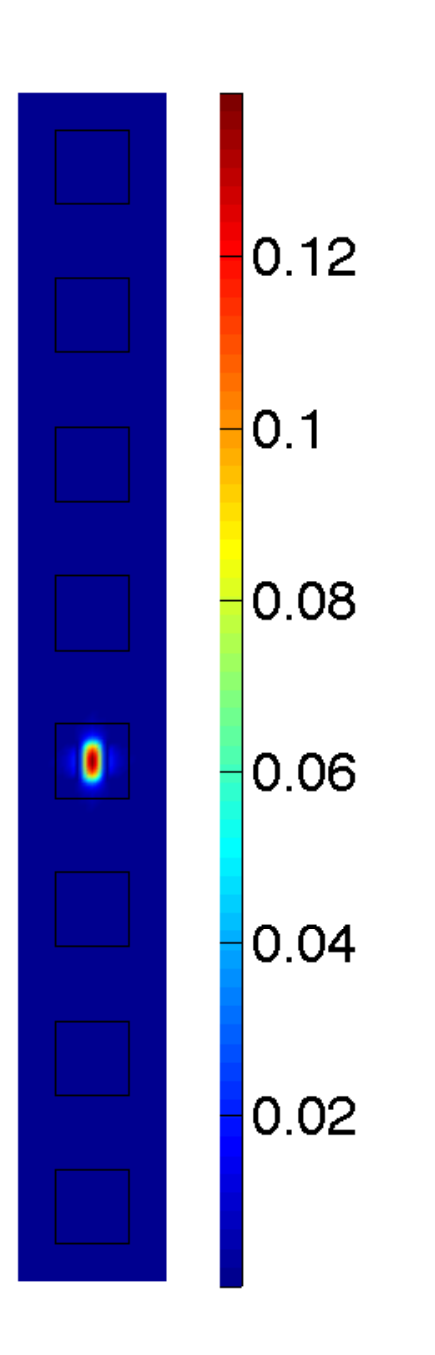}
	\caption*{$\widetilde{S}_7(u)$}%\label{subfig7:sensitivity_indices_multiscale_solution_var_aniso}
\end{subfigure}\hfill
\begin{subfigure}[t]{\figurewidth}
	\centering
	\includegraphics[height=\figureheight]{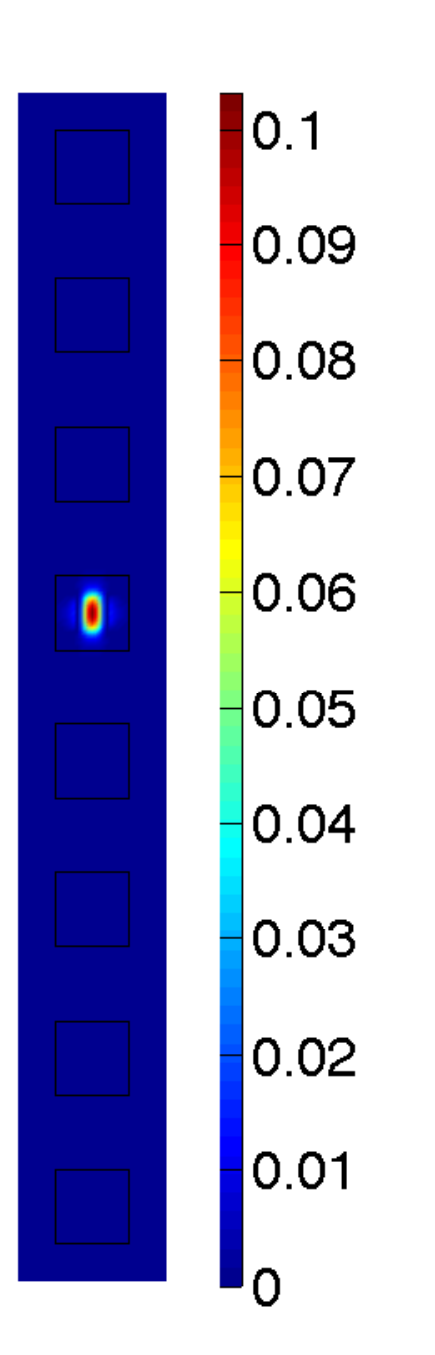}
	\caption*{$\widetilde{S}_9(u)$}%\label{subfig9:sensitivity_indices_multiscale_solution_var_aniso}
\end{subfigure}\hfill
\begin{subfigure}[t]{\figurewidth}
	\centering
	\includegraphics[height=\figureheight]{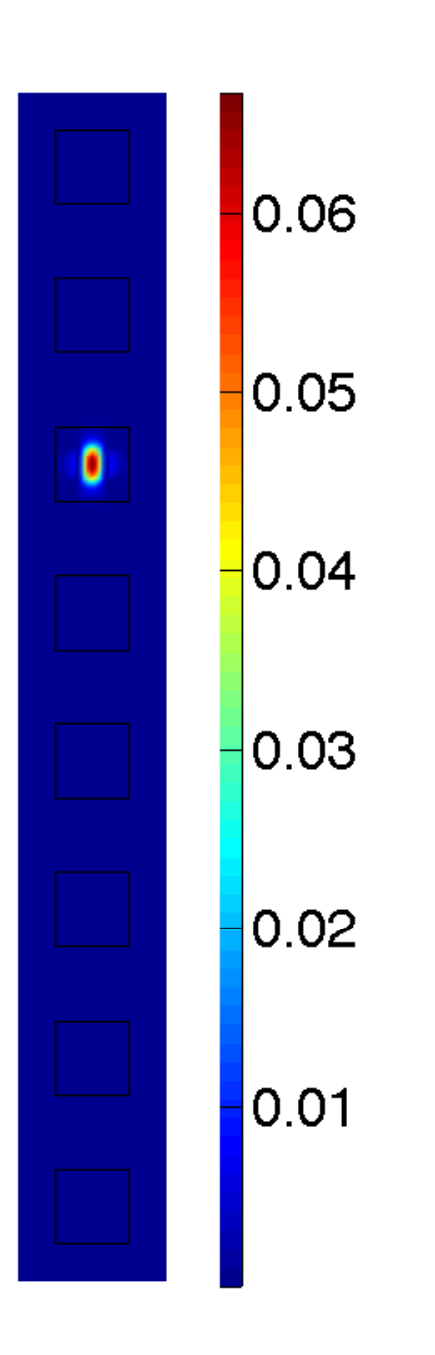}
	\caption*{$\widetilde{S}_{11}(u)$}%\label{subfig11:sensitivity_indices_multiscale_solution_var_aniso}
\end{subfigure}\hfill
\begin{subfigure}[t]{\figurewidth}
	\centering
	\includegraphics[height=\figureheight]{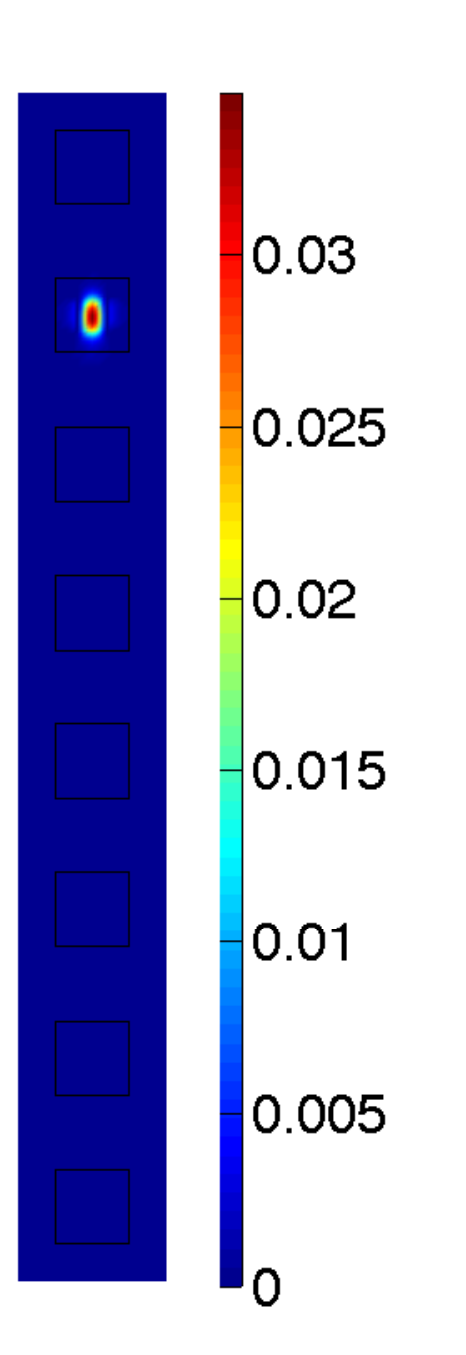}
	\caption*{$\widetilde{S}_{13}(u)$}%\label{subfig13:sensitivity_indices_multiscale_solution_var_aniso}
\end{subfigure}\hfill
\begin{subfigure}[t]{\figurewidth}
	\centering
	\includegraphics[height=\figureheight]{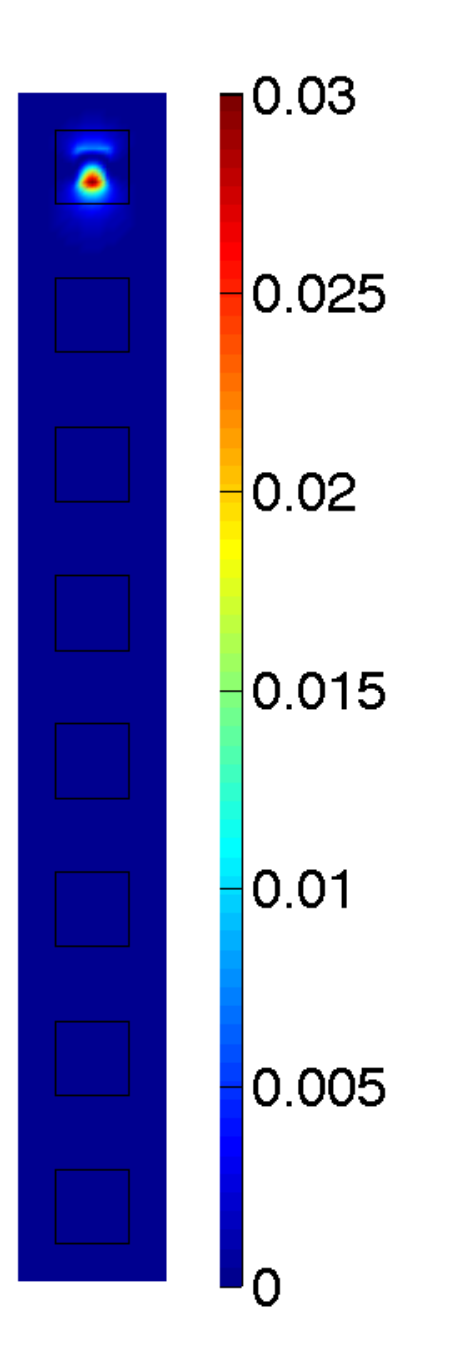}
	\caption*{$\widetilde{S}_{15}(u)$}%\label{subfig15:sensitivity_indices_multiscale_solution_var_aniso}
\end{subfigure}

\begin{subfigure}[t]{\figurewidth}
	\centering
	\includegraphics[height=\figureheight]{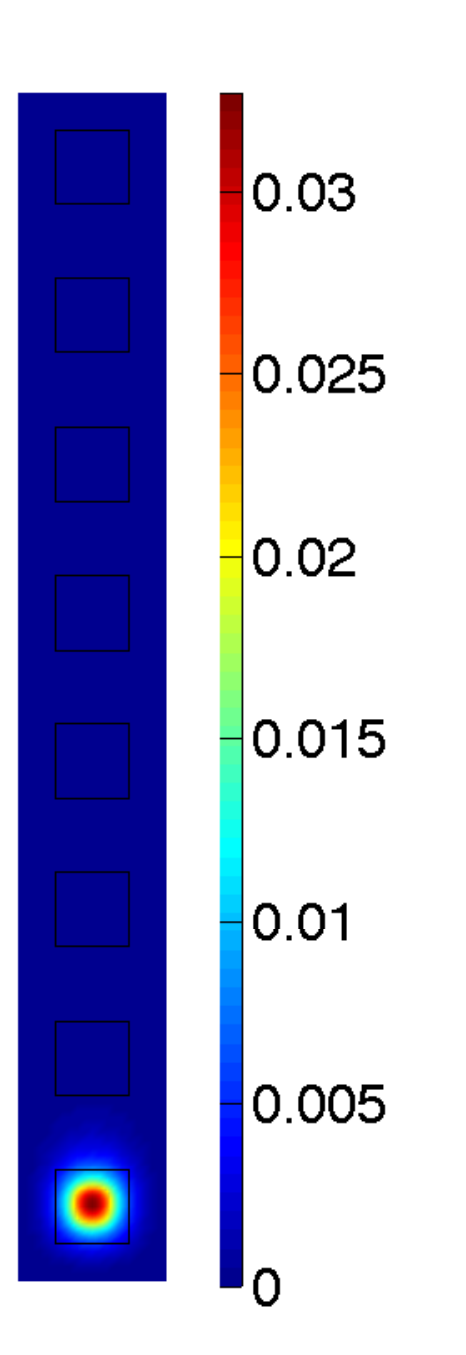}
	\caption*{$\widetilde{S}_2(u)$}%\label{subfig2:sensitivity_indices_multiscale_solution_var_aniso}
\end{subfigure}\hfill
\begin{subfigure}[t]{\figurewidth}
	\centering
	\includegraphics[height=\figureheight]{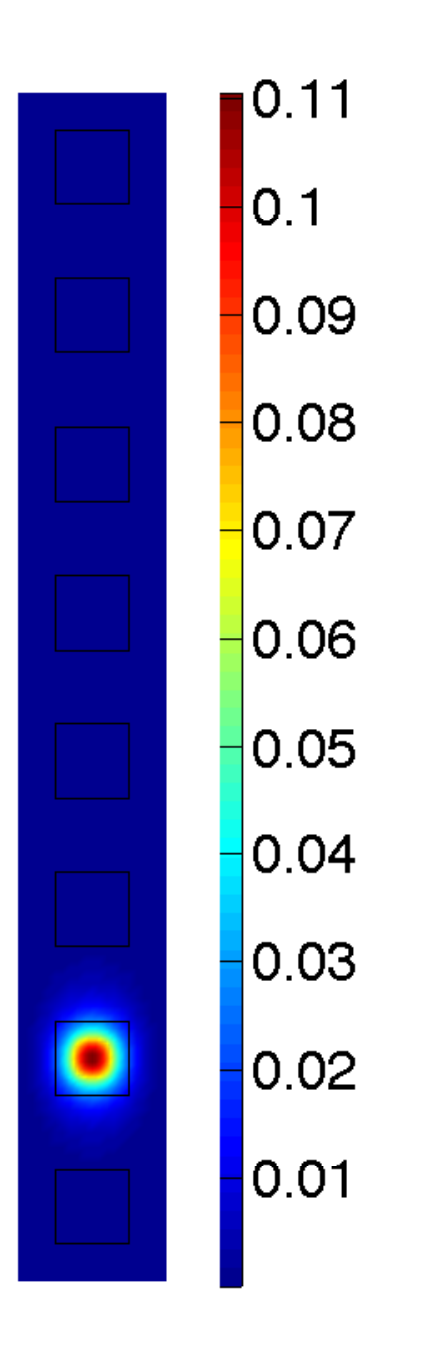}
	\caption*{$\widetilde{S}_4(u)$}%\label{subfig4:sensitivity_indices_multiscale_solution_var_aniso}
\end{subfigure}\hfill
\begin{subfigure}[t]{\figurewidth}
	\centering
	\includegraphics[height=\figureheight]{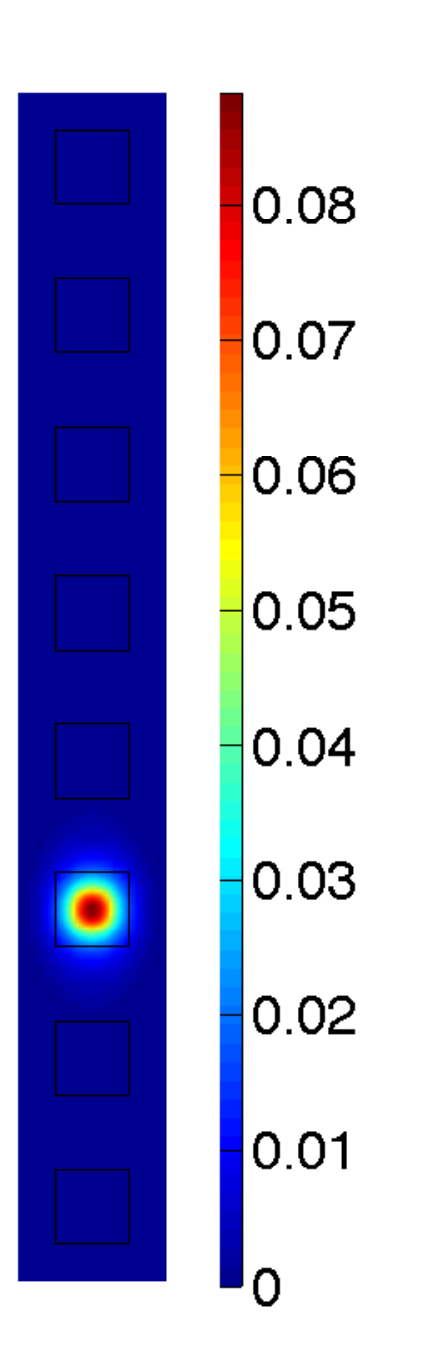}
	\caption*{$\widetilde{S}_6(u)$}%\label{subfig6:sensitivity_indices_multiscale_solution_var_aniso}
\end{subfigure}\hfill
\begin{subfigure}[t]{\figurewidth}
	\centering
	\includegraphics[height=\figureheight]{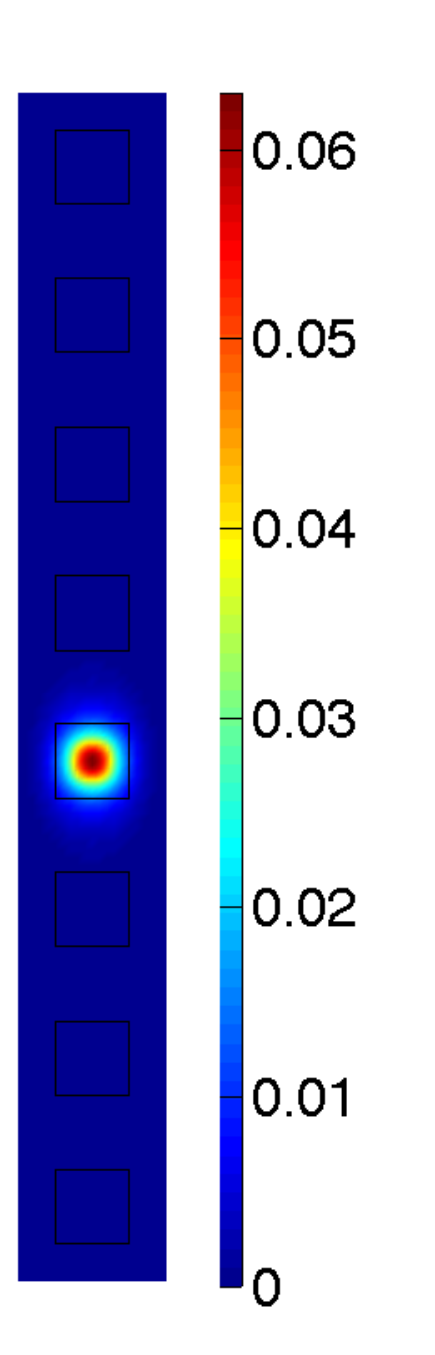}
	\caption*{$\widetilde{S}_8(u)$}%\label{subfig8:sensitivity_indices_multiscale_solution_var_aniso}
\end{subfigure}\hfill
\begin{subfigure}[t]{\figurewidth}
	\centering
	\includegraphics[height=\figureheight]{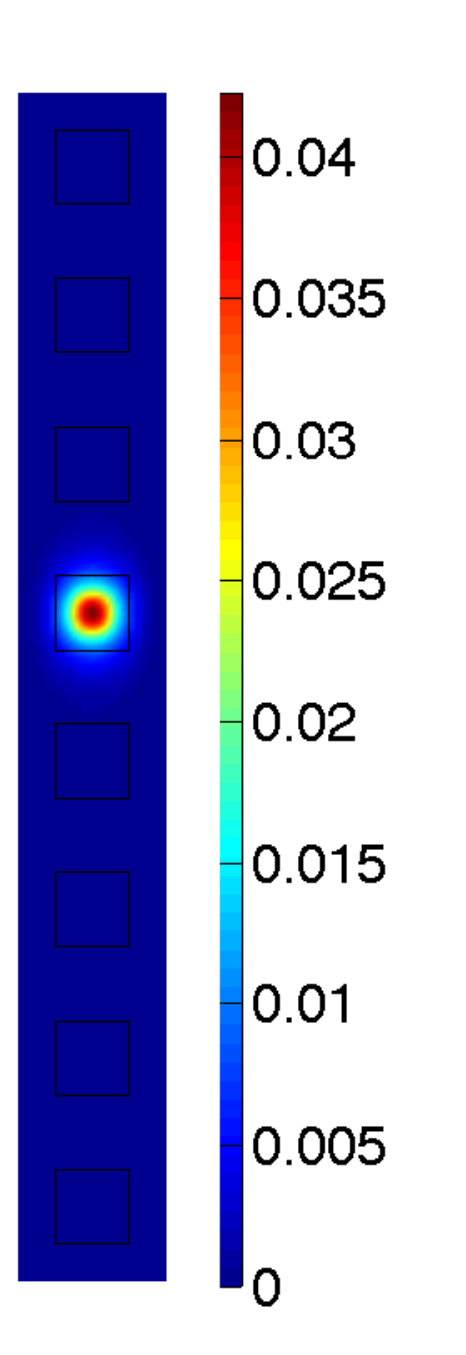}
	\caption*{$\widetilde{S}_{10}(u)$}%\label{subfig10:sensitivity_indices_multiscale_solution_var_aniso}
\end{subfigure}\hfill
\begin{subfigure}[t]{\figurewidth}
	\centering
	\includegraphics[height=\figureheight]{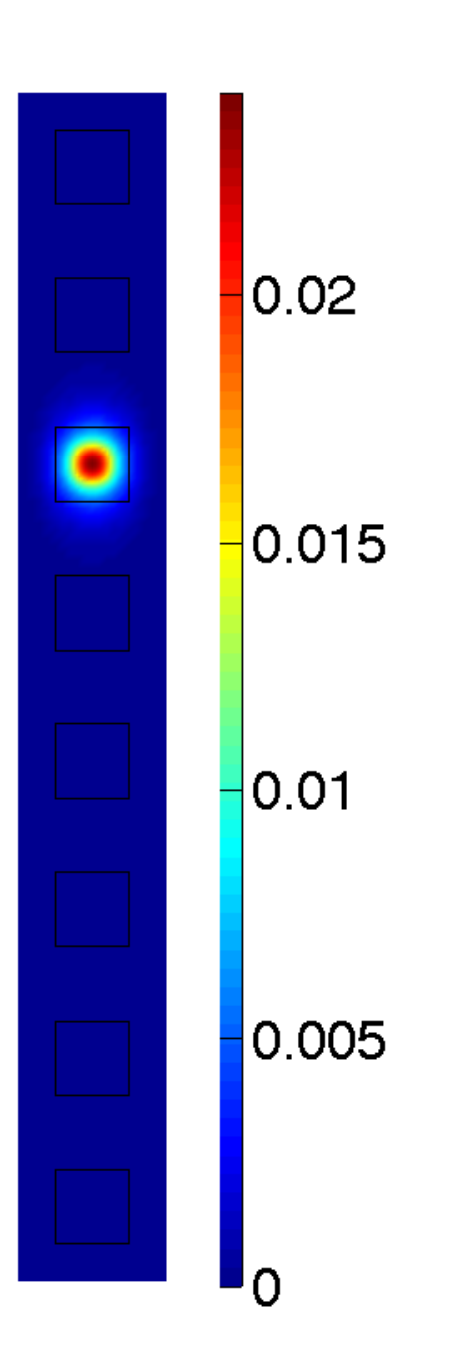}
	\caption*{$\widetilde{S}_{12}(u)$}%\label{subfig12:sensitivity_indices_multiscale_solution_var_aniso}
\end{subfigure}\hfill
\begin{subfigure}[t]{\figurewidth}
	\centering
	\includegraphics[height=\figureheight]{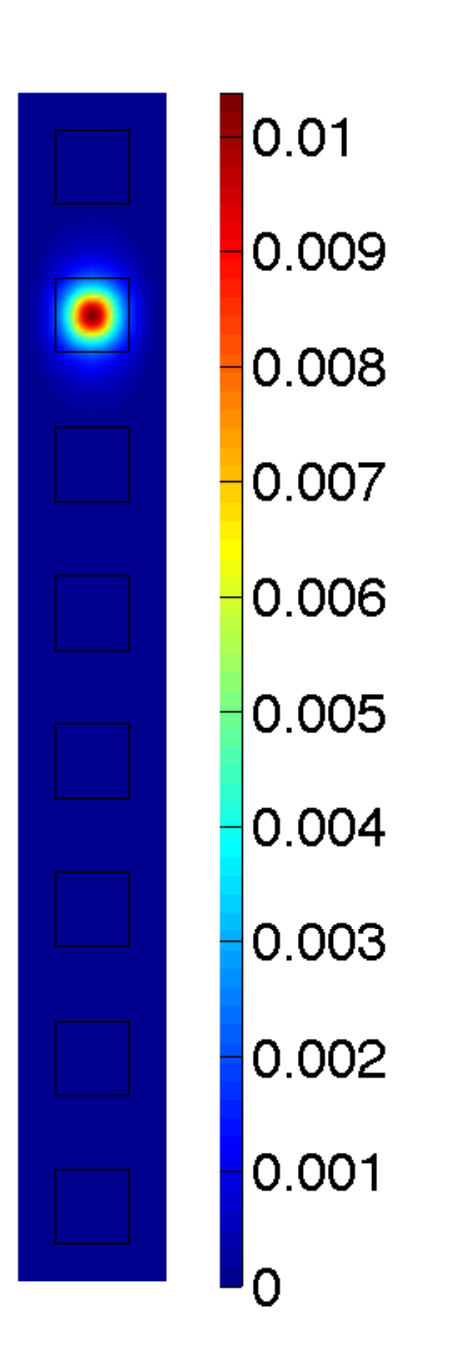}
	\caption*{$\widetilde{S}_{14}(u)$}%\label{subfig14:sensitivity_indices_multiscale_solution_var_aniso}
\end{subfigure}\hfill
\begin{subfigure}[t]{\figurewidth}
	\centering
	\includegraphics[height=\figureheight]{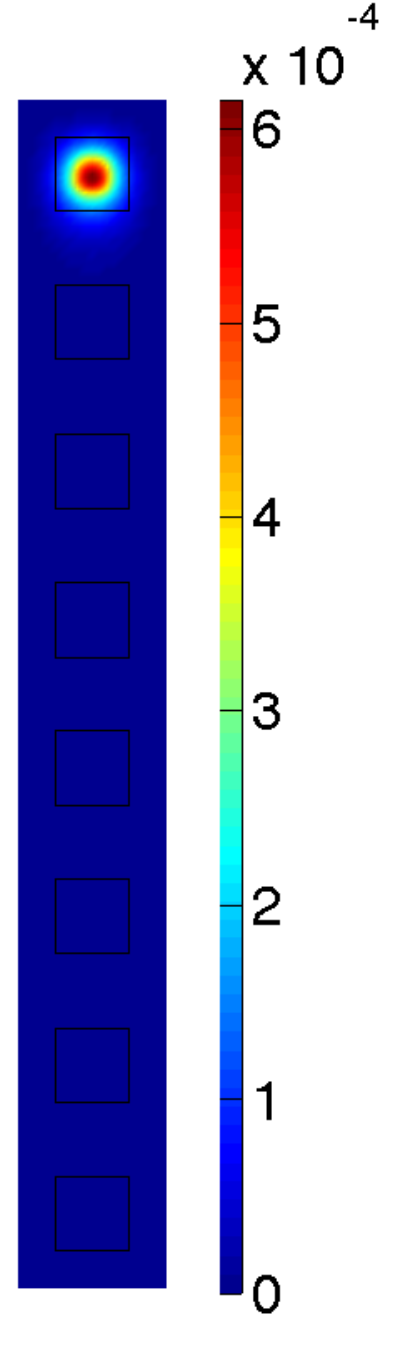}
	\caption*{$\widetilde{S}_{16}(u)$}%\label{subfig16:sensitivity_indices_multiscale_solution_var_aniso}
\end{subfigure}
\caption{Anisotropic case}\label{fig:sensitivity_indices_multiscale_solution_aniso}
\end{subfigure}
\caption{Sensitivity indices $\widetilde{S}_i(u)$ of multiscale solution $u$ with respect to each random variable $\xi_i$, $i \in \set{1,\dots,16}$, where the first row gathers the sensitivity indices $\widetilde{S}_{2q-1}(u)$ associated with random diffusion coefficient $K_q$ (parameterized by $\xi_{2q-1}$) and the second row gathers the sensitivity indices $\widetilde{S}_{2q}(u)$ associated with random reaction parameter $R_q$ (parameterized by $\xi_{2q}$), $q\in\set{1,\dots,8}$}\label{fig:sensitivity_indices_multiscale_solution}
\end{figure}

\section{Conclusion}\label{sec:conclusion}

A global-local iterative method has been proposed for the solution of non-linear stochastic multiscale problems with localized sources of uncertainties and non-linearities.
 The proposed multiscale approach relies on a domain decomposition method with patches.
 The iterative coupling strategy is performed by sequentially solving a linear global problem (with deterministic operator and uncertain right-hand side) and a set of independent non-linear local problems (with uncertain operators and right-hand sides) defined on patches. The global-local iterative coupling algorithm is said non-intrusive in the sense that it does not require any modification of both models and solvers during iterations. The local problems can thus be easily handled using dedicated approximation methods and specific solvers. 
The consistency, convergence and robustness of the proposed algorithm have been analyzed. Numerical results demonstrate the high potential and relevance of the stochastic global-local multiscale approach for dealing with models involving localized uncertainties and possible non-linearities. Several perspectives could be addressed in forthcoming works. First, the modularity of the multiscale approach and in particular its solver coupling capabilities could be exploited in order to take advantage of both commercial software packages available in the industry and in-house research codes, as it was done in recent works \cite{Gen09,Gen11,All11,Duv16}. Second, 
the approach could be extended to more complex non-linear models at multiple scales (\eg{} plasticity, damage or fracture in solid mechanics). Also, in the context of heat and mass transfer in fluid mechanics, localized input uncertainties may result in uncertainties in the solution that are scattered in the whole domain (and not confined within patches) due to high convective terms (even in the linear setting), thus requiring an adaptation of the size of the patches during iterations. Quantifying the effects of localized uncertainties in such multiscale stochastic models is currently one of the appealing engineering and scientific challenges.

\section{Appendix}\label{sec:appendix}

The algorithms and the proofs of some technical results are collected in this appendix.

\subsection{Leave-one-out cross-validation procedure}\label{sec:LOOCValgo}

\begin{lgrthm}[Leave-one-out cross-validation procedure]
\label{LOOCValgo}
$\,$
\begin{algorithmic}[1]
\REQUIRE Coefficients $\mathbf{V} = (\mathbf{v}_{\alpha})_{\alpha \in \Ac}$ of the approximation $\mathbf{v}(\xi) = \sum_{\alpha \in \Ac} \mathbf{v}_{\alpha} \psi_{\alpha}(\xi)$ of $\mathbf{u}(\xi)$, and matrices $\boldsymbol{\Psi} = (\psi_{\alpha}(\xi^l))_{1\leq l\leq N,\alpha\in \Ac}$ and $\mathbf{Y}^T = (\mathbf{u}(\xi^l))_{1\leq l\leq N}$ containing the evaluations of $(\psi_{\alpha}(\xi))_{\alpha \in \Ac}$ and $\mathbf{u}(\xi) = (u_i(\xi))_{i \in \Ic}$
\ENSURE Vector $\boldsymbol{\varepsilon} = (\varepsilon_i)_{i \in \Ic}$, where $\varepsilon_i$ is an estimation of the error $\xE(({{u}_i(\xi) - {v}_i(\xi)})^2)$
\FOR {$i=1,\dots,n$}
	\STATE Compute the set of predicted residuals $\set{\Delta_l^i}_{l=1}^N$ associated with sample set $\set{\xi^l}_{l=1}^N$ using Sherman-Morrison-Woodbury formula: $\Delta^i_l = \frac{r_{li}}{1-h_l}$ with $r_{li}$ the $(l,i)$-th entry of matrix $\mathbf{R} = \boldsymbol{\Psi} \mathbf{V}^T - \mathbf{Y}$, and $h_l$ the $l$-th diagonal term in matrix $\mathbf{H} = \boldsymbol{\Psi} (\boldsymbol{\Psi}^T \boldsymbol{\Psi})^{-1} \boldsymbol{\Psi}^T$
	% standard leave-one-out error
%	\STATE Compute the leave-one-out error $\varepsilon_i = \frac{E_i}{\widehat{m}_2(\mathbf{Y}_i)}$, where $E_i = \frac{1}{N} \sum_{l=1}^N (\Delta^i_l)^2$ and $\widehat{m}_2(\mathbf{Y}_i)$ is the empirical second moment of the $i$-th column $\mathbf{Y}_i = (u_i(\xi^l))_{l=1}^N$ of matrix $\mathbf{Y}$ %containing the evaluations of random variable $u_i$  
	% corrected leave-one-out error
	\STATE Compute the leave-one-out error $e_i = \frac{E_i}{\widehat{m}_2(\mathbf{Y}_i)}$, where $E_i = \frac{1}{N} \sum_{l=1}^N (\Delta^i_l)^2$ and $\widehat{m}_2(\mathbf{Y}_i)$ is the empirical second moment of the $i$-th column $\mathbf{Y}_i = (u_i(\xi^l))_{l=1}^N$ of matrix $\mathbf{Y}$ %containing the evaluations of random variable $u_i$  
	\STATE Compute the %heuristic 
	corrected leave-one-out error $\varepsilon_i = e_i \times T(\Ac,N)$, where $T(\Ac,N) = \(1-\frac{\#\Ac}{N}\)^{-1} \(1 + \frac{\xtr\(\mathbf{\widehat{C}}^{-1}\)}{N}\)$ is a correction factor allowing to reduce the sensitivity %of the error 
	to overfitting \cite{Cha02,Bla11} and where $\mathbf{\widehat{C}} = \frac{1}{N} \boldsymbol{\Psi}^T \boldsymbol{\Psi}% \in \xR^{N \times N}
	$ is the empirical covariance matrix of $(\psi_{\alpha}(\xi))_{\alpha \in \Ac}$
\ENDFOR
\end{algorithmic}
\end{lgrthm}

\subsection{Adaptive sparse least-squares solver with random sampling and working set strategy}\label{sec:LSalgo}

\begin{lgrthm}[Adaptive sparse least-squares solver]
\label{LSalgo}
$\,$
\begin{algorithmic}[1]
\REQUIRE Initial number of samples $N\geq 1$, sampling factor $p_{\text{add}} > 0$, parameter $\theta \in \intervalcc{0}{1}$
\ENSURE Monotone set $\Ac$ and coefficients $\mathbf{V} = (\mathbf{v}_\alpha)_{\alpha\in \Ac}$ of the least-squares approximation $\mathbf{v}(\xi) = \sum_{\alpha \in \Ac} \mathbf{v}_{\alpha} \psi_{\alpha}(\xi)$ of $\mathbf{u}(\xi)$
\STATE Start with null initial set $\Ac=\set{0_{\Fc}}$
\STATE Generate the initial sample set $\set{\xi^l}_{l=1}^{N}$ randomly
\STATE Compute the matrices $\boldsymbol{\Psi}_\Ac = (\psi_{\alpha}(\xi^l))_{1\leq l\leq N,\alpha\in \Ac}$ and $\mathbf{Y}^T = (\mathbf{u}(\xi^l))_{1\leq l\leq N}$
\WHILE {no convergence}
	\STATE \COMMENT{Adaptive random sampling}
	\WHILE {no convergence and no stagnation}
		\STATE Generate the additional sample set $\set{\xi^{N+l}}_{l=1}^{N_{\text{add}}}$ randomly, with $N_{\text{add}} = \ceil(p_{\text{add}} N)$
		\STATE Compute the matrices $\boldsymbol{ \Psi}_{\Ac,\text{add}} = (\psi(\xi^{N+l}))_{1\leq l\leq N_{\text{add}},\alpha \in \Ac}$ and $\mathbf{Y}_{\text{add}}^T = (\mathbf{u}(\xi^{N+l}))_{1\leq l\leq N_{\text{add}}}$
		\STATE Update the number of samples $N \gets N+N_{\text{add}}$ and the matrices $\boldsymbol{\Psi}_\Ac^T \gets (\boldsymbol{\Psi}_\Ac^T,\boldsymbol{\Psi}_{\Ac,\text{add}}^T)$ and $\mathbf{Y}^T \gets (\mathbf{Y}^T,\mathbf{Y}_{\text{add}}^T)$
		\STATE Compute the coefficients $\mathbf{V} = (\mathbf{v}_\alpha)_{\alpha\in \Ac}$ such that $\mathbf{V}^T = (\boldsymbol{\Psi}_{\Ac}^T \boldsymbol{\Psi}_{\Ac})^{-1} \boldsymbol{\Psi}_{\Ac}^T \mathbf{Y}$
	\ENDWHILE
	\STATE \COMMENT{Working set strategy}
	\WHILE {no convergence and no overfitting}
		\STATE Compute the reduced margin $\Mc = \Mc_{\text{red}}(\Ac)$ of monotone set $\Ac$ and the set $\Tc = \Ac \cup \Mc$
		\IF {$\#\Tc > N$} \BREAK \ENDIF
		\STATE Compute the matrices $\boldsymbol{\Psi}_{\Mc} = (\psi_{\alpha}(\xi^l))_{1\leq l\leq N,\alpha\in \Mc}$ and $\boldsymbol{\Psi}_{\Tc} = (\boldsymbol{\Psi}_\Ac,\boldsymbol{\Psi}_{\Mc})$
		\STATE Compute the coefficients $\mathbf{V} = (\mathbf{v}_\alpha)_{\alpha\in \Tc}$ such that $\mathbf{V}^T = (\boldsymbol{\Psi}_{\Tc}^T \boldsymbol{\Psi}_{\Tc})^{-1} \boldsymbol{\Psi}_{\Tc}^T \mathbf{Y}$ 
		\STATE Compute the vector $(\norm{\mathbf{v}_\alpha}_2)_{\alpha \in \Mc}$
		\STATE Define the smallest (monotone) subset $\Nc$ of $\Mc$ such that $e(\Nc) \geq \theta e(\Mc)$, with $e(\Nc) = \sum_{\alpha \in \Nc} \norm{\mathbf{v}_\alpha}^2_2$ and $e(\Mc) = \sum_{\alpha \in \Mc} \norm{\mathbf{v}_\alpha}^2_2$
		\STATE Update the multi-index set $\Ac \gets \Ac \cup \Nc$ and the matrix $\boldsymbol{\Psi}_{\Ac} \gets (\boldsymbol{\Psi}_{\Ac},\boldsymbol{\Psi}_{\Nc})$, where $\boldsymbol{\Psi}_{\Nc}$ is the submatrix of $\boldsymbol{\Psi}_{\Mc}$ whose columns correspond to multi-indices $\alpha \in \Nc$
		\STATE Compute the coefficients $\mathbf{V} = (\mathbf{v}_\alpha)_{\alpha\in \Ac}$ such that $\mathbf{V}^T = (\boldsymbol{\Psi}_{\Ac}^T \boldsymbol{\Psi}_{\Ac})^{-1} \boldsymbol{\Psi}_{\Ac}^T \mathbf{Y}$
	\ENDWHILE
\ENDWHILE
\end{algorithmic}
\end{lgrthm}

\subsection{Proof of Lemma~\ref{lmm:Cv}}\label{sec:lmm_Cv}

The first inequality $\abs{v}_{\Vc} \leq \norm{v}_{\Vc}$ is obvious.
Using \eqref{generalizedpoincare} with $(\Oc,\Ec) = (\Lambda,\Gamma)$, we obtain 
$\norm{v_{\restrictto{\Lambda}}}_{\xHone(\Lambda)} \leq C_{\Lambda,\Gamma} (\abs{v_{\restrictto{\Lambda}}}_{\xHone(\Lambda)} + \norm{v_{\restrictto{\Lambda}}}_{\xHn{1/2}(\Gamma)})$. Then, using the weak continuity of $v$ on $\Gamma$, we have $\norm{v_{\restrictto{\Lambda}}}_{\xHn{1/2}(\Gamma)} = \norm{v_{\restrictto{\Omega \setminus \Lambda}}}_{\xHn{1/2}(\Gamma)} \leq \beta_{\tau} \norm{v_{\restrictto{\Omega \setminus \Lambda}}}_{\xHone(\Omega \setminus \Lambda)}$, where $\beta_{\tau}$ is the norm of the trace operator $\tau \colon \xHone(\Omega \setminus \Lambda) \to \xHn{1/2}(\Gamma)$. Now, using \eqref{generalizedpoincare} with $(\Oc,\Ec) = (\Omega \setminus \Lambda,\Gamma_D \cap \partial (\Omega \setminus \Lambda))$, we obtain 
$\norm{v_{\restrictto{\Omega \setminus \Lambda}}}_{\xHone(\Omega \setminus \Lambda)} \leq \widehat{C} \abs{v_{\restrictto{\Omega \setminus \Lambda}}}_{\xHone(\Omega \setminus \Lambda)}$, with $\widehat{C} = C_{\Omega \setminus \Lambda,\Gamma_D \cap \partial (\Omega \setminus \Lambda)}$. Then, we deduce that
\begin{align*}
\norm{v}_{\Vc}^2 
&\leq \widehat{C}^2 \abs{v_{\restrictto{\Omega \setminus \Lambda}}}_{\xHone(\Omega \setminus \Lambda)}^2 + C_{\Lambda,\Gamma}^2 (\abs{v_{\restrictto{\Lambda}}}_{\xHone(\Lambda)} + \beta_{\tau} \widehat{C} \abs{v_{\restrictto{\Omega \setminus \Lambda}}}_{\xHone(\Omega \setminus \Lambda)})^2 \\
&\leq \widehat{C}^2 ( 1 + 2 C_{\Lambda,\Gamma}^2 \beta_{\tau}^2 ) \abs{v_{\restrictto{\Omega \setminus \Lambda}}}_{\xHone(\Omega \setminus \Lambda)}^2 + 
2 C_{\Lambda,\Gamma}^2 \abs{v_{\restrictto{\Lambda}}}_{\xHone(\Lambda)} ^2 \leq C_{\Vc}^2 \abs{v}_{\Vc}^2,
\end{align*}
where $C_{\Vc}^2 = \max\set{\widehat{C}^2 ( 1 + 2 C_{\Lambda,\Gamma}^2 \beta_{\tau}^2 ),2 C_{\Lambda,\Gamma}^2}$. 
$C_{\Vc}$ is independent of $\xi$ since $\widehat{C}$ and $C_{\Lambda,\Gamma}$ are independent of $\xi$ (assumption~\ref{assgeometry}) and $\beta_{\tau}$ is independent of $\xi$.

\subsection{Proof of Theorem~\ref{thrm:wellposedness_weakformulation}}\label{sec:thrm_wellposedness_weakformulation}

Let $\Vc^{\ast}$ be the topological dual space to $\Vc$ 
and let $\scalproda{\cdot}{\cdot}$ denote the duality pairing between $\Vc$ and $\Vc^{\ast}$. 
For $u \in \Vc$, $v \mapsto d_{\Omega}(u,v;\xi)$ is linear and continuous. Then, there exists a unique non-linear map $\Dc(\xi) \colon \Vc \to \Vc^{\ast}$ such that $d_{\Omega}(u,v;\xi) = \scalproda{\Dc(\xi)(u)}{v}$ for all $u,v \in \Vc$. As $\ell_{\Omega}(\cdot;\xi)$ is linear and continuous on $\Vc$, there exists a unique $L(\xi) \in \Vc^{\ast}$ such that $\ell_{\Omega}(v;\xi) = \scalproda{L(\xi)}{v}$ for all $v \in \Vc$. Problem \eqref{weakformulation} can then be written as $\Dc(\xi)(u(\xi)) = L(\xi)$. First, we have that for all $u,v \in \Vc$, the map $t \mapsto \scalproda{\Dc(\xi)(u+tv)}{v} = a_{\Omega}(u,v;\xi) + t a_{\Omega}(v,v;\xi) + n_{\Lambda}(u+tv,v;\xi)$ is continuous, which implies that $\Dc(\xi)$ is radially continuous. Then, assumption~\eqref{coercivity_a} on $a_{\Omega \setminus \Lambda}$ and $a_{\Lambda}$ and assumption~\eqref{monotonicity_n} on $n_{\Lambda}$ imply that for all $u,v \in \Vc$,
\begin{align*}
\scalproda{\Dc(\xi)(u) - \Dc(\xi)(v)}{u-v} 
&= a_{\Omega \setminus \Lambda}(u-v,u-v) + a_{\Lambda}(u-v,u-v) + n_{\Lambda}(u,u-v) - n_{\Lambda}(v,u-v)\\
&\geq \alpha_a (\abs{u-v}_{\xHone(\Omega \setminus \Lambda)}^2 + \abs{u-v}_{\xHone(\Lambda)}^2) = \alpha_a \abs{u-v}_{\Vc}^2 \geq \frac{\alpha_a}{C_{\Vc}^2} \norm{u-v}_{\Vc}^2 \coloneqq \alpha_{\Dc} \norm{u-v}_{\Vc}^2,
\end{align*}
where the last inequality comes from Lemma~\ref{lmm:Cv}.
Then $\Dc(\xi)$ is strongly monotone with monotonicity constant $\alpha_{\Dc} = \frac{\alpha_a}{C_{\Vc}^2}$.
Also, assumption~\eqref{zerocondition_n} on $n_{\Lambda}$ implies $\Dc(\xi)(0) = 0$. From this latter condition and from the strong monotonicity of $\Dc(\xi)$, we obtain that $\scalproda{\Dc(\xi)(u)}{u} \geq \alpha_{\Dc} \norm{u}_{\Vc}^2$ for all $u \in \Vc$, and therefore $\Dc(\xi)$ is coercive. Accordingly, $\Dc(\xi)$ being radially continuous, monotone and coercive, the Browder-Minty theorem \cite[Theorem~2.18]{Rou05} 
ensures that $\Dc(\xi)$ is surjective, and therefore there exists a solution $u(\xi) \in \Vc$ to problem \eqref{weakformulation}. The strict monotonicity of $\Dc(\xi)$ ensures that this solution is unique, so that we can define an inverse map $\Dc(\xi)^{-1} \colon \Vc^{\ast} \to \Vc$. The strong monotonicity of $\Dc(\xi)$ then implies that $\Dc(\xi)^{-1}$ is Lipschitz continuous, with $\norm{\Dc(\xi)^{-1}(L) - \Dc(\xi)^{-1}(\tilde{L})}_{\Vc} \leq \frac{1}{\alpha_{\Dc}} \norm{L-\tilde{L}}_{\Vc^{\ast}}$ for all $L,\tilde{L} \in \Vc^{\ast}$. Finally, from the strong monotonicity and from assumption~\ref{assdata}, we have that $\norm{u(\xi)}_{\Vc} \leq \frac{1}{\alpha_{\Dc}} \norm{\Dc(\xi)(u(\xi))}_{\Vc^{\ast}} = \frac{1}{\alpha_{\Dc}} \norm{L(\xi)}_{\Vc^{\ast}} = \frac{1}{\alpha_{\Dc}} \sup_{\substack{\norm{v}_{\Vc} =1}} \scalproda{L(\xi)}{v} = \frac{1}{\alpha_{\Dc}} \sup_{\substack{\norm{v}_{\Vc} =1}} \ell_{\Omega}(v;\xi)
\leq \frac{1}{\alpha_{\Dc}} \kappa(\xi)$, with $\kappa \in \xLn{p}_{\mu}(\Xi)$. Since $\alpha_a$ and $C_{\Vc}$ are independent of $\xi$ ({see} Lemma~\ref{lmm:Cv}), $\alpha_{\Dc}$ is independent of $\xi$ and we deduce that $u \in \xLn{p}_{\mu}(\Xi;\Vc)$.

\subsection{Proof of Theorem~\ref{thrm:existence_globallocalformulation}}\label{sec:thrm_existence_globallocalformulation}

Let $u(\xi) \in \widehat{\Vc}$ such that $u(\xi)_{\restrictto{\Omega \setminus \Lambda}} = U(\xi) \in \Uc$ and $u(\xi)_{\restrictto{\Lambda}} = w(\xi) \in \Wc$. Equation 
\eqref{interfacepart} implies that $u(\xi) \in {\Vc}$. Then, considering a test function $\delta u \in \Vc$ and summing \eqref{complementarypart} and \eqref{patchpart}, we obtain that $u(\xi)$ verifies \eqref{weakformulation}. From Theorem~\ref{thrm:wellposedness_weakformulation}, we deduce that problem \eqref{globallocalformulation} admits a unique solution $(U(\xi),w(\xi)) \in \Uc \times \Wc$ which coincides with the solution of \eqref{weakformulation}. Moreover, since $u\in \xLn{p}_{\mu}(\Xi;\Vc)$, we deduce that $U\in \xLn{p}_{\mu}(\Xi;\Uc)$ and $w\in \xLn{p}_{\mu}(\Xi;\Wc)$. Then, let $R \colon \xHn{1/2}(\Gamma) \to \xHone(\Omega \setminus \Lambda)$ denote a linear continuous extension operator with continuity constant $\beta_{R}$. Equation \eqref{complementarypart} yields 
\begin{equation*}
b_{\Gamma}(\lambda(\xi),v) = b_{\Gamma}(\lambda(\xi),R(v)) = \ell_{\Omega \setminus \Lambda}(R(v)) - a_{\Omega \setminus \Lambda}(U(\xi),R(v))
\end{equation*}
for all $v \in \xHn{1/2}(\Gamma)$. The right-hand side being a continuous linear form on $\xHn{1/2}(\Gamma)$, we obtain the existence and uniqueness of a solution $\lambda(\xi) \in \Mc$. Also, 
$\norm{\lambda(\xi)}_{\Mc} = \sup_{\substack{\norm{v}_{\xHn{1/2}(\Gamma)} =1}} b_{\Gamma}(\lambda(\xi),v)
= \sup_{\substack{\norm{v}_{\xHn{1/2}(\Gamma)} =1}} \ell_{\Omega \setminus \Lambda}(R(v)) - a_{\Omega \setminus \Lambda}(U(\xi),R(v)) 
\leq \beta_R (\kappa(\xi) + \beta_a \norm{U(\xi)}_{\Uc})$. Since both $\kappa(\xi)$ and $\norm{U(\xi)}_{\Uc}$ belong to $\xLn{p}_{\mu}(\Xi)$, $\lambda$ belongs to $\xLn{p}_{\mu}(\Xi;\Mc)$.

\subsection{Proof of Lemma~\ref{lmm:equiv_norm_Ustar}}\label{sec:lmm_equiv_norm_Ustar}

Let $\tau \colon \xHone(\Omega \setminus \Lambda) \to \xHn{1/2}(\Gamma)$ denote the trace operator which is linear and continuous with norm $\beta_{\tau}$ independent of $\xi$. Let $R \colon \xHn{1/2}(\Gamma) \to \xHone(\widetilde{\Lambda})$ denote a linear continuous extension operator with norm $\beta_{R}$ independent of $\xi$. For $V \in \widetilde{\Uc}_{\star}$, we write $V_{\restrictto{\widetilde{\Lambda}}} = R(V_{\restrictto{\Gamma}}) + Z$, where $Z \in \xHone_0(\widetilde{\Lambda})$ is such that $c_{\widetilde{\Lambda}}(R(V_{\restrictto{\Gamma}}) + Z,\delta U) = 0$ for all $\delta U \in \xHone_0(\widetilde{\Lambda})$. From assumption~\ref{assc} on $c_{\widetilde{\Lambda}}$ and using \eqref{generalizedpoincare} with $(\Oc,\Ec) = (\widetilde{\Lambda},\Gamma)$, we obtain that 
$\norm{Z}_{\xHone(\widetilde{\Lambda})}^2 \leq C_{\widetilde{\Lambda},\Gamma}^2 \abs{Z}_{\xHone(\widetilde{\Lambda})}^2 \leq \frac{C_{\widetilde{\Lambda},\Gamma}^2}{\alpha_c} c_{\widetilde{\Lambda}}(Z,Z) = \frac{C_{\widetilde{\Lambda},\Gamma}^2}{\alpha_c} c_{\widetilde{\Lambda}}(-R(V_{\restrictto{\Gamma}}),Z) \leq \frac{C_{\widetilde{\Lambda},\Gamma}^2 \beta_c}{\alpha_c} \norm{R(V_{\restrictto{\Gamma}})}_{\xHone(\widetilde{\Lambda})} \norm{Z}_{\xHone(\widetilde{\Lambda})} $. Then, from the continuity of $\tau$ and $R$, we deduce that $\norm{V}_{\xHone(\widetilde{\Lambda})} \leq (1+\frac{\beta_c C_{\widetilde{\Lambda},\Gamma}^2}{\alpha_c}) \beta_{R} \norm{V}_{\xHn{1/2}(\Gamma)} \leq (1+\frac{\beta_c C_{\widetilde{\Lambda},\Gamma}^2}{\alpha_c}) \beta_{R} \beta_{\tau} \norm{V}_{\xHone(\Omega \setminus \Lambda)}$. Finally, since $\norm{V}_{\widetilde{\Uc}}^2 = \norm{V}_{\Uc}^2 + \norm{V}_{\xHone(\widetilde{\Lambda})}^2$, we obtain that $\norm{V}_{\Uc} \leq \norm{V}_{\widetilde{\Uc}} \leq C_{\widetilde{\Uc}} \norm{V}_{\Uc}$ for all $V \in \widetilde{\Uc}_{\star}$, with $C_{\widetilde{\Uc}} = (1+ (1+\frac{\beta_c C_{\widetilde{\Lambda},\Gamma}^2}{\alpha_c})^2 \beta_{R}^2 \beta_{\tau}^2)^{1/2}$ independent of $\xi$.

\subsection{Proof of Lemma~\ref{lmm:continuity_linearmap}}\label{sec:lmm_continuity_linearmap}

First, using property \eqref{coercivity_c} for $c_{\widetilde{\Omega}}$, property \eqref{continuity_c} for 
 $c_{\widetilde{\Lambda}}$ and relation \eqref{generalizedpoincare} for $(\Oc,\Ec) = (\widetilde{\Omega},\Gamma_D \cap \partial \widetilde{\Omega})$ (with constant $\widetilde{C}=C_{\widetilde{\Omega},\Gamma_D \cap \partial \widetilde{\Omega}}$), we obtain that 
 $\norm{\Upsilon(V;\xi)}_{\widetilde{\Uc}} \leq \beta_{\Upsilon} \norm{V}_{\widetilde{\Uc}}$ for all $V \in \widetilde{\Uc}$, with $\beta_{\Upsilon}=\frac{\beta_c \widetilde{C}^2}{\alpha_c}$ independent of $\xi$. Then, using again property \eqref{coercivity_c} for $c_{\widetilde{\Omega}}$ and relation \eqref{generalizedpoincare} for $(\Oc,\Ec) = (\widetilde{\Omega},\Gamma_D \cap \partial \widetilde{\Omega})$, we have that $\norm{\Phi(\beta;\xi)}_{\widetilde{\Uc}}^2 \leq 
\frac{\widetilde{C}^2}{\alpha_c} \abs{b_{\Gamma}(\beta,\Phi(\beta;\xi))} \leq \frac{\widetilde{C}^2}{\alpha_c} \norm{\beta}_{\Mc} \norm{\Phi(\beta;\xi)_{\restrictto{\Gamma}}}_{\xHn{1/2}(\Gamma)} \leq \frac{\beta_{\tau} \widetilde{C}^2}{\alpha_c} \norm{\beta}_{\Mc} \norm{\Phi(\beta;\xi)}_{\widetilde{\Uc}}, $ where $\beta_{\tau}$ is the norm of the trace operator $\tau \colon \xHone(\Omega \setminus \Lambda) \to \xHn{1/2}(\Gamma)$. That proves $\norm{\Phi(\beta;\xi)}_{\widetilde{\Uc}} \leq \beta_{\Phi} \norm{\beta}_{\Mc}$ for all $\beta \in \Mc$, with $\beta_{\Phi} = \frac{\beta_{\tau} \widetilde{C}^2}{\alpha_c}$ independent of $\xi$.

\subsection{Proof of Lemma~\ref{lmm:continuity_nonlinearmap}}\label{sec:lmm_continuity_nonlinearmap}

Given the property \eqref{localization_delta} for the definition of the patch, we can introduce a partition 
$\overline{\Lambda} = \overline{\Lambda_i \cup \Lambda_e}$ with $\Lambda_i \cap \Lambda_e = \emptyset$ and $\dist(\Lambda_i,\Gamma) = \delta$, for which $\Lambda_{\star} \subset \Lambda_i$. That means $\Lambda_e$ is a band of width $\delta$ around $\Gamma$ where the differential operator is linear. Let $\Gamma_e = \partial {\Lambda_e} \cap \partial{\Lambda_i}$ and $\Gamma_e^D = \partial \Lambda_e \cap \Gamma_D$. The restriction $w_e(\xi)$ of the local solution $w(\xi) = \Theta(U(\xi);\xi)$ to $\Lambda_e$ is such that $w_e(\xi) \in \xHone(\Lambda_e)$, $w_e(\xi)_{\restrictto{\Gamma}} = U(\xi)_{\restrictto{\Gamma}}$, $w_e(\xi)_{\restrictto{\Gamma_e}} = w(\xi)_{\restrictto{\Gamma_e}}$, $w_e(\xi)_{\restrictto{\Gamma_e^D}}=0$ and 
\begin{equation}\label{pbLambdae}
a_{\Lambda_e}(w_e(\xi),\delta w_e) = \ell_{\Lambda_e}(\delta w_e;\xi)
\end{equation}
for all $\delta w_e \in \xHone(\Lambda_e)$ such that $\delta w_e=0$ on $\Gamma \cup \Gamma_e \cup \Gamma_e^D$.
Using the linearity of this problem and introducing linear continuous extension operators from $\xHn{1/2}(\Gamma)$ and $\xHn{1/2}(\Gamma_e)$ to $\xHone(\Lambda_e)$, it can be easily proved that $w_e(\xi)$ can be written {as}
\begin{equation}\label{solutionLambdae}
w_e(\xi) = F_{\Gamma}(U(\xi)_{\restrictto{\Gamma}}) + F_{\Gamma_e}(w(\xi)_{\restrictto{\Gamma_e}}) + \overline{w}_e(\xi),
\end{equation}
where $\overline{w}_e(\xi)=0$ on $\Gamma \cup \Gamma_e \cup \Gamma_{e}^D$, and where 
$F_{\Gamma} \colon \xHn{1/2}(\Gamma) \to \xHone(\Lambda_e)$ and $F_{\Gamma_e} \colon \xHn{1/2}(\Gamma_e) \to \xHone(\Lambda_e)$ are linear continuous extension operators with respective norms $\norm{F_{\Gamma}}$ and $\norm{F_{\Gamma_e}}$ independent of $\xi$, such that 
$F_{\Gamma}(v) = v$ on $\Gamma$, $F_{\Gamma}(v)=0$ on $\Gamma_e \cup \Gamma_e^D$, $F_{\Gamma_e}(v) = v$ on $\Gamma_e$, and $F_{\Gamma_e}(v)=0$ on $\Gamma \cup \Gamma_e^D$. From \eqref{pbLambdae}, we {obtain} that the normal flux on $\Gamma_e$, denoted $\lambda_e(\xi) = - {B_L(w(\xi),\nabla w(\xi);\cdot)} \cdot n \in \xHn{1/2}(\Gamma_e)^{\ast}$ (with $n$ the unit normal to $\Gamma_e$ pointing outward $\Lambda_e$), is such that 
$b_{\Gamma_e}(\lambda_e(\xi), \delta w_{e\restrictto{\Gamma_e}}) = \ell_{\Lambda_e}(\delta w_e;\xi) - a_{\Lambda_e}(w_e(\xi),\delta w_e)$
for all $\delta w_e\in \xHone(\Lambda_e)$ such that $\delta w_e=0$ on $\partial \Lambda_e \setminus \Gamma_e$, where $b_{\Gamma_e}$ denotes the duality pairing between $\xHn{1/2}(\Gamma_e)$ and $\xHn{1/2}(\Gamma_e)^{\ast}$. From \eqref{solutionLambdae}, we deduce that 
\begin{equation*}
b_{\Gamma_e}(\lambda_e(\xi), \delta v) = - g_{\Gamma}(U(\xi)_{\restrictto{\Gamma}}, \delta v) - g_{\Gamma_e}(w(\xi)_{\restrictto{\Gamma_e}}, \delta v) + b_{\Gamma_e}(\overline{\lambda}_e(\xi), \delta v)
\end{equation*}
for all $\delta v \in \xHn{1/2}(\Gamma_e)$, where $g_{\Gamma} {\colon \xHn{1/2}(\Gamma) \times \xHn{1/2}(\Gamma_e) \to \xR}$ and $g_{\Gamma_e} {\colon \xHn{1/2}(\Gamma_e) \times \xHn{1/2}(\Gamma_e) \to \xR}$ are continuous bilinear forms with {respective} norms $\norm{g_{\Gamma}}$ and $\norm{g_{\Gamma_e}}$ independent of $\xi$, and $\overline{\lambda}_e(\xi) \in \xHn{1/2}(\Gamma_e)^{\ast}$. By definition, $g_{\Gamma_e}$ is such that $g_{\Gamma_e}(v,v) = a_{\Lambda_e}(F_{\Gamma_e}(v),\delta w)$ for all $v \in \xHn{1/2}(\Gamma_e)$ and $\delta w \in \xHone(\Lambda_e)$ such that $\delta w = v$ on $\Gamma_e$ and $\delta w = 0$ on $\Gamma \cup \Gamma_e^D$. Choosing $\delta w = F_{\Gamma_e}(v)$ and using property \eqref{coercivity_a} for $a_{\Lambda_e}$ and relation \eqref{generalizedpoincare} for $(\Oc,\Ec) = (\Lambda_e,\Gamma)$, we obtain that for all $v \in \xHn{1/2}(\Gamma_e)$,
\begin{equation*}
g_{\Gamma_e}(v,v) = a_{\Lambda_e}(F_{\Gamma_e}(v),F_{\Gamma_e}(v)) \geq \alpha_a \abs{F_{\Gamma_e}(v)}_{\xHone(\Lambda_e)}^2 \geq \frac{\alpha_a}{C^2_{\Lambda_e,\Gamma}} \norm{F_{\Gamma_e}(v)}_{\xHone(\Lambda_e)}^2 \geq 
 \frac{\alpha_a}{C^2_{\Lambda_e,\Gamma} \beta_{\tau_{\Lambda_e,\Gamma_e}}^2} \norm{v}_{\xHn{1/2}(\Gamma_e)}^2,
\end{equation*}
where $\beta_{\tau_{\Lambda_e,\Gamma_e}}$ is the norm of the trace operator from $\xHone(\Lambda_e)$ to $\xHn{1/2}(\Gamma_e)$.
Then, let $\Gamma_i^D = \partial \Lambda_i \cap \Gamma_D$. The restriction $w_i(\xi)$ of the local solution $w(\xi) =\Theta(U(\xi);\xi) $ to $\Lambda_i$ is such that $w_i(\xi) \in \xHone(\Lambda_i)$, $w_i(\xi)_{\restrictto{\Gamma_e}} = w(\xi)_{\restrictto{\Gamma_e}}$, $w_i(\xi)_{\restrictto{\Gamma_i^D}}=0$ and 
\begin{equation*}
d_{\Lambda_i}(w_i(\xi),\delta w_i;\xi) = \ell_{\Lambda_i}(\delta w_i;\xi) - g_{\Gamma}(U(\xi)_{\restrictto{\Gamma}},\delta w_{i\restrictto{\Gamma_e}}) + b_{\Gamma_e}(\overline{\lambda}_e(\xi),\delta w_{i\restrictto{\Gamma_e}}),
\end{equation*}
for all $\delta w_i \in \xHone(\Lambda_i)$ such that $\delta w_i=0$ on $\Gamma_i^D$, where $d_{\Lambda_i}(\cdot,\cdot;\xi)$ is a semi-linear form defined by $d_{\Lambda_i}({u},{v};\xi) = a_{\Lambda_i}({u},{v};\xi) + n_{\Lambda_i}({u},{v};\xi) + g_{\Gamma_e}({u}_{\restrictto{\Gamma_e}},{v}_{\restrictto{\Gamma_e}})$ for $u,v\in \xHone(\Lambda_i)$. 
Now, let $\scalproda{\cdot}{\cdot}$ denote the duality pairing between $\xHone(\Lambda_i)$ and $\xHone(\Lambda_i)^{\ast}$. 
From properties of $a_{\Lambda_i}$, $n_{\Lambda_i}$ and $g_{\Gamma_e}$, we deduce that $d_{\Lambda_i}({u},{v};\xi) = \scalproda{\hat{\Dc}(\xi)(u)}{v}$ {for all $u,v \in \xHone(\Lambda_i)$}, where $\hat{\Dc}(\xi) \colon \xHone(\Lambda_i) \to \xHone(\Lambda_i)^{\ast}$ is {a} radially continuous, coercive and strongly monotone map such that for all $u,v \in \xHone(\Lambda_i)$,
\begin{align*}
\scalproda{\hat{\Dc}(\xi)({u}) - \hat{\Dc}{(\xi)}({v})}{u-v} 
&= a_{\Lambda_i}(u-v,u-v) + n_{\Lambda_i}(u,u-v) - n_{\Lambda_i}(v,u-v) + g_{\Gamma_e}(u_{\restrictto{\Gamma_e}}-v_{\restrictto{\Gamma_e}},u_{\restrictto{\Gamma_e}}-v_{\restrictto{\Gamma_e}}) \\
&\geq \alpha_a \abs{u-v}_{\xHone(\Lambda_i)}^2 + \frac{\alpha_a}{C^2_{\Lambda_e,\Gamma} \beta_{\tau_{\Lambda_e,\Gamma_e}}^2} \norm{u-v}_{\xHn{1/2}(\Gamma_e)}^2 
\geq \alpha_{\hat{\Dc}} \norm{u-v}_{\xHone(\Lambda_i)}^2,
\end{align*}
with $\alpha_{\hat{\Dc}} = \frac{\alpha_a}{2C_{\Lambda_i,\Gamma_{e}}^2} \min\set{1,\frac{1}{C^2_{\Lambda_e,\Gamma} \beta_{\tau_{\Lambda_e,\Gamma_e}}^2}}$ independent of $\xi$. Following the proof of Theorem~\ref{thrm:wellposedness_weakformulation} in Section~\ref{sec:thrm_wellposedness_weakformulation}, we obtain that the solution $w_i{(\xi)}$ is unique and can be written {as}
$w_i (\xi)= G(U(\xi)_{\restrictto{\Gamma}};\xi) {+ \overline{w}_i(\xi)}$, with $\overline{w}_i(\xi)=0$ on $\Gamma \cup \Gamma_{i}^D$, and where $G(\cdot;\xi) \colon \xHn{1/2}(\Gamma) \to \xHone(\Lambda_i)$ {is} a Lipschitz continuous map with $\beta_G = \frac{\norm{{g_{\Gamma}}} \beta_{\tau_{\Lambda_i,\Gamma_e}}}{\alpha_{\hat{\Dc}}}$ 
independent of $\xi$, where $\beta_{\tau_{\Lambda_i,\Gamma_e}}$ is the norm of the trace operator from $\xHone(\Lambda_i)$ to $\xHn{1/2}(\Gamma_e)$. Finally, we deduce that {for all $U,V \in \widetilde{\Uc}$,}
\begin{align*}
&\norm{\Theta(U;\xi) - \Theta({V};\xi)}_{{\Wc}}^2 = 
\norm{\Theta(U;\xi)_{\restrictto{\Lambda_e}} - \Theta({V};\xi)_{\restrictto{\Lambda_e}}}_{\xHone(\Lambda_e)}^2 +
\norm{\Theta(U;\xi)_{\restrictto{\Lambda_i}} - \Theta({V};\xi)_{\restrictto{\Lambda_i}}}_{\xHone(\Lambda_i)}^2 \\
&= \norm{F_{\Gamma}(U_{{\restrictto{\Gamma}}} - {V_{\restrictto{\Gamma}}})}_{\xHone(\Lambda_e)}^2 + \norm{F_{\Gamma_e}(\Theta(U;\xi)_{\restrictto{\Gamma_e}} - \Theta({V};\xi)_{\restrictto{\Gamma_e}})}_{\xHone(\Lambda_e)}^2 + \norm{G(U_{{\restrictto{\Gamma}}};\xi) - G({V_{\restrictto{\Gamma}}};\xi)}_{\xHone(\Lambda_i)}^2 \\
&\leq (\norm{F_{\Gamma}}^2 + \beta_G^2) \norm{U-V}_{\xHn{1/2}(\Gamma)}^2 + 
 \norm{F_{\Gamma_e}}^2 \norm{\Theta(U;\xi)_{\restrictto{\Gamma_e}} - \Theta({V};\xi)_{\restrictto{\Gamma_e}}}_{\xHn{1/2}({\Gamma_e})}^2 \\
& \leq (\norm{F_{\Gamma}}^2 + \beta_G^2 + \norm{F_{\Gamma_e}}^2 \beta_{\tau_{\Lambda_i,\Gamma_e}}^2 \beta_G^2) \norm{U-V}_{\xHn{1/2}(\Gamma)}^2 \\
&\leq \beta_{\Theta}^2 \norm{U-V}_{\widetilde{\Uc}}^2,
\end{align*}
with $\beta_{\Theta}^2 = (\norm{F_{\Gamma}}^2 + \beta_G^2 + \norm{F_{\Gamma_e}}^2 \beta_{\tau_{\Lambda_i,\Gamma_e}}^2 \beta_G^2) \beta_{\tau_{\Omega \setminus \Lambda,\Gamma}}^2$ independent of $\xi$, where $\beta_{\tau_{\Omega \setminus \Lambda,\Gamma}}$ is the norm of the trace operator from $\xHone(\Omega \setminus \Lambda)$ to $\xHn{1/2}(\Gamma)$. 

Using \eqref{localpbTheta} with test functions $\delta w = 0$ on $\Lambda_i$ and introducing {a linear continuous} extension operator $R_{\Gamma}$ from $\xHn{1/2}(\Gamma)$ to $\xHone(\Lambda_e)$ (with norm $\norm{R_{\Gamma}}$ independent of $\xi$) such that $R_{\Gamma}({v})= 0$ on $\Gamma_e \cup \Gamma_e^D$, we have that for all $V \in \widetilde{\Uc}$,
\begin{equation*}
b_{\Gamma}(\Psi(V;\xi),\delta w_e) = a_{\Lambda_e}(\Theta(V;\xi)_{\restrictto{\Lambda_e}},\delta w_e) - \ell_{\Lambda_e} (\delta w_e;\xi)
\end{equation*}
for all $\delta w_e \in \xHone(\Lambda_e) $ such that $\delta w_e=0$ on $\Gamma_e \cup \Gamma_e^D$. Then, we obtain that for all $U,V \in \widetilde{\Uc}$, 
\begin{align*}
\norm{\Psi(U;\xi)-\Psi(V;\xi)}_{{\Mc}} &= \sup_{\substack{v \in \xHn{1/2}(\Gamma) \\ 
\norm{v}_{\xHn{1/2}(\Gamma)}=1}} 
b_{\Gamma}(\Psi(U;\xi) - \Psi({V};\xi),{v}) \\
&= \sup_{\substack{v \in \xHn{1/2}(\Gamma) \\ 
\norm{v}_{\xHn{1/2}(\Gamma)}=1}} 
a_{\Lambda_e}(\Theta(U;\xi)_{\restrictto{\Lambda_e}} - \Theta({V};\xi)_{\restrictto{\Lambda_e}},R_{\Gamma}({v}))\\
&\leq \beta_a \norm{\Theta(U;\xi)_{\restrictto{\Lambda_e}} - \Theta({V};\xi)_{\restrictto{\Lambda_e}}}_{\xHone(\Lambda_e)} \norm{R_{\Gamma}} \\
&\leq \beta_{\Psi} \norm{U - V}_{\widetilde{\Uc}},
\end{align*}
with $\beta_{\Psi} = \beta_a \beta_{\Theta} \norm{R_{\Gamma}}$ independent of $\xi$.

\subsection{Proof of Lemma~\ref{lmm:propD}}\label{sec:lmm_propD}

From the Lipschitz continuity of mappings $\Upsilon$, $\Phi$ and $\Psi$ (see Lemmas~\ref{lmm:continuity_linearmap} and \ref{lmm:continuity_nonlinearmap}) and from the continuity of the linear map $C_{\widetilde{\Omega}}$, we deduce the Lipschitz continuity of $D(\cdot;\xi)$ on $\widetilde{\Uc}$ and \afortiori{} on $\widetilde{\Uc}_{\star}$, with Lipschitz constant $\beta_D = \beta_{c}(1 + \beta_{\Upsilon} + \beta_{\Phi} \beta_{\Psi})$ independent of $\xi$. 
From \eqref{globalpbA}, we have that for all $U,V \in \widetilde{\Uc}$,
\begin{align*}
\scalproda{D(U;\xi) - D(V;\xi)}{U-V}_{\widetilde{\Uc}} 
&= c_{\widetilde{\Omega}}(A(U;\xi) - A(V;\xi),U-V) \\
&= a_{\Omega \setminus \Lambda}(U-V,U-V) + b_{\Gamma}(\Psi(U;\xi)-\Psi(V;\xi),U-V) \\
&= a_{\Omega \setminus \Lambda}(U-V,U-V) + b_{\Gamma}(\Psi(U;\xi)-\Psi(V;\xi),\Theta(U;\xi)-\Theta(V;\xi)) \\
&= a_{\Omega \setminus \Lambda}(U-V,U-V) + a_{\Lambda}(\Theta(U;\xi)-\Theta(V;\xi),\Theta(U;\xi)-\Theta(V;\xi);\xi) \\
&\qquad + n_{\Lambda}(\Theta(U;\xi),\Theta(U;\xi)-\Theta(V;\xi);\xi) - n_{\Lambda}(\Theta(V;\xi),\Theta(U;\xi)-\Theta(V;\xi);\xi).
\end{align*}
From assumptions~\ref{assa} and \ref{assn}, we have $a_{\Lambda}(w,w;\xi) \geq 0$ and $n_{\Lambda}(w,w-w';\xi) - n_{\Lambda}(w',w-w';\xi) \geq 0$ for all $w,w' \in \Wc$. Therefore, using property \eqref{coercivity_a} for $a_{\Omega\setminus \Lambda}$ and relation \eqref{generalizedpoincare} for $(\Oc,\Ec) = (\Omega \setminus \Lambda,\Gamma_D \cap \partial (\Omega \setminus \Lambda))$, we have that for all $U,V \in \widetilde{\Uc}$,
\begin{equation*}
\scalproda{D(U;\xi) - D(V;\xi)}{U-V}_{\widetilde{\Uc}} \geq a_{\Omega \setminus \Lambda}(U-V,U-V) \geq \alpha_a \abs{U-V}^2_{\xHone(\Omega \setminus \Lambda)} \geq \frac{\alpha_a}{\widehat{C}^2} {\norm{U-V}^2_{\Uc}},
\end{equation*}
where $\widehat{C} = C_{\Omega \setminus \Lambda,\Gamma_D \cap \partial(\Omega \setminus \Lambda)}$.
Finally, using Lemma~\ref{lmm:equiv_norm_Ustar}, we obtain that for all $U,V \in \widetilde{\Uc}_{\star}$,
\begin{equation*}
\scalproda{D(U;\xi) - D(V;\xi)}{U-V}_{\widetilde{\Uc}} \geq \alpha_D \norm{U-V}_{\widetilde{\Uc}}^2,
\end{equation*}
with $\alpha_D = \frac{\alpha_a}{\widehat{C}^2 C_{\widetilde{\Uc}}^2}$ independent of $\xi$. That proves the strong monotonicity of $D(\cdot;\xi)$ on $\widetilde{\Uc}_{\star}$.

\subsection{Proof of Theorem~\ref{thrm:convergence}}\label{sec:thrm_convergence}

First recall that the global solution $U(\xi)$ and all global iterates $U^k(\xi)$ are in $\widetilde{\Uc}_{\star}$ (see Lemma~\ref{lmm:Uk_Ustar} and Theorem~\ref{thrm:consistency}). Using relation \eqref{generalizedpoincare} for $(\Oc,\Ec) = (\widetilde{\Omega},\Gamma_D \cap \partial \widetilde{\Omega})$, the symmetry and property \eqref{coercivity_c} for $c_{\widetilde{\Omega}}$, as well as the Lipschitz continuity and the strong monotonicity of map $D(\cdot;\xi)$ on $\widetilde{\Uc}_{\star}$ (see Lemma~\ref{lmm:propD}), we obtain that for all $U,V \in \widetilde{\Uc}_{\star}$,
\begin{align*}
&\norm{B_{\rho_k}(U;\xi) - B_{\rho_k}(V;\xi)}_{\widetilde{\Uc}}^2 \leq \widetilde{C}^2 \abs{B_{\rho_k}(U;\xi) - B_{\rho_k}(V;\xi)}^2_{\xHone(\widetilde{\Omega})} \\
&\leq \frac{\widetilde{C}^2}{\alpha_c} c_{\widetilde{\Omega}}(B_{\rho_k}(U) - B_{\rho_k}(V;\xi),B_{\rho_k}(U;\xi) - B_{\rho_k}(V;\xi);\xi) \\
&= \frac{\widetilde{C}^2}{\alpha_c} \(c_{\widetilde{\Omega}}(U-V,U-V) -2 \rho_k c_{\widetilde{\Omega}}(A(U;\xi) - A(V;\xi),U-V) + \rho_k^2 c_{\widetilde{\Omega}}(A(U;\xi) - A(V;\xi),A(U;\xi) - A(V;\xi))\) \\
&\leq \frac{\widetilde{C}^2}{\alpha_c} \(\scalproda{C_{\widetilde{\Omega}}(U-V)}{U-V}_{\widetilde{\Uc}} -2 \rho_k \scalproda{D(U;\xi) - D(V;\xi)}{U-V}_{\widetilde{\Uc}} + \rho_k^2 \norm{D(U;\xi) - D(V;\xi)}_{\widetilde{\Uc}} \norm{A(U;\xi) - A(V;\xi)}_{\widetilde{\Uc}}\) \\
&\leq \frac{\widetilde{C}^2}{\alpha_c} \(\beta_c - 2 \rho_k \alpha_D + \rho_k^2 \frac{\beta_D^2}{\beta_c}\) \norm{U-V}_{\widetilde{\Uc}}^2,
\end{align*}
where $\widetilde{C}=C_{\widetilde{\Omega},\Gamma_D \cap \partial \widetilde{\Omega}}$. If $\set{\rho_k}_{k \in \xN}$ satisfies \eqref{convergencecondition} with $\rho_{\mathrm{sup}} < \frac{\beta_c}{\beta_D^2}\(2\alpha_D-\frac{1}{\rho_{\mathrm{inf}}} \(\beta_c - \frac{\alpha_c}{\widetilde{C}^2}\)\) \coloneqq \rho_{\mathrm{sup}}^{\ast}$, then the set of mappings $\set{B_{\rho_k}(\cdot;\xi)}_{k \in \xN}$ is uniformly contractive on $\widetilde{\Uc}_{\star}$, with a contractivity constant $\rho_B = \(\frac{\widetilde{C}^2}{\alpha_c} (\beta_c - \rho_{\mathrm{inf}} (2 \alpha_D - \rho_{\mathrm{sup}} \frac{\beta_D^2}{\beta_c}))\)^{1/2} < 1$ independent of $\xi$. Then, we obtain that 
\begin{equation*}
\norm{U^k(\xi) - U(\xi)}_{\widetilde{\Uc}} = \norm{B_{\rho_k}(U^{k-1}(\xi);\xi) - B_{\rho_k}(U(\xi);\xi)}_{\widetilde{\Uc}} \leq \rho_B^k \norm{U(\xi) - U^0(\xi)}_{\widetilde{\Uc}},
\end{equation*}
from which we deduce that the sequence $\set{U^k(\xi)}_{k \in \xN}$ converges almost surely to $U(\xi)$ in $\widetilde{\Uc}_{\star}$. Also, with $U^0=0$, we obtain that $\norm{U^k(\xi) - U(\xi)}_{\widetilde{\Uc}} \leq \norm{U(\xi)}_{\widetilde{\Uc}}$. Since $\norm{U(\xi)}_{\widetilde{\Uc}} \in \xLn{p}_{\mu}(\Xi)$ (see Corollary~\ref{U_Lp}), the dominated convergence theorem gives that the sequence $\set{U^k}_{k \in \xN}$ converges to $U$ in $\xLn{p}_{\mu}(\Xi;\widetilde{\Uc})$.
Finally, from \eqref{local_error} and using the Lipschitz continuity of mappings $\Theta$ and $\Psi$ {(see Lemma~\ref{lmm:continuity_nonlinearmap})}, we directly obtain that the sequence $\set{w^k}_{k \in \xN}$ converges to $w$ almost surely and in $\xLn{p}_{\mu}(\Xi;\Wc)$, and the sequence $\set{\lambda^k}_{k \in \xN}$ converges to $\lambda$ almost surely and in $\xLn{p}_{\mu}(\Xi;\Mc)$.

\subsection{Proof of Theorem~\ref{thrm:robustness}}\label{sec:thrm_robustness}

First, if initial guess $U_{\varepsilon}^0 \in \Vc_{\delta}$ and if $\varepsilon^{\ast} < 1 - \rho_B$ and $\varepsilon \leq \delta (1 - \rho_B - \varepsilon^{\ast})$, we can prove by induction that all approximate global iterates $U_{\varepsilon}^k$ belong to $\Vc_{\delta}$. Indeed, suppose that $U_{\varepsilon}^j \in \Vc_{\delta}$ for all $j<k$. Then, the error at iteration $k$ is such that
\begin{align*}
\norm{U_{\varepsilon}^k - U}_{\xLn{p}_{\mu}(\Xi;\widetilde{\Uc})} &\leq \norm{U_{\varepsilon}^k - U^k}_{\xLn{p}_{\mu}(\Xi;\widetilde{\Uc})}+ \norm{U^k - U}_{\xLn{p}_{\mu}(\Xi;\widetilde{\Uc})}\\
&\leq \norm{B_{\rho_k}^{\varepsilon}(U_{\varepsilon}^{k-1}) - B_{\rho_k}(U_{\varepsilon}^{k-1})}_{\xLn{p}_{\mu}(\Xi;\widetilde{\Uc})} + \norm{B_{\rho_k}(U_{\varepsilon}^{k-1}) - B_{\rho_k}(U)}_{\xLn{p}_{\mu}(\Xi;\widetilde{\Uc})}\\
&\leq \varepsilon \norm{U}_{\xLn{p}_{\mu}(\Xi;\widetilde{\Uc})} + (\rho_B + \varepsilon^{\ast}) \norm{U_{\varepsilon}^{k-1}-U}_{\xLn{p}_{\mu}(\Xi;\widetilde{\Uc})}.
\end{align*}
As $U_{\varepsilon}^{k-1} \in \Vc_{\delta}$ and $\varepsilon \leq \delta (1 - \rho_B - \varepsilon^{\ast})$, $U_{\varepsilon}^k \in \Vc_{\delta}$. By induction, we finally prove that $U_{\varepsilon}^k \in \Vc_{\delta}$ for all $k \in \xN$.
It means that if initial global iterate $U_{\varepsilon}^0$ is contained in the open ball $\Vc_{\delta}$ of radius $\delta \norm{U}_{\xLn{p}_{\mu}(\Xi;\widetilde{\Uc})}$ centered at the exact global solution $U$, then all global iterates $U_{\varepsilon}^k$ remain in this ball. Subsequently, we obtain that
\begin{align*}
\norm{U_{\varepsilon}^k - U}_{\xLn{p}_{\mu}(\Xi;\widetilde{\Uc})} &\leq \varepsilon \norm{U}_{\xLn{p}_{\mu}(\Xi;\widetilde{\Uc})} \sum_{j=0}^{k-1} (\rho_B + \varepsilon^{\ast})^j + (\rho_B + \varepsilon^{\ast})^k \norm{U_{\varepsilon}^0-U}_{\xLn{p}_{\mu}(\Xi;\widetilde{\Uc})}\\
&\leq \dfrac{\varepsilon \(1 - (\rho_B + \varepsilon^{\ast})^k \)}{1 - (\rho_B + \varepsilon^{\ast})} \norm{U}_{\xLn{p}_{\mu}(\Xi;\widetilde{\Uc})} + (\rho_B + \varepsilon^{\ast})^k \norm{U_{\varepsilon}^0-U}_{\xLn{p}_{\mu}(\Xi;\widetilde{\Uc})}\\
&\leq \dfrac{\varepsilon}{1 - (\rho_B + \varepsilon^{\ast})} \norm{U}_{\xLn{p}_{\mu}(\Xi;\widetilde{\Uc})} + (\rho_B + \varepsilon^{\ast})^k \norm{U_{\varepsilon}^0-U}_{\xLn{p}_{\mu}(\Xi;\widetilde{\Uc})}
\end{align*}
and therefore, as $0 < \rho_B + \varepsilon^{\ast} < 1$, the approximate sequence $\set{U_{\varepsilon}^k}_{k \in \xN}$ satisfies \eqref{robustness_Lp}, with $\gamma(\varepsilon,\varepsilon^{\ast}) = \frac{\varepsilon}{1 - (\rho_B + \varepsilon^{\ast})} \to 0$ as $\varepsilon \to 0$.

%%-----------------------------
%%      your bibliography
%%-----------------------------
\bibliographystyle{unsrt}
\bibliography{Biblio}

\end{document}

%% file: partition.tex
\begin{tikzpicture}[scale=0.8,every node/.style={minimum size=1cm},on grid]

% Domain level
\begin{scope}[
	yshift=0,
	every node/.append style={
	%yslant=0.5,xslant=-1.5,
	scale=1.2},
	yslant=0.5,xslant=-1.5
	]
	% Domain coordinates
	\coordinate (A) at (0,0);
	\coordinate (B) at (5,5);
	\coordinate (C) at (0,5);
	\coordinate (D) at (5,0);
	% Domain
	\shadedraw[right color=cyan!10,left color=cyan!50,black,thick,fill=cyan!80,fill opacity=0.75] (A) rectangle (B);
	\node at ($(C)!0.2!(D)$) {$\Omega \setminus \Lambda$};
	%\node[above left] at ($(B)!0.8!(C)$) {$\Omega(\xi) = \(\Omega \setminus \Lambda\) \cup \Lambda(\xi)$};
\end{scope}

% Patch level
\begin{scope}[
	yshift=0,
	every node/.append style={
	%yslant=0.5,xslant=-1.5,
	scale=1.2},
	yslant=0.5,xslant=-1.5,
	decoration=penciline
	]
	% Patch coordinates
	\coordinate (I) at (2.5,2);
	\coordinate (J) at (4,3.5);
	\coordinate (K) at (2.5,3.5);
	\coordinate (L) at (4,2);
	% Patch
	\shadedraw[right color=red!10,left color=red!50,black,%dashed,
	thick,fill=red!60,fill opacity=0.75] (I) rectangle (J);
	\node at ($(I)!0.6!(J)!0.4!(K)$) {$\Lambda$};
	\draw[decorate,thick,black,fill=white] ($(K)!0.5!(L)$) circle (0.2cm);
	\draw[decorate,thick,black,fill=white] ($(I)!0.75!(J)$) circle (0.1cm);
	\draw[decorate,thick] ($(K)!0.8!(L)!0.2!(I)$) -- ($(I)!0.7!(L)!0.4!(J)$);
	\draw[decorate,thick,black] ($(I)!0.25!(J)$) arc (0:-10:1cm);
	% Interface
	\node[below right=-0.5] at ($(I)!0.5!(L)$) {$\Gamma$};
\end{scope}

\end{tikzpicture}

%% file: fictitious_patch.tex
\begin{tikzpicture}[scale=0.8,every node/.style={minimum size=1cm},on grid]

% Fictitious level
\begin{scope}[
	yshift=0,
	every node/.append style={
	%yslant=0.5,xslant=-1.5,
	anchor=center,scale=1.2},
	yslant=0.5,xslant=-1.5
	]
	% Fictitious domain coordinates
	\coordinate (A) at (0,0);
	\coordinate (B) at (5,5);
	\coordinate (C) at (0,5);
	\coordinate (D) at (5,0);
	% Fictitious domain
	\shadedraw[right color=cyan!10,left color=cyan!50,black,thick,fill=cyan!80,fill opacity=0.75] (A) rectangle (B);
	\node at ($(C)!0.2!(D)$) {$\widetilde{\Omega}$};
	%\node at ($(C)!0.2!(D)$) {$\Omega \setminus \Lambda$};
	%\node at ($(A)!0.2!(B)$) {$\widetilde{\Omega}$};
	%\node[left] at ($(B)!0.3!(C)$) {$U$};
	% Fictitious patch coordinates
	\coordinate (E) at (2.5,2);
	\coordinate (F) at (4,3.5);
	\coordinate (G) at (2.5,3.5);
	\coordinate (H) at (4,2);
	% Fictitious patch
	\draw[black,dashed,thick] (E) rectangle (F);
	\node at ($(E)!0.6!(F)!0.2!(G)$) {$\widetilde{\Lambda}$};
	\node[below right=-0.5] at ($(E)!0.5!(H)$) {$\Gamma$};
\end{scope}

% Patch level
\begin{scope}[
	yshift=60,
	every node/.append style={
	%yslant=0.5,xslant=-1.5,
	scale=1.2},
	yslant=0.5,xslant=-1.5
	]
	% Real patch coordinates
	\coordinate (I) at (2.5,2);
	\coordinate (J) at (4,3.5);
	\coordinate (K) at (2.5,3.5);
	\coordinate (L) at (4,2);
\end{scope}

% Vertical lines for linking both levels
\draw[black,dashed,thin] (E) to (I);
\draw[black,dashed,thin] (F) to (J);
\draw[black,dashed,thin] (G) to (K);
\draw[black,dashed,thin] (H) to (L);

% Patch level
\begin{scope}[
	yshift=60,
	every node/.append style={
	%yslant=0.5,xslant=-1.5,
	scale=1.2},
	yslant=0.5,xslant=-1.5,
	decoration=penciline
	]
	% Real patch
	\shadedraw[right color=red!10,left color=red!50,black,thick,fill=red!60,fill opacity=0.75] (I) rectangle (J);
	\node at ($(I)!0.6!(J)!0.4!(K)$) {$\Lambda$};
	\draw[decorate,thick,black,fill=white] ($(K)!0.5!(L)$) circle (0.2cm);
	\draw[decorate,thick,black,fill=white] ($(I)!0.75!(J)$) circle (0.1cm);
	\draw[decorate,thick] ($(K)!0.8!(L)!0.2!(I)$) -- ($(I)!0.7!(L)!0.4!(J)$);
	\draw[decorate,thick,black] ($(I)!0.25!(J)$) arc (0:-10:1cm);
	%\node[left] at ($(J)!0.2!(K)$) {$w$};
\end{scope}

\end{tikzpicture}

%% file: domain_global_patches.tex
\definecolor{mycolor1}{rgb}{1.00000,1.00000,0.00000}%
\definecolor{mycolor2}{rgb}{1.00000,0.00000,1.00000}%
\definecolor{mycolor3}{rgb}{0.00000,1.00000,1.00000}%
\begin{tikzpicture}

\begin{axis}[%
width=0.125\figureheight,
height=\figureheight,
at={(0\figureheight,0\figureheight)},
scale only axis,
every outer x axis line/.append style={black},
every x tick label/.append style={font=\color{black}},
every x tick/.append style={black},
xmin=0,
xmax=2,
every outer y axis line/.append style={black},
every y tick label/.append style={font=\color{black}},
every y tick/.append style={black},
ymin=0,
ymax=16,
axis line style={draw=none},
ticks=none,
axis x line*=bottom,
axis y line*=left,
legend style={legend cell align=left, align=left, draw=black},
legend to name=domainpatches_iso
]

\addplot[area legend, draw=black, fill=mycolor1]
table[row sep=crcr] {%
x	y\\
0	0\\
2	0\\
2	16\\
0	16\\
}--cycle;
\addlegendentry{$\Omega \setminus \Lambda$}

\addplot[area legend, draw=black, fill=mycolor2]
table[row sep=crcr] {%
x	y\\
0.5	0.5\\
1.5	0.5\\
1.5	1.5\\
0.5	1.5\\
}--cycle;
\addlegendentry{$\Lambda_{1}$}

\addplot[area legend, draw=black, fill=mycolor3]
table[row sep=crcr] {%
x	y\\
0.5	2.5\\
1.5	2.5\\
1.5	3.5\\
0.5	3.5\\
}--cycle;
\addlegendentry{$\Lambda_{2}$}

\addplot[area legend, draw=black, fill=red]
table[row sep=crcr] {%
x	y\\
0.5	4.5\\
1.5	4.5\\
1.5	5.5\\
0.5	5.5\\
}--cycle;
\addlegendentry{$\Lambda_{3}$}

\addplot[area legend, draw=black, fill=green]
table[row sep=crcr] {%
x	y\\
0.5	6.5\\
1.5	6.5\\
1.5	7.5\\
0.5	7.5\\
}--cycle;
\addlegendentry{$\Lambda_{4}$}

\addplot[area legend, draw=black, fill=blue]
table[row sep=crcr] {%
x	y\\
0.5	8.5\\
1.5	8.5\\
1.5	9.5\\
0.5	9.5\\
}--cycle;
\addlegendentry{$\Lambda_{5}$}

\addplot[area legend, draw=black, fill=gray]
table[row sep=crcr] {%
x	y\\
0.5	10.5\\
1.5	10.5\\
1.5	11.5\\
0.5	11.5\\
}--cycle;
\addlegendentry{$\Lambda_{6}$}

\addplot[area legend, draw=black, fill=black!50!mycolor1]
table[row sep=crcr] {%
x	y\\
0.5	12.5\\
1.5	12.5\\
1.5	13.5\\
0.5	13.5\\
}--cycle;
\addlegendentry{$\Lambda_{7}$}

\addplot[area legend, draw=black, fill=violet]
table[row sep=crcr] {%
x	y\\
0.5	14.5\\
1.5	14.5\\
1.5	15.5\\
0.5	15.5\\
}--cycle;
\addlegendentry{$\Lambda_{8}$}

\end{axis}
\end{tikzpicture}%

%% file: error_rho.tex
\definecolor{mycolor1}{rgb}{1.00000,0.00000,1.00000}%
\definecolor{mycolor2}{rgb}{0.00000,1.00000,1.00000}%
\begin{tikzpicture}

\begin{axis}[%
width=1.268\figureheight,
height=\figureheight,
at={(0\figureheight,0\figureheight)},
scale only axis,
separate axis lines,
every outer x axis line/.append style={black},
every x tick label/.append style={font=\color{black}},
every x tick/.append style={black},
xmin=0,
xmax=20,
xlabel={Number of iterations},
every outer y axis line/.append style={black},
every y tick label/.append style={font=\color{black}},
every y tick/.append style={black},
ymode=log,
ymin=1e-06,
ymax=1,
yminorticks=true,
ylabel={Error indicator},
axis background/.style={fill=white},
xmajorgrids,
ymajorgrids,
yminorgrids,
legend style={legend cell align=left, align=left, draw=black}
]
\addplot [color=mycolor1, line width=1.0pt]
  table[row sep=crcr]{%
0	1\\
1	0.862214370650809\\
2	0.688998724002749\\
3	0.550483507184351\\
4	0.439742366895042\\
5	0.351227887676224\\
6	0.280495284632231\\
7	0.223983823075717\\
8	0.178842347685651\\
9	0.142788637783672\\
10	0.113996794085151\\
11	0.0910065139700039\\
12	0.0726502443600714\\
13	0.0579948828271758\\
14	0.046294917114044\\
15	0.0369547988951799\\
16	0.0294987409638185\\
17	0.023546806083421\\
18	0.0187957382187601\\
19	0.0150032538138999\\
20	0.0119759706373429\\
};
\addlegendentry{Fixed relaxation $\rho=0.2$};

\addplot [color=mycolor2, line width=1.0pt]
  table[row sep=crcr]{%
0	1\\
1	0.72455841686865\\
2	0.433284749666268\\
3	0.258847268742113\\
4	0.154528621790802\\
5	0.0922105033902803\\
6	0.0550094917637578\\
7	0.0328119599418402\\
8	0.0195702276874062\\
9	0.0116721185372355\\
10	0.00696157757255279\\
11	0.00415226310246421\\
12	0.00247649875794521\\
13	0.00147714283417654\\
14	0.000881133700519013\\
15	0.000525688223966771\\
16	0.000313652659846821\\
17	0.000187197449023181\\
18	0.000111739075324327\\
19	6.65898325169523e-05\\
20	3.9793045065413e-05\\
};
\addlegendentry{Fixed relaxation $\rho=0.4$};

\addplot [color=red, line width=1.0pt]
  table[row sep=crcr]{%
0	1\\
1	0.587123356284947\\
2	0.232866635245057\\
3	0.0920619220553043\\
4	0.0363445807398275\\
5	0.014341741597183\\
6	0.0056590525414231\\
7	0.00223318668387288\\
8	0.000881559738457987\\
9	0.000348191319529526\\
10	0.000137514515097248\\
11	5.42796737563882e-05\\
12	2.16708439928134e-05\\
13	9.00858542237297e-06\\
14	4.18643155596775e-06\\
15	2.89633025014918e-06\\
16	2.28557650411268e-06\\
17	2.51042235200969e-06\\
18	2.41690350995381e-06\\
19	2.66766567626404e-06\\
20	2.87269272565215e-06\\
};
\addlegendentry{Fixed relaxation $\rho=0.6$};

\addplot [color=blue, line width=1.0pt]
  table[row sep=crcr]{%
0	1\\
1	0.450111574031453\\
2	0.0877250363735664\\
3	0.0169254488053321\\
4	0.00326247937492746\\
5	0.000629598157999198\\
6	0.000121811715156641\\
7	2.37959139104408e-05\\
8	5.61264210653049e-06\\
9	2.9592958121934e-06\\
10	2.67547725274817e-06\\
11	2.54344906296527e-06\\
12	2.75582004945523e-06\\
13	2.91197801295022e-06\\
14	2.447135566112e-06\\
15	2.7245783062489e-06\\
16	2.46626177228876e-06\\
17	2.88686036459867e-06\\
18	2.82658474767089e-06\\
19	2.80983043177813e-06\\
20	2.60981560345364e-06\\
};
\addlegendentry{Fixed relaxation $\rho=0.8$};

\addplot [color=green, line width=1.0pt]
  table[row sep=crcr]{%
0	1\\
1	0.314077506677132\\
2	0.00265652161044692\\
3	6.18648261019734e-05\\
4	4.94166410189848e-06\\
5	3.13034331575921e-06\\
6	3.0087583221405e-06\\
7	2.89968893965285e-06\\
8	2.9526425651444e-06\\
9	2.87200503281611e-06\\
10	3.10079914944099e-06\\
11	2.63276118875778e-06\\
12	2.87242704181046e-06\\
13	2.75878350434004e-06\\
14	3.62994158292854e-06\\
15	3.1851714169681e-06\\
16	3.07847475126319e-06\\
17	3.12009721317716e-06\\
18	3.08457645603045e-06\\
19	3.02317077536545e-06\\
20	2.85520074248011e-06\\
};
\addlegendentry{Fixed relaxation $\rho=1$};

\addplot [color=orange, line width=1.0pt]
  table[row sep=crcr]{%
0	1\\
1	0.181236279110016\\
2	0.0381409911141259\\
3	0.0081920317256996\\
4	0.00175777881193758\\
5	0.000378600897495022\\
6	8.17265476088537e-05\\
7	1.82776921543423e-05\\
8	5.10279057780435e-06\\
9	3.50540379374344e-06\\
10	3.21913540269491e-06\\
11	3.79283768644457e-06\\
12	3.2334728198983e-06\\
13	3.33702286152444e-06\\
14	3.07164147744175e-06\\
15	3.75751320507889e-06\\
16	3.4447391718777e-06\\
17	3.46793351886072e-06\\
18	3.376420210437e-06\\
19	3.07951262156194e-06\\
20	3.16406882804136e-06\\
};
\addlegendentry{Fixed relaxation $\rho=1.2$};

\addplot [color=red!50!green, line width=1.0pt]
  table[row sep=crcr]{%
0	1\\
1	0.0721536906362558\\
2	0.0293545495773964\\
3	0.0120540315796554\\
4	0.00495637090062224\\
5	0.00204679074727727\\
6	0.00084750487704641\\
7	0.0003520905784077\\
8	0.000146574029711823\\
9	6.13218179238547e-05\\
10	2.59769017091049e-05\\
11	1.12428796320697e-05\\
12	6.1654358582859e-06\\
13	4.18347862199472e-06\\
14	4.32928650669389e-06\\
15	3.81329777061698e-06\\
16	4.03165594112894e-06\\
17	3.87535581793133e-06\\
18	4.10881545742817e-06\\
19	4.37589797201692e-06\\
20	3.75591226975599e-06\\
};
\addlegendentry{Fixed relaxation $\rho=1.4$};

\addplot [color=violet, line width=1.0pt]
  table[row sep=crcr]{%
0	1\\
1	0.125392897808761\\
2	0.075572655205153\\
3	0.0455814825860962\\
4	0.0275187370462643\\
5	0.016630812915022\\
6	0.0100630400989888\\
7	0.0060977552794876\\
8	0.00370086147167178\\
9	0.00225073854777016\\
10	0.00137202485371628\\
11	0.00083846331976764\\
12	0.000514050683004162\\
13	0.000316306974392692\\
14	0.000195555561534125\\
15	0.000121655413333122\\
16	7.59249531898025e-05\\
17	4.79895797326116e-05\\
18	3.04817602979768e-05\\
19	1.98550438055019e-05\\
20	1.29642562968396e-05\\
};
\addlegendentry{Fixed relaxation $\rho=1.6$};

\addplot [color=teal, line width=1.0pt]
  table[row sep=crcr]{%
0	1\\
1	0.253768122724245\\
2	0.205688431328953\\
3	0.166560483805079\\
4	0.135131676373545\\
5	0.10958212043056\\
6	0.089011282785836\\
7	0.0722974516956036\\
8	0.0588130376607904\\
9	0.0478583662086816\\
10	0.039003836185522\\
11	0.0318086735668558\\
12	0.0259827874199739\\
13	0.0212464866675467\\
14	0.0174048885545437\\
15	0.0142790098829337\\
16	0.0117393895665503\\
17	0.00967060862322355\\
18	0.00798641432940317\\
19	0.00661248265550844\\
20	0.00549173147560807\\
};
\addlegendentry{Fixed relaxation $\rho=1.8$};

\addplot [color=black!50!red, line width=1.0pt]
  table[row sep=crcr]{%
0	1\\
1	0.314077506677132\\
2	0.00265659420899902\\
3	4.79927261748233e-05\\
4	4.87680461021245e-06\\
5	3.45782419364074e-06\\
6	2.83654088534255e-06\\
7	2.62728552231246e-06\\
8	2.57642623027019e-06\\
9	2.47678789610921e-06\\
10	2.4593297403098e-06\\
11	2.45114707472263e-06\\
12	2.44287719788855e-06\\
13	2.44072823964957e-06\\
14	2.43749192401308e-06\\
15	2.43662260274237e-06\\
16	2.43649336029344e-06\\
17	2.43640461446543e-06\\
18	2.4362731247986e-06\\
19	2.43623578707055e-06\\
20	2.43620773743519e-06\\
};
\addlegendentry{Aitken's dynamic relaxation $\rho_k$};

\end{axis}
\end{tikzpicture}%

%% file: error_tol_no_legend.tex
\definecolor{mycolor1}{rgb}{1.00000,0.00000,1.00000}%
\definecolor{mycolor2}{rgb}{0.00000,1.00000,1.00000}%
\begin{tikzpicture}

\begin{axis}[%
width=1.268\figureheight,
height=\figureheight,
at={(0\figureheight,0\figureheight)},
scale only axis,
separate axis lines,
every outer x axis line/.append style={black},
every x tick label/.append style={font=\color{black}},
every x tick/.append style={black},
xmin=0,
xmax=20,
xlabel={Number of iterations},
every outer y axis line/.append style={black},
every y tick label/.append style={font=\color{black}},
every y tick/.append style={black},
ymode=log,
ymin=1e-08,
ymax=1,
yminorticks=true,
ylabel={Error indicator},
ytick={1,1e-1,1e-2,1e-3,1e-4,1e-5,1e-6,1e-7,1e-8},
axis background/.style={fill=white},
xmajorgrids,
ymajorgrids,
yminorgrids,
legend style={legend cell align=left, align=left, draw=black}
]
\addplot [color=mycolor1, line width=1.0pt]
  table[row sep=crcr]{%
0	1\\
1	0.314077506677132\\
2	0.00265609546186943\\
3	5.75946150772207e-05\\
4	2.52629323652626e-05\\
5	2.12082347837431e-05\\
6	2.11781838501828e-05\\
7	2.09340223899649e-05\\
8	2.04044550357652e-05\\
9	2.0106644912903e-05\\
10	2.00951180394087e-05\\
11	2.01001038773484e-05\\
12	2.00764977537806e-05\\
13	2.00687438519246e-05\\
14	2.00611874619156e-05\\
15	2.00577073299354e-05\\
16	2.0054780260824e-05\\
17	2.00534499401982e-05\\
18	2.00532798211468e-05\\
19	2.00529437293365e-05\\
20	2.00527913043235e-05\\
};
%\addlegendentry{$\varepsilon_{\text{cv}}=10^{-2}$};

\addplot [color=mycolor2, line width=1.0pt]
  table[row sep=crcr]{%
0	1\\
1	0.314077506677132\\
2	0.00265659420899902\\
3	4.79927261748233e-05\\
4	4.87680461021245e-06\\
5	3.45782419364074e-06\\
6	2.83654088534255e-06\\
7	2.62728552231246e-06\\
8	2.57642623027019e-06\\
9	2.47678789610921e-06\\
10	2.4593297403098e-06\\
11	2.45114707472263e-06\\
12	2.44287719788855e-06\\
13	2.44072823964957e-06\\
14	2.43749192401308e-06\\
15	2.43662260274237e-06\\
16	2.43649336029344e-06\\
17	2.43640461446543e-06\\
18	2.4362731247986e-06\\
19	2.43623578707055e-06\\
20	2.43620773743519e-06\\
};
%\addlegendentry{$\varepsilon_{\text{cv}}=10^{-3}$};

\addplot [color=red, line width=1.0pt]
  table[row sep=crcr]{%
0	1\\
1	0.314077506677132\\
2	0.00265640365640891\\
3	4.78523833599226e-05\\
4	4.04348137084533e-06\\
5	6.49135785322948e-07\\
6	2.51136681477065e-07\\
7	2.22256925413535e-07\\
8	1.89121667861771e-07\\
9	1.76020666832373e-07\\
10	1.70828631733564e-07\\
11	1.68735168629954e-07\\
12	1.67338691264245e-07\\
13	1.66621370019213e-07\\
14	1.66343842445902e-07\\
15	1.66260062953711e-07\\
16	1.6622704449577e-07\\
17	1.6619348906532e-07\\
18	1.66177161325065e-07\\
19	1.66174409987178e-07\\
20	1.66167531975762e-07\\
};
%\addlegendentry{$\varepsilon_{\text{cv}}=10^{-4}$};

\addplot [color=green, line width=1.0pt]
  table[row sep=crcr]{%
0	1\\
1	0.314077506677132\\
2	0.00265639414400724\\
3	4.78528128971129e-05\\
4	4.03244365966293e-06\\
5	5.98714535551216e-07\\
6	4.84095959443581e-08\\
7	3.25758766624387e-08\\
8	2.35100965180228e-08\\
9	2.20991076990192e-08\\
10	2.04759102347364e-08\\
11	2.02002625554687e-08\\
12	1.9926554218181e-08\\
13	1.98003998512144e-08\\
14	1.97336741648814e-08\\
15	1.97071359023912e-08\\
16	1.96833814385193e-08\\
17	1.96782473160647e-08\\
18	1.96737778950716e-08\\
19	1.96711395238354e-08\\
20	1.96698306629676e-08\\
};
%\addlegendentry{$\varepsilon_{\text{cv}}=10^{-5}$};

\end{axis}
\end{tikzpicture}%

%% file: error_tol.tex
\definecolor{mycolor1}{rgb}{1.00000,0.00000,1.00000}%
\definecolor{mycolor2}{rgb}{0.00000,1.00000,1.00000}%
\begin{tikzpicture}

\begin{axis}[%
width=1.268\figureheight,
height=\figureheight,
at={(0\figureheight,0\figureheight)},
scale only axis,
separate axis lines,
every outer x axis line/.append style={black},
every x tick label/.append style={font=\color{black}},
every x tick/.append style={black},
xmin=0,
xmax=20,
xlabel={Number of iterations},
every outer y axis line/.append style={black},
every y tick label/.append style={font=\color{black}},
every y tick/.append style={black},
ymode=log,
ymin=1e-08,
ymax=1,
yminorticks=true,
ylabel={Error indicator},
ytick={1,1e-1,1e-2,1e-3,1e-4,1e-5,1e-6,1e-7,1e-8},
axis background/.style={fill=white},
xmajorgrids,
ymajorgrids,
yminorgrids,
legend style={legend cell align=left, align=left, draw=black}
]
\addplot [color=mycolor1, line width=1.0pt]
  table[row sep=crcr]{%
0	1\\
1	0.315720459145447\\
2	0.00130323635212895\\
3	3.81947212926923e-05\\
4	1.85121582672015e-05\\
5	1.52709703155066e-05\\
6	1.51223042457269e-05\\
7	1.39490748079489e-05\\
8	1.38760393985822e-05\\
9	1.38651193457628e-05\\
10	1.38645885055502e-05\\
11	1.38323642216946e-05\\
12	1.38045388708349e-05\\
13	1.38030990240004e-05\\
14	1.38011693812321e-05\\
15	1.37988488861961e-05\\
16	1.37994167848322e-05\\
17	1.37985809232917e-05\\
18	1.37983488614551e-05\\
19	1.37983924133059e-05\\
20	1.37983802744092e-05\\
};
\addlegendentry{$\varepsilon_{\text{cv}}=10^{-2}$};

\addplot [color=mycolor2, line width=1.0pt]
  table[row sep=crcr]{%
0	1\\
1	0.315720459145447\\
2	0.00130275262856929\\
3	3.3479643875952e-05\\
4	4.58312350535212e-06\\
5	2.30563091174599e-06\\
6	2.21890321792615e-06\\
7	2.1651394739217e-06\\
8	2.09092247356788e-06\\
9	2.05086241096813e-06\\
10	2.02880569968408e-06\\
11	2.01379219890618e-06\\
12	2.00868315035051e-06\\
13	2.00759174554654e-06\\
14	2.00580441547434e-06\\
15	2.00499071316946e-06\\
16	2.00485551682053e-06\\
17	2.0048010248984e-06\\
18	2.00471991944203e-06\\
19	2.00467711550699e-06\\
20	2.0046583752109e-06\\
};
\addlegendentry{$\varepsilon_{\text{cv}}=10^{-3}$};

\addplot [color=red, line width=1.0pt]
  table[row sep=crcr]{%
0	1\\
1	0.315720459145447\\
2	0.00130335339065967\\
3	3.30413422448043e-05\\
4	3.7862251816073e-06\\
5	5.58819685413194e-07\\
6	2.9037291260901e-07\\
7	2.68471517007941e-07\\
8	2.45908781087157e-07\\
9	2.3663520713563e-07\\
10	2.29224017205511e-07\\
11	2.24344207176785e-07\\
12	2.23545517745786e-07\\
13	2.22735394526003e-07\\
14	2.22341947148948e-07\\
15	2.2210596991645e-07\\
16	2.22036695914144e-07\\
17	2.21995000582594e-07\\
18	2.21959046709272e-07\\
19	2.21950753366795e-07\\
20	2.21945383812017e-07\\
};
\addlegendentry{$\varepsilon_{\text{cv}}=10^{-4}$};

\addplot [color=green, line width=1.0pt]
  table[row sep=crcr]{%
0	1\\
1	0.315720459145447\\
2	0.00130331826397226\\
3	3.30299292444694e-05\\
4	3.76976119909163e-06\\
5	4.33770988506494e-07\\
6	3.48081236107827e-08\\
7	2.98194720441073e-08\\
8	2.44520461757631e-08\\
9	2.07469829075955e-08\\
10	2.02060434378388e-08\\
11	1.99433247020555e-08\\
12	1.96866589081263e-08\\
13	1.95875359439635e-08\\
14	1.95427210135112e-08\\
15	1.95178048217075e-08\\
16	1.95034632517949e-08\\
17	1.94976315641165e-08\\
18	1.94948396235848e-08\\
19	1.94932795076467e-08\\
20	1.94922994487458e-08\\
};
\addlegendentry{$\varepsilon_{\text{cv}}=10^{-5}$};

\end{axis}
\end{tikzpicture}%

%% file: cpu_time_no_legend.tex
\definecolor{mycolor1}{rgb}{1.00000,0.00000,1.00000}%
\definecolor{mycolor2}{rgb}{0.00000,1.00000,1.00000}%
\begin{tikzpicture}

\begin{axis}[%
width=1.268\figureheight,
height=\figureheight,
at={(0\figureheight,0\figureheight)},
scale only axis,
separate axis lines,
every outer x axis line/.append style={black},
every x tick label/.append style={font=\color{black}},
every x tick/.append style={black},
xmin=0,
xmax=20,
xlabel={Number of iterations},
every outer y axis line/.append style={black},
every y tick label/.append style={font=\color{black}},
every y tick/.append style={black},
ymin=0,
ymax=450,
ylabel={CPU time (s)},
axis background/.style={fill=white},
xmajorgrids,
ymajorgrids
]
\addplot [color=mycolor1, line width=1.0pt]
  table[row sep=crcr]{%
1	106.186672\\
2	110.982003\\
3	116.042764\\
4	123.657356\\
5	108.587467\\
6	120.329155\\
7	116.744587\\
8	119.44611\\
9	119.517767\\
10	117.624156\\
11	114.181662\\
12	131.16139\\
13	131.21148\\
14	116.439152\\
15	127.410593\\
16	119.402162\\
17	109.915048\\
18	112.993006\\
19	112.763079\\
20	123.90161\\
};
%\addlegendentry{$\varepsilon_{\text{cv}}=10^{-2}$};

\addplot [color=mycolor2,solid,line width=1.0pt]
  table[row sep=crcr]{%
1	102.698832\\
2	137.795294\\
3	146.412154\\
4	163.898746\\
5	138.560114\\
6	152.81248\\
7	143.159454\\
8	135.255021\\
9	147.872734\\
10	148.472218\\
11	132.005092\\
12	133.589049\\
13	143.476591\\
14	144.653401\\
15	146.721547\\
16	144.445801\\
17	148.481932\\
18	138.645096\\
19	147.351889\\
20	144.662357\\
};
%\addlegendentry{$\varepsilon_{\text{cv}}=10^{-3}$};

\addplot [color=red,solid,line width=1.0pt]
  table[row sep=crcr]{%
1	121.747708\\
2	169.38097\\
3	191.482439\\
4	241.841167\\
5	196.583574\\
6	214.711833\\
7	226.49788\\
8	263.998732\\
9	223.366553\\
10	211.046026\\
11	206.121747\\
12	214.724002\\
13	226.915471\\
14	219.424097\\
15	240.609912\\
16	218.412067\\
17	207.718098\\
18	231.868949\\
19	267.103662\\
20	211.914743\\
};
%\addlegendentry{$\varepsilon_{\text{cv}}=10^{-4}$};

\addplot [color=green,solid,line width=1.0pt]
  table[row sep=crcr]{%
1	126.940582\\
2	259.406203\\
3	309.819547\\
4	377.527374\\
5	336.736105\\
6	387.213362\\
7	392.50116\\
8	354.842674\\
9	443.437456\\
10	367.7231\\
11	378.599385\\
12	367.965085\\
13	381.145219\\
14	342.505101\\
15	391.175307\\
16	361.615104\\
17	359.257049\\
18	339.095303\\
19	333.958773\\
20	377.271626\\
};
%\addlegendentry{$\varepsilon_{\text{cv}}=10^{-5}$};

\end{axis}
\end{tikzpicture}%

%% file: cpu_time.tex
\definecolor{mycolor1}{rgb}{1.00000,0.00000,1.00000}%
\definecolor{mycolor2}{rgb}{0.00000,1.00000,1.00000}%
\begin{tikzpicture}

\begin{axis}[%
width=1.268\figureheight,
height=\figureheight,
at={(0\figureheight,0\figureheight)},
scale only axis,
separate axis lines,
every outer x axis line/.append style={black},
every x tick label/.append style={font=\color{black}},
every x tick/.append style={black},
xmin=0,
xmax=20,
xlabel={Number of iterations},
every outer y axis line/.append style={black},
every y tick label/.append style={font=\color{black}},
every y tick/.append style={black},
ymin=0,
ymax=350,
ylabel={CPU time (s)},
axis background/.style={fill=white},
xmajorgrids,
ymajorgrids,
legend style={legend cell align=left, align=left, draw=black}
]
\addplot [color=mycolor1, line width=1.0pt]
  table[row sep=crcr]{%
1	94.961037\\
2	96.486102\\
3	100.429646\\
4	101.398912\\
5	99.247174\\
6	98.447348\\
7	99.77097\\
8	101.296769\\
9	101.91025\\
10	103.75204\\
11	98.823006\\
12	100.218087\\
13	94.446287\\
14	95.448088\\
15	95.410049\\
16	102.96006\\
17	96.480508\\
18	96.220784\\
19	95.930688\\
20	99.37199\\
};
\addlegendentry{$\varepsilon_{\text{cv}}=10^{-2}$};

\addplot [color=mycolor2, line width=1.0pt]
  table[row sep=crcr]{%
1	99.415563\\
2	103.882396\\
3	109.880284\\
4	123.398258\\
5	124.947046\\
6	109.830554\\
7	110.793023\\
8	114.35586\\
9	141.507042\\
10	118.595862\\
11	110.617532\\
12	113.201358\\
13	138.694745\\
14	125.268379\\
15	121.626299\\
16	113.724981\\
17	108.772941\\
18	115.137177\\
19	118.078092\\
20	112.516932\\
};
\addlegendentry{$\varepsilon_{\text{cv}}=10^{-3}$};

\addplot [color=red, line width=1.0pt]
  table[row sep=crcr]{%
1	97.408212\\
2	138.154411\\
3	148.865786\\
4	169.046004\\
5	175.393302\\
6	160.360997\\
7	151.120532\\
8	150.363304\\
9	154.130747\\
10	187.055159\\
11	162.64174\\
12	145.999911\\
13	150.756426\\
14	156.888162\\
15	153.82795\\
16	159.04433\\
17	150.367372\\
18	156.774043\\
19	163.390833\\
20	160.354122\\
};
\addlegendentry{$\varepsilon_{\text{cv}}=10^{-4}$};

\addplot [color=green, line width=1.0pt]
  table[row sep=crcr]{%
1	108.819139\\
2	186.444431\\
3	211.696693\\
4	235.522006\\
5	236.429948\\
6	265.444377\\
7	240.165142\\
8	249.446569\\
9	247.172435\\
10	237.091644\\
11	224.400816\\
12	239.89302\\
13	239.898396\\
14	235.706383\\
15	236.794332\\
16	266.696836\\
17	301.331136\\
18	266.009226\\
19	243.211218\\
20	224.990554\\
};
\addlegendentry{$\varepsilon_{\text{cv}}=10^{-5}$};

\end{axis}
\end{tikzpicture}%

%% file: relaxation_parameter_no_legend.tex
\definecolor{mycolor1}{rgb}{1.00000,0.00000,1.00000}%
\definecolor{mycolor2}{rgb}{0.00000,1.00000,1.00000}%
\begin{tikzpicture}

\begin{axis}[%
width=1.268\figureheight,
height=\figureheight,
at={(0\figureheight,0\figureheight)},
scale only axis,
separate axis lines,
every outer x axis line/.append style={black},
every x tick label/.append style={font=\color{black}},
every x tick/.append style={black},
xmin=0,
xmax=20,
xlabel={Number of iterations},
every outer y axis line/.append style={black},
every y tick label/.append style={font=\color{black}},
every y tick/.append style={black},
ymin=0,
ymax=1.4,
ylabel={Relaxation parameter},
axis background/.style={fill=white},
xmajorgrids,
ymajorgrids,
yminorgrids,
legend style={legend cell align=left, align=left, draw=black}
]
\addplot [color=mycolor1, line width=1.0pt]
  table[row sep=crcr]{%
1	1\\
2	1\\
3	0.993082060898308\\
4	0.979268170882619\\
5	0.780381951133088\\
6	0.421889886746609\\
7	0.217831174512013\\
8	0.115986047671824\\
9	0.0594352718631697\\
10	0.0318474064655007\\
11	0.0143515046370355\\
12	0.00720860581143631\\
13	0.00357805588908727\\
14	0.00195595858146103\\
15	0.000954962531083144\\
16	0.000457018973614378\\
17	0.000215568915025792\\
18	0.000110016180389838\\
19	6.8034498415679e-05\\
20	2.5763836568284e-05\\
};
%\addlegendentry{$\varepsilon_{\text{cv}} = 10^{-2}$};

\addplot [color=mycolor2, line width=1.0pt]
  table[row sep=crcr]{%
1	1\\
2	1\\
3	0.993087674925416\\
4	0.980253192503484\\
5	0.994999499917481\\
6	0.600714171484318\\
7	0.352068897085812\\
8	0.175868578222567\\
9	0.0985983142422735\\
10	0.0517082902759818\\
11	0.0197199645403382\\
12	0.0128766685204325\\
13	0.0059746221999437\\
14	0.00290290467000426\\
15	0.00150378333122332\\
16	0.000796821310841657\\
17	0.000376045226538225\\
18	0.000175322705331062\\
19	9.48845150402699e-05\\
20	4.85587721414944e-05\\
};
%\addlegendentry{$\varepsilon_{\text{cv}} = 10^{-3}$};

\addplot [color=red, line width=1.0pt]
  table[row sep=crcr]{%
1	1\\
2	1\\
3	0.993088484858308\\
4	0.980197893637968\\
5	0.998856423179244\\
6	1.11650676812512\\
7	0.823190365679039\\
8	0.457567697175187\\
9	0.250888489208065\\
10	0.130676286057878\\
11	0.0647396191022343\\
12	0.0361966590431485\\
13	0.0174798168888589\\
14	0.00703710891042577\\
15	0.00403389885463086\\
16	0.00168828039515254\\
17	0.00109898637552216\\
18	0.000510826763593734\\
19	0.000209036440370484\\
20	0.000120799355452498\\
};
%\addlegendentry{$\varepsilon_{\text{cv}} = 10^{-4}$};

\addplot [color=green, line width=1.0pt]
  table[row sep=crcr]{%
1	1\\
2	1\\
3	0.993088510413416\\
4	0.980201182169025\\
5	0.998964570448125\\
6	1.13222848912368\\
7	1.16563971030489\\
8	0.663561000270022\\
9	0.419624266955341\\
10	0.216853575446323\\
11	0.106959333725091\\
12	0.0614577688482941\\
13	0.0279896810470092\\
14	0.0122677295974082\\
15	0.00695670885238706\\
16	0.00359445059761837\\
17	0.00177763189759459\\
18	0.000903686059517854\\
19	0.000468012561386868\\
20	0.000243855781376796\\
};
%\addlegendentry{$\varepsilon_{\text{cv}} = 10^{-5}$};

\end{axis}
\end{tikzpicture}%

%% file: relaxation_parameter.tex
\definecolor{mycolor1}{rgb}{1.00000,0.00000,1.00000}%
\definecolor{mycolor2}{rgb}{0.00000,1.00000,1.00000}%
\begin{tikzpicture}

\begin{axis}[%
width=1.268\figureheight,
height=\figureheight,
at={(0\figureheight,0\figureheight)},
scale only axis,
separate axis lines,
every outer x axis line/.append style={black},
every x tick label/.append style={font=\color{black}},
every x tick/.append style={black},
xmin=0,
xmax=20,
xlabel={Number of iterations},
every outer y axis line/.append style={black},
every y tick label/.append style={font=\color{black}},
every y tick/.append style={black},
ymin=0,
ymax=1.4,
ylabel={Relaxation parameter},
axis background/.style={fill=white},
xmajorgrids,
ymajorgrids,
yminorgrids,
legend style={legend cell align=left, align=left, draw=black}
]
\addplot [color=mycolor1, line width=1.0pt]
  table[row sep=crcr]{%
1	1\\
2	1\\
3	0.997485624047853\\
4	0.984098573479415\\
5	0.714370850490764\\
6	0.429310888117022\\
7	0.219553144745064\\
8	0.123823278669317\\
9	0.0619953711922955\\
10	0.0276065517575781\\
11	0.0169383084311074\\
12	0.00930260560846952\\
13	0.00440565360417066\\
14	0.00222408460042551\\
15	0.000782605909642561\\
16	0.000393445779102838\\
17	0.000237597367011339\\
18	0.000114991353694213\\
19	5.2403641060682e-05\\
20	2.69536826851834e-05\\
};
\addlegendentry{$\varepsilon_{\text{cv}} = 10^{-2}$};

\addplot [color=mycolor2, line width=1.0pt]
  table[row sep=crcr]{%
1	1\\
2	1\\
3	0.997485809169807\\
4	0.985677422558444\\
5	1.03339119666131\\
6	0.680758604918496\\
7	0.408549665101266\\
8	0.139313347488461\\
9	0.0560760140413451\\
10	0.0458258981084213\\
11	0.0246324922596268\\
12	0.0123851382674805\\
13	0.00670506723867599\\
14	0.00343533886207484\\
15	0.00138111893240651\\
16	0.000493484933503815\\
17	0.000214041830128205\\
18	0.00013196140482117\\
19	8.14685978419161e-05\\
20	4.18186831317547e-05\\
};
\addlegendentry{$\varepsilon_{\text{cv}} = 10^{-3}$};

\addplot [color=red, line width=1.0pt]
  table[row sep=crcr]{%
1	1\\
2	1\\
3	0.997483296383369\\
4	0.985539307811957\\
5	1.04481442524357\\
6	1.12214879220255\\
7	0.649657122131399\\
8	0.365796550228566\\
9	0.230774400568042\\
10	0.114293657154837\\
11	0.0603470581476444\\
12	0.0332195897078455\\
13	0.0148248212566307\\
14	0.00777646463147022\\
15	0.00406843679532172\\
16	0.00199736024732236\\
17	0.00101088136167746\\
18	0.000525504488544177\\
19	0.000196699112160435\\
20	0.000118458825692555\\
};
\addlegendentry{$\varepsilon_{\text{cv}} = 10^{-4}$};

\addplot [color=green, line width=1.0pt]
  table[row sep=crcr]{%
1	1\\
2	1\\
3	0.997483501095794\\
4	0.985565848520018\\
5	1.04510792627575\\
6	1.15356031400754\\
7	1.16789653618465\\
8	0.633349322028259\\
9	0.367721715749131\\
10	0.202384362243854\\
11	0.104742594624667\\
12	0.0546662211618187\\
13	0.0299802025134833\\
14	0.0143237799605059\\
15	0.00737712307636672\\
16	0.00355577615818079\\
17	0.00189888891368969\\
18	0.000853602660004076\\
19	0.000421330503151731\\
20	0.000216299374018625\\
};
\addlegendentry{$\varepsilon_{\text{cv}} = 10^{-5}$};

\end{axis}
\end{tikzpicture}%

%% file: nb_samples_4_tol_no_legend.tex
\definecolor{mycolor1}{rgb}{1.00000,0.00000,1.00000}%
\definecolor{mycolor2}{rgb}{0.00000,1.00000,1.00000}%
\begin{tikzpicture}

\begin{axis}[%
width=1.268\figureheight,
height=\figureheight,
at={(0\figureheight,0\figureheight)},
scale only axis,
separate axis lines,
every outer x axis line/.append style={black},
every x tick label/.append style={font=\color{black}},
every x tick/.append style={black},
xmin=0,
xmax=20,
xlabel={Number of iterations},
every outer y axis line/.append style={black},
every y tick label/.append style={font=\color{black}},
every y tick/.append style={black},
ymin=0,
ymax=500,
ylabel={Number of samples},
axis background/.style={fill=white},
xmajorgrids,
ymajorgrids,
legend style={legend cell align=left, align=left, draw=black}
]
\addplot [color=mycolor1, line width=1.0pt]
  table[row sep=crcr]{%
1	21\\
2	24\\
3	27\\
4	30\\
5	41\\
6	30\\
7	37\\
8	57\\
9	30\\
10	30\\
11	27\\
12	41\\
13	57\\
14	41\\
15	37\\
16	57\\
17	51\\
18	46\\
19	27\\
20	41\\
};
%\addlegendentry{$\varepsilon_{\text{cv}}=10^{-2}$};

\addplot [color=mycolor2, line width=1.0pt]
  table[row sep=crcr]{%
1	30\\
2	63\\
3	70\\
4	70\\
5	85\\
6	63\\
7	85\\
8	77\\
9	77\\
10	85\\
11	85\\
12	70\\
13	70\\
14	85\\
15	77\\
16	70\\
17	70\\
18	63\\
19	63\\
20	94\\
};
%\addlegendentry{$\varepsilon_{\text{cv}}=10^{-3}$};

\addplot [color=red, line width=1.0pt]
  table[row sep=crcr]{%
1	30\\
2	127\\
3	170\\
4	170\\
5	170\\
6	206\\
7	154\\
8	250\\
9	170\\
10	187\\
11	187\\
12	170\\
13	206\\
14	170\\
15	187\\
16	206\\
17	187\\
18	187\\
19	154\\
20	170\\
};
%\addlegendentry{$\varepsilon_{\text{cv}}=10^{-4}$};

\addplot [color=green, line width=1.0pt]
  table[row sep=crcr]{%
1	41\\
2	275\\
3	334\\
4	405\\
5	368\\
6	405\\
7	405\\
8	368\\
9	405\\
10	405\\
11	405\\
12	368\\
13	446\\
14	334\\
15	405\\
16	405\\
17	368\\
18	368\\
19	368\\
20	446\\
};
%\addlegendentry{$\varepsilon_{\text{cv}}=10^{-5}$};

\end{axis}
\end{tikzpicture}%

%% file: nb_samples_4_tol.tex
\definecolor{mycolor1}{rgb}{1.00000,0.00000,1.00000}%
\definecolor{mycolor2}{rgb}{0.00000,1.00000,1.00000}%
\begin{tikzpicture}

\begin{axis}[%
width=1.268\figureheight,
height=\figureheight,
at={(0\figureheight,0\figureheight)},
scale only axis,
separate axis lines,
every outer x axis line/.append style={black},
every x tick label/.append style={font=\color{black}},
every x tick/.append style={black},
xmin=0,
xmax=20,
xlabel={Number of iterations},
every outer y axis line/.append style={black},
every y tick label/.append style={font=\color{black}},
every y tick/.append style={black},
ymin=0,
ymax=300,
ylabel={Number of samples},
axis background/.style={fill=white},
xmajorgrids,
ymajorgrids,
legend style={legend cell align=left, align=left, draw=black}
]
\addplot [color=mycolor1, line width=1.0pt]
  table[row sep=crcr]{%
1	21\\
2	24\\
3	24\\
4	27\\
5	27\\
6	41\\
7	27\\
8	21\\
9	30\\
10	24\\
11	27\\
12	24\\
13	30\\
14	27\\
15	21\\
16	27\\
17	30\\
18	21\\
19	27\\
20	27\\
};
\addlegendentry{$\varepsilon_{\text{cv}}=10^{-2}$};

\addplot [color=mycolor2, line width=1.0pt]
  table[row sep=crcr]{%
1	21\\
2	33\\
3	41\\
4	57\\
5	63\\
6	51\\
7	46\\
8	57\\
9	51\\
10	51\\
11	46\\
12	51\\
13	46\\
14	51\\
15	63\\
16	51\\
17	41\\
18	51\\
19	70\\
20	63\\
};
\addlegendentry{$\varepsilon_{\text{cv}}=10^{-3}$};

\addplot [color=red, line width=1.0pt]
  table[row sep=crcr]{%
1	33\\
2	94\\
3	115\\
4	140\\
5	104\\
6	115\\
7	115\\
8	115\\
9	115\\
10	127\\
11	115\\
12	104\\
13	115\\
14	127\\
15	115\\
16	104\\
17	115\\
18	127\\
19	127\\
20	127\\
};
\addlegendentry{$\varepsilon_{\text{cv}}=10^{-4}$};

\addplot [color=green, line width=1.0pt]
  table[row sep=crcr]{%
1	33\\
2	154\\
3	187\\
4	227\\
5	206\\
6	206\\
7	227\\
8	227\\
9	250\\
10	227\\
11	206\\
12	227\\
13	227\\
14	227\\
15	227\\
16	250\\
17	206\\
18	250\\
19	227\\
20	227\\
};
\addlegendentry{$\varepsilon_{\text{cv}}=10^{-5}$};

\end{axis}
\end{tikzpicture}%

%% file: dim_stochastic_basis_4_tol_no_legend.tex
\definecolor{mycolor1}{rgb}{1.00000,0.00000,1.00000}%
\definecolor{mycolor2}{rgb}{0.00000,1.00000,1.00000}%
\begin{tikzpicture}

\begin{axis}[%
width=1.268\figureheight,
height=\figureheight,
at={(0\figureheight,0\figureheight)},
scale only axis,
separate axis lines,
every outer x axis line/.append style={black},
every x tick label/.append style={font=\color{black}},
every x tick/.append style={black},
xmin=0,
xmax=20,
xlabel={Number of iterations},
every outer y axis line/.append style={black},
every y tick label/.append style={font=\color{black}},
every y tick/.append style={black},
ymin=0,
ymax=160,
ylabel={Dimension of stochastic basis},
axis background/.style={fill=white},
xmajorgrids,
ymajorgrids,
legend style={legend cell align=left, align=left, draw=black}
]
\addplot [color=mycolor1, line width=1.0pt]
  table[row sep=crcr]{%
1	3\\
2	3\\
3	5\\
4	5\\
5	5\\
6	3\\
7	6\\
8	3\\
9	5\\
10	3\\
11	3\\
12	5\\
13	9\\
14	4\\
15	3\\
16	13\\
17	5\\
18	5\\
19	5\\
20	5\\
};
%\addlegendentry{$\varepsilon_{\text{cv}}=10^{-2}$};

\addplot [color=mycolor2, line width=1.0pt]
  table[row sep=crcr]{%
1	6\\
2	9\\
3	11\\
4	11\\
5	16\\
6	14\\
7	10\\
8	13\\
9	11\\
10	11\\
11	12\\
12	11\\
13	11\\
14	13\\
15	11\\
16	11\\
17	11\\
18	11\\
19	11\\
20	11\\
};
%\addlegendentry{$\varepsilon_{\text{cv}}=10^{-3}$};

\addplot [color=red, line width=1.0pt]
  table[row sep=crcr]{%
1	8\\
2	19\\
3	34\\
4	32\\
5	29\\
6	30\\
7	28\\
8	29\\
9	30\\
10	36\\
11	34\\
12	29\\
13	31\\
14	31\\
15	29\\
16	32\\
17	28\\
18	28\\
19	36\\
20	33\\
};
%\addlegendentry{$\varepsilon_{\text{cv}}=10^{-4}$};

\addplot [color=green, line width=1.0pt]
  table[row sep=crcr]{%
1	12\\
2	77\\
3	59\\
4	70\\
5	73\\
6	75\\
7	70\\
8	76\\
9	66\\
10	77\\
11	90\\
12	84\\
13	69\\
14	76\\
15	74\\
16	85\\
17	73\\
18	75\\
19	72\\
20	67\\
};
%\addlegendentry{$\varepsilon_{\text{cv}}=10^{-5}$};

\addplot [color=mycolor1, dashed, line width=1.0pt]
  table[row sep=crcr]{%
1	4\\
2	6\\
3	7\\
4	7\\
5	11\\
6	7\\
7	9\\
8	16\\
9	8\\
10	7\\
11	8\\
12	16\\
13	10\\
14	11\\
15	12\\
16	18\\
17	15\\
18	12\\
19	7\\
20	12\\
};
%\addlegendentry{$\varepsilon_{\text{cv}}=10^{-2}$};

\addplot [color=mycolor2, dashed, line width=1.0pt]
  table[row sep=crcr]{%
1	6\\
2	15\\
3	19\\
4	18\\
5	20\\
6	18\\
7	19\\
8	23\\
9	27\\
10	28\\
11	24\\
12	19\\
13	25\\
14	33\\
15	22\\
16	19\\
17	24\\
18	20\\
19	21\\
20	25\\
};
%\addlegendentry{$\varepsilon_{\text{cv}}=10^{-3}$};

\addplot [color=red, dashed, line width=1.0pt]
  table[row sep=crcr]{%
1	9\\
2	31\\
3	46\\
4	64\\
5	52\\
6	64\\
7	52\\
8	97\\
9	65\\
10	57\\
11	60\\
12	67\\
13	76\\
14	57\\
15	65\\
16	70\\
17	73\\
18	51\\
19	62\\
20	67\\
};
%\addlegendentry{$\varepsilon_{\text{cv}}=10^{-4}$};

\addplot [color=green, dashed, line width=1.0pt]
  table[row sep=crcr]{%
1	12\\
2	78\\
3	114\\
4	138\\
5	117\\
6	136\\
7	137\\
8	117\\
9	143\\
10	131\\
11	146\\
12	146\\
13	151\\
14	126\\
15	157\\
16	142\\
17	143\\
18	118\\
19	138\\
20	153\\
};
%\addlegendentry{$\varepsilon_{\text{cv}}=10^{-5}$};

\end{axis}
\end{tikzpicture}%

%% file: dim_stochastic_basis_4_tol.tex
\definecolor{mycolor1}{rgb}{1.00000,0.00000,1.00000}%
\definecolor{mycolor2}{rgb}{0.00000,1.00000,1.00000}%
\begin{tikzpicture}

\begin{axis}[%
width=1.268\figureheight,
height=\figureheight,
at={(0\figureheight,0\figureheight)},
scale only axis,
separate axis lines,
every outer x axis line/.append style={black},
every x tick label/.append style={font=\color{black}},
every x tick/.append style={black},
xmin=0,
xmax=20,
xlabel={Number of iterations},
every outer y axis line/.append style={black},
every y tick label/.append style={font=\color{black}},
every y tick/.append style={black},
ymin=0,
ymax=100,
ylabel={Dimension of stochastic basis},
axis background/.style={fill=white},
xmajorgrids,
ymajorgrids,
legend style={legend cell align=left, align=left, draw=black}
]
\addplot [color=mycolor1, line width=1.0pt]
  table[row sep=crcr]{%
1	4\\
2	4\\
3	3\\
4	3\\
5	3\\
6	3\\
7	3\\
8	3\\
9	3\\
10	3\\
11	3\\
12	3\\
13	3\\
14	3\\
15	3\\
16	3\\
17	4\\
18	3\\
19	4\\
20	3\\
};
\addlegendentry{$\varepsilon_{\text{cv}}=10^{-2}$};

\addplot [color=mycolor2, line width=1.0pt]
  table[row sep=crcr]{%
1	4\\
2	7\\
3	10\\
4	8\\
5	12\\
6	7\\
7	9\\
8	11\\
9	10\\
10	8\\
11	8\\
12	8\\
13	9\\
14	7\\
15	9\\
16	8\\
17	8\\
18	8\\
19	12\\
20	8\\
};
\addlegendentry{$\varepsilon_{\text{cv}}=10^{-3}$};

\addplot [color=red, line width=1.0pt]
  table[row sep=crcr]{%
1	6\\
2	13\\
3	16\\
4	17\\
5	17\\
6	15\\
7	19\\
8	15\\
9	15\\
10	15\\
11	15\\
12	17\\
13	19\\
14	18\\
15	26\\
16	16\\
17	15\\
18	15\\
19	18\\
20	23\\
};
\addlegendentry{$\varepsilon_{\text{cv}}=10^{-4}$};

\addplot [color=green, line width=1.0pt]
  table[row sep=crcr]{%
1	11\\
2	35\\
3	36\\
4	37\\
5	40\\
6	50\\
7	41\\
8	36\\
9	40\\
10	44\\
11	39\\
12	39\\
13	38\\
14	37\\
15	39\\
16	38\\
17	43\\
18	37\\
19	37\\
20	39\\
};
\addlegendentry{$\varepsilon_{\text{cv}}=10^{-5}$};

\addplot [color=mycolor1, dashed, line width=1.0pt]
  table[row sep=crcr]{%
1	4\\
2	6\\
3	6\\
4	8\\
5	5\\
6	11\\
7	6\\
8	5\\
9	6\\
10	6\\
11	7\\
12	5\\
13	7\\
14	7\\
15	5\\
16	6\\
17	6\\
18	5\\
19	7\\
20	5\\
};
%\addlegendentry{$\varepsilon_{\text{cv}}=10^{-2}$};

\addplot [color=mycolor2, dashed, line width=1.0pt]
  table[row sep=crcr]{%
1	4\\
2	8\\
3	12\\
4	14\\
5	15\\
6	15\\
7	12\\
8	13\\
9	11\\
10	16\\
11	13\\
12	11\\
13	14\\
14	16\\
15	14\\
16	15\\
17	12\\
18	17\\
19	16\\
20	16\\
};
%\addlegendentry{$\varepsilon_{\text{cv}}=10^{-3}$};

\addplot [color=red, dashed, line width=1.0pt]
  table[row sep=crcr]{%
1	6\\
2	24\\
3	34\\
4	51\\
5	28\\
6	33\\
7	34\\
8	36\\
9	32\\
10	29\\
11	29\\
12	31\\
13	33\\
14	29\\
15	30\\
16	31\\
17	29\\
18	35\\
19	36\\
20	38\\
};
%\addlegendentry{$\varepsilon_{\text{cv}}=10^{-4}$};

\addplot [color=green, dashed, line width=1.0pt]
  table[row sep=crcr]{%
1	9\\
2	52\\
3	56\\
4	76\\
5	64\\
6	76\\
7	63\\
8	68\\
9	74\\
10	71\\
11	66\\
12	74\\
13	78\\
14	82\\
15	79\\
16	85\\
17	75\\
18	86\\
19	76\\
20	73\\
};
%\addlegendentry{$\varepsilon_{\text{cv}}=10^{-5}$};

\end{axis}
\end{tikzpicture}%

%% file: m2an170020.bbl
\begin{thebibliography}{10}

\bibitem{Nou09}
Anthony Nouy.
\newblock Recent developments in spectral stochastic methods for the numerical
  solution of stochastic partial differential equations.
\newblock {\em Archives of Computational Methods in Engineering}, 16:251--285,
  2009.

\bibitem{Xiu09}
D.~Xiu.
\newblock {Fast numerical methods for stochastic computations: a review}.
\newblock {\em Communications in computational physics}, 5(2-4):242--272, 2009.

\bibitem{LeMai10}
O.~P. Le~Ma{\^\i}tre and O.~M. Knio.
\newblock {\em {Spectral Methods for Uncertainty Quantification With
  Applications to Computational Fluid Dynamics}}.
\newblock Springer Netherlands, 2010.

\bibitem{Hou97}
Thomas~Y. Hou and Xiao-Hui Wu.
\newblock {A Multiscale Finite Element Method for Elliptic Problems in
  Composite Materials and Porous Media}.
\newblock {\em Journal of Computational Physics}, 134(1):169--189, 1997.

\bibitem{Hug98}
Thomas~J.R. Hughes, Gonzalo~R. Feij\'oo, Luca Mazzei, and Jean-Baptiste Quincy.
\newblock {The variational multiscale method---a paradigm for computational
  mechanics}.
\newblock {\em Computer Methods in Applied Mechanics and Engineering},
  166(1-2):3--24, 1998.
\newblock Advances in Stabilized Methods in Computational Mechanics.

\bibitem{E03}
Weinan E and Bj{\"o}rn Engquist.
\newblock The heterogeneous multiscale methods.
\newblock {\em Communications in Mathematical Sciences}, 1(1):87--132, 2003.

\bibitem{Xu07}
X.~Frank Xu.
\newblock {A multiscale stochastic finite element method on elliptic problems
  involving uncertainties}.
\newblock {\em Computer Methods in Applied Mechanics and Engineering},
  196(25-28):2723--2736, 2007.

\bibitem{Nar05}
Velamur Asokan~Badri Narayanan and Nicholas Zabaras.
\newblock {Variational multiscale stabilized FEM formulations for transport
  equations: stochastic advection--diffusion and incompressible stochastic
  Navier--Stokes equations}.
\newblock {\em Journal of Computational Physics}, 202(1):94--133, 2005.

\bibitem{Aso06}
Badrinarayanan~Velamur Asokan and Nicholas Zabaras.
\newblock {A stochastic variational multiscale method for diffusion in
  heterogeneous random media}.
\newblock {\em Journal of Computational Physics}, 218(2):654--676, 2006.

\bibitem{Gan07bis}
Baskar Ganapathysubramanian and Nicholas Zabaras.
\newblock {Modeling diffusion in random heterogeneous media: Data-driven
  models, stochastic collocation and the variational multiscale method}.
\newblock {\em Journal of Computational Physics}, 226(1):326--353, 2007.

\bibitem{Dos08}
P.~Dostert, Y.~Efendiev, and T.Y. Hou.
\newblock {Multiscale finite element methods for stochastic porous media flow
  equations and application to uncertainty quantification}.
\newblock {\em Computer Methods in Applied Mechanics and Engineering},
  197(43-44):3445--3455, 2008.

\bibitem{Xu09}
X.~Frank Xu, Xi~Chen, and Lihua Shen.
\newblock {A Green-function-based multiscale method for uncertainty
  quantification of finite body random heterogeneous materials}.
\newblock {\em Computers \& Structures}, 87(21-22):1416--1426, 2009.

\bibitem{Gan09}
B.~Ganapathysubramanian and N.~Zabaras.
\newblock {A stochastic multiscale framework for modeling flow through random
  heterogeneous porous media}.
\newblock {\em Journal of Computational Physics}, 228(2):591--618, 2009.

\bibitem{Gin10}
V.~Ginting, A.~M\r{a}lqvist, and M.~Presho.
\newblock {A Novel Method for Solving Multiscale Elliptic Problems with
  Randomly Perturbed Data}.
\newblock {\em Multiscale Modeling \& Simulation}, 8(3):977--996, 2010.

\bibitem{LeBri14}
{Le Bris, Claude}, {Legoll, Fr{\'e}d{\'e}ric}, and {Thomines, Florian}.
\newblock Multiscale finite element approach for ``weakly'' random problems and
  related issues.
\newblock {\em ESAIM: M2AN}, 48(3):815--858, 2014.

\bibitem{Jin07}
C.~Jin, X.~Cai, and C.~Li.
\newblock {Parallel Domain Decomposition Methods for Stochastic Elliptic
  Equations}.
\newblock {\em SIAM Journal on Scientific Computing}, 29(5):2096--2114, 2007.

\bibitem{Zha08}
Kai Zhang, Ran Zhang, Yunguang Yin, and Shi Yu.
\newblock {Domain decomposition methods for linear and semilinear elliptic
  stochastic partial differential equations}.
\newblock {\em Applied Mathematics and Computation}, 195(2):630--640, 2008.

\bibitem{Sar09}
Abhijit Sarkar, Nabil Benabbou, and Roger Ghanem.
\newblock {Domain decomposition of stochastic PDEs: Theoretical formulations}.
\newblock {\em International Journal for Numerical Methods in Engineering},
  77(5):689--701, 2009.

\bibitem{Gan11}
B.~Ganis, I.~Yotov, and M.~Zhong.
\newblock {A Stochastic Mortar Mixed Finite Element Method for Flow in Porous
  Media with Multiple Rock Types}.
\newblock {\em SIAM Journal on Scientific Computing}, 33(3):1439--1474, 2011.

\bibitem{Whe11}
Mary~F. Wheeler, Tim Wildey, and Ivan Yotov.
\newblock {A multiscale preconditioner for stochastic mortar mixed finite
  elements}.
\newblock {\em Computer Methods in Applied Mechanics and Engineering},
  200(9-12):1251--1262, 2011.

\bibitem{Ver96}
R.~Verf{\"u}rth.
\newblock {\em {A Review of A Posteriori Error Estimation and Adaptive
  Mesh-refinement Techniques}}.
\newblock Wiley-Teubner, Stuttgart, 1996.

\bibitem{Ste97}
E.~Stein and S.~Ohnimus.
\newblock {Coupled model- and solution-adaptivity in the finite-element
  method}.
\newblock {\em Computer Methods in Applied Mechanics and Engineering},
  150(1-4):327--350, 1997.

\bibitem{Bel99bis}
T.~Belytschko and T.~Black.
\newblock {Elastic crack growth in finite elements with minimal remeshing}.
\newblock {\em International Journal for Numerical Methods in Engineering},
  45(5):601--620, 1999.

\bibitem{Moe99}
N.~Mo{\"e}s, J.~Dolbow, and T.~Belytschko.
\newblock {A finite element method for crack growth without remeshing}.
\newblock {\em International Journal for Numerical Methods in Engineering},
  46(1):131--150, 1999.

\bibitem{Str00bis}
T.~Strouboulis, I.~Babu\v{s}ka, and K.~Copps.
\newblock {The design and analysis of the Generalized Finite Element Method}.
\newblock {\em Computer Methods in Applied Mechanics and Engineering},
  181(1-3):43--69, 2000.

\bibitem{Lio99}
Jacques-Louis Lions and Olivier Pironneau.
\newblock {Domain decomposition methods for CAD}.
\newblock {\em Comptes Rendus de l'Acad\'emie des Sciences - Series I -
  Mathematics}, 328(1):73--80, 1999.

\bibitem{Glo05}
Roland Glowinski, Jiwen He, Alexei Lozinski, Jacques Rappaz, and Jo{\"e}l
  Wagner.
\newblock {Finite element approximation of multi-scale elliptic problems using
  patches of elements}.
\newblock {\em Numerische Mathematik}, 101(4):663--687, 2005.

\bibitem{He07}
Jiwen He, Alexei Lozinski, and Jacques Rappaz.
\newblock {Accelerating the method of finite element patches using
  approximately harmonic functions}.
\newblock {\em Comptes Rendus Mathematique}, 345(2):107--112, 2007.

\bibitem{Ste83}
Joseph~L Steger, F~Carroll Dougherty, and John~A Benek.
\newblock {A Chimera grid scheme}.
\newblock In K.N. Ghia and U.~Ghia, editors, {\em Advances in Grid Generation},
  volume~5, pages 59--69, New York, 1983. American Society of Mechanical
  Engineers, FED.

\bibitem{Bre01}
Franco Brezzi, Jacques-Louis Lions, and Olivier Pironneau.
\newblock {Analysis of a Chimera method}.
\newblock {\em Comptes Rendus de l'Acad\'emie des Sciences - Series I -
  Mathematics}, 332(7):655--660, 2001.

\bibitem{Pir09}
Olivier Pironneau, Alexei Lozinski, Alain Perronnet, and Fr\'ed\'eric Hecht.
\newblock {Numerical zoom for multiscale problems with an application to flows
  through porous media}.
\newblock {\em Discrete \& Continuous Dynamical Systems - A}, 23(1/2):265--280,
  2009.

\bibitem{Gen09}
Lionel Gendre, Olivier Allix, Pierre Gosselet, and Fran{\c c}ois Comte.
\newblock {Non-intrusive and exact global/local techniques for structural
  problems with local plasticity}.
\newblock {\em Computational Mechanics}, 44:233--245, 2009.

\bibitem{Loz10}
Alexei Lozinski.
\newblock {\em {M\'ethodes num\'eriques et mod\'elisation pour certains
  probl\`emes multi-\'echelles}}.
\newblock Habilitation \`a diriger des recherches, Universit\'e Paul Sabatier,
  Toulouse 3, France, 2010.

\bibitem{Gen11}
L.~Gendre, O.~Allix, and P.~Gosselet.
\newblock {A two-scale approximation of the Schur complement and its use for
  non-intrusive coupling}.
\newblock {\em International Journal for Numerical Methods in Engineering},
  87(9):889--905, 2011.

\bibitem{Hag12}
Corinna Hager, Patrice Hauret, Patrick~Le Tallec, and Barbara~I. Wohlmuth.
\newblock Solving dynamic contact problems with local refinement in space and
  time.
\newblock {\em Computer Methods in Applied Mechanics and Engineering},
  201-204(0):25--41, 2012.

\bibitem{Gra08}
A.~Gravouil, J.~Rannou, and M.-C. Ba{\"\i}etto.
\newblock {A Local Multi-grid X-FEM approach for 3D fatigue crack growth}.
\newblock {\em International Journal of Material Forming}, 1(1):1103--1106,
  2008.

\bibitem{Ran09}
J.~Rannou, A.~Gravouil, and M.~C. Ba{\"\i}etto-Dubourg.
\newblock {A local multigrid X-FEM strategy for 3-D crack propagation}.
\newblock {\em International Journal for Numerical Methods in Engineering},
  77(4):581--600, 2009.

\bibitem{Pas11}
J.~C. Passieux, A.~Gravouil, J.~R\'ethor\'e, and M.~C. Ba{\"\i}etto.
\newblock {Direct estimation of generalized stress intensity factors using a
  three-scale concurrent multigrid X-FEM}.
\newblock {\em International Journal for Numerical Methods in Engineering},
  85(13):1648--1666, 2011.

\bibitem{Pas13}
Jean-Charles Passieux, Julien R{\'e}thor{\'e}, Anthony Gravouil, and
  Marie-Christine Ba{\"\i}etto.
\newblock {Local/global non-intrusive crack propagation simulation using a
  multigrid X-FEM solver}.
\newblock {\em Computational Mechanics}, 52(6):1381--1393, 2013.

\bibitem{Dhia98}
Hachmi~Ben Dhia.
\newblock {Probl\`emes m\'ecaniques multi-\'echelles: la m\'ethode Arlequin -
  Multiscale mechanical problems: the Arlequin method}.
\newblock {\em {Comptes Rendus de l'Acad\'emie des Sciences - Series IIB -
  Mechanics-Physics-Astronomy}}, 326(12):899--904, 1998.

\bibitem{Cha08b}
Ludovic Chamoin, J.T. Oden, and Serge Prudhomme.
\newblock {A stochastic coupling method for atomic-to-continuum Monte-Carlo
  simulations}.
\newblock {\em Computer Methods in Applied Mechanics and Engineering},
  197(43-44):3530--3546, 2008.

\bibitem{Cot11}
R.~Cottereau, D.~Clouteau, H.~Ben Dhia, and C.~Zaccardi.
\newblock {A stochastic-deterministic coupling method for continuum mechanics}.
\newblock {\em Computer Methods in Applied Mechanics and Engineering},
  200(47-48):3280--3288, 2011.

\bibitem{Che13a}
M.~Chevreuil, A.~Nouy, and E.~Safatly.
\newblock {A multiscale method with patch for the solution of stochastic
  partial differential equations with localized uncertainties}.
\newblock {\em Computer Methods in Applied Mechanics and Engineering},
  255:255--274, 2013.

\bibitem{Hen14}
{Henning, Patrick}, {M\r{a}lqvist, Axel}, and {Peterseim, Daniel}.
\newblock {A localized orthogonal decomposition method for semi-linear elliptic
  problems}.
\newblock {\em ESAIM: M2AN}, 48(5):1331--1349, 2014.

\bibitem{Efe04}
Y.~Efendiev, T.~Y. Hou, and V.~Ginting.
\newblock {Multiscale Finite Element Methods for Nonlinear Problems and Their
  Applications}.
\newblock {\em Commun. Math. Sci.}, 2(4):553--589, 12 2004.

\bibitem{Efe09}
Yalchin Efendiev and Thomas~Y Hou.
\newblock {\em {Multiscale Finite Element Methods: Theory and Applications}},
  volume~4 of {\em Surveys and Tutorials in the Applied Mathematical Sciences}.
\newblock Springer-Verlag New York, 2009.

\bibitem{Nor10}
Jan~Martin Nordbotten.
\newblock {\em {Variational and Heterogeneous Multiscale Methods}}, pages
  713--720.
\newblock Springer Berlin Heidelberg, Berlin, Heidelberg, 2010.

\bibitem{Hen15}
Patrick Henning and Mario Ohlberger.
\newblock Error control and adaptivity for heterogeneous multiscale
  approximations of nonlinear monotone problems.
\newblock {\em Discrete and Continuous Dynamical Systems - Series S},
  8(1):119--150, 2015.

\bibitem{All11}
Olivier Allix, Lionel Gendre, Pierre Gosselet, and Guillaume Guguin.
\newblock {Non-intrusive Coupling: An Attempt to Merge Industrial and Research
  Software Capabilities}.
\newblock In Dana Mueller-Hoeppe, Stefan Loehnert, and Stefanie Reese, editors,
  {\em Recent Developments and Innovative Applications in Computational
  Mechanics}, chapter~15, pages 125--133. Springer Berlin Heidelberg, 2011.

\bibitem{Chk13}
Abdellah Chkifa, Albert Cohen, Ronald DeVore, and Christoph Schwab.
\newblock {Sparse adaptive Taylor approximation algorithms for parametric and
  stochastic elliptic PDEs}.
\newblock {\em ESAIM: M2AN}, 47(1):253--280, 2013.

\bibitem{Baz69}
A.~D. Bazykin.
\newblock {Hypothetical Mechanism of Speciaton}.
\newblock {\em Evolution}, 23(4):685--687, 1969.

\bibitem{Aro75}
D.~G. Aronson and H.~F. Weinberger.
\newblock {\em {Nonlinear diffusion in population genetics, combustion, and
  nerve pulse propagation}}, pages 5--49.
\newblock Springer Berlin Heidelberg, Berlin, Heidelberg, 1975.

\bibitem{Bra04}
B.H Bradshaw-Hajek and P~Broadbridge.
\newblock {A robust cubic reaction-diffusion system for gene propagation}.
\newblock {\em Mathematical and Computer Modelling}, 39(9):1151--1163, 2004.

\bibitem{Han82}
Adel Hanna, Alan Saul, and Kenneth Showalter.
\newblock {Detailed studies of propagating fronts in the iodate oxidation of
  arsenous acid}.
\newblock {\em Journal of the American Chemical Society}, 104(14):3838--3844,
  1982.

\bibitem{Hua96}
Jie Huang and Boyd~F. Edwards.
\newblock {Pattern formation and evolution near autocatalytic reaction fronts
  in a narrow vertical slab}.
\newblock {\em Phys. Rev. E}, 54:2620--2627, Sep 1996.

\bibitem{Edw02}
Boyd~F. Edwards.
\newblock {Poiseuille Advection of Chemical Reaction Fronts}.
\newblock {\em Phys. Rev. Lett.}, 89:104501, Aug 2002.

\bibitem{Kae02}
Mads K{\ae}rn and Michael Menzinger.
\newblock {Propagation of Excitation Pulses and Autocatalytic Fronts in
  Packed-Bed Reactors}.
\newblock {\em The Journal of Physical Chemistry B}, 106(14):3751--3758, 2002.

\bibitem{Spa03}
Robert~S. Spangler and Boyd~F. Edwards.
\newblock {Poiseuille advection of chemical reaction fronts: Eikonal
  approximation}.
\newblock {\em The Journal of Chemical Physics}, 118(13):5911--5915, 2003.

\bibitem{Lec03}
M.~Leconte, J.~Martin, N.~Rakotomalala, and D.~Salin.
\newblock {Pattern of Reaction Diffusion Fronts in Laminar Flows}.
\newblock {\em Phys. Rev. Lett.}, 90:128302, Mar 2003.

\bibitem{Kop08}
Igor~V. Koptyug, Vladimir~V. Zhivonitko, and Renad~Z. Sagdeev.
\newblock {Advection of Chemical Reaction Fronts in a Porous Medium}.
\newblock {\em The Journal of Physical Chemistry B}, 112(4):1170--1176, 2008.
\newblock PMID: 18173259.

\bibitem{Sah13}
Sandeep Saha, Severine Atis, Dominique Salin, and Laurent Talon.
\newblock {Phase diagram of sustained wave fronts opposing the flow in
  disordered porous media}.
\newblock {\em EPL (Europhysics Letters)}, 101(3):38003, 2013.

\bibitem{Atk09}
Kendall Atkinson and Weimin Han.
\newblock {\em {Theoretical Numerical Analysis: A Functional Analysis
  Framework}}, volume~39.
\newblock Springer, 2009.

\bibitem{Xiu06}
Dongbin Xiu and Daniel~M. Tartakovsky.
\newblock {Numerical Methods for Differential Equations in Random Domains}.
\newblock {\em SIAM Journal on Scientific Computing}, 28(3):1167--1185, 2006.

\bibitem{Tar06}
Daniel~M. Tartakovsky and Dongbin Xiu.
\newblock {Stochastic analysis of transport in tubes with rough walls}.
\newblock {\em Journal of Computational Physics}, 217(1):248--259, 2006.
\newblock Uncertainty Quantification in Simulation Science.

\bibitem{Canu07}
Claudio Canuto and Tomas Kozubek.
\newblock {A fictitious domain approach to the numerical solution of PDEs in
  stochastic domains}.
\newblock {\em Numerische Mathematik}, 107(2):257, May 2007.

\bibitem{Nou08a}
A.~Nouy, A.~Cl\'ement, F.~Schoefs, and N.~Mo{\"e}s.
\newblock {An extended stochastic finite element method for solving stochastic
  partial differential equations on random domains}.
\newblock {\em Computer Methods in Applied Mechanics and Engineering},
  197(51):4663--4682, 2008.

\bibitem{Nou11b}
Anthony Nouy, Mathilde Chevreuil, and Elias Safatly.
\newblock Fictitious domain method and separated representations for the
  solution of boundary value problems on uncertain parameterized domains.
\newblock {\em Computer Methods in Applied Mechanics and Engineering},
  200(45--46):3066--3082, 2011.

\bibitem{Bel99}
Faker~Ben Belgacem.
\newblock {The Mortar finite element method with Lagrange multipliers}.
\newblock {\em Numerische Mathematik}, 84(2):173--197, 1999.

\bibitem{Woh01}
Barbara~I. Wohlmuth.
\newblock {\em {Discretization Methods and Iterative Solvers Based on Domain
  Decomposition}}, volume~17 of {\em Lecture Notes in Computational Science and
  Engineering}.
\newblock Springer, Berlin, New York, 2001.

\bibitem{Kim01}
C.~Kim, R.~Lazarov, J.~Pasciak, and P.~Vassilevski.
\newblock {Multiplier Spaces for the Mortar Finite Element Method in Three
  Dimensions}.
\newblock {\em SIAM Journal on Numerical Analysis}, 39(2):519--538, 2001.

\bibitem{Duv16}
Micka{\"e}l Duval, Jean-Charles Passieux, Michel Sala{\"u}n, and St{\'e}phane
  Guinard.
\newblock {Non-intrusive Coupling: Recent Advances and Scalable Nonlinear
  Domain Decomposition}.
\newblock {\em Archives of Computational Methods in Engineering}, 23(1):17--38,
  Mar 2016.

\bibitem{Mat05}
Hermann~G. Matthies and Andreas Keese.
\newblock {Galerkin methods for linear and nonlinear elliptic stochastic
  partial differential equations}.
\newblock {\em Computer Methods in Applied Mechanics and Engineering},
  194(12-16):1295--1331, 2005.

\bibitem{Gir14}
L.~Giraldi, A.~Litvinenko, D.~Liu, H.~Matthies, and A.~Nouy.
\newblock {To Be or Not to Be Intrusive? The Solution of Parametric and
  Stochastic Equations---the ``Plain Vanilla'' Galerkin Case}.
\newblock {\em SIAM Journal on Scientific Computing}, 36(6):A2720--A2744, 2014.

\bibitem{Gir15}
L.~Giraldi, D.~Liu, H.~G. Matthies, and A.~Nouy.
\newblock {To Be or Not to be Intrusive? The Solution of Parametric and
  Stochastic Equations---Proper Generalized Decomposition}.
\newblock {\em SIAM Journal on Scientific Computing}, 37(1):A347--A368, 2015.

\bibitem{Chk15b}
Abdellah Chkifa, Albert Cohen, Giovanni Migliorati, Fabio Nobile, and Ra\'ul
  Tempone.
\newblock {Discrete least squares polynomial approximation with random
  evaluations - application to parametric and stochastic elliptic PDEs}.
\newblock {\em ESAIM: M2AN}, 49(3):815--837, 2015.

\bibitem{Bar51}
Maurice~S. Bartlett.
\newblock {An Inverse Matrix Adjustment Arising in Discriminant Analysis}.
\newblock {\em The Annals of Mathematical Statistics}, 22(1):107--111, 1951.

\bibitem{Caw04}
Gavin~C. Cawley and Nicola~L.C. Talbot.
\newblock {Fast exact leave-one-out cross-validation of sparse least-squares
  support vector machines}.
\newblock {\em Neural Networks}, 17(10):1467--1475, 2004.

\bibitem{Iro69}
Bruce~M Irons and Robert~C Tuck.
\newblock {A version of the Aitken accelerator for computer iteration}.
\newblock {\em International Journal for Numerical Methods in Engineering},
  1(3):275--277, 1969.

\bibitem{Mac86}
Allan~J. Macleod.
\newblock {Acceleration of vector sequences by multi-dimensional $\Delta^2$
  methods}.
\newblock {\em Communications in Applied Numerical Methods}, 2(4):385--392,
  1986.

\bibitem{Kut08}
Ulrich K{\"u}ttler and Wolfgang~A. Wall.
\newblock {Fixed-point fluid--structure interaction solvers with dynamic
  relaxation}.
\newblock {\em Computational Mechanics}, 43(1):61--72, 2008.

\bibitem{Liu14}
Y.J. Liu, Q.~Sun, and X.L. Fan.
\newblock {A non-intrusive global/local algorithm with non-matching interface:
  Derivation and numerical validation}.
\newblock {\em Computer Methods in Applied Mechanics and Engineering},
  277(0):81--103, 2014.

\bibitem{Sud08}
Bruno Sudret.
\newblock {Global sensitivity analysis using polynomial chaos expansions}.
\newblock {\em Reliability Engineering \& System Safety}, 93(7):964--979, 2008.

\bibitem{Cha02}
Olivier Chapelle, Vladimir Vapnik, and Yoshua Bengio.
\newblock {Model Selection for Small Sample Regression}.
\newblock {\em Machine Learning}, 48(1-3):9--23, 2002.

\bibitem{Bla11}
G{\'e}raud Blatman and Bruno Sudret.
\newblock {Adaptive sparse polynomial chaos expansion based on least angle
  regression}.
\newblock {\em Journal of Computational Physics}, 230(6):2345--2367, 2011.

\bibitem{Rou05}
Tom\'a{\v{s}} Roub\'i{\v{c}}ek.
\newblock {\em {Nonlinear Partial Differential Equations with Applications}},
  volume 153.
\newblock Springer, 2005.

\end{thebibliography}
